\documentclass[aos]{imsart}

\RequirePackage{amsthm,amsmath,amsfonts,amssymb,{mathtools}}
\RequirePackage[numbers]{natbib}
\usepackage{bibentry}
\usepackage{bbm}
\usepackage[colorlinks=true,citecolor=blue]{hyperref}
\RequirePackage{epsfig, graphicx, color}
\RequirePackage{latexsym}
\RequirePackage{psfrag}

\allowdisplaybreaks

\usepackage{epsf,epsfig,subfigure,caption}
\usepackage{graphics,graphicx}
\usepackage{epstopdf}
\usepackage{comment}

\startlocaldefs
\newtheorem{theorem}{Theorem}
\newtheorem{lemma}{Lemma}

\newtheorem{corollary}{Corollary}
\newtheorem{remark}{Remark}

\newtheorem{assump}{Assumption}

\usepackage{enumitem}

\newcommand{\bfm}[1]{\ensuremath{\mathbf{#1}}}

   \def\bA{\bfm A}  %\def\AA{\mathbb{A}}

     \def\EE{\mathbb{E}}

     \def\II{\mathbb{I}}

     \def\OO{\mathbb{O}}
\def\bp{\bfm p}     \def\PP{\mathbb{P}}
     
\def\br{\bfm r}     \def\RR{\mathbb{R}}
\def\bs{\bfm s}     \def\SS{\mathbb{S}}

\def\calA{{\cal  A}} 
\def\calB{{\cal  B}}

\def\calE{{\cal  E}} 
 
\def\calG{{\cal  G}}

\def\calM{{\cal  M}} 
\def\calN{{\cal  N}} 
 
\def\calP{{\cal  P}}

\def\calS{{\cal  S}} 
\def\calT{{\cal  T}} 
\def\calU{{\cal  U}}

\def\calX{{\cal  X}} 
 
\def\calZ{{\cal  Z}}

\def\frakJ{\mathfrak{J}}
\def\frakE{\mathfrak{E}}
\def\frakP{\mathfrak{P}}
\def\frakM{\mathfrak{M}}
\makeatletter
\newcommand*{\rom}[1]{\expandafter\@slowromancap\romannumeral #1@}
\makeatother
\newcommand{\bfsym}[1]{\ensuremath{\boldsymbol{#1}}}

        \def\bLambda {\bfsym {\Lambda}}

\usepackage[ruled,vlined]{algorithm2e}

\DeclareMathOperator{\argmin}{argmin}

\DeclareMathOperator{\diag}{diag}

\DeclareMathOperator{\rank}{rank}

\DeclareMathOperator{\sgn}{sgn}

\DeclareMathOperator{\Var}{Var}

\DeclareMathOperator{\tr}{tr}

\def\newpage{\vfill\eject}

\def\sgn{\mbox{sgn}}

\newcommand{\beq}{\begin{equation}}
\newcommand{\eeq}{\end{equation}}
\newcommand{\beqn}{\begin{eqnarray}}
\newcommand{\eeqn}{\end{eqnarray}}
\newcommand{\beqnn}{\begin{eqnarray*}}
	\newcommand{\eeqnn}{\end{eqnarray*}}

\numberwithin{equation}{section}

\usepackage {tikz}
\usetikzlibrary {positioning}

\newcounter{CondCounter}
\newcommand{\F}{{\rm F}}

\SetKwInput{KwInput}{Input}                % Set the Input
\SetKwInput{KwOutput}{Output}              % set the Output

\def\submin{\scriptscriptstyle \sf min}
\def\submax{\scriptscriptstyle \sf max}

\usepackage{mathtools}

\endlocaldefs

\begin{document}
\begin{frontmatter}
\title{Inference for Low-rank Tensors -- No Need to Debias$^{\dagger, \ddagger}$}
%\title{A sample article title with some additional note\thanksref{T1}}
\runtitle{Inference for Low-rank Tensors}
%\thankstext{T1}{A sample of additional note to the title.}

\begin{aug}
%%%%%%%%%%%%%%%%%%%%%%%%%%%%%%%%%%%%%%%%%%%%%%
%%Only one address is permitted per author. %%
%%Only division, organization and e-mail is %%
%%included in the address.                  %%
%%Additional information can be included in %%
%%the Acknowledgments section if necessary. %%
%%%%%%%%%%%%%%%%%%%%%%%%%%%%%%%%%%%%%%%%%%%%%%
\author[A]{\fnms{Dong} \snm{Xia}\ead[label=e1]{}},
\author[B]{\fnms{Anru R.} \snm{Zhang}\ead[label=e2]{}}
\and
\author[C]{\fnms{Yuchen} \snm{Zhou}\ead[label=e3]{}}
%%%%%%%%%%%%%%%%%%%%%%%%%%%%%%%%%%%%%%%%%%%%%%
%% Addresses                                %%
%%%%%%%%%%%%%%%%%%%%%%%%%%%%%%%%%%%%%%%%%%%%%%
\address[A]{Hong Kong University of Science and Technology \printead{e1}}
\address[B]{Duke University}
%\address[C]{University of Wisconsin-Madison}
\address[C]{Princeton University}
\end{aug}

\begin{abstract}
In this paper, we consider the statistical inference for several low-rank tensor models. Specifically, in the Tucker low-rank tensor PCA or regression model, provided with any estimates achieving some attainable error rate, we develop the data-driven confidence regions for the singular subspace of the parameter tensor based on the asymptotic distribution of an updated estimate by two-iteration alternating minimization. The asymptotic distributions are established under some essential conditions on the signal-to-noise ratio (in PCA model) or sample size (in regression model). If the parameter tensor is further orthogonally decomposable, we develop the methods and non-asymptotic theory for inference on each individual singular vector. For the rank-one tensor PCA model, we establish the asymptotic distribution for general linear forms of principal components and confidence interval for each entry of the parameter tensor. Finally, numerical simulations are presented to corroborate our theoretical discoveries. 

In all these models, we observe that different from many matrix/vector settings in existing work, debiasing is not required to establish the asymptotic distribution of estimates or to make statistical inference on low-rank tensors. In fact, due to the widely observed statistical-computational-gap for low-rank tensor estimation, one usually requires stronger conditions than the statistical (or information-theoretic) limit to ensure the computationally feasible estimation is achievable. Surprisingly, such conditions ``incidentally" render a feasible low-rank tensor inference without debiasing. 
\end{abstract}



\end{frontmatter}
%%%%%%%%%%%%%%%%%%%%%%%%%%%%%%%%%%%%%%%%%%%%%%
%% Please use \tableofcontents for articles %%
%% with 50 pages and more                   %%
%%%%%%%%%%%%%%%%%%%%%%%%%%%%%%%%%%%%%%%%%%%%%%
%\tableofcontents

%%%%%%%%%%%%%%%%%%%%%%%%%%%%%%%%%%%%%%%%%%%%%%
%%%% Main text entry area:

\begin{sloppypar}
%%%%%%%%%%%%%%%%%%%%%%
\section{Introduction}\label{sec:intro}
%%%%%%%%%%%%%%%%%%%%%%

\footnotetext[2]{Dong Xia's research was partially supported by Hong Kong RGC Grant ECS 26302019 and GRF 16303320. Anru R. Zhang and Yuchen Zhou's research was partially supported by NSF Grants CAREER-1944904, NSF DMS-1811868, and grants from Wisconsin Alumni Research Foundation (WARF). }
\footnotetext[3]{The authors are listed alphabetically. This work was done while Anru R. Zhang and Yuchen Zhou were at the University of Wisconsin-Madison.}

An $m$th order tensor is a multiway array along $m$ directions. Recent years have witnessed a fast growing demand for the collection, processing, and analysis of data in the form of tensors. These tensor data commonly arise, to name a few, when features are collected from different domains, or when multiple data copies are provided by various agents or sources. For instances, the worldwide food trading flows \citep{de2015structural,jing2020community} produce a fourth order tensor (countries $\times$ countries $\times$ food $\times$ years); the online click-through data \citep{han2020optimal,sun2017provable} in e-commerce form a third order tensor (users $\times$ categories $\times$ periods); Berkeley human mortality data \citep{wilmoth2016human,zhang2019optimal} yield a third order tensor (ages $\times$ years $\times$ countries). In addition, the applications of tensor also include collaborative filtering \citep{karatzoglou2010multiverse,shah2019iterative}, recommender system design \citep{bi2018multilayer}, computational imaging \citep{zhang2020denoising}, and neuroimaging \citep{zhou2013tensor}. Researchers have made tremendous efforts to innovate effective methods for the analysis of tensor data. 

Low-rank models have rendered fundamental toolkits to analyze tensor data. A tensor $\calT\in\RR^{p_1\times\cdots\times p_m}$ has low Tucker rank (or multilinear rank) if all fibers\footnote{Here, the tensor fibers are the counterpart of matrix columns and rows for tensors. See \cite{kolda2009tensor} for a review.} of $\calT$ along different ways lie in rank-reduced subspaces of high-dimension, say $\{U_j\}_{j=1}^m$, respectively \citep{tucker1966some}. The core assumption of low-rank tensor models is that the observed data is driven by an {\it unknown} low-rank tensor $\calT$, while the Tucker low-rank conditions can significantly reduce the model complexity. Consequently, the analysis of tensor data often boils down to the estimation and inference of the low-rank tensor $\calT$ or its principal components based on the given datasets.

In the literature, a rich list of methods have been developed for the {\it estimation} of low-rank tensor $\calT$ and the associated subspace $U_j$, such as alternating minimization \citep{anandkumar2014tensor}, convex regularization \citep{tomioka2013convex,yuan2016tensor}, power iterations \citep{anandkumar2014tensor}, orthogonal iteration \citep{de2000best,zhang2018tensor}, vanilla gradient descent with spectral initialization \citep{cai2019nonconvex}, projected gradient descent \citep{chen2019non}, simultaneous gradient descent \citep{han2020optimal}, etc. However, in many practical scenarios, to enable more reliable decision making and prediction, it is important to quantify the estimation error in addition to point estimations. This task, referred to as {\it uncertainty quantification} or {\it statistical inference}, usually involves the construction of confidence intervals/regions for the unknown parameters through the development of the (approximate) distributions of the estimators. The statistical inference or uncertainty quantification for low-rank tensor models remains largely unexplored. In this paper, we aim to make an attempt to this fundamental and challenging problem. Our focus is on two basic yet important settings: low-rank tensor PCA and tensor regression, which we briefly summarize as follows. 

{\it Tensor principal component analysis (PCA)} is among the most basic problem of unsupervised inference for low-rank tensors. We consider the tensor PCA model \citep{anandkumar2014tensor,chen2019phase,liu2017characterizing,perry2020statistical,richard2014statistical,zhang2018tensor}, which assumes 
\begin{equation}\label{eq:PCA_model}
\calA=\calT+\calZ,
\end{equation}
where the signal $\calT$ admits a low-rank decomposition (\ref{eq:tucker}) and the noise $\calZ$ contains i.i.d. entries with mean zero and variance $\sigma^2$. A central goal of tensor PCA is on the estimation and inference of $\calT$ and/or $\{U_j\}_j$, i.e. the low-rank structure from $\calA$. Tensor PCA has been proven effective for learning hidden components in Gaussian mixture models \cite{anandkumar2014tensor}, where $\{U_j\}_j$ represent the hidden components. By constructing confidence regions of $\{U_j\}_j$, we are able to make uncertainty quantification for the hidden components of Gaussian mixture models. In addition, confidence regions of $\{U_j\}_j$ can be useful for the inference of spatial and temporal patterns of gene regulation during brain development \cite{liu2017characterizing}. When applying the tensor PCA model to community detection in hypergraph networks \citep{ke2019community} or multilayer networks \citep{jing2020community}, $U_j$ is directly related to the estimated community structures and the confidence region of $U_j$ is an important tool to quantify the uncertainty of community detection. This also applies to the uncertainty quantification for tensor/high-order clustering \citep{han2020exact,luo2020tensor}.

{\it Low-rank tensor regression} can be seen as one of the most basic setting of {\it supervised inference} for low-rank tensors. Specifically, suppose we observe a set of random pairs $\{\calX_i, Y_i\}_{i=1}^n$ associated as
\begin{equation}\label{eq:tr_model}
Y_i=\langle \calT,\calX_i\rangle+\xi_i.
\end{equation}
Here, the main point of interest is $\calT$, a low-rank tensor that characterizes the association between response $Y$ and covariate $\calX$, and $\xi_i$ is the noise term. 
When the tensor order is $m=2$, this problem is reduced to the widely studied {\it trace matrix regression model} in the literature \citep{candes2010matrix,cai2013sparse,chen2019non,fan2019generalized,koltchinskii2011nuclear,koltchinskii2015optimal,raskutti2019convex,rauhut2017low,tomioka2013convex}. This model can also be used as the prototype of many problems in high-dimensional statistics and machine learning, including phase retrieval \citep{candes2013phaselift} and blind deconvolution \citep{li2019rapid}. When $m\geq 3$, this problem has been studied under the scenario of high-order interaction pursuit \citep{hao2020sparse} and large-scale linear system from partial differential equations \citep{lynch1964tensor}. In applications of tensor regression to neuroimaging analysis, the principal components of $\calT$ are useful in the understanding of the association between disease outcomes and brain image patterns \cite{zhou2013tensor}. In addition, the principal components determine the cluster memberships of neuroimaging data \citep{sun2019dynamic}. Confidence regions of $\{U_j\}_j$ in the aforementioned applications allow us to make significance test for the detected regions of interest, and to make uncertainty quantification for clustering outcomes, respectively. 

In addition to tensor PCA and regression, there is a broad range of low-rank tensor models, such as tensor completion \citep{montanari2018spectral,xia2019polynomial,yuan2016tensor,zhang2019cross}, generalized tensor estimation \citep{han2020optimal}, and tensor high-order clustering \citep{chi2018provable,feizi2017tensor,han2020exact,luo2020tensor,sun2019dynamic,wu2016general}. A common goal of these problems is to accurately estimate and make inference on some type of low-rank structures.

%%%%%%%%%%%%%%%%%%%%%%%%%%
\subsection{Summary of the Main Results}\label{sec:contribution}
%%%%%%%%%%%%%%%%%%%%%%%%%%

In this paper, we aim to develop the methods and non-asymptotic theory for statistical inference under the low-rank tensor PCA and regression models. First, suppose the target tensor $\calT$ is Tucker low-rank with singular subspace $U_j$ as the point of interest. Given any estimator $\hat{U}^{(0)}_j$ that achieves some reasonable estimation error, we introduce a straightforward two-iteration alternating minimization scheme (Algorithms~\ref{algo:am_optimal_PCA} and \ref{algo:am_optimal_regression} in Section~\ref{sec:est}) and obtain $\hat{U}_j$. Surprisingly, we are able to derive an asymptotic distribution of $\|\sin\Theta(\hat U_j , U_j)\|_{\rm F}^2$ (definition of sin-theta distance is postponed to Section \ref{sec:notation}) even though $\hat U_j$ is from non-convex iterations. Under the tensor PCA model with some essential conditions on SNR, we prove that
\begin{equation}\label{eq:asymptotic-PCA-hatU_j}
\frac{\|\sin\Theta(\hat U_j, U_j)\|_{\rm F}^2-p_j\sigma^2\|\Lambda_j^{-1}\|_{\rm F}^2}{\sqrt{2p_j}\sigma^2\|\Lambda_j^{-2}\|_{\rm F}}\stackrel{{\rm d.}}{\longrightarrow} N(0,1)\quad {\rm as}\quad p_j\to\infty.
\end{equation}
Here, $\Lambda_j$ is the diagonal matrix containing all non-zero singular values of the $j$th matricization of $\calG$ (see definition of matricization in Section \ref{sec:notation}). Under the tensor regression model with some essential conditions on sample size and SNR, we prove that
\begin{equation}\label{eq:asymptotic-regression-hatU_j}
\frac{\|\sin\Theta(\hat U_j, U_j)\|_{\rm F}^2-p_jn^{-1}\sigma^2\|\Lambda_j^{-1}\|_{\rm F}^2}{\sqrt{2p_j}n^{-1}\sigma^2\|\Lambda_j^{-2}\|_{\rm F}}\stackrel{{\rm d.}}{\longrightarrow} N(0,1)\quad {\rm as}\quad p_j\to\infty.
\end{equation}
We also develop the non-asymptotic Berry-Essen-type bounds for the limiting distributions in \eqref{eq:asymptotic-PCA-hatU_j} and \eqref{eq:asymptotic-regression-hatU_j}.

Then, we consider a special class of \emph{orthogonally decomposable tensors} $\calT$ in the sense that $\calT=\sum_{j=1}^r\lambda_{j} \cdot u_j\otimes v_j\otimes w_j \in \mathbb{R}^{p_1\times p_2\times p_3}$ for orthonormal vectors $\{u_j\}_j$, $\{v_j\}_j$, and $\{w_j\}_j$. The orthogonally decomposable tensor has been widely studied as a benchmark setting for tensor decomposition in the literature \citep{auddy2020perturbation,belkin2018eigenvectors,chen2009tensor,kolda2001orthogonal,robeva2016orthogonal}. In addition, the (near-)orthogonally decomposable tensors have been used in various applications of statistics and machine learning, such as latent variable model \citep{anandkumar2014tensor}, hidden Markov models \citep{anandkumar2012method}, etc. Under the tensor PCA model, we prove that
\begin{equation}\label{eq:uuj_tsvd}
\frac{\langle \hat u_j, u_j\rangle^2-(1-p_j\sigma^2\lambda_{j}^{-2})}{\sqrt{2p_j}\sigma^2\lambda_{j}^{-2}}\stackrel{{\rm d.}}{\longrightarrow} N(0,1)\quad {\rm as}\quad p_1\to\infty
\end{equation}
for $j=1,\cdots,r$ when some essential SNR condition holds. Here, $\{\hat u_j, \hat v_j, \hat w_j\}_{j}$ are the estimates of $\{u_j, v_j, w_j\}_{j}$ (up to some permutation of index $j$) based on a two-step power iteration  (Algorithm \ref{algo:power_iter}).  Similar results can also be obtained for $\langle \hat v_j, v_j\rangle^2$ and $\langle \hat w_j, w_j\rangle^2$. 

Next, we propose the estimates of $\Lambda_j, \lambda_{j}, \sigma^2$ that are involved in the asymptotic distributions of $\|\sin\Theta(\hat{U}_j, U_j)\|_\F^2$ in \eqref{eq:asymptotic-PCA-hatU_j}\eqref{eq:asymptotic-regression-hatU_j} and $\langle \hat u_j, u_j\rangle^2$ in \eqref{eq:uuj_tsvd}. We prove that the asymptotic normality in \eqref{eq:asymptotic-PCA-hatU_j}\eqref{eq:asymptotic-regression-hatU_j}\eqref{eq:uuj_tsvd} still hold after plugging in these estimates. These results immediately yield the data-driven confidence regions for $U_j$ (Tucker low-rank settings) or $\{u_j\}_j$ (orthogonally decomposable settings).

If $\calA$ is a rank-1 tensor, the low-rank tensor PCA model reduces to the widely studied \emph{rank-1 tensor PCA} (see a literature survey in Section \ref{sec:related-prior-work}). Under this model, we establish the asymptotic normality of any linear functionals for the power iteration estimators $\hat u, \hat v, \hat w$: for all unit vectors $q_i \in \mathbb{R}^{p_i}$, under regularity conditions, we have
$$\bigg(\frac{\langle q_1, \hat u - u\rangle + \frac{p_1\langle q_1, u\rangle}{2(\lambda/\sigma)^2}}{\sqrt{\frac{p_1\langle q_1, u\rangle^2}{2(\lambda/\sigma)^4} + \frac{1 - \langle q_1, u\rangle^2}{(\lambda/\sigma)^2}}}, \frac{\langle q_2, \hat v - v\rangle + \frac{p_2\langle q_2, v\rangle}{2(\lambda/\sigma)^2}}{\sqrt{\frac{p_2\langle q_2, v\rangle^2}{2(\lambda/\sigma)^4} + \frac{1 - \langle q_2, v\rangle^2}{(\lambda/\sigma)^2}}}, \frac{\langle q_3, \hat w - w\rangle + \frac{p_3\langle q_3, w\rangle}{2(\lambda/\sigma)^2}}{\sqrt{\frac{p_3\langle q_3, w\rangle^2}{2(\lambda/\sigma)^4} + \frac{1 - \langle q_3, w\rangle^2}{(\lambda/\sigma)^2}}}\bigg)^\top \stackrel{{\rm d.}}{\to} N(0, I_3)$$ as $p_1, p_2, p_3 \to \infty$. 
We further derive the entrywise asymptotic distribution for each entry of the estimator $\hat{\calT}$, and propose a thresholding procedure to construct the asymptotic $1-\alpha$ entrywise confidence interval for $\calT$, which is the first of such work to our best knowledge. 

Our theoretical results reveal a key message: under the tensor PCA and regression model, {\it the inference of principal components can be efficiently done when a computationally feasible optimal estimate is achievable.} In recent literature, it is widely observed in many low-rank tensor models (See \ref{sec:related-prior-work} for a review of literature) that in order to achieve an accurate estimation in polynomial time, one often requires a more stringent condition than what is needed in the statistical (or information-theoretic) limit. Such a statistical and computational gap becomes a ``blessing" to the statistical inference of low-rank tensor models, as debiasing can become unnecessary if those strong but essential conditions for computational feasibility are met! 

%%%%%%%%%%%%%%%%%%%%%%
\subsection{Related Prior Work}\label{sec:related-prior-work}
%%%%%%%%%%%%%%%%%%%%%%

This paper is related to a broad range of literature in high-dimensional statistics and matrix/tensor analysis. First, a variety of methods have been proposed for tensor PCA in the literature. A non-exhaustive list include high-order orthogonal iteration \citep{de2000best}; sequential-HOSVD \citep{vannieuwenhoven2012new}, inference for low-rank matrix completion \citep{chen2019inference,foucart2017biasing}, (truncated) power iteration \citep{anandkumar2014guaranteed,liu2017characterizing,sun2017provable}, STAT-SVD \citep{zhang2019optimal}. In addition, the computational hardness was widely considered for tensor PCA. Particularly in the worse case scenario, the best low-rank approximation of tensors can be NP hard \citep{de2008tensor,hillar2013most}. The average-case computational complexity for tensor PCA model has also been widely studied under various computational models, including the Sum-of-Squares \citep{hopkins2015tensor}, optimization landscape \citep{arous2019landscape}, average-case reduction \citep{brennan2020reducibility,luo2020tensor,luo2020open,zhang2018tensor}, and statistical query \citep{dudeja2021statistical}. It has now been widely justified that the SNR condition $\lambda_{\min}/\sigma \geq Cp^{3/4}$ is essential to ensure tensor PCA is solvable in polynomial time.

Regression of low-rank tensor has attracted enormous attention recently. Various methods, such as the (regularized) alternating minimization \citep{li2018tucker,sun2017store,zhou2013tensor}, convex regularization \citep{raskutti2019convex,tomioka2013convex}, projected gradient descent \citep{chen2019non,rauhut2017low}, importance sketching \citep{zhang2020islet} were studied. Recently, \cite{han2020optimal} proved that a gradient descent algorithm can recover a low-rank third order tensor $\calT$ with statistically optimal convergence rate when the sample size $n$ is much greater than the tensor dimension $p^{3/2}$. It was widely conjectured that $n\geq Cp^{3/2}$ is essential for the problem being solvable in polynomial time (see \cite{barak2016noisy} for the evidence). 

While the statistical inference for low-rank tensor models remain largely unexplored, there have been several recent results demystifying the statistical inference for low-rank {\it matrix} models. For matrix PCA, \cite{xia2019normal} introduced an explicit representation formula for $\hat U_j\hat U_j^{\top}$. A more precise characterization of the distribution of $\|\sin\Theta(\hat U_j, U_j)\|_{\rm F}^2$ was established in \cite{bao2018singular} by random matrix theory. On the other hand, the estimators of tensor PCA are often calculated from iterative optimization algorithms (e.g., power iterations or gradient descent) in existing literature, while the estimator of matrix PCA is based on non-iterative schemes. Due to the complex statistical dependence involved in iterative optimization algorithms, it is significantly more challenging to analyze the asymptotic distribution of the estimator in tensor PCA than the one in matrix PCA. We also note that, when studying the asymptotic distributions of individual eigenvectors, an eigengap condition is often crucial for matrix PCA but not required for tensor PCA.

The inference and uncertainty quantification were also considered for low-rank matrix regression. For example, \cite{carpentier2019uncertainty} introduced a debiased estimator based on the nuclear norm penalized low-rank estimator. \cite{cai2016geometric} introduced another debiasing technique and characterize the entrywise distribution of the debiased estimator under the restricted isometry property. \cite{xia2019confidence} studied a debiased estimator for matrix regression under the isotropic Gaussian design and established the distribution of $\|\sin\Theta(\hat U_j, U_j)\|_{\rm F}^2$ under nearly optimal sample size conditions. All these approaches rely on suitable debiasing of certain initial estimates. In addition to low-rank estimation, an appropriate debias was found crucial for high-dimensional sparse regression \citep{zhang2008sparsity}, and various debiasing schemes were introduced \citep{javanmard2014confidence,van2014asymptotically,zhang2014confidence}. Interestingly, as will be shown in Section \ref{sec:normal_app}, our estimating and inference procedure for low-rank tensor regression does not involve debiasing.

Statistical inference for low-rank models are particularly challenging for tensor problems. In a concurrent work, \cite{huang2020power} studied the statistical inference and power iteration for tensor PCA. Recently, \cite{cai2020uncertainty} studied the entrywise statistical inference for noisy low-rank tensor completion based on a nearly unbiased estimator and an incoherence condition on $U_j$s, i.e., all the rows of $U_j$ have comparable magnitudes. In comparison, our results do not require further conditions on $U_j$s or debiasing. 

%%%%%%%%%%%%%%%%%%%%%%%%
\subsection{Organizations}\label{sec:organization}
%%%%%%%%%%%%%%%%%%%%%%%%

The rest of the paper is organized as follows. After an introduction on notation and preliminaries in Section \ref{sec:notation}, we discuss the inference for principal components under the Tucker low-rank models in Section \ref{sec:normal_app}. Specifically, a general two-iteration alternating minimization procedure, inference for tensor PCA, inference for tensor regression, and a proof sketch are given in Sections \ref{sec:est}, \ref{sec:inference-PCA}, \ref{asymptotic-distribution-regression}, and \ref{sec:proof-sketch}, respectively. In Section~\ref{sec:T-orth}, we focus on the inference for individual singular vectors of orthogonally decomposable tensors. The asymptotic distribution and entrywise confidence interval are discussed for rank-1 tensor PCA model in Section \ref{sec:entrywise-rank1}. In the supplementary materials, Section \ref{sec:additional-procedure} includes some algorithms for tensor PCA and regression in the literature. All proofs of the main technical results are collected in Section \ref{sec:proof}. 

%%%%%%%%%%%%%%%%%%%%%%%
\section{Notation and Preliminaries}\label{sec:notation}
%%%%%%%%%%%%%%%%%%%%%%%

Let $\{a_k\}$ and $\{b_k\}$ be two sequences of non-negative numbers. We denote $a_k\ll b_k$ if $\lim_{k\to \infty} a_k/b_k=0$ and $a_k\gg b_k$ if $\lim_{k\to \infty} a_k/b_k=\infty$. We use calligraphic letters $\calT,\calG$ to denote tensors, upper-case letters $U, W$ to denote matrices, and lower-case letters $u,w$ to denote vectors or scalars. For a random variable $X$ and $\alpha > 0$, the Orlicz $\psi_{\alpha}$-norm of $X$ is defined as
\begin{equation*}
\left\|X\right\|_{\psi_{\alpha}} = \inf\{K > 0: \EE\{\exp(|X|/K)^{\alpha}\} \leq 2\}.
\end{equation*}
Specifically, a random variable with finite $\psi_2$-norm or $\psi_1$-norm is called the sub-Gaussian or sub-exponential random variable, respectively. Let $e_j$ denote the $j$th canonical basis vector whose dimension varies at different places. Let $\rank(\calT)$ be the Tucker rank of $\calT$ and write $(a_1, \ldots, a_m)\leq (b_1, \ldots, b_m)$ if $a_j\leq b_j$ for all $j\in[m]$. We use $\|\cdot\|_{\rm F}$ for Frobenius norm, $\|\cdot\|$ for matrix spectral norm and $\|\cdot\|_q$ for vector $\ell_q$-norm. Denote $\mathbb{S}^{p-1} = \{v\in\mathbb{R}^{p}: \|v\|_2\leq 1\}$ as the set of $p$-dimensional unit vectors. Define $\mathbb{O}_{p, r} = \{U\in \mathbb{R}^{p\times r}: U^\top U = I_r\}$ as the set of all $p$-by-$r$ matrices with orthonormal columns. In particular, $\OO_r$ is the set of all $r\times r$ orthogonal matrices. 

We denote $\times_j$ the $j$th multi-linear product between a tensor and matrix. For instance, if $\calG\in\RR^{r_1\times r_2\times r_3}$ and $V_1\in\RR^{p_1\times r_1}$, then  
$$
\calG\times_1 V_1=\Big(\sum_{j_1=1}^{r_1}\calG(j_1,i_2,i_3)V(i_1,j_1)\Big)_{i_1\in[p_1],i_2\in[r_2],i_3\in[r_3]}.
$$
We write $(U_1,\cdots,U_m)\cdot\calG$ in short for $\calG\times_1 U_1\times_2\cdots\times_m U_m$. Let $\calM_j$ be the $j$th tensor matricization that rearranges each mode-$j$ fiber of $\calT\in \mathbb{R}^{p_1\times \cdots \times p_d}$ to a column of $\calM_j(\calT)\in \mathbb{R}^{p_j\times(p_1\cdots p_d/p_j)}$. 

We say $\calT$ has Tucker rank $(r_1,\cdots,r_m)$ if it admits a Tucker decomposition 
\begin{equation}
\label{eq:tucker}
\calT=(U_1,\cdots, U_m)\cdot \calG,
\end{equation}
where $\calG \in \RR^{r_1 \times \cdots \times r_m}$ and $U_i \in \OO_{p_i,r_i}$ for $i \in [m]$. The Tucker decomposition (\ref{eq:tucker}) can be roughly seen as a generalization of matrix singular value decomposition (SVD) to higher-order tensors, where $U_j$ can be viewed as principal components of the $j$th matricization of $\calT$, and $\calG$ contains the singular values. In the case that $r_1=\cdots=r_m= r$ and $\calG$ is diagonalizable, we say $\calT$ is \emph{orthogonally decomposable}. %See formal definition of Tucker rank, Tucker decomposition, and orthogonally decomposable tensors in Section \ref{sec:notation}. 
If $\calT$ satisfies Tucker decomposition (\ref{eq:tucker}), one has
$$
\calM_j(\calT)=U_j\calM_j(\calG)\big(U_1\otimes\cdots\otimes U_{j-1}\otimes U_{j+1}\otimes\cdots\otimes U_{m}\big)^{\top}\in\RR^{p_j\times (p_1\cdots p_m/p_j)}.
$$
Here $\otimes$ stands for Kronecker product so that $U\otimes W\in\RR^{(p_1p_2)\times(r_1r_2)}$ if $U\in\RR^{p_1\times r_1}$ and $W\in\RR^{p_2\times r_2}$. The readers are referred to \cite{kolda2009tensor} for a comprehensive survey on tensor algebra. 

Let $\sigma_r(\cdot)$ be the $r$th largest singular value of a matrix. If $\calT$ has Tucker ranks $(r_1,\cdots,r_m)$, the signal strength of $\calT$ is defined by 
$$
\lambda_{\submin}:=\lambda_{\submin}(\calT)=\min\big\{\sigma_{r_1}\big(\calM_1(\calT)\big), \sigma_{r_2}\big(\calM_2(\calT)\big),\cdots, \sigma_{r_m}\big(\calM_m(\calT)\big)\big\},
$$
i.e., the smallest positive singular value of all matricizations. Similarly, define $\lambda_{\submax}:=\lambda_{\submax}(\calT)=\max_j \sigma_1\big(\calM_j(\calT)\big)$. The condition number of $\calT$ is defined by $\kappa(\calT):=\lambda_{\submax}(\calT)\lambda_{\submin}^{-1}(\calT)$. We let $\Lambda_j$ be the $r_j\times r_j$ diagonal matrix containing the singular values of $\calM_j(\calG)$ (or equivalently the singular values of $\calM_j(\calT)$). Note that $\Lambda_j$s are not necessarily equal for different $j$, although $\|\Lambda_1\|_{\rm F}=\cdots=\|\Lambda_m\|_{\rm F} = \|\calT\|_{\rm F}$. 

We define the principle angles between $U, \widehat{U} \in \mathbb{O}_{p ,r}$ as an $r$-by-$r$ diagonal matrix: $\Theta(U, \widehat{U}) = \diag(\text{arccos}(\sigma_1), \dots, \text{arccos}(\sigma_r))$, where $\sigma_1 \geq \cdots \geq \sigma_r \geq 0$ are the singular values of $U^\top \widehat{U}$. Then the $\sin\Theta$ distances between $\hat{U}$ and $U$ are defined as
\begin{equation*}
\|\sin\Theta(U, \widehat{U})\| = \|\diag\left(\sin(\text{arccos}(\sigma_1)), \dots, \sin(\text{arccos}(\sigma_r))\right)\| = \sqrt{1 - \sigma_r^2},
\end{equation*}
\begin{equation*}
\|\sin\Theta(U, \widehat{U})\|_{\F} = \left(\sum_{i=1}^{r}\sin^2(\text{arccos}(\sigma_i))\right)^{1/2} = \left(r - \sum_{i=1}^{r}\sigma_i^2\right)^{1/2}.
\end{equation*}

%%%%%%%%%%%%%%%%%%%%%
\section{Inference for Principal Components of Tucker Low-rank Tensor}\label{sec:normal_app}
%%%%%%%%%%%%%%%%%%%%%

For notational simplicity, we focus on the inference for third-order tensors, i.e., $m=3$, while the results for general $m$th order tensor essentially follows and will be briefly discussed in Section \ref{sec:diss}. 

%%%%%%%%%%%%%%%%%%%%%
\subsection{Estimating Procedure}\label{sec:est}
%%%%%%%%%%%%%%%%%%%%%

An accurate estimation is often the starting point for statistical inference and uncertainty quantification. In this section, we briefly discuss the estimation procedure for both tensor regression and PCA models. First, we summarize both models as follows:
$$Y_i = \langle \calX_i, \calT\rangle + \xi_i, \quad i=1,\ldots, n.$$
Here, $\calX_i$ can be the covariate in tensor regression; $n=p_1p_2p_3$, $Y_i=\calA(j_1,j_2,j_3)$, and $\calX_i=(e_{j_1},e_{j_2},e_{j_3})\cdot 1$ with $i=(j_1-1)p_2p_3+(j_2-1)p_3+j_3$, $j_1\in[p_1], j_2\in[p_2], j_3\in[p_3]$ in tensor PCA. Let $l_n(\calT) = \sum_{i=1}^n(Y_i - \langle \calX_i, \calT\rangle)^2$ be the loss function in both settings. Then a straightforward solution to both problems is via the following Tucker rank constrained least squares estimator:
\begin{equation}\label{eq:MLE}
\begin{split}
\underset{\rank(\calT)\leq (r_1,r_2,r_3)}{\min}\ \ell_n(\calT) & :=\frac{1}{n}\sum_{i=1}^n\big(Y_i-\langle \calX_i, \calT\rangle\big)^2,\\
\text{or equivalently} \quad (\hat \calG, \hat U_1, \hat U_2, \hat U_3) & :=\underset{\calG\in\RR^{r_1\times r_2\times r_3},\ U_j\in\OO_{p_j,r_j}}{\arg\min}\ \ell_n\big((U_1,U_2,U_3)\cdot \calG\big).
\end{split}
\end{equation}
Since the objective function (\ref{eq:MLE}) is highly non-convex, an efficient algorithm with provable guarantees is crucial for both tensor PCA and regression. As discussed earlier, various computationally feasible procedures have been proposed in the literature. For tensor regression, \cite{han2020optimal} recently introduced a simultaneous gradient descent algorithm and proved their proposed procedure achieves the minimax optimal estimation error; for tensor PCA, a simpler and more direct approach, higher-order orthogonal iteration (HOOI), was introduced by \cite{de2000best}. The implementation details of both algorithms are provided in Section \ref{sec:additional-procedure} in the supplementary materials. 

Moreover, the primary interest of this paper is on the statistical inference for $\calT$ or $U_j$, far beyond deriving estimators achieving optimal estimation error.  In general, even estimators achieving minimax optimal estimation error rate may not enjoy a proper asymptotic distribution. For example, the true parameter $\calT$ or $U_j$ plus a small enough perturbation can achieve optimal estimation error but does not satisfy any tractable distribution.

To this end, we introduce a \emph{two-iteration alternating minimization} algorithm for both Tucker low-rank tensor PCA and tensor regression in Algorithms~\ref{algo:am_optimal_PCA} and \ref{algo:am_optimal_regression}, respectively. 
Our theory in later this section reveals a surprising fact: if any estimator $\tilde \calT=(\hat U_1^{(0)},\hat U_2^{(0)}, \hat U_3^{(0)})\cdot \hat\calG^{(0)}$ achieving some attainable estimation error is provided as the input, the two-iteration alternating minimization in Algorithms \ref{algo:am_optimal_PCA} and \ref{algo:am_optimal_regression} will provide an estimator enjoying asymptotic normality and being ready to use for confidence region construction. 
\begin{algorithm}
	\SetAlgoLined
	\KwInput{$\ell_n(\cdot)$: Objective function  (\ref{eq:MLE}); Initializations $(\hat U_1^{(0)}, \hat U_2^{(0)}, \hat U_3^{(0)})$;}
	\For{$t=0,1$}{
		$\hat{U}_1^{(t+1)} = \text{leading $r_1$ left singular vectors of } \mathcal{M}_1(\calA\times_2 \hat U_2^{(t)\top} \times_3 \hat U_3^{(t)\top});$\\
		$\hat{U}_2^{(t+1)} = \text{leading $r_2$ left singular vectors of } \mathcal{M}_2(\calA\times_1 \hat U_1^{(t)\top} \times_3 \hat U_3^{(t)\top});$\\
		$\hat{U}_3^{(t+1)} = \text{leading $r_3$ left singular vectors of } \mathcal{M}_3(\calA\times_1 \hat U_1^{(t)\top} \times_2 \hat U_2^{(t)\top});$\\
	}
	\KwOutput{Test statistic $\hat U_1:=\hat U_1^{(2)}, \hat U_2:=\hat U_2^{(2)}, \hat U_3:=\hat U_3^{(2)}$, and $\hat\calG=(\hat U_1^{(2)\top}, \hat U_2^{(2)\top}, \hat U_3^{(2)\top})\cdot \calA$. 
	}
	\caption{Power Iteration for Tensor PCA}
	\label{algo:am_optimal_PCA}
\end{algorithm}

\begin{algorithm}
	\SetAlgoLined
	\KwInput{$\ell_n(\cdot)$: Objective function (\ref{eq:MLE}); Initializations $(\hat U_1^{(0)}, \hat U_2^{(0)}, \hat U_3^{(0)})$, and $\hat\calG^{(0)}$ is the solution of $\argmin_{\calG}\ell_n\big((\hat U_1^{(0)},\hat U_2^{(0)},\hat U_3^{(0)})\cdot \calG\big)$ for tensor regression model\;}
	\For{$t = 0, 1$}{
		Solve $\nabla_{U_1}\ell_n\big((\hat{U}_1^{(t+0.5)}, \hat{U}_2^{(t)}, \hat{U}_3^{(t)})\cdot\hat\calG^{(t)}\big) = 0$ to obtain $\hat U_1^{(t+0.5)}$\;
		Update by $\hat U_1^{(t+1)}={\rm SVD}_{r_1}\big(\hat{U}_1^{(t+0.5)}\big)$\;
		Solve $\nabla_{U_2}\ell_n\big((\hat{U}_1^{(t)}, \hat{U}_2^{(t+0.5)}, \hat{U}_3^{(t)})\cdot\hat\calG^{(t)}\big) = 0$ to obtain $\hat U_2^{(t+0.5)}$\;
		Update by $\hat U_2^{(t+1)}={\rm SVD}_{r_2}\big(\hat{U}_2^{(t+0.5)}\big)$\;
		Solve $\nabla_{U_3}\ell_n\big((\hat{U}_1^{(t)}, \hat{U}_2^{(t)}, \hat{U}_3^{(t+0.5)})\cdot\hat\calG^{(t)}\big) = 0$ to obtain $\hat U_3^{(t+0.5)}$\;
		Update by $\hat U_3^{(t+1)}={\rm SVD}_{r_3}\big(\hat{U}_3^{(t+0.5)}\big)$\;
		Solve $\nabla_{\calG}\ell_n\big((\hat{U}_1^{(t+1)}, \hat{U}_2^{(t+1)}, \hat{U}_3^{(t+1)})\cdot \hat\calG^{(t+1)}\big) = 0$ to obtain $\hat\calG^{(t+1)}$\;
	}
	\KwOutput{Test statistic $\hat U_1:=\hat U_1^{(2)}, \hat U_2:=\hat U_2^{(2)}, \hat U_3:=\hat U_3^{(2)}$, and $\hat\calG:=\hat\calG^{(2)}$.}
	\caption{Alternating Minimization for Tensor Regression}
	\label{algo:am_optimal_regression}
\end{algorithm}

\begin{remark}[Interpretation of Alternating Minimization Update in Tensor PCA]
	A key observation by \citep[Theorems 4.1, 4.2]{de2000best} shows minimizing $\min_{\rank(\calT)\leq (r_1,r_2,r_3)}\|\calT-\calA\|_{\rm F}^2$ is equivalent to maximizing $\max_{U_j\in\OO_{p_j,r_j}}\ \|(U_1^{\top},U_2^{\top},U_3^{\top})\cdot \calA\|_{\rm F}^2$.  Therefore, the optimization in tensor PCA is equivalent to
	\begin{equation*}
	\begin{split}
	(\hat U_1, \hat U_2, \hat U_3):= & \underset{U_j\in\OO_{p_j,r_j}}{\arg\min}\ \ell_n((U_1, U_2, U_3)\cdot \calG) :=  \underset{U_j\in\OO_{p_j,r_j}}{\arg\max}\ \|(U_1^{\top},U_2^{\top},U_3^{\top})\cdot \calA\|_{\rm F}^2\\
	= & \underset{U_j\in \mathbb{O}_{p_j, r_j}}{\arg\max} \left\|U_j \mathcal{M}_j(\calA\times_{j+1} U_{j+1} \times_{j+2} U_{j+2})\right\|_\F^2. 
	\end{split}
	\end{equation*}
	Here, for convenience of notation, $U_4=U_1, U_5=U_2, r_4=r_1, r_5=r_2$. Note that, given fixed $\hat U_{j+1}^{(t)}$ and $\hat U_{j+2}^{(t)}$, Eckart-Young-Mirsky Theorem \citep{eckart1936approximation} implies the optimal solution to $\max_{U_j\in\OO_{p_j,r_j}}\ \|(U_j^{\top},\hat U_{j+1}^{(t)\top}, \hat U_{j+2}^{(t)\top})\cdot \calA\|_{\rm F}^2$ is attainable via singular value decomposition: 
	$$\hat{U}_j^{(t+1)} = \text{leading $r_j$ left singular vectors of } \mathcal{M}_j\left(\calA\times_{j+1} \hat U_{j+1}^{(t)\top} \times_{j+2} \hat U_{j+2}^{(t)\top}\right).$$
	This explains the alternating minimization update steps for tensor PCA in Algorithm \ref{algo:am_optimal_PCA}.
\end{remark}

Hereinafter, we denote $\hat U_j$ the output of Algorithms~\ref{algo:am_optimal_PCA} and \ref{algo:am_optimal_regression}, $p= \max\{p_1,p_2,p_3\}$ and $r_{\submax}=\max\{r_1,r_2,r_3\}$. Next, we establish the asymptotic distribution and develop the inference procedure for $\|\sin\Theta(\hat U_j, U_j)\|_\F^2$ in tensor PCA and tensor regression models when $\calT$ admits the Tucker decomposition (\ref{eq:tucker}). 

%%%%%%%%%%%%%%%%%%%%%%%%%%%%%%
\subsection{Inference for Tucker Low-rank Tensor PCA}\label{sec:inference-PCA}
%%%%%%%%%%%%%%%%%%%%%%%%%%%%%%

We assume the following condition on initialization $(\hat U_1^{(0)},\hat U_2^{(0)},\hat U_3^{(0)})$ of Algorithm~\ref{algo:am_optimal_PCA} holds.
\begin{assump}\label{assump:tpca} 
	Under tensor PCA model (\ref{eq:PCA_model}) with $\calZ_{i_1,i_2,i_3} \stackrel{i.i.d.}{\sim} N(0,\sigma^2)$, there is an event $\calE_0$ with $\PP(\calE_0)\geq 1-C_1e^{-c_1p}$ for some absolute constants $c_1,C_1>0$ so that, under $\calE_0$, the initialization $(\hat U_1^{(0)},\hat U_2^{(0)},\hat U_3^{(0)})$ satisfy $\max_{j=1,2,3}\|\sin\Theta(\hat U_j^{(0)}, U_j)\| \leq C_2\sqrt{p}\sigma/\lambda_{\submin}$ for some absolute constant $C_2>0$.
\end{assump}
The claimed error rates in Assumption~\ref{assump:tpca} are attainable by the algorithm HOOI under the SNR condition $\lambda_{\submin}/\sigma\geq Cp^{3/4}$ \citep[Theorem 1]{zhang2018tensor}. 
Such the SNR condition is essential to ensure a consistent estimator is achievable in polynomial time as illustrated by the literature reviewed in Section \ref{sec:related-prior-work}. Note that \cite[Theorem 1]{zhang2018tensor} presented an expectation error bound $\EE\|\sin\Theta(\hat U_j^{(0)}, U_j)\|$, while its proof indeed involved a desired probabilistic bound as claimed by Assumption~\ref{assump:tpca}. If a given initialization estimation error upper bound is in a metric other than the $\sin\Theta$ distance described in Assumption \ref{assump:tpca}, we may apply Lemma \ref{lm:equivalence-different-error} in the supplementary materials to ``translate" the upper bound in another metric to the desired $\sin\Theta$ distance.

Suppose $\hat{U}_j$ is the output of Algorithm \ref{algo:am_optimal_PCA}. Built on Assumption~\ref{assump:tpca}, we characterize the distribution of $\|\sin\Theta(\hat U_j, U_j)\|_{\rm F}^2$ by the following theorem. 
\begin{theorem}[Asymptotic normality of principal components in tensor PCA]\label{thm:na_tsvd}
	Suppose Assumption~\ref{assump:tpca} holds for tensor PCA model (\ref{eq:PCA_model}), $\calZ(i_1,i_2,i_3) \stackrel{i.i.d.}{\sim} N(0,\sigma^2)$, $p_j\asymp p$ for $j=1,2,3$, and $\kappa(\calT)\leq \kappa_0$. Let $\hat U_j$s be the output of Algorithm~\ref{algo:am_optimal_PCA} for tensor PCA model. There exist absolute constants $c_1,C_0,C_1,C_2,C_3>0$ such that if $\lambda_{\submin}/\sigma \geq C_0(p^{3/4}+\kappa_0^2p^{1/2})$, then
	\begin{equation*}
	\begin{split}
	\sup_{x\in\RR}&\left|\PP\left(\frac{\|\sin\Theta(\hat U_j, U_j)\|_{\rm F}^2 - p_j\sigma^2\|\Lambda_j^{-1}\|_{\rm F}^2}{\sqrt{2p_j }\sigma^2\|\Lambda_j^{-2}\|_{\rm F}} \leq x\right) - 
	\Phi(x)\right|\\
	&\hspace{2cm} \leq C_1e^{-c_1p}+C_2\left(\frac{\kappa_0^{6}(pr_{\submax})^{3/2}}{(\lambda_{\submin}/\sigma)^2} + \frac{\kappa_0^2(p\log p)^{1/2}}{\lambda_{\submin}/\sigma}\right)  + C_3\frac{r_{\submax}^{3/2}}{\sqrt{p}},
	\end{split}
	\end{equation*}
	where $\Lambda_j = \diag(\lambda_1^{(j)}, \dots, \lambda_{r_j}^{(j)})$ is the diagonal matrix containing the singular values of $\calM_j(\calG)$, and $\Phi(x)$ is the cumulative distribution function of $N(0, 1)$. 
\end{theorem}
If the condition number $\kappa_0=O(1)$, $(pr_{\submax})^{3/4}(\lambda_{\submin}/\sigma)^{-1}\to 0$ and $r_{\submax}^3/p\to 0$ as $p\to\infty$, Theorem~\ref{thm:na_tsvd} yields
\begin{equation*}
\frac{\|\sin\Theta(\hat U_j, U_j)\|_{\rm F}^2 - p_j\sigma^2\|\Lambda_j^{-1}\|_{\rm F}^2}{\sqrt{2p_j}\sigma^2\|\Lambda_j^{-2}\|_{\rm F}} \stackrel{\rm d.}{\longrightarrow} N(0, 1)\quad {\rm as}\quad p\to\infty. 
\end{equation*}
By the proof of Theorem \ref{thm:na_tsvd}, we can further establish the following joint distribution of all $U_j$s:
\begin{align*}
\bigg( & \frac{\|\sin\Theta(\hat U_1, U_1)\|_{\rm F}^2 - p_1\sigma^2\|\Lambda_1^{-1}\|_{\rm F}^2}{\sqrt{2p_1}\sigma^2\|\Lambda_1^{-2}\|_{\rm F}}, \frac{\|\sin\Theta(\hat U_2, U_2)\|_{\rm F}^2 - p_2\sigma^2\|\Lambda_2^{-1}\|_{\rm F}^2}{\sqrt{2p_2}\sigma^2\|\Lambda_2^{-2}\|_{\rm F}},\\ 
& \quad \frac{\|\sin\Theta(\hat U_3, U_3)\|_{\rm F}^2 - p_3\sigma^2\|\Lambda_3^{-1}\|_{\rm F}^2}{\sqrt{2p_3}\sigma^2\|\Lambda_3^{-2}\|_{\rm F}}\bigg) \stackrel{\rm d.}{\longrightarrow} N(0, I_3)\quad {\rm as}\quad p\to\infty.
\end{align*}

\begin{remark}
	We briefly compare Theorem \ref{thm:na_tsvd} with the existing results in the literature. The asymptotic normality of $\|\sin\Theta(\hat U_j, U_j)\|_{\rm F}^2$ in Theorem \ref{thm:na_tsvd} requires SNR condition $\lambda_{\submin}/\sigma \gg (r_{\submax}p)^{3/4}$, which is slightly stronger than the optimal SNR condition $\lambda_{\submin}/\sigma \geq C_0p^{3/4}$ for achieving the consistent estimation in \citep[Theorem 1]{zhang2018tensor} (if $r\geq 1$), matches the condition in \cite[Theorem 1]{zheng2015interpolating} (if $r=1$), and weaker than the condition in \cite[Theorem 4]{richard2014statistical} (if $r=1$). Second, note that Theorem~\ref{thm:na_tsvd} implies $\EE\|\sin\Theta(\hat U_j, U_j)\|_{\rm F}^2=\big(1+o(1)\big)p_j\sigma^2\|\Lambda^{-1}_j\|_{\rm F}^2$. To the best of our knowledge, this is the first result with a precise constant characterization of the estimation error in tensor PCA.  

	Compared with the conditions for consistent estimation in \cite{zhang2018tensor}, our Theorem \ref{thm:na_tsvd} is for valid statistical inference, which requires the additional $\kappa_0^2p^{1/2}$ term in SNR and a rank condition $r_{\submax}^3/p\to 0$. These terms emerge from  technical issues, in particular from the way we bound the higher-order terms in the empirical spectral projector $\hat U_j\hat U_j^{\top}$ to better cope with the dependence across iterations. We note \cite{chen2019phase} proves that a rank-one planted tensor is distinguishable from the pure noise tensor if SNR $\lambda_{\submin}/\sigma\geq C_0p^{1/2}$ holds for a certain positive constant threshold $C_0$. In comparison, our Theorem \ref{thm:na_tsvd} requires a stronger condition $\lambda_{\submin}/\sigma\gg p^{3/4}$. In fact, the gap between $p^{1/2}$ and $p^{3/4}$ is fundamental. Without considering the computational feasibility, the SNR threshold $p^{1/2}$ is sufficient not only for detection but also for estimation (see \cite[Theorem 2]{zhang2018tensor}). However, when SNR falls below the threshold $p^{3/4}$, various pieces of evidence were established in the literature, as described in the first paragraph of Section \ref{sec:related-prior-work}, that show no polynomial time algorithms can reliably estimate the principal components. 
\end{remark}

While Theorem~\ref{thm:na_tsvd} characterizes the asymptotic distribution of $\|\sin\Theta(\hat U_j, U_j)\|_{\rm F}^2$ for tensor PCA model, the result is not immediately applicable to uncertainty quantification of $\hat U_j$ since $\|\Lambda_j^{-1}\|_{\rm F}^2$,  $\|\Lambda_j^{-2}\|_{\rm F}$, and $\sigma^2$ are often unknown in practice. We thus propose an estimate for $\Lambda_j$, $\sigma$:
\begin{equation}\label{eq:lambda_j-sigma_estimate}
\begin{split}
& \hat{\Lambda}_j =\text{ diagonal matrix with the top $r_j$ singular values of } \calM_j\big(\calA \times_{j+1}\hat{U}_{j+1}^{\top} \times_{j+2}\hat{U}_{j+2}^{\top}\big),\\
& \hat\sigma = \big\|\calA-\calA\times_1 \hat U_1\hat U_1^\top\times_2 \hat U_2\hat U_2^\top\times_3 \hat U_3\hat U_3^\top\big\|_{\rm F}/\sqrt{p_1p_2p_3}.
\end{split}
\end{equation}
We can prove a deviation bound for $\hat{\sigma}$ and the normal approximation for $\|\sin\Theta(\hat U_j, U_j)\|_{\rm F}^2$ with the proposed plug-in estimators.
\begin{lemma}\label{lm:variance_estimation}
	Under conditions of Theorem~\ref{thm:na_tsvd}, there exist two constants $C_1, C_2 > 0$ such that 
	\begin{equation*}
	\PP\left\{|\hat\sigma^2/\sigma^2 - 1| \leq C_2(\kappa_0\sqrt{r_{\submax}}p^{-1} + p^{-3/4}\sqrt{\log(p)})\right\} \geq 1-C_1p^{-3}.
	\end{equation*}
\end{lemma}
\begin{theorem}[Inference for Tucker Low-rank Tensor PCA]\label{thm:adaptive_tsvd}
	Suppose the conditions in Theorem \ref{thm:na_tsvd} hold. Let $\hat{\Lambda}_1 \in \RR^{r_1 \times r_1}$ and $\hat\sigma$ be defined as \eqref{eq:lambda_j-sigma_estimate}. There exist absolute constants $c_1,C_0,C_1,C_2,C_3>0$ such that if $\lambda_{\submin}/\sigma \geq C_0(p^{3/4} + \kappa_0^2p^{1/2})$, then for $j=1,2,3$,
	\begin{equation*}
	\begin{split}
	& \sup_{x\in\RR} \left|\PP\left(\frac{\|\sin\Theta(\hat U_j, U_j)\|_{\rm F}^2 - p_j\hat\sigma^2\|\hat{\Lambda}_j^{-1}\|_{\rm F}^2}{\sqrt{2p_j}\hat\sigma^2\|\hat{\Lambda}_j^{-2}\|_{\rm F}} \leq x\right) - \Phi(x)\right|\\
	\leq & C_1e^{-c_1p}+C_2\bigg(\frac{r_{\submax}^{3/2}\kappa_0^{6}p^{3/2}}{(\lambda_{\submin}/\sigma)^2} +  \frac{\kappa_0^3\sqrt{pr_{\submax}(r_{\submax}^2+\log p)}}{\lambda_{\submin}/\sigma} + \frac{\sqrt{\log(p)}}{p^{1/4}} + \frac{\kappa_0\sqrt{r_{\submax}}}{\sqrt{p}}\bigg) + C_3\frac{r_{\submax}^{3/2}}{\sqrt{p}}.
	\end{split}
	\end{equation*}
\end{theorem}
When the condition number $\kappa_0=O(1)$, $(pr_{\submax})^{3/4}(\lambda_{\submin}/\sigma)^{-1}\to 0$, and $r_{\submax}^3/p\to 0$ as $p\to\infty$, Theorem~\ref{thm:adaptive_tsvd} implies
\begin{equation}\label{eq:tsvd_data_na}
\frac{\|\sin\Theta(\hat U_j, U_j)\|_\F^2 - p_j\hat\sigma^2\|\hat{\Lambda}_j^{-1}\|_{\rm F}^2}{\sqrt{2p_j}\hat{\sigma}^2\|\hat{\Lambda}_j^{-2}\|_{\rm F}}\stackrel{{\rm d.}}{\longrightarrow} N(0,1)\quad {\rm as}\quad p\to\infty.
\end{equation}
Equation (\ref{eq:tsvd_data_na}) is readily applicable to statistical inference for $U_j$. After getting $\hat U_j$ by Algorithm \ref{algo:am_optimal_PCA}, we propose a $(1-\alpha)$-level confidence region for $U_j$ as
\begin{align}
{\rm CR}_{\alpha}(\hat U_j):=\Big\{V\in\OO_{p_j,r_j}: \|\sin\Theta(\hat U_j,V)\|_{\rm F}^2 \leq p_j\hat{\sigma}^2\|\hat\Lambda_j^{-1}\|_{\rm F}^2 + z_{\alpha}\sqrt{2p_j}\hat{\sigma}^2\|\hat\Lambda_j^{-2}\|_{\rm F} \Big\},\label{eq:tsvd_CI}
\end{align}
where $z_{\alpha}=\Phi^{-1}(1-\alpha)$ is the $(1-\alpha)$ quantile of the standard normal distribution. The following corollary is an immediate result of Theorem~\ref{thm:adaptive_tsvd}, which confirms that the confidence region ${\rm CR}_{\alpha}(\hat U_1)$ is indeed asymptotically accurate.
\begin{corollary}[Confidence region for tensor PCA]
	Suppose the conditions of Theorem~\ref{thm:adaptive_tsvd} hold and the confidence region ${\rm CR}_{\alpha}(\hat U_j)$ is defined in (\ref{eq:tsvd_CI}). If $\kappa_0^{6}(r_{\submax}^{3/2}p^{3/2}+r_{\submax}p\log p)(\lambda_{\submin}/\sigma)^{-2} \to 0$ and $r_{\submax}^3/p\to 0$ as $p\to\infty$, then
	\begin{align*}
	\lim_{p\to\infty}\PP\big(U_j\in {\rm CR}_{\alpha}(\hat U_j)\big)= 1-\alpha. 
	\end{align*}
\end{corollary}

We note that, through a more sophisticated analysis,  Theorem \ref{thm:na_tsvd} can be generalized to the setting with sub-Gaussian noise. For simplicity, we only prove the following rank-one case, which has been the focus of many papers on tensor PCA.
\begin{theorem}[Rank-one tensor PCA under sub-Gaussian noise]\label{thm:PCA-sub-Gaussian}
Suppose Assumption \ref{assump:tpca} holds for tensor PCA model (\ref{eq:PCA_model}) with $r=1$, $p_j\asymp p$ for $j=1,2,3$, $\calZ$ has i.i.d. entries with $\EE\calZ(i_1,i_2,i_3)=0, \EE[\calZ(i_1,i_2,i_3)^2]=\sigma^2$, $\EE[\calZ(i_1,i_2,i_3)^4]/\sigma^4=\nu$ and $\|\calZ(i_1,i_2,i_3)\|_{\psi_2}\leq C\sigma$ for some constant $C>0$.  There exist absolute constants $c_1, C_0, C_1, C_2, C_3>0$ such that if $\lambda/\sigma\geq C_0p^{3/4}$, then 
\begin{align*}
\sup_{x\in\RR}&\left|\PP\left(\frac{\|\sin\Theta(\hat U_1, U_1)\|_{\rm F}^2-p_1\sigma^2\lambda^{-2}}{\sigma^2\lambda^{-2}\sqrt{p_1\big(2+(\nu-3)\|U_2\|_{4}^4\|U_3\|_{4}^4\big)}}\leq x\right) -\Phi(x)\right|\leq \frac{C_2}{p^{1/2}\big(2+(\nu-3)\|U_2\|_{4}^4\|U_3\|_{4}^4\big)^{3/2}}\\
&\quad+C_3\Big(\frac{\sqrt{p\log p}}{\lambda/\sigma}+\frac{p^{3/2}}{(\lambda/\sigma)^2}+\frac{\log p}{\sqrt{p}}\Big)\cdot \frac{1}{\sqrt{2+(\nu-3)\|U_2\|_{4}^4\|U_3\|_{4}^4}}+C_1e^{-c_1p}.
\end{align*}
Similar results can be derived for $\|\sin\Theta(\hat U_2, U_2)\|_{\rm F}^2$ and $\|\sin\Theta(\hat U_3, U_3)\|_{\rm F}^2$. 
\end{theorem}

If $\nu-1\geq c_0$ for some absolute constant $c_0>0$ and $\lambda^{-1}\sigma p^{3/4}\to 0$ as $p\to \infty$,  Theorem \ref{thm:PCA-sub-Gaussian} implies 
$$
\frac{\|\sin\Theta(\hat U_1, U_1)\|_{\rm F}^2-p_1\sigma^2\lambda^{-2}}{\sigma^2\lambda^{-2}\sqrt{p_1\big(2+(\nu-3)\|U_2\|_{4}^4\|U_3\|_{4}^4\big)}} \stackrel{{\rm d.}}{\longrightarrow} N(0,1)\quad {\rm as}\quad p\to\infty.
$$
The SNR condition in Theorem \ref{thm:PCA-sub-Gaussian} is the same as that in Theorem \ref{thm:na_tsvd}. Moreover, the asymptotic variance of $\|\sin\Theta(\hat U_1, U_1)\|_{\rm F}^2$ includes the kurtosis of the noise distribution $\nu$, which can be challenging to estimate data-drivenly. We leave the estimation of the kurtosis and the data-driven inference for $\hat{U}_k$ as future research. 

%%%%%%%%%%%%%%%%%%%%%%
\subsection{Inference for Tucker Low-rank Tensor Regression}\label{asymptotic-distribution-regression}
%%%%%%%%%%%%%%%%%%%%%%

This section is devoted to the asymptotic distribution and inference in low-rank tensor regression. We first introduce the following assumption on the initialization for Algorithm \ref{algo:am_optimal_regression}.
\begin{assump}\label{assump:tr}
	Under tensor regression model (\ref{eq:tr_model}) with $\calX(i_1,i_2,i_3)\stackrel{i.i.d.}{\sim} N(0,1)$, $\Var(\xi_i) = \sigma^2$ and $\|\xi_i\|_{\psi_2} \leq C\sigma$ for some constant $C > 0$, there is an event $\calE_0$ with $\PP(\calE_0)\geq 1-C_1e^{-c_1p}$ for some absolute constants $c_1,C_1>0$ so that, under $\calE_0$, the initialization $\tilde\calT=(\hat U_1^{(0)},\hat U_2^{(0)},\hat U_3^{(0)})\cdot \hat \calG^{(0)}$ satisfy $\|\tilde\calT-\calT\|_{\rm F}^2\leq C_2pr_{\submax}\sigma^2/n$ or $\max_j\|\sin\Theta(\hat U_j^{(0)},U_j)\| \leq C_2\sqrt{p/n}\sigma/\lambda_{\submin}$ for some absolute constant $C_2>0$. 
\end{assump}

The claimed bound of $\|\tilde\calT-\calT\|_{\rm F}^2$ in Assumption~\ref{assump:tr} is attainable, for instance, by the gradient descent algorithm developed in \cite{han2020optimal} and the importance sketching algorithm developed in \cite{zhang2020islet} under the SNR condition $n(\lambda_{\submin}/\sigma)^2\geq  Cp^{3/2}$ and the sample size condition $n\geq  Cp^{3/2}r_{\submax}$. The theoretical guarantees for this claim can be found in \cite[Theorem 4.2]{han2020optimal} and \cite[Theorem 4]{zhang2020islet}.

Based on Assumption \ref{assump:tr}, we establish the following asymptotic results for tensor regression.
\begin{theorem}\label{thm:regression_alternating}
	Suppose Assumption~\ref{assump:tr} holds for tensor regression model (\ref{eq:tr_model}), $\calX(i_1,i_2,i_3)\stackrel{i.i.d.}{\sim} N(0,1)$, $\Var(\xi_i) = \sigma^2$, and $\|\xi_i\|_{\psi_2} \leq C\sigma$ for some constant $C > 0$, $p_j\asymp p$ for $j=1,2,3$, and $\kappa(\calT)\leq \kappa_0$. Let $\hat U_j$s be the output of two-iteration alternating minimization (Algorithm~\ref{algo:am_optimal_regression}). There exist absolute constants $c_1,C_0,C_1,C_2,C_3,C_4>0$ such that if $n(\lambda_{\submin}/\sigma)^2\geq C_0(p^{3/2} \vee \kappa_0^4pr_{\submax}^2)$ and $n\geq C_2(p^{3/2} \vee \kappa_0^2pr_{\submax}^3)$, then
	\begin{equation*}
	\begin{split}
	&\sup_{x \in \RR}\left|\PP\left(\frac{\|\sin\Theta(\hat U_j,U_j)\|_{\rm F}^2 -p_jn^{-1}\sigma^2\|\Lambda_j^{-1}\|_{\rm F}^2}{\sqrt{2p_j}n^{-1}\sigma^2\|\Lambda_j^{-2}\|_{\rm F}} \leq x\right) - \Phi(x)\right|\\ 
	\leq & \frac{C_3\kappa_0^{4}r_{\submax}^{5/2}p^{3/2}}{n} + C_3\kappa_0^3\Big(\frac{r_{\submax}^5p\log^2n}{n}\Big)^{1/2} + \frac{C_3p^{3/2}}{n}\Big(\frac{\kappa_0^5r_{\submax}^2}{\lambda_{\submin}/\sigma}+\frac{\kappa_0^5r_{\submax}^{3/2}}{(\lambda_{\submin}/\sigma)^2}\Big)\\
	& \hspace{1cm} + C_3\kappa_0^4\Big(\frac{pr_{\submax}^3+r_{\submax}p\log p}{n(\lambda_{\submin}/\sigma)^2}\Big)^{1/2} 
	+C_1e^{-c_1p}+ C_4\frac{r_{\submax}^{3/2}}{\sqrt{p}},
	\end{split}
	\end{equation*}
	where $\Lambda_j$ is the $r_j\times r_j$ diagonal matrix containing the singular values of $\calM_j(\calT)$.
\end{theorem}
If the condition number $\kappa_0=O(1)$, $r_{\submax}^3/p\to 0$, $(r_{\submax}^{5/2}p^{3/2}+r_{\submax}^5p\log^2n)/n\to 0$ and $r_{\submax}^{3/2}p^{3/2}/(n(\lambda_{\submin}/\sigma)^2)\to 0$
as $p\to\infty$, Theorem~\ref{thm:regression_alternating} implies 
$$
\frac{\|\sin\Theta(\hat U_j,U_j)\|_{\rm F}^2 -p_jn^{-1}\sigma^2\|\Lambda_j^{-1}\|_{\rm F}^2}{\sqrt{2p_j}n^{-1}\sigma^2\|\Lambda_j^{-2}\|_{\rm F}} \stackrel{\rm d.}{\longrightarrow} N(0, 1)\quad {\rm as}\quad  p\to\infty. 
$$

To make inference for tensor regression, we develop the following asymptotic normal distribution for $\|\sin\Theta(\hat{U}_j, U_j)\|_\F$ with the plug-in estimates of $\Lambda_j$.
\begin{theorem}[Tensor regression]\label{thm:regression_adaptive}
	Suppose the conditions in Theorem \ref{thm:regression_alternating} hold. Let $\hat\Lambda_j = \diag(\hat\lambda_1, \dots, \hat\lambda_{r_j})$ be a diagonal matrix containing the singular values of $\calM_1(\hat \calG)$, where $\hat\calG$ is the output of Algorithm~\ref{algo:am_optimal_regression}. There exist absolute constants $c_1,C_0,C_1,C_2,C_3,C_4>0$ such that if $n(\lambda_{\submin}/\sigma)^2\geq C_0(p^{3/2} \vee \kappa_0^6pr_{\submax}^2)$ and $n\geq C_2(p^{3/2} \vee \kappa_0^8pr_{\submax}^3)$, then for $j=1,2,3$, 
	\begin{equation*}
	\begin{split}
	&\sup_{x \in \RR}\left|\PP\left(\frac{\|\sin\Theta(\hat U_j,U_j)\|_{\rm F}^2 - p_jn^{-1}\sigma^2\|\hat\Lambda_j^{-1}\|_{\rm F}^2}{\sqrt{2p_j}n^{-1}\sigma^2\|\hat\Lambda_j^{-2}\|_{\rm F}} \leq x\right) - \Phi(x)\right|\\ 
	& \leq C_3\frac{\kappa_0^{4}r_{\submax}^{5/2}p^{3/2}}{n} + C_3\kappa_0^3\Big(\frac{r_{\submax}^5p\log^2n}{n}\Big)^{1/2} + \frac{C_3p^{3/2}}{n}\Big(\frac{\kappa_0^5r_{\submax}^2}{\lambda_{\submin}/\sigma}+\frac{\kappa_0^5r_{\submax}^{3/2}}{(\lambda_{\submin}/\sigma)^2}\Big)\\
	& \hspace{1cm} + C_3\kappa_0^4\Big(\frac{pr_{\submax}^3+r_{\submax}p\log p}{n(\lambda_{\submin}/\sigma)^2}\Big)^{1/2}+C_1e^{-c_1p}+ C_4\frac{r_{\submax}^{3/2}}{\sqrt{p}}.
	\end{split}
	\end{equation*}
\end{theorem}
We propose the following $(1-\alpha)$-level confidence region for $U_j$:
\begin{equation}\label{eq:tr_CI}
\widetilde{\rm CR}_{\alpha}(\hat U_j):=\bigg\{V\in\OO_{p_j,r_j}: \|\sin\Theta(\hat U_j,V)\|_{\rm F}^2 \leq \frac{p_j\sigma^2\|\hat\Lambda_j^{-1}\|_{\rm F}^2}{n}+z_\alpha\frac{\sqrt{2p_j}\sigma^2\|\hat\Lambda_j^{-2}\|_{\rm F}}{n}\bigg\}.
\end{equation}
The following corollary establishes the coverage probability of the proposed confidence region. 
\begin{corollary}[Confidence region for tensor regression]
	Suppose the conditions of Theorem~\ref{thm:regression_adaptive} hold and the confidence region $\widetilde{\rm CR}_{\alpha}(\hat U_j)$ is defined by (\ref{eq:tr_CI}). If $(\kappa_0^{5}r_{\submax}^{5/2}p^{3/2}+\kappa_0^6r_{\submax}^5p\log^2n)/n \to 0$, $\kappa_0^5r_{\submax}^{3/2}p^{3/2}/(n(\lambda_{\submin}/\sigma)^2)\to 0$ and $r_{\submax}^3/p\to 0$ as $p\to\infty$, then
	\begin{align*}
	\lim_{p\to\infty}\PP\big(U_j\in \widetilde{\rm CR}_{\alpha}(\hat U_j)\big)= 1-\alpha. 
	\end{align*}
\end{corollary}

\begin{remark}[Selection of $\sigma$]
	When $\sigma$ is unknown, we can estimate it by a sample splitting scheme as follows. First, we retain a part of sample $\{(\calX_k, Y_k)\}_{k=1}^{\lceil p^{3/2}\rceil}$ and use the other samples to compute the estimator $\tilde \calT$. Define
	$$
	\hat\sigma^2:= \sum\nolimits_{k=1}^{\lceil p^{3/2}\rceil} \big(Y_k-\langle \tilde\calT, \calX_k\rangle\big)^2 / \lceil p^{3/2}\rceil.
	$$
	Under Assumption~\ref{assump:tr} and conditions of Theorem~\ref{thm:regression_alternating}, we can show with probability at least $1-p^{-3}$, 
	$\big|\hat \sigma^2/\sigma^2-1\big|=O\big(p^{-3/4}\sqrt{\log p}+r_{\submax}pn^{-1}\big)$. By plugging in $\hat\sigma$ to \eqref{eq:tr_CI}, we obtain a data-driven $(1-\alpha)$ asymptotic confidence region for $U_j$.
\end{remark}

%%%%%%%%%%%%%%%%%%%%%
\subsection{Proof Sketch}\label{sec:proof-sketch}
%%%%%%%%%%%%%%%%%%%%%

In this section, we briefly explain the proof strategy for tensor PCA model, i.e., Theorem~\ref{thm:na_tsvd}. The proof for tensor regression model is more complicated but shares similar spirits. Without loss of generality, we assume $\sigma=1$. First,
\begin{equation*}
\begin{split}
2\|\sin\Theta(\hat U_1,U_1)\|_{\rm F}^2= & \|\hat U_1\hat U_1^{\top}-U_1U_1^{\top}\|_{\rm F}^2=2r_1-2\big<\hat U_1\hat U_1^{\top},U_1U_1^{\top}\big>\\
= & -2\big<U_1U_1^{\top}, \hat U_1\hat U_1^{\top}-U_1U_1^{\top}\big>.
\end{split}
\end{equation*}
It thus suffices to investigate the distribution of $\big<U_1U_1^{\top}, \hat U_1\hat U_1^{\top}-U_1U_1^{\top}\big>$. By Algorithm~\ref{algo:am_optimal_PCA}, $\hat U_1$ are the top-$r_1$ left singular vectors of $\calM_1\big(\calA\times_2 \hat U_2^{(1)\top}\times_3 \hat U_3^{(1)\top}\big)$. As a result, $\hat U_1\hat U_1^{\top}$ is the spectral projector and can be decomposed as
$$
\calM_1(\calA)\Big(\hat U_2^{(1)}\hat U_2^{(1)\top}\otimes \hat U_3^{(1)}\hat U_3^{(1)\top}\Big)\calM_1^{\top}(\calA)=:\calM_1(\calT)\calM_1^{\top}(\calT)+D_1^{(1)}.
$$
The high-level ideas of the proof include the following steps.

{\it Step 1}: We apply the spectral representation formula (\cite{xia2019normal}; also see the statement in Lemma \ref{lem:spectral} from the supplementary materials) and expand
$$
\hat U_1\hat U_1^{\top}=U_1U_1^{\top}+S_1(D_1^{(1)})+S_2(D_1^{(1)})+S_3(D_1^{(1)})+\sum_{k\geq 4}S_k(D_1^{(1)}),
$$
where $S_k(\cdot)$ denotes the $k$th order perturbation term:
$$
S_{k}(D_1^{(1)})=\sum_{s_1+\cdots+s_{k+1}=k}(-1)^{1+\tau(\bs)}\cdot B_1^{-s_1}D_1^{(1)}B_{1}^{-s_2}D_1^{(1)}B_1^{-s_3}\cdots B_1^{-s_k}D_1^{(1)}B_1^{-s_{k+1}},
$$
where $B_1^{-k}=U_1\Lambda_1^{-2k}U_1^{\top}$ for each positive integer $k$, $B_1^{0}:= I_{p_1} - U_{1}U_{1}^\top$, $s_1,\cdots,s_{k+1}$ are non-negative integers, and $\tau(\bs)=\sum_{j=1}^{k+1}\II(s_j>0)$. 

{\it Step 2}: Since $\big<U_1U_1^{\top}, S_1(D_1^{(1)})\big>=0$ and 
$\|S_k(D_1^{(1)}) \|\leq (C_1\kappa_0^2\sqrt{p}/\lambda_{\submin})^k$ with high probability, we can write
$$
\big<\hat U_1\hat U_1^{\top}-U_1U_1^{\top}, U_1U_1^{\top}\big>=\big<S_2(D_1^{(1)}),U_1U_1^{\top}\big>+\big<S_3(D_1^{(1)}),U_1U_1^{\top}\big>+O\Big(\frac{r_{\submax}\kappa_0^8p^2}{\lambda_{\submin}^4}\Big).
$$
In other words, the higher order terms ($k\geq 4$) can be bounded with high probability, which becomes small order terms. 

{\it Step 3}: We show, with high probability, the third order term can be bounded by
$$
\big|\big<S_3(D_1^{(1)}),U_1U_1^{\top}\big>\big|=O\Big(\frac{\kappa_0^3p\sqrt{r_{\submax}\log p}}{\lambda_{\submin}^3}+\frac{\kappa_0^3p^2r_{\submax}^{3/2}}{\lambda_{\submin}^4}\Big)
$$
and becomes small order term. Now, it suffices to only investigate the second order term carefully. 

{\it Step 4}: We decompose the second order term $\big<S_2(D_1^{(1)}),U_1U_1^{\top}\big>$ into a leading term and remainder terms. Similarly to {\it Step 2} and {\it Step 3}, we show that the remainder terms are, with high probability, bounded by $O(\kappa_0^3p\sqrt{r_{\submax}\log p}\lambda_{\submin}^{-3}+\kappa_0^3p^2r_{\submax}^{3/2}\lambda_{\submin}^{-4})$. 

{\it Step 5}: We prove that the leading term of $\big<S_2(D_1^{(1)}),U_1U_1^{\top}\big>$ can be written as a sum of independent random variables, which yields a normal approximation by Berry-Essen Theorem. Finally, combining all these steps, we get the normal approximation for $\|\hat U_1\hat U_1^{\top}-U_1U_1^{\top}\|_{\rm F}^2$.

Among these steps, Steps 4 and 5 are the most technically involved. Throughout the proof, we apply the spectral representation formula at multiple stages to prove sharp upper bounds for higher-order terms, and establish central limit theorem for the second-order term. 

The following lemmas are used in our proof and could be of independent interest. First, Lemma \ref{lm:orlicz} is used to establish the concentration inequalities for the sum of random variables that have heavier tails than Gaussian.
\begin{lemma}[Orlicz $\psi_{\alpha}$-norm for product of random variables]\label{lm:orlicz}
	Suppose $X_1, \dots, X_n$ are $n$ random variables (not necessarily independent) satisfying $\|X_i\|_{\psi_{\alpha_i}} \leq K_i$. Define $\bar{\alpha} = \left(\sum_{i=1}^{n}\alpha_i^{-1}\right)^{-1}$. Then 
	\begin{equation*}
	\left\|\prod_{i=1}^{n}X_i\right\|_{\psi_{\bar{\alpha}}} \leq \prod_{i=1}^n K_i.
	\end{equation*}
\end{lemma}
Next, Lemma \ref{lm:third_moment_gaussian} provides a tight probabilistic upper bound for sum of third moments of Gaussian random matrices.
\begin{lemma}\label{lm:third_moment_gaussian}
	Suppose $Z_1, \dots, Z_n \in \RR^{p \times r}$ are independent random matrices satisfying $Z_i(j,k) \stackrel{i.i.d.}{\sim} N(0, 1)$. Then there exist two universal constants $C, C_1 > 0$ such that for fixed $M_1, \dots, M_n \in \RR^{p \times r}$, 
	\begin{equation*}
	\PP\left(\left|\sum_{i=1}^{n}\|Z_i\|_{\rm F}^2\langle Z_i, M_i\rangle\right| \geq Cpr\left(\sum_{i=1}^{n}\|M_i\|_{\rm F}^2\right)^{1/2}\sqrt{\log(p)}\right) \leq p^{-C_1}.
	\end{equation*}
\end{lemma}
%%%%%%%%%%%%%%%%%%%%%%%

\section{PCA for Orthogonally Decomposable Tensors}\label{sec:T-orth}
%%%%%%%%%%%%%%%%%%%%%%%

In this section, we specifically focus on the tensor PCA model (\ref{eq:PCA_model}) with orthogonally decomposable signal tensor $\calT$:
\begin{equation}\label{eq:T-orth}
\calT=\sum_{i=1}^r \lambda_i \cdot u_i\otimes v_i\otimes w_i,
\end{equation}
where $U=(u_1,\cdots, u_r)\in \OO_{p_1, r}$, $V=(v_1,\cdots,v_{r})\in\OO_{p_2, r}$, and $W=(w_1,\cdots,w_{r})\in\OO_{p_3, r}$ all have orthonormal columns; the singular values satisfy $\lambda_{\submin} = \min\{\lambda_1, \ldots, \lambda_r\}>0$. Here, for any $u \in \RR^{p_1}, v \in \RR^{p_2}, w \in \RR^{p_3}$, $u \otimes v \otimes w$ is a $p_1 \times p_2 \times p_3$ tensor whose $(i,j,k)$th entry is $u(i)v(j)w(k)$.

Our goal is to make inference on the principal components based on a noisy observation $\calA=\calT+\calZ$. Different from the inference for Tucker low-rank tensor discussed in Section~\ref{sec:normal_app}, where an accurate estimation is hopeful only for the joint column space of $U_j$ due to the non-identifiability of Tucker decomposition, we can make inference for each individual vector $\{u_j, v_j, w_j\}$ if $\calT$ is orthogonally decomposable as (\ref{eq:T-orth}). Given some estimates $\{\hat u_j^{(0)}, \hat v_j^{(0)}, \hat w_j^{(0)}\}_{j=1}^r$, we propose to pass them to a post-processing step by two-iteration procedure in Algorithm~\ref{algo:power_iter} to obtain the test statistics $\{\hat{u}_j, \hat{v}_j, \hat{w}_j\}_{j=1}^r$.
\begin{algorithm}
	\SetAlgoLined
	\KwInput{$\calA$, initialization $\{\hat u_j^{(0)}, \hat v_j^{(0)}, \hat w_j^{(0)}\}_{j=1}^r$;}
	\For{$t = 0, 1$}{
		\For{$j=1,2,\cdots,r$}{
			Compute $\hat u_j^{(t+0.5)}=\calA\times_2 \hat v_j^{(t)\top}\times_3 \hat w_j^{(t)\top}$; Update $\hat u_j^{(t+1)}=\hat u_j^{(t+0.5)}\|\hat u_j^{(t+0.5)}\|_2^{-1}$\;
			Compute $\hat v_j^{(t+0.5)}=\calA\times_1 \hat u_j^{(t)\top}\times_3 \hat w_j^{(t)\top}$; Update $\hat v_j^{(t+1)}=\hat v_j^{(t+0.5)}\|\hat v_j^{(t+0.5)}\|_2^{-1}$\;
			Compute $\hat w_j^{(t+0.5)}=\calA\times_1 \hat u_j^{(t)\top}\times_2 \hat v_j^{(t)\top}$; Update $\hat w_j^{(t+1)}=\hat w_j^{(t+0.5)}\|\hat w_j^{(t+0.5)}\|_2^{-1}$\;
		}	
	}
	\KwOutput{$\hat u_j=\hat u_j^{(2)}, \hat v_j=\hat v_j^{(2)}$ and $\hat w_j =\hat w_j^{(2)}$ for all $j=1,\cdots,r$.}
	\caption{Power Iterations for Orthogonally decomposable $\calT$}
	\label{algo:power_iter}
\end{algorithm}

Since our primary interest is about the statistical inference for $\{u_j, v_j, w_j\}$, we assume that the initializations of Algorithm~\ref{algo:power_iter} satisfies the following Assumption \ref{assump:T-orth}. Such an assumption is achievable by the power iteration method with $k$-means initialization introduced in \cite{anandkumar2014tensor} along with the theoretical guarantees developed in \cite{liu2017characterizing} when $\lambda/\sigma\geq Cp^{3/4}$. 
\begin{assump}\label{assump:T-orth}
	Under the tensor PCA model (\ref{eq:PCA_model}) with $\calT$ being orthogonally decomposable as (\ref{eq:T-orth}), there is an event $\calE_0$ with $\PP(\calE_0)\geq 1-C_1e^{-c_1p}$ for some absolute constants $c_1,C_1>0$ such that, under $\calE_0$, the initializations $\{\hat u_j^{(0)}, \hat v_j^{(0)},\hat w_j^{(0)}\}_j$ satisfy $\max\big\{\|\hat u_{\pi(j)}^{(0)}-u_j\|_2, \|\hat v_{\pi(j)}^{(0)}-v_j\|_2, \|\hat w_{\pi(j)}^{(0)}-w_j\|_2\big\}\leq C_2\sigma\sqrt{p}/\lambda_j$ for some permutation  $\pi: [r]\to [r]$, all $1\leq j\leq r$, and some absolute constant $C_2>0$. 
\end{assump}

We establish the asymptotic normality for the outcome of Algorithm \ref{algo:power_iter} as follows. 

\begin{theorem}[PCA for orthogonally decomposable tensors]\label{thm:T-orth}
	Suppose Assumption~\ref{assump:T-orth} holds for tensor PCA model (\ref{eq:PCA_model}) with an orthogonally decomposable $\calT$ as (\ref{eq:T-orth}), $\calZ(i_1,i_2,i_3)\stackrel{i.i.d.}{\sim} N(0,\sigma^2)$, $p_j\asymp p$ for $j=1,2,3$, and $\kappa(\calT)\leq \kappa_0$. Let $\{\hat u_j, \hat v_j, \hat w_j\}_{j=1}^r$ be the output of Algorithm~\ref{algo:power_iter}. There exist absolute constants $c_1,C_0,C_1,C_2,C_3>0$ such that if $\lambda_{\submin}/\sigma\geq C_0(p^{3/4}+\kappa_0^2p^{1/2})$, then
	\begin{align}
	\sup_{x\in\RR}\bigg|\PP\bigg(&\frac{\langle \hat u_{\pi(j)}, u_j\rangle^2-(1-p_j\sigma^2\lambda_j^{-2})}{\sqrt{2p_j}\sigma^2\lambda_j^{-2}}\leq x\bigg)-\Phi(x) \bigg|\notag\\
	&\hspace{2cm}\leq C_1e^{-c_1p}+C_2\left(\frac{\kappa_0^{6}\sigma^2(pr)^{3/2}}{\lambda_{\submin}^2} + \frac{\kappa_0^2\sigma(p\log p)^{1/2}}{\lambda_{\submin}}\right)  + C_3\frac{r^{3/2}}{\sqrt{p}} \label{eq:T-orth-na}
	\end{align}
	for all $j=1,\cdots,r$. Here, $\pi(\cdot)$ is the permutation introduced in Assumption \ref{assump:T-orth}. Moreover, let $\hat \lambda_j = \|\calA\times_2 \hat v_j^{\top}\times_3 \hat w_j^{\top}\|_2$. Then, (\ref{eq:T-orth-na}) also holds if $\lambda_j$ is replaced by $\hat \lambda_j$ and $\kappa_0^2\sigma(p\log p)^{1/2}\lambda_{\submin}^{-1}$ is replaced by $\kappa_0^3\sigma\sqrt{pr(r^2+\log p)}\lambda_{\submin}^{-1}$. 
	Similar results also hold for $\langle\hat v_{\pi(j)}, v_j \rangle^2$ and $\langle\hat w_{\pi(j)}, w_j \rangle^2$. 
\end{theorem}
By Theorem~\ref{thm:T-orth}, if $\lambda_{\submin}/\sigma \gg \kappa_0^3(pr)^{3/4}+\kappa_0^2(p\log p)^{1/2}$ and $r\ll p^{1/3}$, then for each $j=1,\cdots,r$, 
$$
\frac{\langle \hat u_{\pi(j)}, u_j\rangle^2-(1-p_j\sigma^2\lambda_j^{-2})}{\sqrt{2p_j}\sigma^2\lambda_j^{-2}}\stackrel{{\rm d.}}{\longrightarrow} N(0,1)\quad {\rm as }\quad p\to\infty.
$$
Similarly to Section~\ref{sec:inference-PCA}, we plug in data-driven estimates of $\lambda_j$ and $\sigma^2$ and construct a $(1-\alpha)$ confidence region for $u_j$ as
\begin{equation}\label{eq:CR-u_j}
{\rm CR}_{\alpha}(\hat u_{\pi(j)}):=\Big\{v\in\RR^{p_j}: \|v\|_2=1\ {\rm and}\ \langle \hat u_{\pi(j)}, v\rangle^2\geq (1-p_j\hat\sigma^2\hat\lambda_{\pi(j)}^{-2})-z_{\alpha}\sqrt{2p_j}\hat\sigma^2\hat\lambda_{\pi(j)}^{-2}\Big\}.
\end{equation}
The confidence region for $v_j, w_j$ can be constructed similarly.

\section{Entry-wise Inference for Rank-$1$ Tensors}\label{sec:entrywise-rank1}
%%%%%%%%%%%%%%%%%%

In this section, we consider the statistical inference for tensor PCA model with a rank-$1$ signal tensor:
\begin{equation}
\calA = \calT + \calZ, \quad \calT = \lambda \cdot u \otimes v \otimes w.
\end{equation}
Here, $u \in \mathbb{S}^{p_1-1}, v \in \mathbb{S}^{p_2-1}, w \in \mathbb{S}^{p_3-1}$, the singular value $\lambda > 0$, and $\calZ\overset{i.i.d.}{\sim}N(0, \sigma^2)$. We specifically aim to study the inference for any linear form of $u, v, w$, i.e., $\langle q_1, u\rangle, \langle q_2, v\rangle$, and $\langle q_3, w\rangle$, with arbitrary deterministic unit vectors $\{q_1, q_2, q_3\}$. We also aim to study the inference for each entry $\calT_{ijk}, i \in [p_1], j \in [p_2], k \in [p_3]$. To this end, we first apply the rank-1 power iteration in Algorithm \ref{algo:hooi} \citep{richard2014statistical,zhang2001rank}. Algorithm \ref{algo:hooi} can be roughly seen as a rank-1 special case of Algorithm \ref{algo:power_iter} for the Tucker low-rank tensor PCA and Algorithm \ref{al:tensor-PCA} for the orthogonally decomposable tensor PCA. 
\begin{algorithm}
	\SetAlgoLined
	\KwInput{$\calA$}
	Initialize $\hat u^{(0)}={\rm SVD}_{1}(\calM_1(\calA))$, $\hat v^{(0)}={\rm SVD}_{1}(\calM_2(\calA))$, $\hat w^{(0)}={\rm SVD}_{1}(\calM_3(\calA))$, $t=1$\;
	\While{$t < t_{\submax}$}{
		Compute $\hat u^{(t+0.5)}=\calA\times_2 \hat v^{(t)\top}\times_3 \hat w^{(t)\top}$; Update $\hat u^{(t+1)}=\hat u^{(t+0.5)}\|\hat u^{(t+0.5)}\|_2^{-1}$\;
		Compute $\hat v^{(t+0.5)}=\calA\times_1 \hat u^{(t)\top}\times_3 \hat w^{(t)\top}$; Update $\hat v^{(t+1)}=\hat v^{(t+0.5)}\|\hat v^{(t+0.5)}\|_2^{-1}$\;
		Compute $\hat w^{(t+0.5)}=\calA\times_1 \hat u^{(t)\top}\times_2 \hat v^{(t)\top}$; Update $\hat w^{(t+1)}=\hat w^{(t+0.5)}\|\hat w^{(t+0.5)}\|_2^{-1}$\;
		$t=t+1$;
	}
	$\hat\lambda=\calA \times_1 \hat u^{(t_{\submax})\top} \times_2 \hat v^{(t_{\submax})\top} \times_3 \hat w^{(t_{\submax})\top}$\;
	$\hat\calT=\hat\lambda(\hat u^{(t_{\submax})} \otimes \hat v^{(t_{\submax})} \otimes \hat w^{(t_{\submax})})$\;
	\KwOutput{$\hat u=\hat u^{(t_{\submax})}, \hat v=\hat v^{(t_{\submax})}$ $\hat w =\hat w^{(t_{\submax})}$, $\hat\lambda$ and $\hat\calT$.}
	\caption{Power iterations for rank-$1$ tensor $\calT$}
	\label{algo:hooi}
\end{algorithm}

Next, we establish the asymptotic normality for the output of Algorithm \ref{algo:hooi}, $\hat{u}, \hat{v}, \hat{w}$, under the essential SNR condition that ensures tensor PCA is solvable in polynomial time. Without loss of generality, we assume that the signs of $\hat u, \hat v, \hat w$ satisfy $\langle \hat u, u\rangle\geq 0, \langle \hat v, v\rangle\geq 0$ and $\langle \hat w, w\rangle\geq 0$ (otherwise one can flip the sign of $\hat u, \hat{v}, \hat{w}$ without changing the problem essentially). With a slight abuse of notation, let $u_i$, $v_j$, and $w_k$ be the $i$th entry of $u$, the $j$th entry of $v$, and the $k$th entry of $w$, respectively. 

\begin{theorem}\label{thm:entry_inference}
	Consider the tensor PCA model \eqref{eq:PCA_model} with Gaussian noise $\calZ(i_1,i_2,i_3) \stackrel{i.i.d.}{\sim} N(0,\sigma^2)$ and rank$(\calT) = 1$, $p_j\asymp p$ for $j=1,2,3$. Let $\hat{\lambda}, \hat u, \hat v, \hat w, \hat \calT$ be the outputs of  Algorithm \ref{algo:hooi} with iteration number $t_{\submax} \geq C_1\log(p)$ for constant $C_1 > 0$. Suppose $\lambda/\sigma \gg p^{3/4}$. 
	For any deterministic array $\{q_1^{(k)}, q_2^{(k)}, q_3^{(k)}\}_{k=1}^\infty$ satisfying $q_i^{(k)} \in \SS^{k-1}$, denote
	\begin{equation*}
	\begin{split}
	&T_{q_1^{(p_1)}, q_2^{(p_2)}, q_3^{(p_3)}} =\\ 
	& \left(\frac{\langle q_1^{(p_1)}, \hat u - u\rangle + \frac{p_1\langle q_1^{(p_1)}, u\rangle}{2(\lambda/\sigma)^2}}{\sqrt{\frac{p_1\langle q_1^{(p_1)}, u\rangle^2}{2(\lambda/\sigma)^4} + \frac{1 - \langle q_1^{(p_1)}, u\rangle^2}{(\lambda/\sigma)^2}}}, \frac{\langle q_2^{(p_2)}, \hat v - v\rangle + \frac{p_2\langle q_2^{(p_2)}, v\rangle}{2(\lambda/\sigma)^2}}{\sqrt{\frac{p_2\langle q_2^{(p_2)}, v\rangle^2}{2(\lambda/\sigma)^4} + \frac{1 - \langle q_2^{(p_2)}, v\rangle^2}{(\lambda/\sigma)^2}}}, \frac{\langle q_3^{(p_3)}, \hat w - w\rangle + \frac{p_3\langle q_3^{(p_3)}, w\rangle}{2(\lambda/\sigma)^2}}{\sqrt{\frac{p_3\langle q_3^{(p_3)}, w\rangle^2}{2(\lambda/\sigma)^4} + \frac{1 - \langle q_3^{(p_3)}, w\rangle^2}{(\lambda/\sigma)^2}}}\right)^\top.
	\end{split}
	\end{equation*}
	Then 
	\begin{equation}\label{ineq:q_inner_inference}
	T_{q_1^{(p_1)}, q_2^{(p_2)}, q_3^{(p_3)}} \stackrel{\rm d.}{\to} N(0, I_3) \quad \text{\rm as} \quad p \to \infty.
	\end{equation}
	Specifically, if $|u_i|, |v_j|, |w_k| \ll \min\{\lambda/(\sigma p), 1\}$ for some $i \in [p_1], j \in [p_2], k \in [p_3]$, then
	\begin{equation}\label{ineq:hat_u_asymp}
	\begin{pmatrix}
	\frac{\lambda}{\sigma}(\hat u_i-u_i), ~ 
	\frac{\lambda}{\sigma}(\hat v_j-v_j), ~ 
	\frac{\lambda}{\sigma}(\hat w_k-w_k)
	\end{pmatrix}^\top \stackrel{{\rm d.}}{\to} N(0, I_3)\quad {\rm as}\quad p\to\infty.
	\end{equation}
	If, furthermore, $\sigma/\lambda \ll |u_i|, |v_j|, |w_k| \ll \min\{\lambda/(\sigma p), 1/\sqrt{\log(p)}\}$, then
	\begin{equation}\label{ineq:entrywise_asymp}
	\frac{\hat{\calT}_{ijk} - \calT_{ijk}}{\sigma\sqrt{\hat u_i^2\hat v_j^2 + \hat v_j^2\hat w_k^2 + \hat w_k^2\hat u_i^2}} \stackrel{{\rm d.}}{\to} N(0, 1) \quad {\rm as}\quad p\to\infty.
	\end{equation}
\end{theorem}
Theorem \ref{thm:entry_inference} establishes the asymptotic distribution for any linear functional $q_1^\top \hat u, q_2^\top \hat v, q_3^\top \hat w$. Theorem \ref{thm:entry_inference} also implies that $[\hat{\calT}_{ijk} - z_{\alpha/2}\sigma\sqrt{\hat u_i^2\hat v_j^2 + \hat v_j^2\hat w_k^2 + \hat w_k^2\hat u_i^2}, \hat{\calT}_{ijk} + z_{\alpha/2}\sigma\sqrt{\hat u_i^2\hat v_j^2 + \hat v_j^2\hat w_k^2 + \hat w_k^2\hat u_i^2}]$ is an asymptotic $(1-\alpha)$ confidence interval for $T_{ijk}$ under some boundedness conditions of $|u_i|, |v_j|, |w_k|$. Here, the upper bound $|u_i|, |v_j|, |w_k| \ll \min\{\lambda/(\sigma p), 1/\sqrt{\log(p)}\}$ is significantly weaker than the incoherence condition commonly used in the matrix/tensor estimation/inference literature.  On the other hand, the lower bound condition, $|u_i|, |v_j|, |w_k| \gg \sigma/\lambda$, is essential to ensure the asymptotic normality of $\hat\calT$. To see this, consider a special case that $u_i=v_j=w_k=0$, then \eqref{ineq:hat_u_asymp} implies
$$\frac{\lambda^2\hat\calT_{ijk}}{\sigma^3} \overset{\rm d.}{\to} G_1 G_2 G_3 ~~ \text{as $p \to \infty$,} \quad  (G_1, G_2, G_3)^\top \sim N(0, I_3).$$
In other words, $\hat\calT_{ijk}$ satisfies a third moment Gaussian, not a Gaussian distribution. 

To cover the broader scenarios that the lower bound conditions are absent, we consider the following lower-thresholding procedure. Let $s(t) = \max\{t, \log(p)\sigma^2\hat\lambda^{-2}\}$\footnote{Here, $\log(p)$ can be replaced by any value that grows to infinity as $p$ grows.} for $t \geq 0$ and define the confidence interval for $\calT_{ijk}$ as
\begin{align}
\widetilde{CI}_{\alpha}(\hat \calT_{ijk})\notag := & \Big[\hat{\calT}_{ijk} - z_{\alpha/2}\sigma\sqrt{s(\hat u_i^2)s(\hat v_j^2)+s(\hat v_j^2)s(\hat w_k^2)+s(\hat w_k^2)s(\hat u_i^2)}, \\
& \hat{\calT}_{ijk} + z_{\alpha/2}\sigma\sqrt{s(\hat u_i^2)s(\hat v_j^2)+s(\hat v_j^2)s(\hat w_k^2)+s(\hat w_k^2)s(\hat u_i^2)}\Big].\label{eq:CI}
\end{align}
We can prove $\widetilde{CI}_{\alpha}(\hat \calT_{ijk})$ is a valid $(1-\alpha)$-level asymptotic confidence interval.
\begin{theorem}\label{th:CI}
	Suppose the conditions in Theorem \ref{thm:entry_inference} hold. If $\lambda/\sigma \gg p^{3/4}$ and $|u_i|, |v_j|, |w_k| \ll \min\{\lambda/(\sigma p), 1/\sqrt{\log(p)}\}$ for $i \in [p_1], j \in [p_2], k \in [p_3]$, then
	\begin{equation}\label{ineq:CI}
	\liminf_{p \to \infty}\PP\left(\calT_{ijk} \in \widetilde{CI}_{\alpha}(\hat \calT_{ijk})\right) \geq 1-\alpha.
	\end{equation}
\end{theorem}

\begin{remark}[Proof sketch of Theorem \ref{thm:entry_inference}]\label{rmk:proof-sketch}
	The proof scheme for Theorem \ref{thm:entry_inference} is essentially different from many recent literature on the entrywise inference \citep{cai2020uncertainty,chen2019inference,xia2019statistical} and we provide a proof sketch here. Without loss generality, we assume $\sigma = 1$ and $\langle u, \hat u\rangle, \langle v, \hat v\rangle, \langle w, \hat w\rangle \geq 0$. First, we can decompose $\langle \hat u, q_1\rangle$ into two terms:
	\begin{equation}\label{eq:hatu-decomposition}
	\langle q_1, \hat u \rangle = \langle \hat u, uu^\top q_1\rangle + \langle\hat u, (I - uu^\top) q_1\rangle = (q_1^\top u)\hat u^\top u + (U_{\perp}^\top q_1)^\top U_{\perp}^\top\hat{u}.
	\end{equation}
	Similar decompositions hold for $\langle \hat{v}, q_2\rangle$ and $\langle \hat{w}, q_3\rangle$. For any $O_i \in \OO_{p_i-1}$, we construct three rotation matrices as
	$$\tilde{O}_1 = uu^\top + U_{\perp}O_1U_{\perp}^\top \in \OO_{p_1}, \quad \tilde{O}_2 = vv^\top + V_{\perp}O_2V_{\perp}^\top \in \OO_{p_2}, \quad \tilde{O}_3 = ww^\top + W_{\perp}O_3W_{\perp}^\top \in \OO_{p_3},$$ 
	where $U_{\perp} \in \mathbb{O}_{p_1, p_1-1}, V_{\perp} \in \mathbb{O}_{p_2, p_2-1}, W_{\perp} \in \mathbb{O}_{p_3, p_3-1}$ are the orthogonal complement of $u, v, w$, respectively. A key observation is that $\tilde{\calA} = \calA \times_1 \tilde{O}_1^\top \times_2 \tilde{O}_2^\top \times_3 \tilde{O}_3^\top$ and $\calA$ share the same distribution. Suppose $\tilde{u}, \tilde{v}, \tilde{w}$ are the outputs of Algorithm \ref{algo:hooi}. Then we have $\tilde{u} = \tilde{O}_1^\top\hat{u}, \tilde{v} = \tilde{O}_2^\top\hat{v}, \tilde{w} = \tilde{O}_3^\top\hat{w}$ and can further prove that given $\langle u, \hat{u}\rangle, \langle v, \hat{v}\rangle$ and $\langle w, \hat{w}\rangle$, $\bigg(\frac{\hat{u}^\top U_{\perp}}{\|U_{\perp}^\top\hat{u}\|_2}, \frac{\hat{v}^\top V_{\perp}}{\|V_{\perp}^\top\hat{v}\|_2}, \frac{\hat{w}^\top W_{\perp}}{\|W_{\perp}^\top\hat{w}\|_2}\bigg)$ and $\bigg(\frac{\hat{u}^\top U_{\perp}}{\|U_{\perp}^\top\hat{u}\|_2}O_1, \frac{\hat{v}^\top V_{\perp}}{\|V_{\perp}^\top\hat{v}\|_2}O_2, \frac{\hat{w}^\top W_{\perp}}{\|W_{\perp}^\top\hat{w}\|_2}O_3\bigg)$ have the same distribution. By the uniqueness of the Haar measure \citep{neumann1936uniqueness,weil1940integration} and Theorem \ref{thm:na_tsvd}, we can further prove for any fixed vectors $f_1 \in \mathbb{S}^{p_1-2}, f_2 \in \mathbb{S}^{p_2-2}, f_3 \in \mathbb{S}^{p_3-2}$, we have
	\begin{align*}
	&\bigg(\lambda\hat{u}^\top U_{\perp}f_1, ~ \lambda\hat{v}^\top V_{\perp}f_2, ~ \lambda\hat{w}^\top W_{\perp}f_3,\\
	& \quad \frac{\hat{u}^\top u - (1 - p_1\lambda^{-2}/2)}{\sqrt{p_1/2}\lambda^{-2}},~ \frac{\hat v^\top v - (1 - p_2\lambda^{-2}/2)}{\sqrt{p_2/2}\lambda^{-2}},~ \frac{\hat w^\top w - (1 - p_3\lambda^{-2}/2)}{\sqrt{p_3/2}\lambda^{-2}}\bigg)^\top \stackrel{\rm d.}{\to} N(0, I_6).
	\end{align*}
	This inequality and \eqref{eq:hatu-decomposition} result in \eqref{ineq:q_inner_inference}\eqref{ineq:hat_u_asymp}\eqref{ineq:entrywise_asymp}.
\end{remark}
\begin{remark}
	The entrywise inference for Tucker low-rank or orthogonal decomposable tensor PCA can be significantly more challenging due to the dependence among different factors. We leave it as future research.
\end{remark}

%%%%%%%%%%%%%%%%%%%%%%
\section{Numerical Simulations}\label{sec:numerical-simulation}
%%%%%%%%%%%%%%%%%%%%%%

We now conduct numerical studies to support our theoretical findings in previous sections. Each experiment is repeated for $2000$ times, from which we obtain $2000$ realizations of the respective statistics. Then we draw histograms or boxplots, and compare with the corresponding baselines. In each histogram, the red line is the density of the standard normal distribution.

\begin{figure}[h!]
	\centering
	\subfigure{\includegraphics[height=40mm,width=68mm]{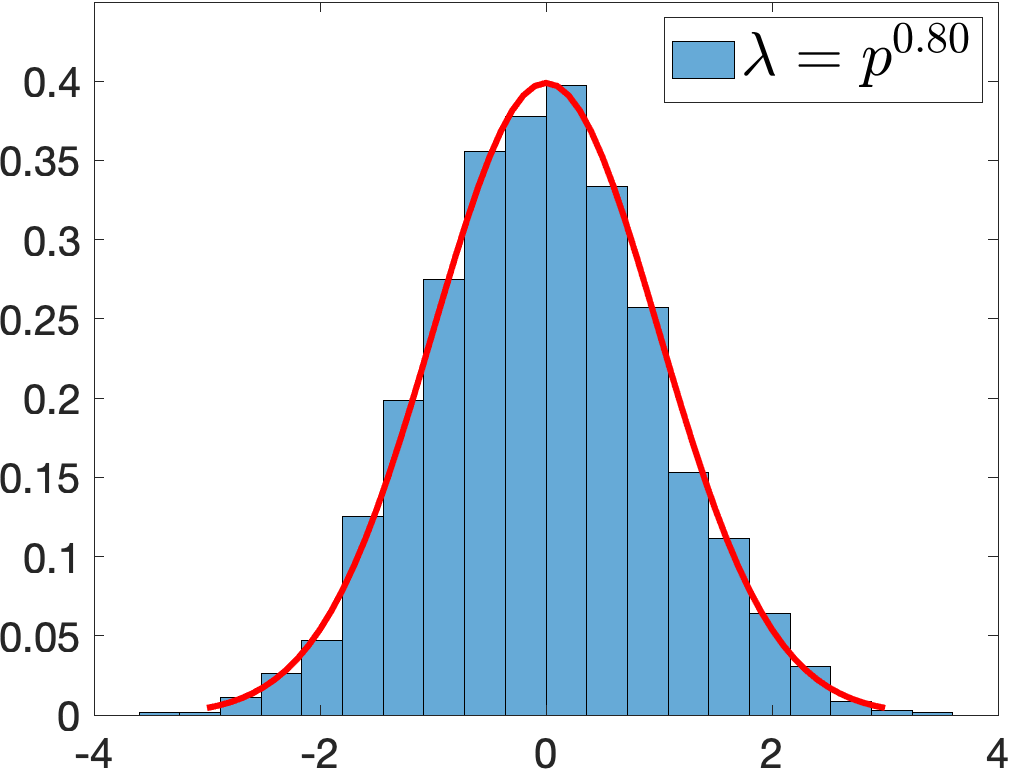}}
	\hspace{0.3cm}
	\subfigure{\includegraphics[height=40mm,width=68mm]{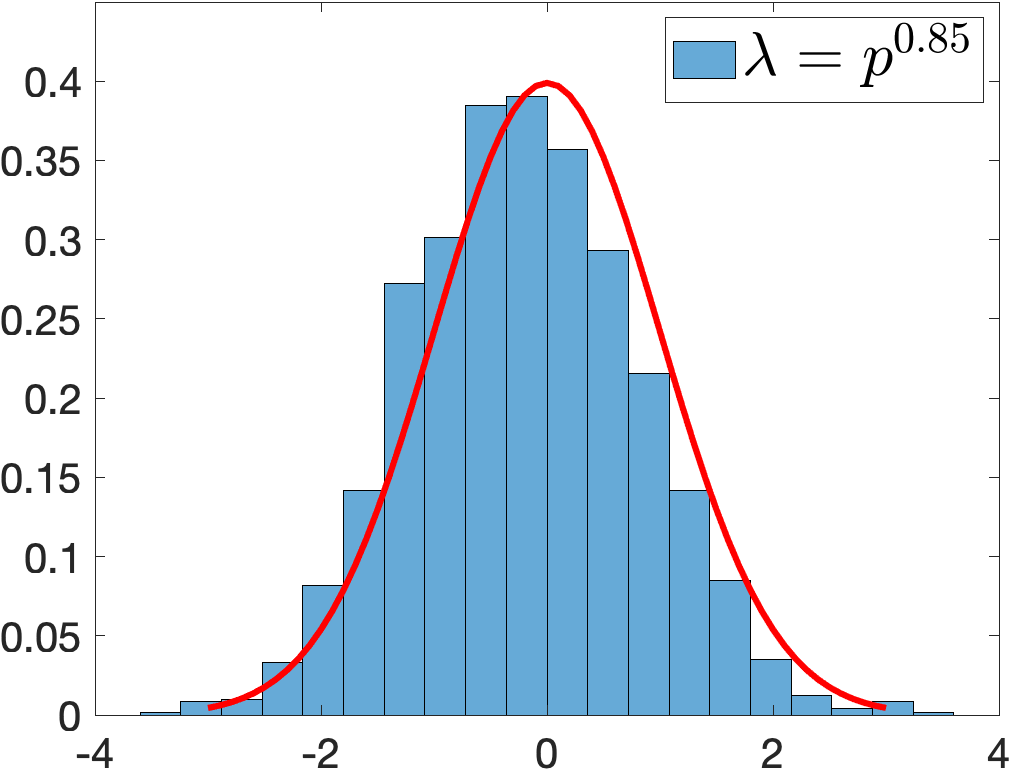}}
	\subfigure{\includegraphics[height=40mm,width=68mm]{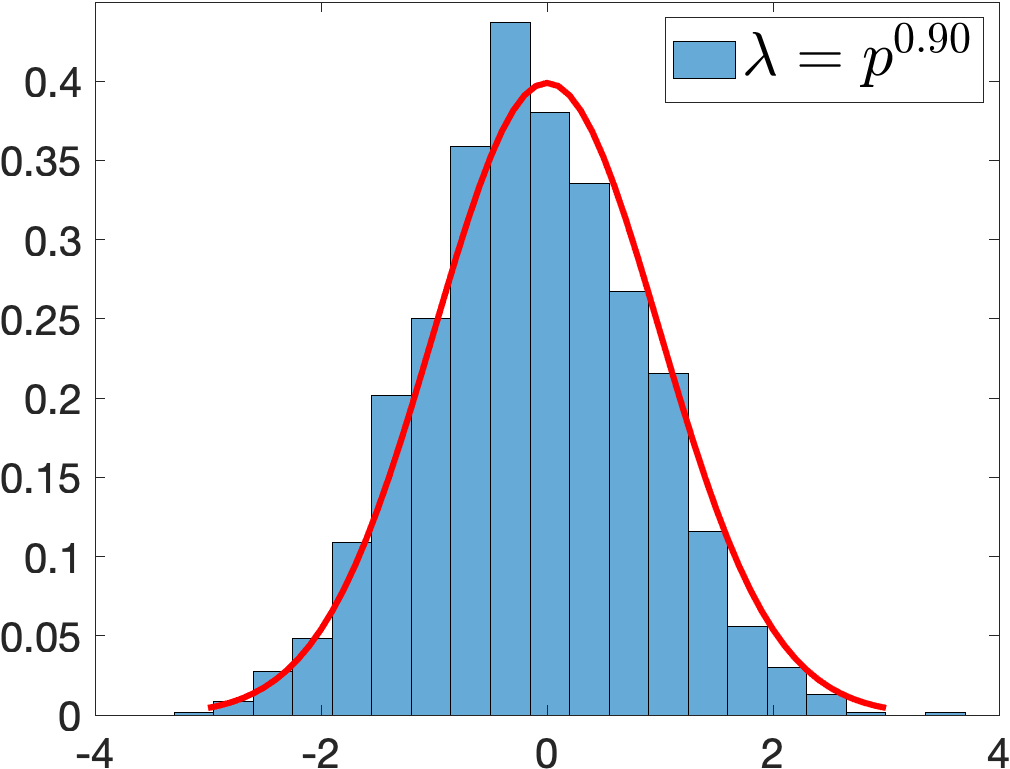}}
	\hspace{0.3cm}
	\subfigure{\includegraphics[height=40mm,width=68mm]{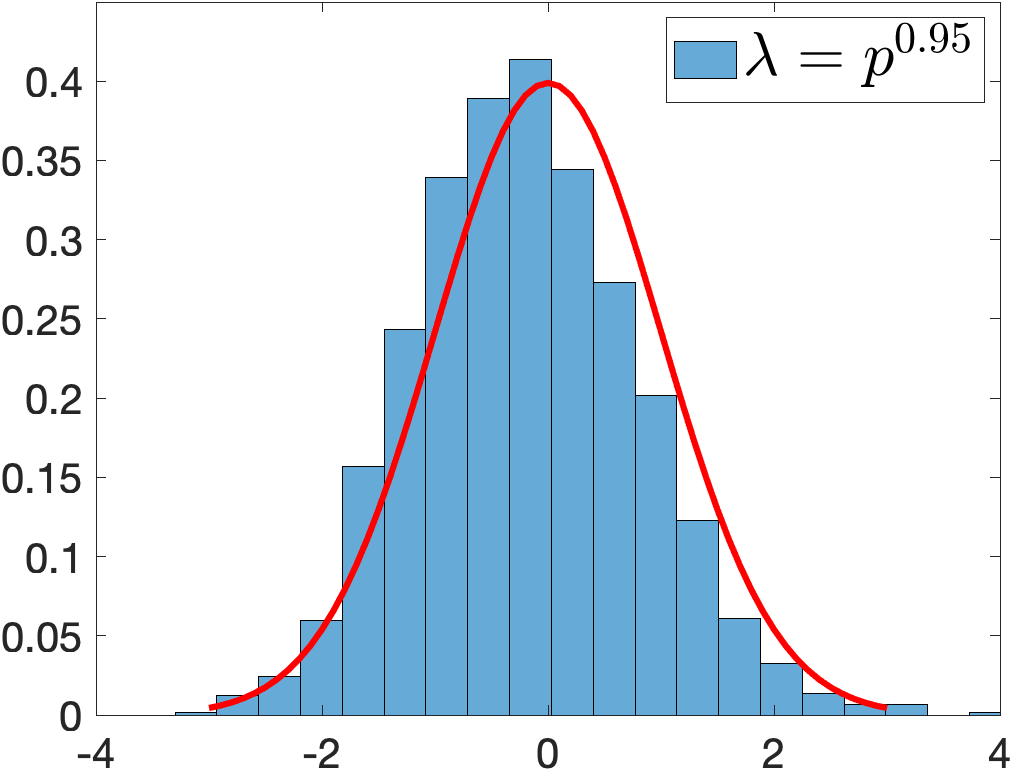}}
	\caption{Normal approximation of $\frac{\|\sin\Theta(\hat U_1, U_1)\|_{\rm F}^2-p\|\Lambda_1^{-1}\|_{\rm F}^2}{\sqrt{2p}\|\Lambda_1^{-2}\|_{\rm F}}$ for order-3 Tucker low-rank tensor PCA model (\ref{eq:PCA_model}). Here, $p_1=p_2=p_3=p=200$, $r=3$, $\sigma=1$.}
	\label{fig:na2}
\end{figure}

\begin{figure}[h!]
	\centering
	\subfigure{\includegraphics[height=40mm,width=68mm]{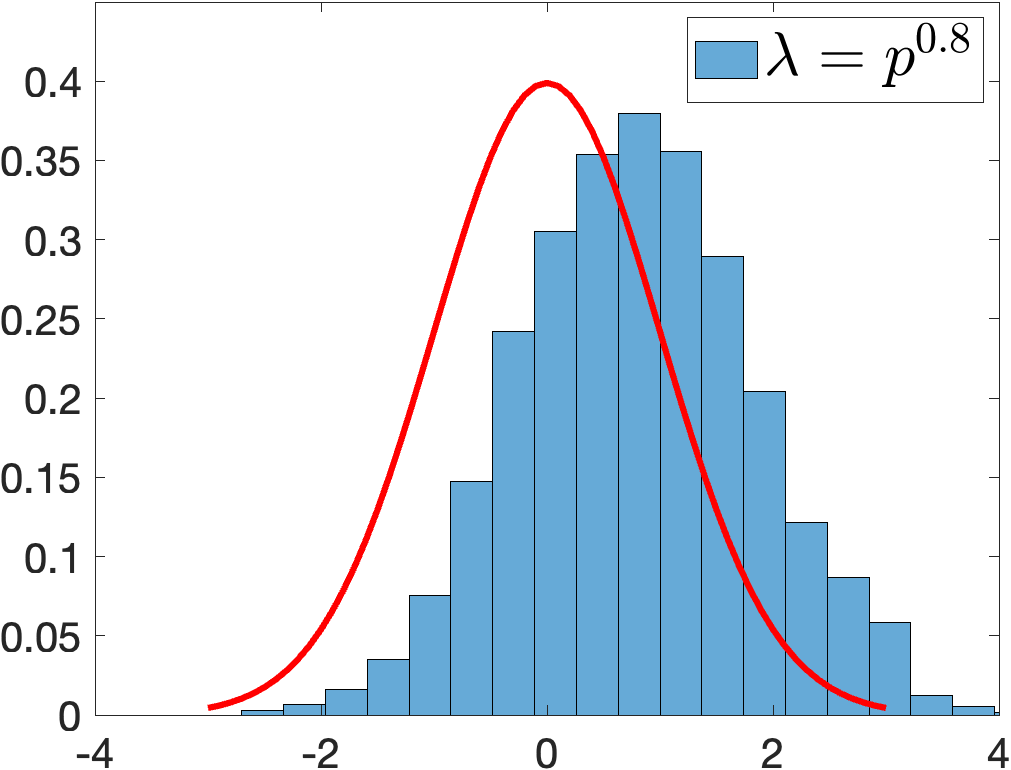}}
	\hspace{0.3cm}
	\subfigure{\includegraphics[height=40mm,width=68mm]{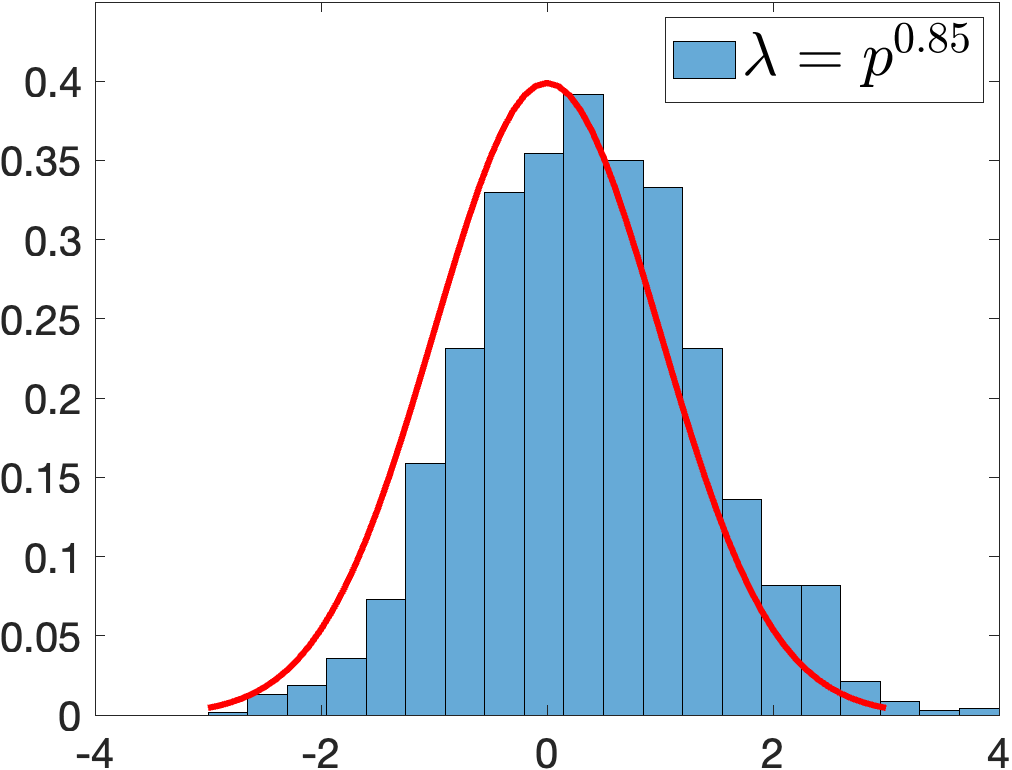}}
	\subfigure{\includegraphics[height=40mm,width=68mm]{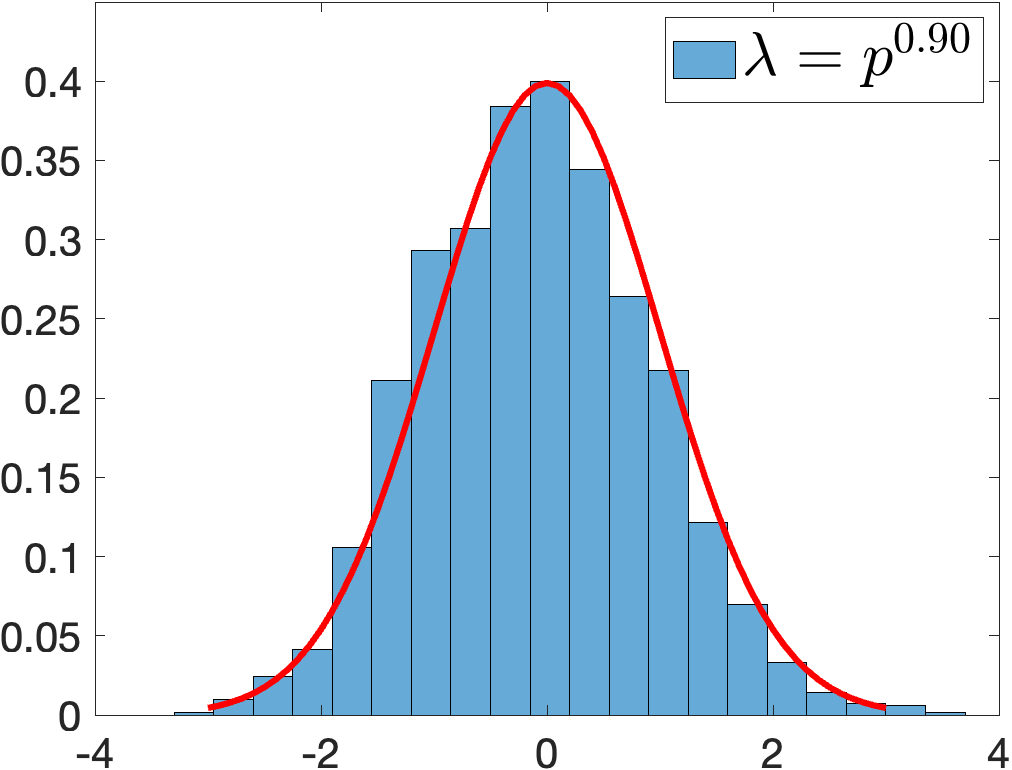}}
	\hspace{0.3cm}
	\subfigure{\includegraphics[height=40mm,width=68mm]{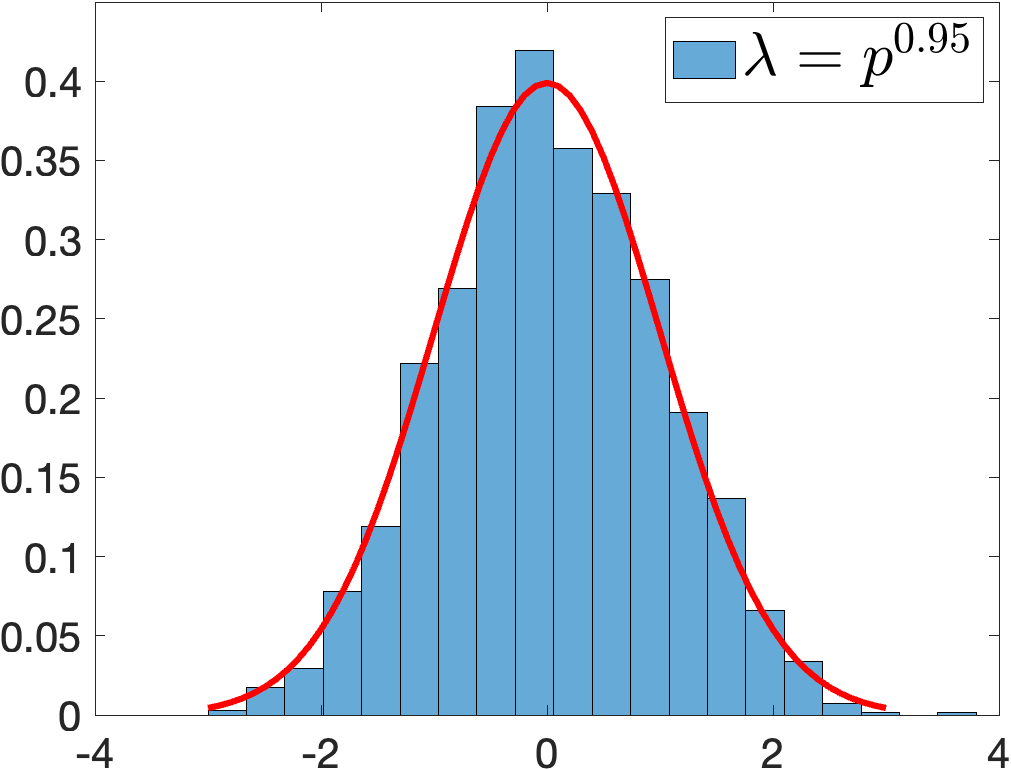}}
	\caption{Normal approximation of $\frac{\|\sin\Theta(\hat U_1, U_1)\|_{\rm F}^2-p\hat{\sigma}^2\|\hat\Lambda_1^{-1}\|_{\rm F}^2}{\sqrt{2p}\hat{\sigma}^2\|\hat\Lambda_1^{-2}\|_{\rm F}}$ for order-3 Tucker low-rank tensor PCA model (\ref{eq:PCA_model}). Here, $p_1=p_2=p_3=p=200$, $r=3$, $\sigma=1$.}
	\label{fig:a2}
\end{figure}

We begin with the inference for principal components of Tucker low-rank tensors. Specifically, we randomly draw $\check{U}_j\in\mathbb{R}^{p_j\times r_j}$ with i.i.d. standard normal entries and normalize to $U_j = \text{QR}(\check{U}_j)$. We then draw core tensor $\check{\calG}\in \mathbb{R}^{r\times r\times r}$ with i.i.d. standard normal entries and rescale to $\calG = \check{\calG}\cdot p^{\gamma}/ (\lambda_{\submin}(\check\calG))$. Consequently, $U_j$ is uniform randomly selected from $\mathbb{O}_{p_j, r_j}$ and $\lambda_{\submin}(\calG)=\lambda=p^{\gamma}$. For $p_1=p_2=p_3=200$, $r=3$, and $\sigma=1$, each value of $\gamma\in \{0.80, 0.85, 0.90, 0.95\}$, we observe $\calA$ under tensor PCA model (\ref{eq:PCA_model}) and apply Algorithm~\ref{algo:am_optimal_PCA} to obtain realizations of 
$$
T_1= \frac{\|\sin\Theta(\hat U_1, U_1)\|_{\rm F}^2-p\|\Lambda_1^{-1}\|_{\rm F}^2}{\sqrt{2p_1}\|\Lambda_1^{-2}\|_{\rm F}} \quad \text{and} \quad T_2=\frac{\|\sin\Theta(\hat U_1, U_1)\|_{\rm F}^2-p\hat{\sigma}^2\|\hat\Lambda_1^{-1}\|_{\rm F}^2}{\sqrt{2p}\hat{\sigma}^2\|\hat\Lambda_1^{-2}\|_{\rm F}}. 
$$
We repeat this procedure for $2000$ times, from which we obtain $2000$ realizations of the respective statistics and plot the density histograms in Figures~\ref{fig:na2} and \ref{fig:a2}, respectively. We can see $T_1$ and $T_2$ both achieve good normal approximation in these settings. 

\begin{figure}[h!]
	\centering
	\subfigure{\includegraphics[height=40mm,width=68mm]{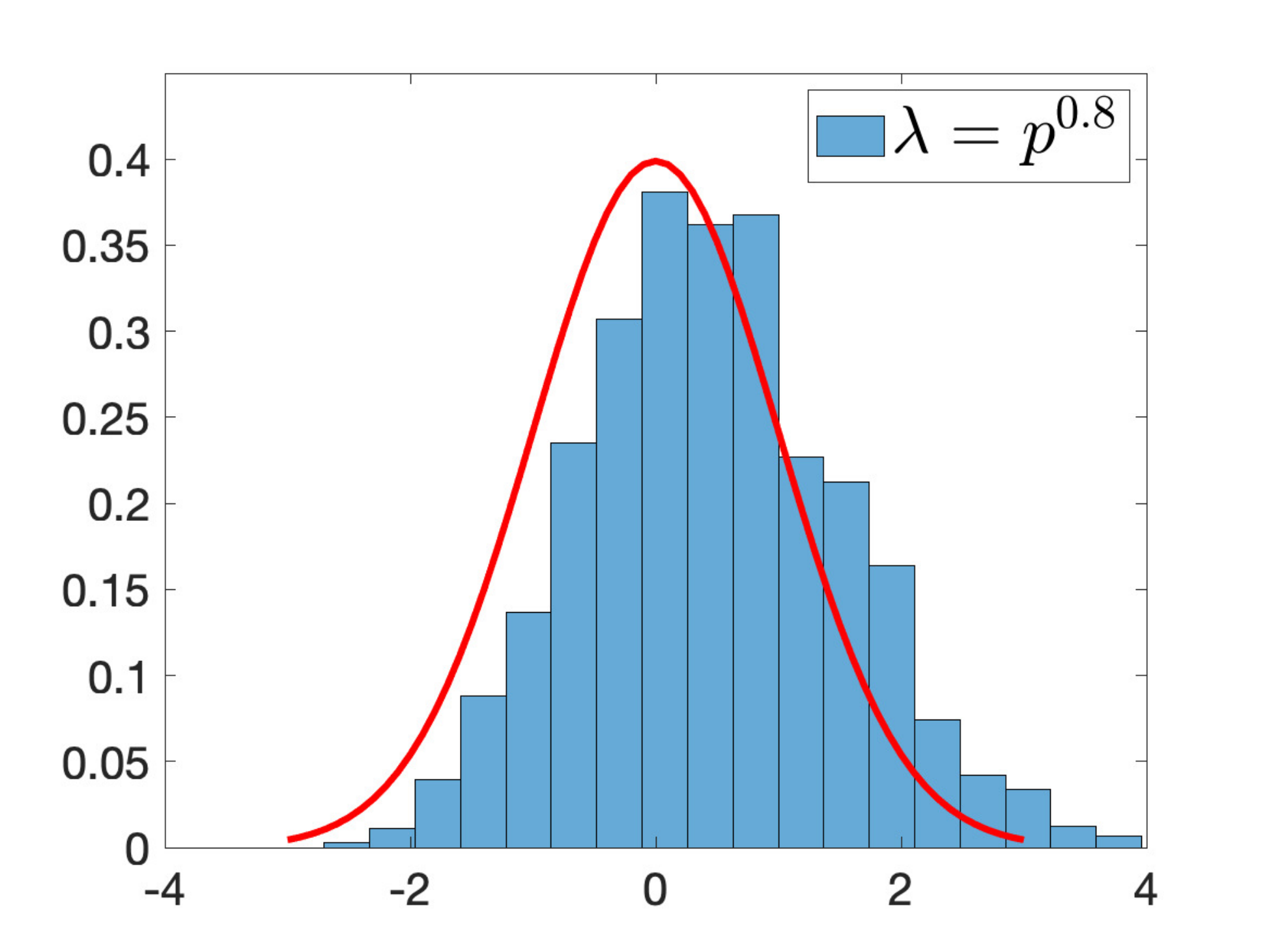}}
	\hspace{0.3cm}
	\subfigure{\includegraphics[height=40mm,width=68mm]{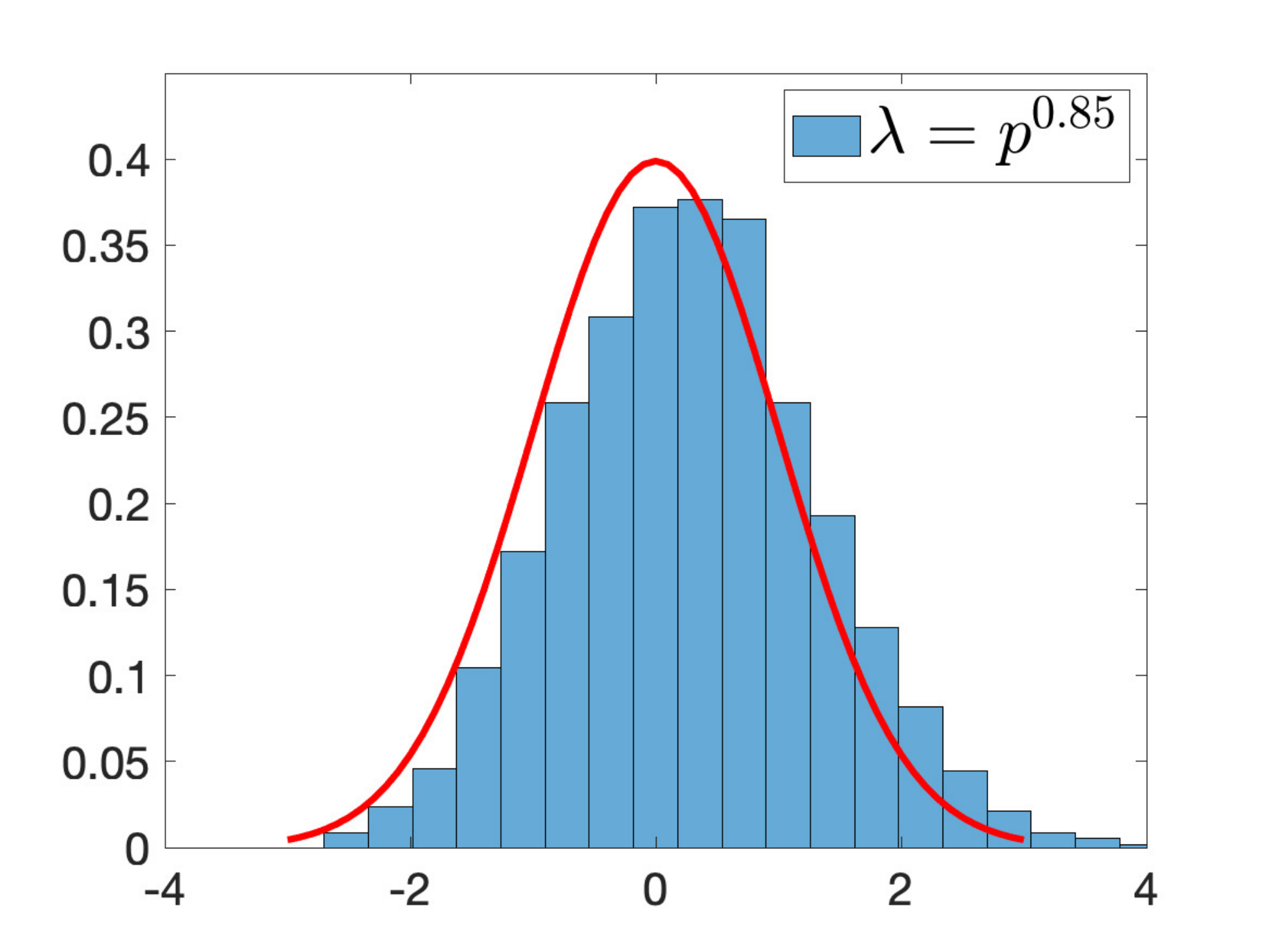}}
	\subfigure{\includegraphics[height=40mm,width=68mm]{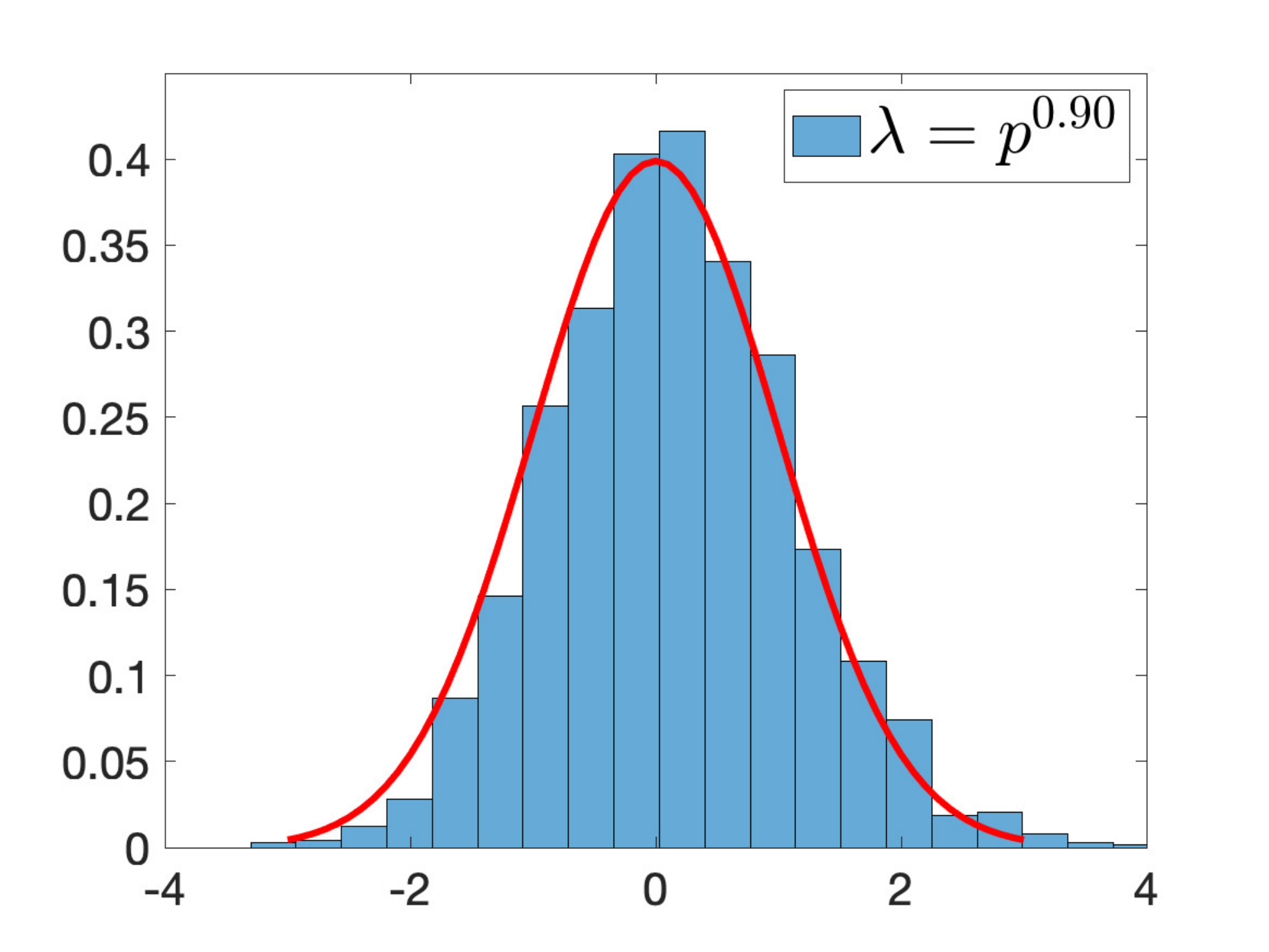}}
	\hspace{0.3cm}
	\subfigure{\includegraphics[height=40mm,width=68mm]{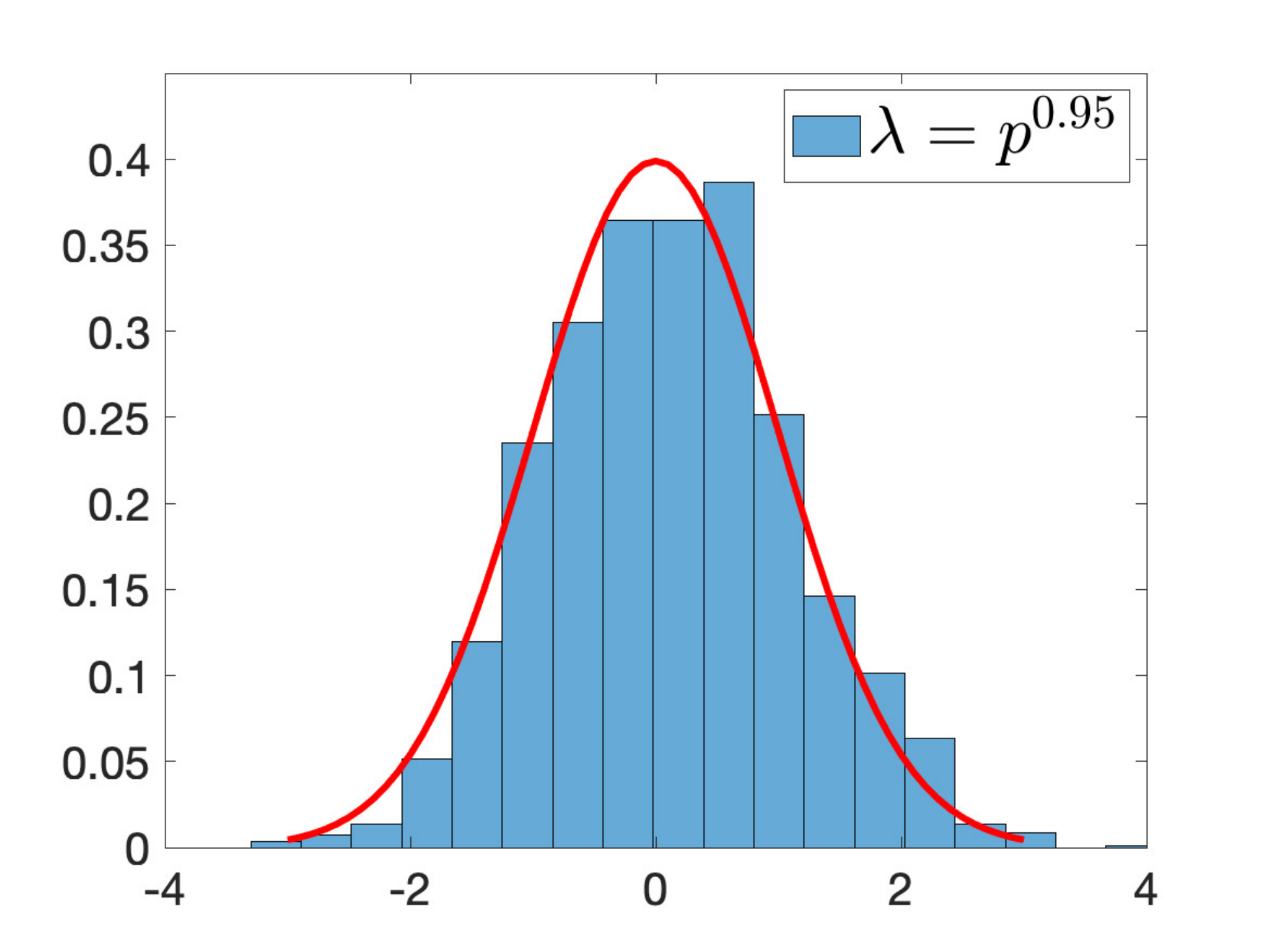}}
	\caption{Normal approximation of $\frac{\langle \hat u_3, u_3\rangle^2-(1-p\lambda^{-2})}{\sqrt{2p}\lambda^{-2}}$ for tensor PCA model (\ref{eq:PCA_model}) when $\calT$ is a third-order orthogonally decomposable tensor and $\sigma=1$. Here, $p_1=p_2=p_3=p=200, r=3, \lambda_{\submin}=\lambda$.}
	\label{fig:na1}
\end{figure}

We then consider the asymptotic normality in orthogonally decomposable tensors under the tensor PCA model. Similarly, we fix $p=200$, $r=3$, and construct the orthogonally decomposable tensor as $\calT=\sum_{i=1}^r (r+1-i)\lambda\cdot (u_i\otimes v_i\otimes w_i)$, where $[u_1,\ldots, u_r], [v_1,\ldots,  v_r], [w_1,\ldots, w_r]$ are drawn uniform randomly from $\mathbb{O}_{p, r}$ similarly to the previous setting and $\lambda=p^{\gamma}$ with $\gamma=0.80, 0.85, 0.90, 0.95$. For each $\gamma$, we obtain $2000$ replicates of $T = \frac{\langle \hat u_3, u_3\rangle^2-(1-p\lambda^{-2})}{\sqrt{2p}\lambda^{-2}}$, draw the density histogram, and plot the results in Figure~\ref{fig:na1}. We can see the normal approximation of $T$ becomes more accurate as the signal strength $\lambda$ grows. 

\begin{figure}[h!]
	\centering
	\subfigure{\includegraphics[height=40mm,width=68mm]{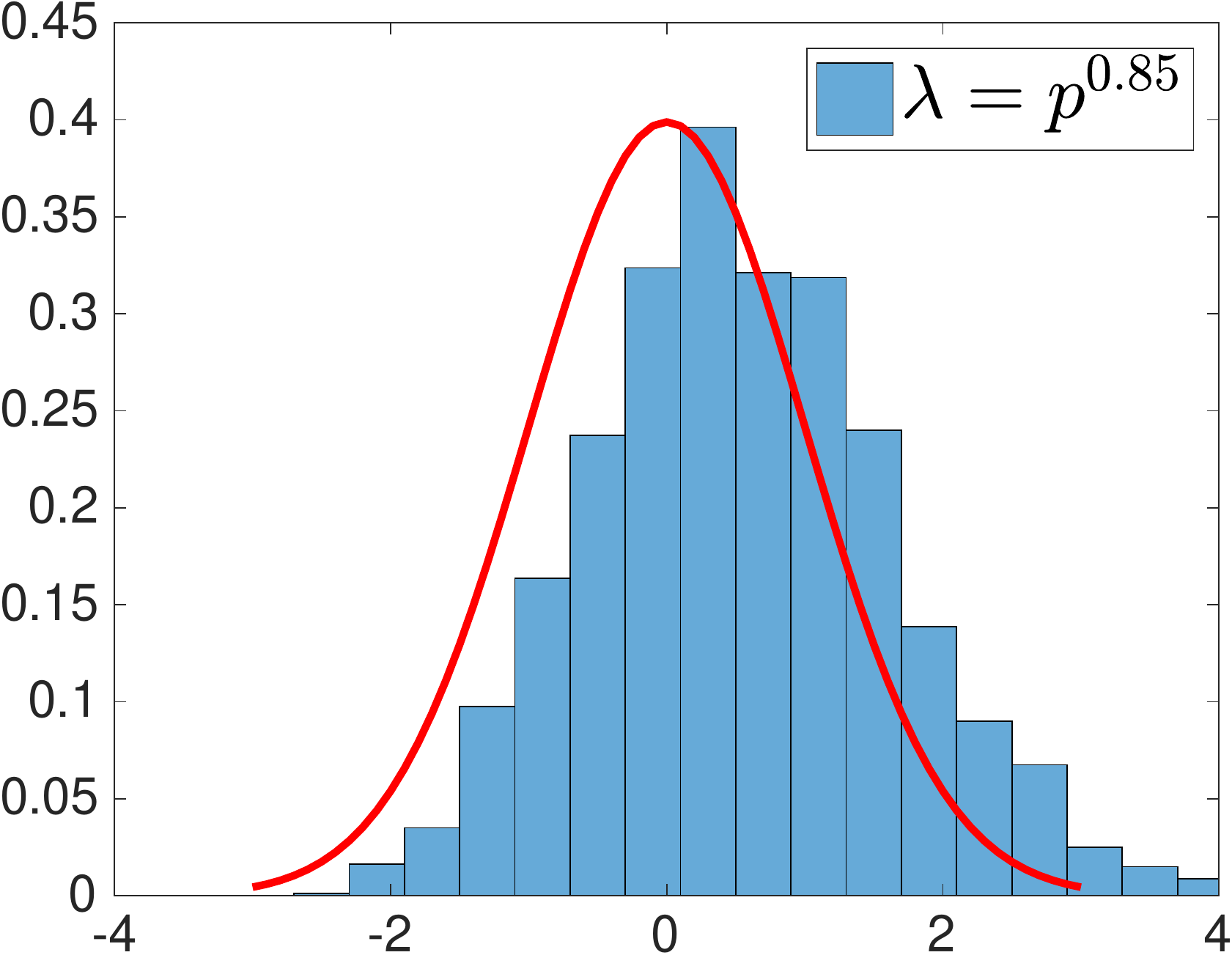}}
	\hspace{0.3cm}
	\subfigure{\includegraphics[height=40mm,width=68mm]{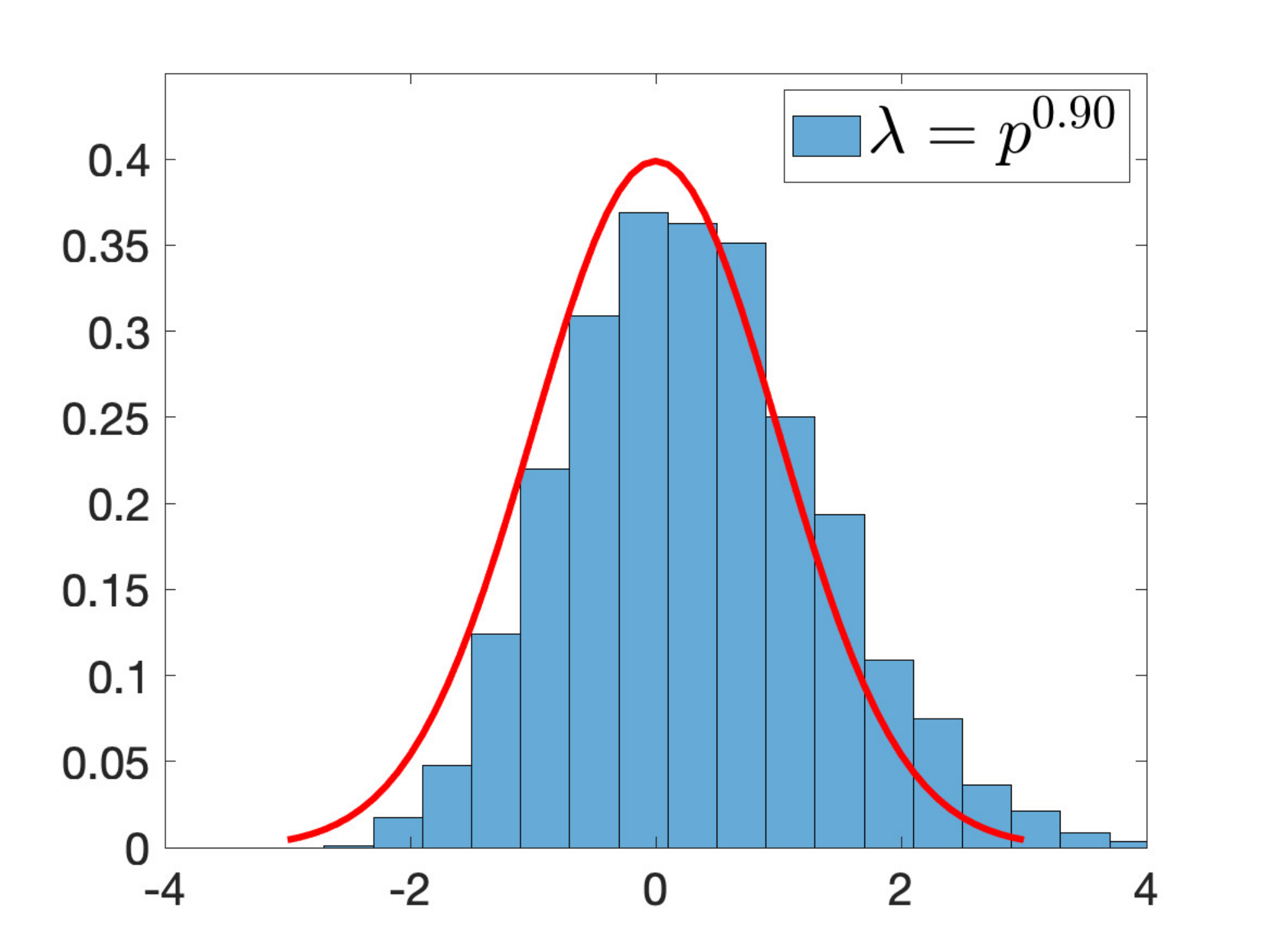}}
	\subfigure{\includegraphics[height=40mm,width=68mm]{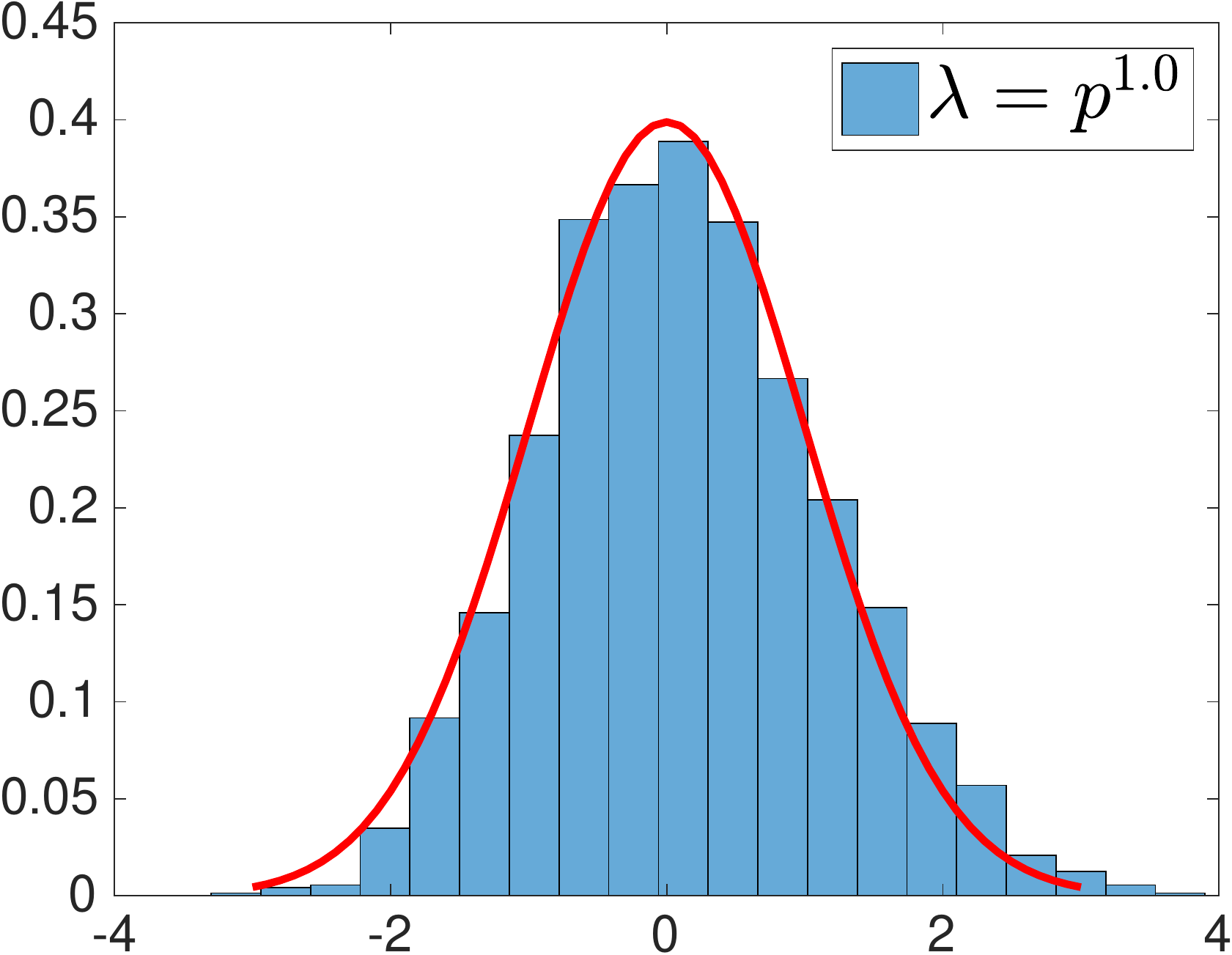}}
	\hspace{0.3cm}
	\subfigure{\includegraphics[height=40mm,width=68mm]{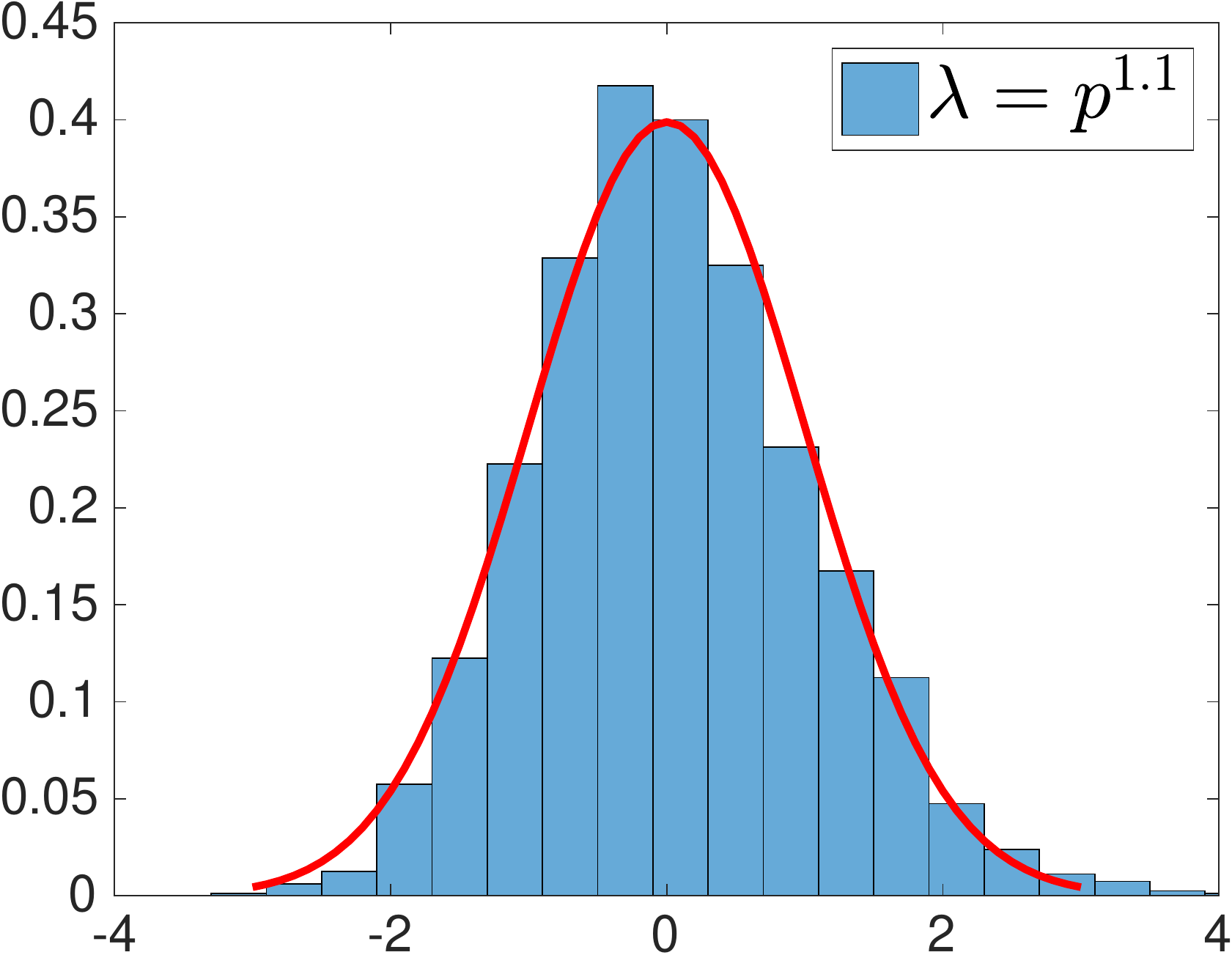}}
	\caption{Normal approximation of $\frac{\langle \hat u_1, u_1\rangle^2-(1-p\lambda^{-2})}{\sqrt{2p}\lambda^{-2}}$ for tensor PCA model (\ref{eq:PCA_model}) when $\calT=\lambda\cdot (u_1\otimes v_1\otimes w_1\otimes q_1)$ is a fourth-order tensor and $\sigma=1$. Here, $p_1=p_2=p_3=p_4=p=100, r=1$ and $\lambda_{\submin}=\lambda$.  
	}
	\label{fig:na3}
\end{figure}

Though the focus of this paper is on third-order tensors, we will explain later in Section~\ref{sec:diss} that the results can be generalized to higher-order ones. Next, we conduct simulation study on tensor PCA model for fourth-order orthogonally decomposable tensors when $p=100$ and $r=1$. With a few modifications on the proof, we can show $\big(\langle \hat u_1, u_1\rangle^2-(1-p\lambda^{-2})\big)(\sqrt{2p}\lambda^{-2})^{-1}$ is asymptotically normal under the required SNR assumption for efficient computation: SNR $\geq Cp$. The simulation results in Figure~\ref{fig:na3} show that equipped with a warm initialization, the two-iteration alternating minimization yields an estimator achieving good normal approximation even if SNR $\approx p^{0.9}$, which is strictly weaker than the required SNR assumption for efficient computation. See more discussions in Section~\ref{sec:diss}. 

\begin{figure}[h!]
	\centering
	\subfigure{\includegraphics[height=40mm,width=68mm]{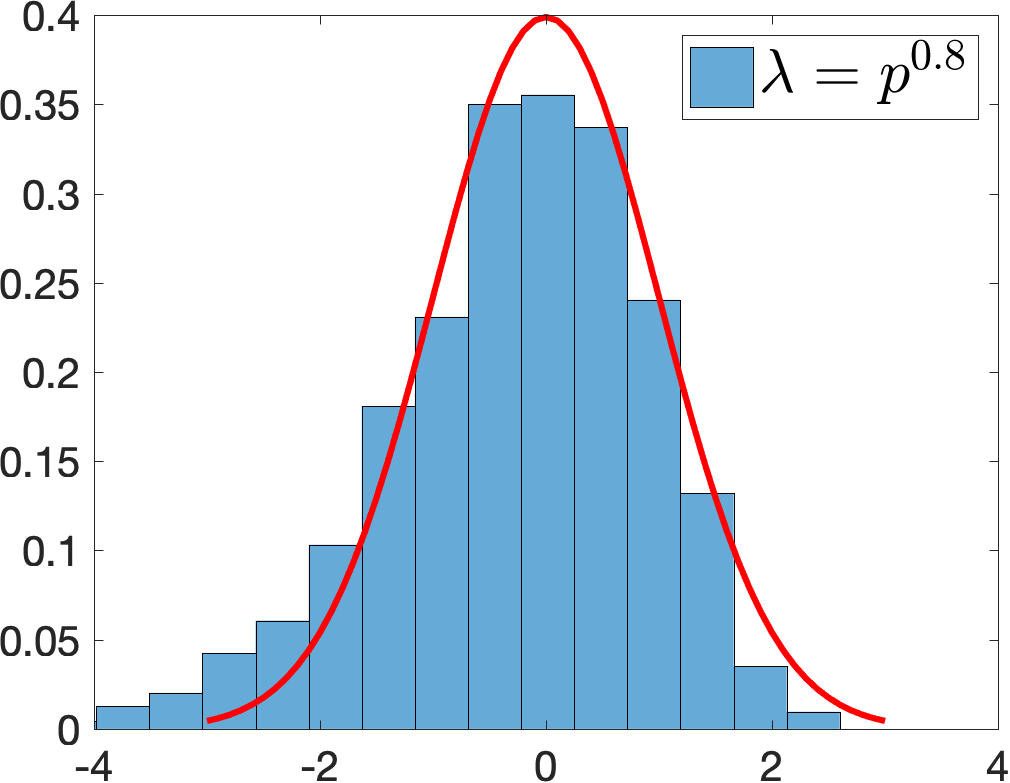}}
	\hspace{0.3cm}
	\subfigure{\includegraphics[height=40mm,width=68mm]{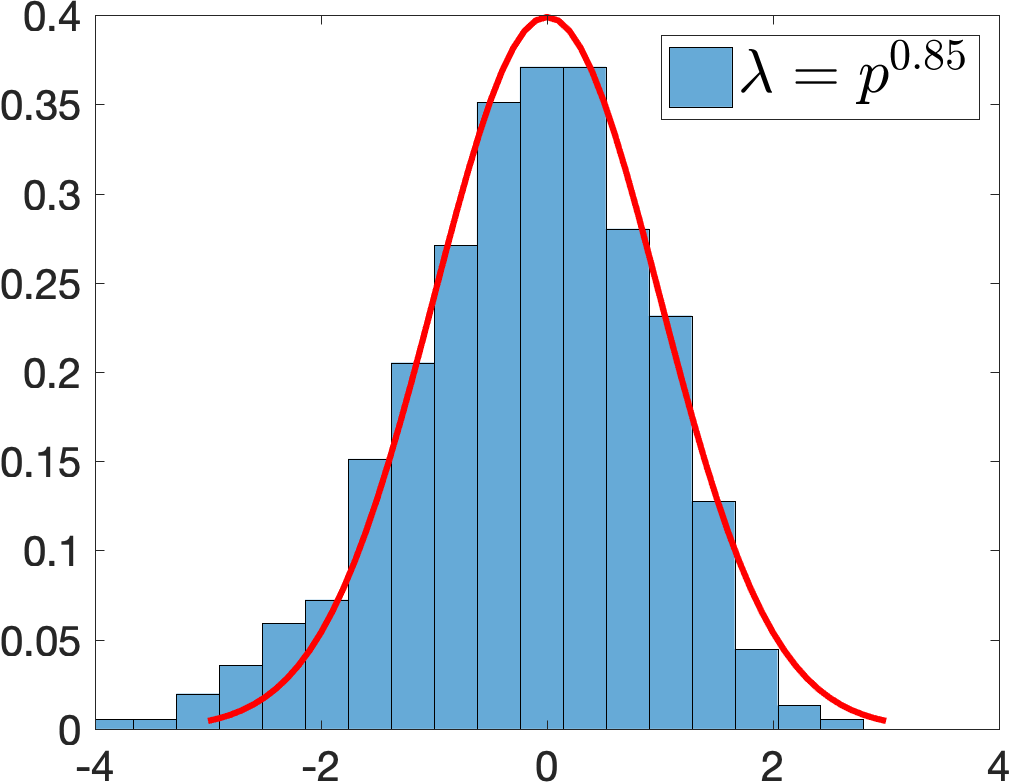}}
	\subfigure{\includegraphics[height=40mm,width=68mm]{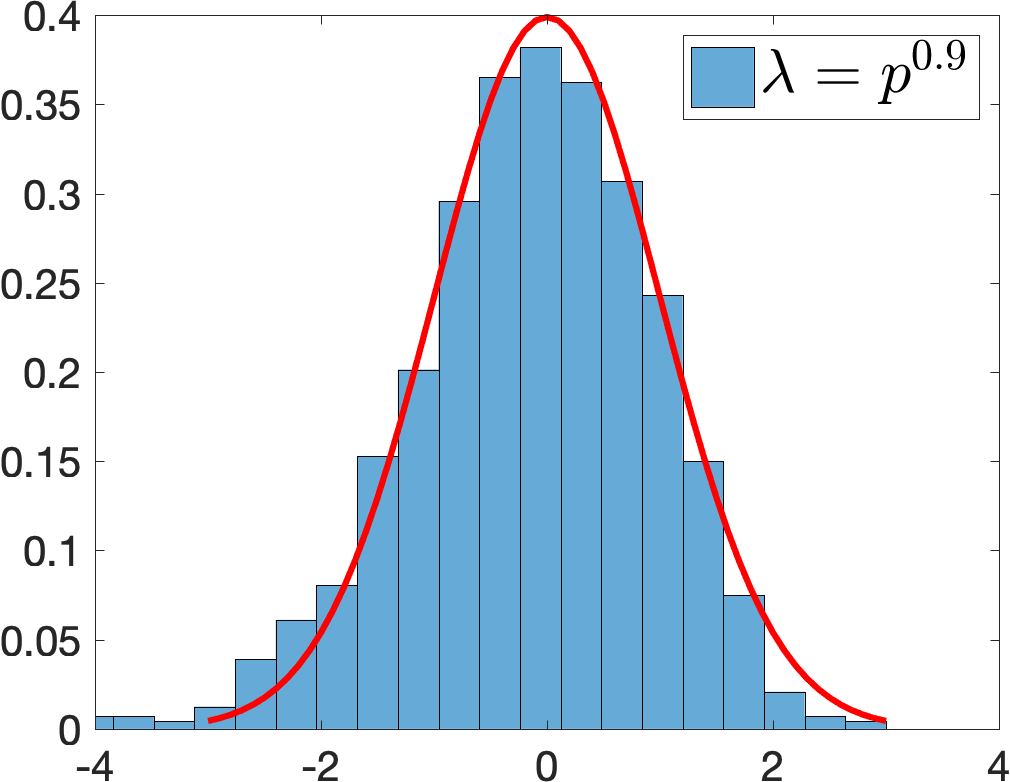}}
	\hspace{0.3cm}
	\subfigure{\includegraphics[height=40mm,width=68mm]{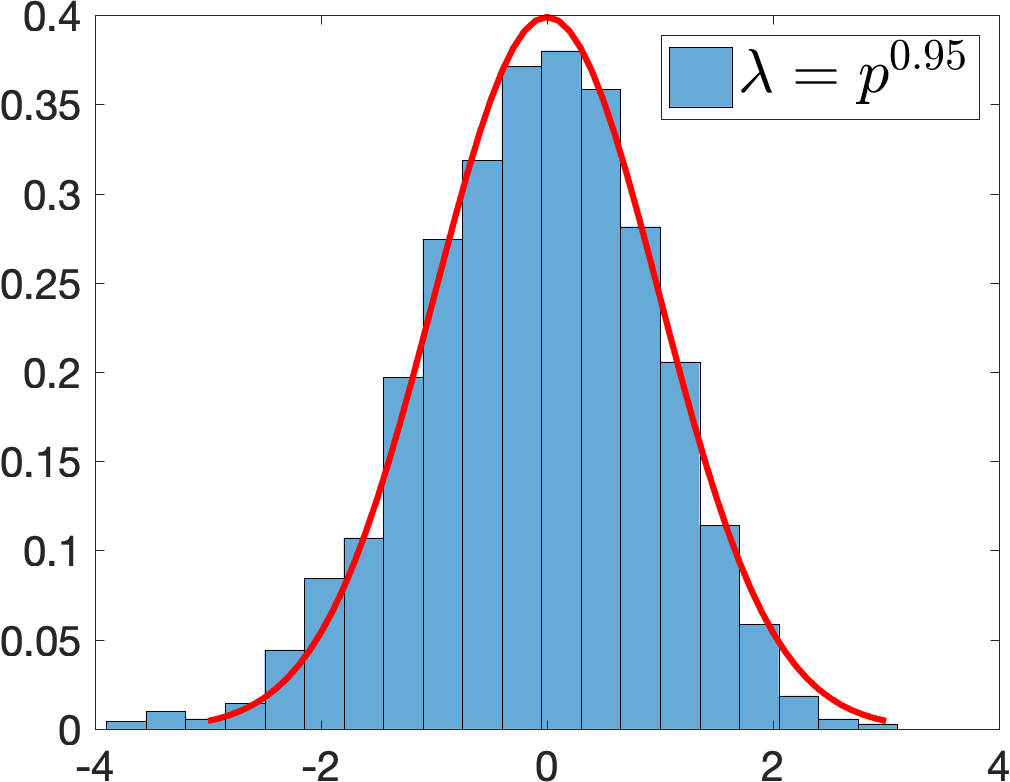}}
	\caption{Normal approximation of $\frac{\hat{\calT}_{1,1,1} - \calT_{1,1,1}}{\sqrt{\hat u_1^2\hat v_1^2 + \hat v_1^2\hat w_1^2 + \hat w_1^2\hat u_1^2}}$ for tensor PCA model (\ref{eq:PCA_model}) when $\calT$ is a rank-$1$ tensor and $\sigma=1$. The parameters are $p_1=p_2=p_3=p=200$ with signal strength $\lambda$.}
	\label{fig:entrywise}
\end{figure}

Then, we consider the entrywise inference under the rank-1 tensor PCA model. We construct $\calT = \lambda\cdot u \otimes v \otimes w \in \mathbb{R}^{p\times p\times p}$, where $u = v = w = (1/\sqrt{p}, \dots, 1/\sqrt{p})^\top$ and $\lambda = p^{\gamma}$ with $\gamma\in \{0.80, 0.85, 0.90, 0.95\}$. For each value of $\gamma$, we draw a random observation $\calA$ under the tensor PCA model (\ref{eq:PCA_model}) and apply Algorithm \ref{algo:hooi} with $t_{\submax}=10$. We present the histogram  in Figure \ref{fig:entrywise} based on 2000 replicate values of $\frac{\hat{\calT}_{1,1,1} - \calT_{1,1,1}}{\sqrt{\hat u_1^2\hat v_1^2 + \hat v_1^2\hat w_1^2 + \hat w_1^2\hat u_1^2}}$. The simulation results validate the asymptotic normality of $\frac{\hat{\calT}_{ijk} - \calT_{ijk}}{\sqrt{\hat u_i^2\hat v_j^2 + \hat v_j^2\hat w_k^2 + \hat w_k^2\hat u_i^2}}$ when $u, v, w$ have balanced entry values, which are in line with the theory in Theorem \ref{thm:entry_inference}. 

\begin{figure}[h!]
	\centering
	\includegraphics[height=80mm]{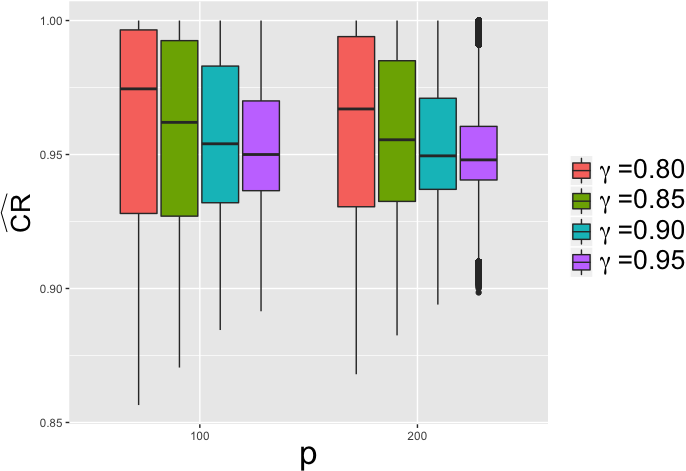}
	\caption{Boxplots for empirical coverage of entrywise confidence interval $\widehat{\text{CR}}_{ijk}$}
	\label{fig:boxplot}
\end{figure}

Finally, we consider the accuracy of the asymptotic entrywise confidence interval proposed in \eqref{eq:CI} under the tensor PCA model. Let $\calT = \lambda\cdot u\otimes v\otimes w$ be a rank-$1$ tensor, where $u, v, w$ are uniform randomly drawn from $\mathbb{S}^{p-1}$ for $p\in\{100, 200\}$ and $\lambda = p^{\gamma}$ for $\gamma\in\{0.80, 0.85, 0.90, 0.95\}$. For each combination of $(p, \gamma)$, we report the empirical coverage rates for the $0.95$-confidence interval $\widehat{\text{CR}}_{ijk}$ by boxplots in Figure \ref{fig:boxplot}. The results show the empirical coverage rates are close to 0.95 in all settings and larger values of $(\gamma, p)$ lead to more accurate coverage.

%%%%%%%%%%%%%%%%%%%%%%%%%%%%
\section{Discussion}\label{sec:diss}
%%%%%%%%%%%%%%%%%%%%%%%%%%%%

In this paper, we investigate the inference for low-rank tensors under two basic and fundamentally important tensor models: tensor PCA and regression. Based on an initial estimator achieving a reasonable estimation error, we propose to update by a two-iteration alternating minimization algorithm then establish the asymptotic distribution for the singular subspace outcomes. Distributions of general linear forms of the singular vectors are also established for rank-one tensor PCA model, which further enables the entrywise inference on the parameter tensor. 

Although our main focus is on third-order tensors, the results in this paper can be extended to higher-order tensors. For example, suppose $m\geq4$ and $\calT=\sum_{j=1}^r \lambda_j \cdot u_j^{(1)} \otimes\cdots\otimes u_j^{(m)}$ is orthogonally decomposable. Given $\calA$ from the tensor PCA model (\ref{eq:PCA_model}) and Assumption~\ref{assump:T-orth} holds, we can refine by two power iterations similarly to Algorithm \ref{algo:power_iter}, then obtain $\{\hat u_j^{(1)}, \hat u_j^{(2)}, \cdots, \hat u_j^{(m)}\}_{j=1}^r$. Similarly to Theorem \ref{thm:T-orth}, we can prove 
$$\frac{\langle u_j^{(k)}, \hat u_j^{(k)}\rangle^2 - (1-p_k\lambda_j^{-2}\sigma^2)}{\sqrt{2p_k}\lambda_j^{-2}\sigma^2} \overset{d.}{\to} N(0, 1), \quad k=1,\ldots, m,$$ 
if $\lambda_{\submin}/\sigma\gg p^{3/4}$ and other regularity conditions holds. %when $r=O(1)$ and $p_k\asymp p$ for all $k=1,\cdots,m$. 
If $m\geq 4$, the SNR condition $\lambda_{\submin}\gg p^{3/4}$ is weaker than the condition that ensure a computationally feasible estimator exists, i.e., $\lambda_{\submin}/\sigma\gg p^{m/4}$ \citep{zhang2018tensor}.  
In other words, if an sufficiently good initial estimate is already available, a weaker SNR condition $\lambda_{\submin}/\sigma\gg p^{3/4}$ is sufficient to guarantee the asymptotic normality of our final estimates. This phenomenon is further justified by the simulation results in Figure~\ref{fig:na3}.

\section*{Acknowledgments}
The authors thank the Editor, the Associate Editor, and three anonymous referees for their comments that help improve the presentation of this paper. 

\bibliographystyle{imsart-number} % Style BST file (imsart-number.bst or imsart-nameyear.bst)
\bibliography{refer.bib}       % Bibliography file (usually '*.bib')
%%%%%%%%%%%%%%%%%%%%%%%%%%%%%%%%%%%%%%%%%%%%%%
%% Supplementary Material, if any, should   %%
%% be provided in {supplement} environment  %%
%% with title and short description.        %%
%%%%%%%%%%%%%%%%%%%%%%%%%%%%%%%%%%%%%%%%%%%%%%
\appendix
\newpage
\setcounter{page}{1}

\begin{center}
	{\LARGE Supplementary Materials to ``Inference for Low-rank Tensors}
	
	\bigskip	
	{\LARGE   -- No Need to Debias"}
	
	\bigskip
	
	{Dong Xia$^\ast$, ~ Anru Zhang$^{\dagger}$, ~ and ~ Yuchen Zhou$^\ddagger$\\
		{\small Hong Kong University of Science and Technology}$^\ast$\\
		{\small  Duke University$^\dagger$}\\}
		{\small  Princeton University}$^\ddagger$\\
\end{center}

\footnotetext[1]{
	Dong Xia's research was partially supported by Hong Kong RGC Grant ECS 26302019 and GRF 16303320.}

\footnotetext[2]{Anru R. Zhang and Yuchen Zhou's research was partially supported by NSF Grants CAREER-1944904, NSF DMS-1811868, and grants from Wisconsin Alumni Research Foundation (WARF). This work was done while Anru R. Zhang and Yuchen Zhou were at the University of Wisconsin-Madison.}

\begin{abstract}
	In this supplement, we provide the optimal estimation procedures for Tucker low-rank tensor PCA and tensor regression, and proofs of technical results in the main content. 
\end{abstract}

%%%%%%%%%%%%%%%%%%%%%%
\section{Optimal Estimation Procedure of Tucker Low-rank Tensor PCA and Tensor Regression}\label{sec:additional-procedure}
%%%%%%%%%%%%%%%%%%%%%%

We collect the estimation procedures for Tucker low-rank tensor PCA and tensor regression in this section. Consider the Tucker low-rank tensor PCA: $\calA = \calX + \calZ$, where $\calX = (U_1, U_2, U_3)\calG$. As proved by \cite{zhang2018tensor}, the following Algorithm \ref{al:tensor-PCA} achieves the optimal rate in estimation error.
\begin{algorithm}
	\SetAlgoLined
	\KwInput{$\calA$, $r_1, r_2, r_3$, iteration $t_{\submax}$;}
	Initialize $\hat U_1^{(0)}={\rm SVD}_{r_1}(\calM_1(\calA))$, $\hat U_2^{(0)}={\rm SVD}_{r_2}(\calM_2(\calA))$, $\hat U_3^{(0)}={\rm SVD}_{r_3}(\calM_3(\bA))$, $t=1$\;
	\While{$t\leq t_{\submax}$}{
		$\hat U_1^{(t)}=\text{leading $r_1$ left singular vectors of }\calM_1(\calA)\times_2 \hat U_2^{(t-1)\top}\times_3\hat U_3^{(t-1)\top}$\;
		$\hat U_2^{(t)}=\text{leading $r_2$ left singular vectors of }\calM_2(\calA)\times_1 \hat U_1^{(t-1)\top}\times_3\hat U_3^{(t-1)\top}$\;
		$\hat U_3^{(t)}=\text{leading $r_3$ left singular vectors of }\calM_3(\calA)\times_1\hat U_1^{(t-1)\top}\times_2 \hat U_2^{(t-1)\top}$\;
		$t=t+1$\;
	}
	$\hat\calG = \calA \times_1 \hat U_1^{(t_{\submax})\top} \times_2 \hat U_2^{(t_{\submax})\top} \times_3 \hat U_3^{(t_{\submax})\top}$\;
	\KwOutput{$\hat U_1=\hat U_1^{(t_{\submax})}, \hat U_2=\hat U_2^{(t_{\submax})}, \hat U_3=\hat U_3^{(t_{\submax})}, \hat\calG$}
	\caption{Higher Order Orthogonal Iteration (HOOI) \citep{de2000best,richard2014statistical,zhang2001rank}}
	\label{algo:power}
	\label{al:tensor-PCA}
\end{algorithm}

Next, we introduce the simultaneous gradient descent in Algorithm \ref{al:regression} for Tucker low-rank regression. \cite[Theorem 4.2]{han2020optimal} proved that Algorithm \ref{al:regression} achieves the optimal rate of estimation error for $\calT$.
\begin{algorithm}
	\SetAlgoLined
	\KwInput{$\ell_n(\cdot)$: the objective function (\ref{eq:MLE}) for tensor regression, $\{(\calX_i, Y_i)\}_{i=1}^n$, $r_1, r_2, r_3$, tuning parameters $a, b > 0$, step size $\eta$;}
	$\tilde U_1, \tilde U_2, \tilde U_3, \tilde\calG =$ HOOI($\sum_{i=1}^{n}Y_i\calX_i$);
	Initialize $\hat U_j^{(0)} = b\tilde U_j$ for $j \in [3]$, $\hat\calG^{(0)}=\tilde\calG/b^3$\;
	\For{$t = 0, \dots, t_{\submax}-1$}{
		\For{$j = 1,2,3$}{
			$\hat U_j^{(t+1)} = U^{(t)} - \eta\big(\nabla_{U_j}\ell_n\big((\hat{U}_1^{(t)}, \hat{U}_2^{(t)}, \hat{U}_3^{(t)})\cdot\hat\calG^{(t)}\big) + a\hat U_j^{(t)}(\hat U_j^{(t)\top}\hat U_j^{(t)} - b^2I_{r_j})\big)$\;
		}
		$\hat\calG^{(t+1)}=\hat\calG^{(t)} - \eta\nabla_{\calG}\ell_n\big((\hat{U}_1^{(t)}, \hat{U}_2^{(t)}, \hat{U}_3^{(t)})\cdot\hat\calG^{(t)}\big)$\;
		$t=t+1$;
	}
	$\hat\calT=(\hat U_1^{(t_{\submax})}, \hat U_2^{(t_{\submax})}, \hat U_3^{(t_{\submax})})\cdot\hat\calG^{(t_{\submax})}$\;
	$\hat U_j = {\rm SVD}_{r_j}(\calM_j(\hat\calT))$ for $j \in [3]$; $\hat\calG=\hat\calT \times_1 \hat U_1^\top \times_2 \hat U_2^\top \times_3 \hat U_3^\top$\;
	\KwOutput{$\hat U_1, \hat U_2, \hat U_3, \hat\calG$}
	\caption{Simultaneous Gradient Descent \citep{han2020optimal}}
	\label{al:regression}
\end{algorithm}

%%%%%%%%%%%%%%%%%%%%%%%
\section{Proofs}\label{sec:proof}
%%%%%%%%%%%%%%%%%%%%%%%

We collect all proofs for the technical results in this section. Without loss of generality, we assume that $r_j \asymp r_{\submax} \asymp r$ for $j \in [3]$.
\subsection{Proof of Theorem~\ref{thm:na_tsvd}}
    Note that $\calA/\sigma = \calT/\sigma + \calZ/\sigma.$
    We can replace $\calA, \calT, \calZ$ by $\calA/\sigma, \calT/\sigma, \calZ/\sigma$ without essentially changing the problem. Thus, we assume that $\sigma = 1$ without loss of generality.
	To simplify the notations, we write $\calP_{U}=UU^{\top}$ as the spectral projector for any orthonormal columns $U$, i.e.,  $U^{\top}U$ being an identity matrix. Then, write $\calP_U^{\perp}=I-\calP_U$. 
	Denote $A_j=\calM_j(\calA)$, $T_j=\calM_j(\calT)$, $G_j=\calM_j(\calG)$, and $Z_j=\calM_j(\calZ)$ the corresponding matricizations for all $j=1,2,3$. 
	
	Without loss of generality, we only consider $j = 1$ and prove the theorem for $\|\hat U_1\hat U_1^{\top}-U_1U_1^{\top}\|_{\rm F}^2$. Notice that Theorem~\ref{thm:na_tsvd} automatically holds if $p \leq r^{1/3}$, we only need to consider the case $p \geq r^{1/3}$. By Algorithm~\ref{algo:am_optimal_PCA}, $\hat U_1=\hat U_1^{(2)}$ contains the top-$r_1$ eigenvectors of $A_1(\calP_{\hat U_2^{(1)}}\otimes \calP_{\hat U_3^{(1)}})A_1^{\top}$. As a result, $\hat U_1^{(2)}\hat U_1^{(2)\top}$ is the spectral projector for the top-$r_1$ eigenvectors of 
	\begin{align*}
	A_1(\calP_{\hat U_2^{(1)}}\otimes \calP_{\hat U_3^{(1)}})A_1^{\top}=T_1(\calP_{U_2}\otimes \calP_{U_3})T_1^{\top} +\frakJ_1+\frakJ_2+\frakJ_3+\frakJ_4
	\end{align*}
	where $\frakJ_1=T_1(\calP_{\hat U_2^{(1)}}\otimes \calP_{\hat U_3^{(1)}})Z_1^{\top}$, $\frakJ_2=Z_1(\calP_{\hat U_2^{(1)}}\otimes \calP_{\hat U_3^{(1)}})T_1^{\top}$, $\frakJ_3=Z_1(\calP_{\hat U_2^{(1)}}\otimes \calP_{\hat U_3^{(1)}})Z_1^{\top}$, and $\frakJ_4 = T_1((\calP_{\hat U_2^{(1)}} - \calP_{U_{2}})\otimes \calP_{\hat U_3^{(1)}})T_1^{\top} + T_1(\calP_{U_2} \otimes (\calP_{\hat U_3^{(1)}} - \calP_{U_3})T_1^{\top}$.
	%Denote 
	%$$
	%\frakJ_0=T_1(\calP_{U_2}\otimes\calP_{U_3})Z_1^{\top}+Z_1(\calP_{U_2}\otimes\calP_{U_3})T_1^{\top}+Z_1(\calP_{U_2}\otimes\calP_{U_3})Z_1^{\top}
	%$$
	
	To this end, we write 
	\begin{equation*}
	\begin{split}
	A_1(\calP_{\hat U_2^{(1)}}\otimes \calP_{\hat U_3^{(1)}})A_1^{\top}=&U_1G_1G_1^{\top}U_1^{\top}+\frakJ_1+\frakJ_2+\frakJ_3+\frakJ_4\\=:&U_1G_1G_1^{\top}U_1^{\top}+\frakE_1.
	\end{split}
	\end{equation*}
	
	\begin{lemma}\label{lem:Ebound}Under Assumption~\ref{assump:tpca} and conditions of Theorem~\ref{thm:na_tsvd}, 
		there exist absolute constants $c_1,C_1,C_2>0$ so that with probability at least $1 - C_1e^{-c_1p}$,
		$$\|\frakJ_1\| = \|\frakJ_2\|\leq C_2\kappa_0\lambda_{\submin}\sqrt{p},\quad \|\frakJ_3\|\leq C_2p,\quad \|\frakJ_4\|\leq C_2\kappa_0^2p, \quad \|\frakE_1\|\leq C_2\kappa_0\lambda_{\submin}\sqrt{p}.$$
	\end{lemma}
    Moreover, by \cite[Theorem 1]{zhang2018tensor}, the following bounds hold:
    \begin{equation}\label{ineq16}
    	\max\{\|\hat U_k^{(0)}\hat U_k^{(0)\top}-U_kU_k^{\top}\|, \|\hat U_k^{(1)}\hat U_k^{(1)\top}-U_kU_k^{\top}\|, \|\hat U_k^{(2)}\hat U_k^{(2)\top}-U_kU_k^{\top}\|\} \leq C_2\sqrt{p}\lambda_{\submin}^{-1}
    \end{equation}
    for all $k \in [3]$.
	Denote $\calE_0$ the event of Lemma~\ref{lem:Ebound} and \eqref{ineq16} so that $\PP(\calE_0)\geq 1-C_1e^{-c_1p}$. 
	By definition, $\Lambda_j^2$ is a diagonal matrix containing the eigenvalues of $G_jG_j^{\top}$.  Without loss of generality, we assume that $G_jG_j^\top = \Lambda_j^2$ is a diagonal matrix. Then immediately we have
	\begin{equation}\label{eq4}
		\|\Lambda_j^{-1}G_j\| = 1, \forall j \in [3].
	\end{equation}
	
	\paragraph*{Step 1: representation of spectral projector $\hat U_1\hat U_1^{\top}$.}
	
	We write
	\begin{align*}
	\|\hat U_1\hat U_1^{\top}-U_1U_1^{\top}\|_{\rm F}^2=2r_1-2\big<\hat U_1\hat U_1^{\top}, U_1U_1^{\top}\big>=-2\big<\hat U_1\hat U_1^{\top}-U_1U_1^{\top}, U_1U_1^{\top}\big>.
	\end{align*}
	Define, for a positive integer $k$, $\frakP_j^{-k}=U_j\Lambda_j^{-2k}U_j^{\top}$. With a little abuse of notations, denote $\mathfrak{P}_j^{0}:=\mathfrak{P}_j^{\perp}:=\calP_{U_j}^{\perp}$. 
	Note that, under the event $\calE_0$ of Lemma~\ref{lem:Ebound} 
	$$
	\|\frakE_1\| \leq C_2\kappa_0\lambda_{\submin}\sqrt{p}  < \frac{\lambda_{\submin}^2}{2}.
	$$
	implying that the condition of \cite[Theorem~1]{xia2019normal} is satisfied. 
	\begin{lemma}\label{lem:spectral}\citep[Theorem~1]{xia2019normal}
		If $\|\frakE_1\| \leq \frac{\lambda_{\submin}^2}{2}$, the following equation holds 
		\begin{equation}\label{eq2}
			\hat U_1\hat U_1^{\top}-U_1U_1^{\top}=\sum\nolimits_{k\geq 1}\calS_{G_1,k}(\frakE_1)
		\end{equation}
		where for each positive integer $k$
		$$
		\calS_{G_1,k}(\frakE_1)=\sum_{s_1+\cdots+s_{k+1}=k}(-1)^{1+\tau(\bs)}\cdot \frakP_1^{-s_1}\frakE_1\frakP_{1}^{-s_2}\frakE_1\frakP_1^{-s_3}\cdots\frakP_1^{-s_k}\frakE_1\frakP_1^{-s_{k+1}}
		$$
		where $s_1,\cdots,s_{k+1}$ are non-negative integers and $\tau(\bs)=\sum_{j=1}^{k+1}\II(s_j>0)$.
	\end{lemma}
    By Lemma \ref{lem:Ebound} and \ref{lem:spectral}, eq.\eqref{eq2} holds under event $\calE_0$ of Lemma~\ref{lem:Ebound}.
	Since $\mathfrak{P}_j^{0}U_jU_j^{\top} = U_jU_j^{\top}\mathfrak{P}_j^{0} = 0$, we have $$\big<\calS_{G_1,1}(\frakE_1), U_1U_1^{\top}\big>= \big<\frakP_1^{-1}\frakE_1\frakP_{1}^{\perp}+\frakP_{1}^{\perp}\frakE_1\frakP_1^{-1}, U_1U_1^{\top}\big> = 0.$$ 
	Similarly, $\big<\calS_{G_1,2}(\frakE_1), U_1U_1^{\top}\big>=-\big<\frakP_1^{-1}\frakE_1\frakP_1^{\perp}\frakE_1\frakP_1^{-1}, U_1U_1^{\top}\big>$ and 
	$$
	\big<\calS_{G_1,3}(\frakE_1), U_1U_1^{\top}\big>=-2\tr\big(\frakP_1^{-1}\frakE_1\frakP_1^{\perp}\frakE_1\frakP_1^{\perp}\frakE_1\frakP_1^{-2}\big)+2\tr\big(\frakP_1^{-1}\frakE_1\frakP_1^{\perp}\frakE_1\frakP_1^{-1}\frakE_1\frakP_1^{-1}\big).
	$$
	Note that 
	\begin{equation}\label{ineq:bound_S}
	\begin{split}
	\|\calS_{G_1,k}(\frakE_1)\| 
	\leq& \sum_{s_1+\cdots+s_{k+1}=k}\left\|(-1)^{1+\tau(\bs)}\cdot \frakP_1^{-s_1}\frakE_1\frakP_{1}^{-s_2}\frakE_1\frakP_1^{-s_3}\cdots\frakP_1^{-s_k}\frakE_1\frakP_1^{-s_{k+1}}\right\|\\
	\leq& \binom{2k}{k}\frac{\|\frakE_1\|^k}{\lambda_{\submin}^{2k}} \leq \left(\frac{4\|\frakE_1\|}{\lambda_{\submin}^{2}}\right)^k,
	\end{split}
	\end{equation}
	implying that 
	\begin{equation*}
	\begin{split}
	&\Big|\sum\nolimits_{k\geq 4}\big<\calS_{G_1,k}(\frakE_1),U_1U_1^{\top}\big>\Big|\leq r_1\sum\nolimits_{k\geq 4}\Big(\frac{4\|\frakE_1\|}{\lambda_{\submin}^2}\Big)^k\leq C_2r_1\frac{\kappa_0^{4}p^2}{\lambda_{\submin}^4}
	\end{split}
	\end{equation*}
	where the last inequality holds under event $\calE_0$ by Lemma~\ref{lem:Ebound}. 
	
	Therefore, under event $\calE_0$, we write 
	$$
	\|\hat U_1\hat U_1^{\top}-U_1U_1^{\top}\|_{\rm F}^2=-2\big<\calS_{G_1,2}(\frakE_1), U_1U_1^{\top}\big>-2\big<\calS_{G_1,3}(\frakE_1), U_1U_1^{\top}\big>+O\Big(\frac{r_1\kappa_0^{4}p^2}{\lambda_{\submin}^4}\Big).
	$$
	Now, it suffices to investigate the first two terms on RHS of above equation.
	
	\paragraph*{Step 2: bounding $\big<\calS_{G_1,3}(\frakE_1), U_1U_1^{\top}\big>$.}
	Since $T_1^\top \frakP_1^{\perp} = 0$ and $\frakP_1^{\perp}T_1 = 0$,
	\begin{align}
	\big<\calS_{G_1,3}(\frakE_1), U_1U_1^{\top}\big>=& -2\tr\big(\frakP_1^{-1}\frakE_1\frakP_1^{\perp}\frakE_1\frakP_1^{\perp}\frakE_1\frakP_1^{-2}\big)+2\tr\big(\frakP_1^{-1}\frakE_1\frakP_1^{\perp}\frakE_1\frakP_1^{-1}\frakE_1\frakP_1^{-1}\big)\notag\\
	=&-2\tr\big(\frakP_1^{-1}(\frakJ_1+\frakJ_3)\frakP_1^{\perp}\frakJ_3\frakP_1^{\perp}(\frakJ_2+\frakJ_3)\frakP_1^{-2}\big)\notag\\
	&+2\tr\big(\frakP_1^{-1}(\frakJ_1+\frakJ_3)\frakP_1^{\perp}(\frakJ_2+\frakJ_3)\frakP_1^{-1}(\frakJ_1+\frakJ_2+\frakJ_3+\frakJ_4)\frakP_1^{-1}\big)\notag\\
	=&-2\tr\big(\frakP_1^{-1}(\frakJ_1+\frakJ_3)\frakP_1^{\perp}\frakJ_3\frakP_1^{\perp}(\frakJ_2+\frakJ_3)\frakP_1^{-2}\big)\label{eq:S_{G_1,3}}\\
	&+2\tr\big(\frakP_1^{-1}\frakJ_1\frakP_1^{\perp}\frakJ_2\frakP_1^{-1}(\frakJ_1+\frakJ_2)\frakP_1^{-1}\big)+2\tr\big(\frakP_1^{-1}\frakJ_1\frakP_1^{\perp}\frakJ_2\frakP_1^{-1}(\frakJ_3+\frakJ_4)\frakP_1^{-1}\big)\notag\\
	&+2\tr\big(\frakP_1^{-1}\frakJ_3\frakP_1^{\perp}(\frakJ_2+\frakJ_3)\frakP_1^{-1}(\frakJ_1+\frakJ_2+\frakJ_3+\frakJ_4)\frakP_1^{-1}\big)\notag\\
	&+2\tr\big(\frakP_1^{-1}\frakJ_1\frakP_1^{\perp}\frakJ_3\frakP_1^{-1}(\frakJ_1+\frakJ_2+\frakJ_3+\frakJ_4)\frakP_1^{-1}\big)\notag
	\end{align}
	Define the term
	\begin{equation*}
	\begin{split}
	\frakM =& \big<\calS_{G_1,3}(\frakE_1), U_1U_1^{\top}\big> - 2\tr\big(\frakP_1^{-1}\frakJ_1\frakP_1^{\perp}\frakJ_2\frakP_1^{-1}(\frakJ_1+\frakJ_2)\frakP_1^{-1}\big)\\
	=& -2\tr\big(\frakP_1^{-1}(\frakJ_1+\frakJ_3)\frakP_1^{\perp}\frakJ_3\frakP_1^{\perp}(\frakJ_2+\frakJ_3)\frakP_1^{-2}\big)+2\tr\big(\frakP_1^{-1}\frakJ_1\frakP_1^{\perp}\frakJ_2\frakP_1^{-1}(\frakJ_3+\frakJ_4)\frakP_1^{-1}\big)\\
	&+2\tr\big(\frakP_1^{-1}\frakJ_3\frakP_1^{\perp}(\frakJ_2+\frakJ_3)\frakP_1^{-1}(\frakJ_1+\frakJ_2+\frakJ_3+\frakJ_4)\frakP_1^{-1}\big)\\
	&+2\tr\big(\frakP_1^{-1}\frakJ_1\frakP_1^{\perp}\frakJ_3\frakP_1^{-1}(\frakJ_1+\frakJ_2+\frakJ_3+\frakJ_4)\frakP_1^{-1}\big).
	\end{split}
	\end{equation*}
	By Lemma \ref*{lem:Ebound}, under event $\calE_0$,
	\begin{equation}\label{ineq:trace_J3}
	\begin{split}
	\frakM \leq C_2r_1\frac{\kappa_0^2p\cdot (\kappa_0\sqrt{p}\lambda_{\submin})(\kappa_0\sqrt{p}\lambda_{\submin})}{\lambda_{\submin}^6} \leq C_2r_1\frac{\kappa_0^4p^2}{\lambda_{\submin}^4}.
	\end{split}
	\end{equation}
	Therefore, we conclude on event $\calE_0$ that
	\begin{equation}\label{ineq:second_order}
	\Big|\|\hat U_1\hat U_1^{\top}-U_1U_1^{\top}\|_{\rm F}^2-2\tr\big(\frakP_1^{-1}\frakE_1\frakP_1^{\perp}\frakE_1\frakP_1^{-1}\big)+4\tr\big(\frakP_1^{-1}\frakJ_1\frakP_1^{\perp}\frakJ_2\frakP_1^{-1}(\frakJ_1+\frakJ_2)\frakP_1^{-1}\big)\Big|\leq C_2\frac{r_1\kappa_0^{4}p^2}{\lambda_{\submin}^4}.
	\end{equation}
	We begin with considering $\tr\big(\frakP_1^{-1}\frakJ_1\frakP_1^{\perp}\frakJ_2\frakP_1^{-1}(\frakJ_1+\frakJ_2)\frakP_1^{-1}\big)$. Clearly, 
	\begin{equation}\label{ineq8}
	\begin{split}
	&\left|\tr\big(\frakP_1^{-1}\frakJ_1\frakP_1^{\perp}\frakJ_2\frakP_1^{-1}(\frakJ_1+\frakJ_2)\frakP_1^{-1}\big)\right|\\
	\leq& \left|\tr\big(\frakP_1^{-1}\frakJ_1\frakP_1^{\perp}\frakJ_2\frakP_1^{-1}\frakJ_1\frakP_1^{-1}\big)\right|
	+ \left|\tr\big(\frakP_1^{-1}\frakJ_1\frakP_1^{\perp}\frakJ_2\frakP_1^{-1}\frakJ_2\frakP_1^{-1}\big)\right|.
	\end{split}
	\end{equation}
	It suffices to bound $\left|\tr\big(\frakP_1^{-1}\frakJ_1\frakP_1^{\perp}\frakJ_2\frakP_1^{-1}\frakJ_1\frakP_1^{-1}\big)\right|$ and $\left|\tr\big(\frakP_1^{-1}\frakJ_1\frakP_1^{\perp}\frakJ_2\frakP_1^{-1}\frakJ_2\frakP_1^{-1}\big)\right|$, respectively. By the proof of Lemma~\ref{lem:Ebound}, on event $\calE_0$, there exist two (random) matrices $R_2 \in \OO_{r_2}$ and $R_3\in \OO_{r_3}$ such that $\|\hat U_2^{(1)} - U_2R_2\|, \|\hat U_3^{(1)} - U_3R_3\| \leq C_2\sqrt{p}/\lambda_{\submin}$. Therefore, on event $\calE_0$,  \eqref{ineq:projection}, 
	\begin{align}
	&\left\|\frakJ_1 - T_1(\calP_{U_2}\otimes\calP_{U_3})Z_1^{\top}\right\| = \left\|T_1(\calP_{\hat U_2^{(1)}}\otimes \calP_{\hat U_3^{(1)}})Z_1^{\top} - T_1(\calP_{U_2}\otimes\calP_{U_3})Z_1^{\top}\right\|\notag\\
	\leq& \left\|\left[T_1(\hat U_2^{(1)}\otimes \hat U_3^{(1)})\right]\left[Z_1(\hat U_2^{(1)}\otimes \hat U_3^{(1)})\right]^\top - \left[T_1((U_2R_2)\otimes (U_3R_3))\right]\left[Z_1((U_2R_2)\otimes (U_3R_3))\right]^\top\right\|\notag\\
	\leq &\left\|Z_1\left[(\hat U_2^{(1)}\otimes \hat U_3^{(1)})-((U_2R_2)\otimes (U_3R_3)) \right]\right\|\left\|T_1(\hat U_2^{(1)}\otimes \hat U_3^{(1)})\right\|\notag\\
	&+ \left\|Z_1((U_2R_2)\otimes (U_3R_3))\right\|\left\|T_1\left[(\hat U_2^{(1)}\otimes \hat U_3^{(1)})-((U_2R_2)\otimes (U_3R_3)) \right]\right\|\notag\\
	\leq& \kappa_0\lambda_{\submin}\left\|Z_1\left[(\hat U_2^{(1)}\otimes \hat U_3^{(1)})-((U_2R_2)\otimes (U_3R_3))\right]\right\|\notag\\
	&+C_0\sqrt{p}\left\|T_1\left[(\hat U_2^{(1)}\otimes \hat U_3^{(1)})-((U_2R_2)\otimes (U_3R_3)) \right]\right\|\notag\\
	\leq& \kappa_0\lambda_{\submin}\left(\left\|Z_1\left((\hat U_2^{(1)} - U_2R_2) \otimes \hat U_3^{(1)}\right)\right\|+\left\|Z_1\left((U_2R_2) \otimes (\hat U_3^{(1)} - U_3R_3)\right)\right\|\right)\notag\\
	&+C_0\sqrt{p}\Big(\big\|T_1\big((\hat U_2^{(1)}-U_2R_2)\otimes \hat U_3^{(1)}\big)\big\|+\big\|T_1\big((U_2R_2)\otimes (\hat U_3^{(1)}-U_3R_3)\big) \big\|\Big)\notag\\
	&\stackrel{\textrm{proof of lemma~\ref{lem:Ebound}, \ref{lm:best_rotation}}}{\leq} C_2\kappa_0\lambda_{\submin}\cdot \sqrt{pr}\left(\left\|\hat U_2^{(1)} - U_2R_2\right\| + \left\|\hat U_3^{(1)} - U_3R_3\right\|\right)+C_3\kappa_0p^{3/2}/\lambda_{\submin}\notag\\
	\leq& C_2\kappa_0 p\sqrt{r}.\label{ineq:diff_J_1}
	\end{align}
	By Lemma \ref{lem:Ebound} and \eqref{ineq:diff_J_1}, on event $\calE_0$,
	\begin{align}
	&\left|\tr\big(\frakP_1^{-1}\frakJ_1\frakP_1^{\perp}\frakJ_2\frakP_1^{-1}\frakJ_1\frakP_1^{-1}\big)\right|\notag\\
	\leq& \left|\tr\big(\frakP_1^{-1}T_1(\calP_{U_2}\otimes\calP_{U_3})Z_1^{\top}\frakP_1^{\perp}Z_1(\calP_{U_2}\otimes\calP_{U_3})T_1^{\top}\frakP_1^{-1}T_1(\calP_{U_2}\otimes\calP_{U_3})Z_1^{\top}\frakP_1^{-1}\big)\right|\notag\\
	&+\left|\tr\big(\frakP_1^{-1}(\frakJ_1 - T_1(\calP_{U_2}\otimes\calP_{U_3})Z_1^{\top})\frakP_1^{\perp}Z_1(\calP_{U_2}\otimes\calP_{U_3})T_1^{\top}\frakP_1^{-1}T_1(\calP_{U_2}\otimes\calP_{U_3})Z_1^{\top}\frakP_1^{-1}\big)\right|\notag\\
	&+ \left|\tr\big(\frakP_1^{-1}\frakJ_1\frakP_1^{\perp}(\frakJ_1 - T_1(\calP_{U_2}\otimes\calP_{U_3})Z_1^{\top})^\top\frakP_1^{-1}T_1(\calP_{U_2}\otimes\calP_{U_3})Z_1^{\top}\frakP_1^{-1}\big)\right|\notag\\
	&+ \left|\tr\big(\frakP_1^{-1}\frakJ_1\frakP_1^{\perp}\frakJ_2\frakP_1^{-1}(\frakJ_1 - T_1(\calP_{U_2}\otimes\calP_{U_3})Z_1^{\top})\frakP_1^{-1}\big)\right|\notag\\
	\leq& \left|\tr\big(\frakP_1^{-1}T_1(\calP_{U_2}\otimes\calP_{U_3})Z_1^{\top}\frakP_1^{\perp}Z_1(\calP_{U_2}\otimes\calP_{U_3})T_1^{\top}\frakP_1^{-1}T_1(\calP_{U_2}\otimes\calP_{U_3})Z_1^{\top}\frakP_1^{-1}\big)\right|\notag\\
	&+ C_2r_1\frac{\kappa_0 p\sqrt{r}\cdot \kappa_0\sqrt{p}\lambda_{\submin}\cdot \kappa_0\sqrt{p}\lambda_{\submin}}{\lambda_{\submin}^6}\notag\\
	=&\left|\tr\big(U_1\Lambda_1^{-2}U_1^\top T_1(\calP_{U_2}\otimes\calP_{U_3})Z_1^{\top}U_{1\perp}U_{1\perp}^\top Z_1(\calP_{U_2}\otimes\calP_{U_3})T_1^{\top}U_1\Lambda_1^{-2}U_1^\top T_1(\calP_{U_2}\otimes\calP_{U_3})Z_1^{\top}U_1\Lambda_1^{-2}U_1^\top\big)\right|\notag\\
	&+ C_2\kappa_0^3r^{3/2}p^2\lambda_{\submin}^{-4}\notag\\
	=&\left|\tr\big(U_1\Lambda_1^{-2}G_1W_2^\top W_2 G_1^{\top}\Lambda_1^{-2} G_1W_1^\top \Lambda_1^{-2}U_1^\top\big)\right|+ C_2\kappa_0^3r^{3/2}p^2\lambda_{\submin}^{-4}\label{ineq7}
	\end{align}
	where we denote
	\begin{equation}\label{eq:Gaussian}
	W_1 = U_1^\top Z_1(U_2 \otimes U_3) \in \RR^{r_1 \times (r_2r_3)}, \quad W_2 = U_{1\perp}^\top Z_1(U_2 \otimes U_3) \in \RR^{(p_1 - r_1) \times (r_2r_3)}.
	\end{equation}
	Due to the property of Gaussian matrices, $[W_1\ W_2] = [U_1^\top\ U_{1\perp}^\top]Z_1(U_2 \otimes U_3) \stackrel{i.i.d.}{\sim} N(0, 1)$. Therefore, $W_1 \stackrel{i.i.d.}{\sim} N(0, 1)$, $W_{1\perp} \stackrel{i.i.d.}{\sim} N(0, 1)$, and $W_1, W_2$ are independent. Conditioning on $W_2$, we have
	\begin{align*}
	&\tr\big(U_1\Lambda_1^{-2}G_1W_2^\top W_2 G_1^{\top}\Lambda_1^{-2} G_1W_1^\top \Lambda_1^{-2}U_1^\top\big)\bigg|W_2
	=  \tr\big(\Lambda_1^{-4}G_1W_2^\top W_2 G_1^{\top}\Lambda_1^{-2} G_1W_1^\top\big)\bigg|W_2\\
	\sim& N(0, \|\Lambda_1^{-4}G_1W_2^\top W_2 G_1^{\top}\Lambda_1^{-2} G_1\|_{\rm F}^2)\bigg| W_2.
	\end{align*}
	By the Gaussian concentration inequality, we get
	\begin{equation*}
	\begin{split}
	& \PP\left(\left|\tr\big(U_1\Lambda_1^{-2}G_1W_2^\top W_2 G_1^{\top}\Lambda_1^{-2} G_1W_1^\top \Lambda_1^{-2}U_1^\top\big)\bigg|W_2\right| \leq C_2\sqrt{\log(p)}\|\Lambda_1^{-4}G_1W_2^\top W_2 G_1^{\top}\Lambda_1^{-2} G_1\|_{\rm F}\bigg|W_2\right)\\ 
	\geq& 1- p^{-3}
	\end{split}
	\end{equation*}
	for some absolute constant $C_2>0$. Denote the above event $\calE_1$ so that $\PP(\calE_1)\geq 1-p^{-3}$. 

In addition, by \eqref{eq4}, on event $\calE_0$, 
	\begin{equation*}
	\|\Lambda_1^{-4}G_1W_2^\top W_2 G_1^{\top}\Lambda_1^{-2} G_1\|_{\rm F} \leq C_2\sqrt{r_1}\frac{\|W_2\|^2}{\lambda_{\submin}^3} \leq C_2\sqrt{r_1}\frac{p}{\lambda_{\submin}^3}.
	\end{equation*}
By the previous two inequalities, on event $\calE_0\cap \calE_1$, 
	\begin{equation*}
	\left|\tr\big(U_1\Lambda_1^{-2}G_1W_2^\top W_2 G_1^{\top}\Lambda_1^{-2} G_1W_1^\top \Lambda_1^{-2}U_1^\top\big)\bigg|W_2\right| \leq C_2\sqrt{r_1}\frac{p\sqrt{\log(p)}}{\lambda_{\submin}^3}.
	\end{equation*}
By combining eq. \eqref{ineq7} and the above inequality, we conclude on event $\calE_0\cap \calE_1$ that  
	\begin{equation*}
	\left|\tr\big(\frakP_1^{-1}\frakJ_1\frakP_1^{\perp}\frakJ_2\frakP_1^{-1}\frakJ_1\frakP_1^{-1}\big)\right| \leq C_2\left(r^{1/2}p\sqrt{\log(p)}\lambda_{\submin}^{-3} + r^{3/2}\kappa_0^3p^{2}\lambda_{\submin}^{-4}\right).
	\end{equation*}
	Similarly, on event $\calE_0\cap \calE_1$, 
	\begin{equation*}
	\left|\tr\big(\frakP_1^{-1}\frakJ_1\frakP_1^{\perp}\frakJ_2\frakP_1^{-1}\frakJ_2\frakP_1^{-1}\big)\right| \leq C_2\left(r^{1/2}p\sqrt{\log(p)}\lambda_{\submin}^{-3} + r^{3/2}\kappa_0^3p^{2}\lambda_{\submin}^{-4}\right).
	\end{equation*}
	Combining e.q. \eqref{ineq8} and the above two inequalities, with probability at least $1 - C_1p^{-3}$, 
	\begin{equation}\label{ineq:second_term}
	\left|\tr\big(\frakP_1^{-1}\frakJ_1\frakP_1^{\perp}\frakJ_2\frakP_1^{-1}(\frakJ_1+\frakJ_2)\frakP_1^{-1}\big)\right| \leq C_2\left(r^{1/2}p\sqrt{\log(p)}\lambda_{\submin}^{-3} + r^{3/2}\kappa_0^3p^{2}\lambda_{\submin}^{-4}\right).
	\end{equation}
	
	\paragraph*{Step 3: bounding smaller terms of $\langle \calS_{G_1,2}(\frakE_1),U_1U_1^{\top}\rangle$}
	Recall that $-\langle \calS_{G_1,2}(\frakE_1),U_1U_1^{\top}\rangle=\tr\big(\frakP_1^{-1}\frakE_1\frakP_1^{\perp}\frakE_1\frakP_1^{-1}\big)$. Since $T_1^\top \frakP_1^{\perp} = 0$ and $\frakP_1^{\perp}T_1 = 0$, we write
	\begin{align}
	&\tr\big(\frakP_1^{-1}\frakE_1\frakP_1^{\perp}\frakE_1\frakP_1^{-1}\big)= \tr\big(\frakP_1^{-1}(\frakJ_1 + \frakJ_3)\frakP_1^{\perp}(\frakJ_2 + \frakJ_3)\frakP_1^{-1}\big)\notag\\
	=& \tr\big(\frakP_1^{-1}\frakJ_1\frakP_1^{\perp}\frakJ_2\frakP_1^{-1}\big) + \tr\big(\frakP_1^{-1}\frakJ_1\frakP_1^{\perp}\frakJ_3\frakP_1^{-1}\big) + \tr\big(\frakP_1^{-1}\frakJ_3\frakP_1^{\perp}\frakJ_2\frakP_1^{-1}\big) +
	\tr\big(\frakP_1^{-1}\frakJ_3\frakP_1^{\perp}\frakJ_3\frakP_1^{-1}\big)\notag\\
	=&: \text{\rom{1}} + \text{\rom{2}} + \text{\rom{3}} + \text{\rom{4}}.\label{eq:main_term}
	\end{align}
	By Lemma~\ref{lem:Ebound}, on event $\calE_0$, 
	\begin{equation}\label{ineq:rom4}
	|\text{\rom{4}}| \leq r_1\|\frakP_1^{-1}\|^2\|\frakJ_3\|^2 \leq C_2r_1\frac{p^2}{\lambda_{\submin}^4}.
	\end{equation}
Next, we show that $\text{\rom{2}} = \text{\rom{3}} \leq C_2\left(r^{1/2} p\sqrt{\log(p)}\lambda_{\submin}^{-3} + r^{3/2}\kappa_0p^{2}\lambda_{\submin}^{-4}\right)$ with probability at least $1 - C_1p^{-3}$. Similarly to \eqref{ineq:diff_J_1}, with probability at least $1 - C_1e^{-c_1p}$,
	\begin{align}
	&\left\|\frakJ_3 - Z_1(\calP_{U_2}\otimes\calP_{U_3})Z_1^{\top}\right\| = \left\|Z_1(\calP_{\hat U_2^{(1)}}\otimes \calP_{\hat U_3^{(1)}})Z_1^{\top} - Z_1(\calP_{U_2}\otimes\calP_{U_3})Z_1^{\top}\right\|\notag\\
	\leq& \left\|\left[Z_1(\hat U_2^{(1)}\otimes \hat U_3^{(1)})\right]\left[Z_1(\hat U_2^{(1)}\otimes \hat U_3^{(1)})\right]^\top - \left[Z_1((U_2R_2)\otimes (U_3R_3))\right]\left[Z_1((U_2R_2)\otimes (U_3R_3))\right]^\top\right\|\notag\\
	\leq &\left\|Z_1\left[(\hat U_2^{(1)}\otimes \hat U_3^{(1)})-((U_2R_2)\otimes (U_3R_3)) \right]\right\|\left\|Z_1(\hat U_2^{(1)}\otimes \hat U_3^{(1)})\right\|\notag\\
	&+ \left\|Z_1((U_2R_2)\otimes (U_3R_3))\right\|\left\|Z_1\left[(\hat U_2^{(1)}\otimes \hat U_3^{(1)})-((U_2R_2)\otimes (U_3R_3)) \right]\right\|\notag\\
	=& \left\|Z_1\left[(\hat U_2^{(1)}\otimes \hat U_3^{(1)})-((U_2R_2)\otimes (U_3R_3))\right]\right\|\left(C_2\sqrt{p} +\left\|Z_1(U_2\otimes U_3)\right\| \right)\notag\\
	\leq& C_2\sqrt{p}\left(\left\|Z_1\left((\hat U_2^{(1)} - U_2R_2) \otimes \hat U_3^{(1)}\right)\right\|+\left\|Z_1\left((U_2R_2) \otimes (\hat U_3^{(1)} - U_3R_3)\right)\right\|\right)\notag\\
	&\stackrel{\textrm{proof of Lemma~\ref{lem:Ebound}}}{\leq} C_2\sqrt{p}\cdot \sqrt{pr}\left\|\hat U_2^{(1)} - U_2R_2\right\| + C_2\sqrt{p}\cdot \sqrt{pr}\left\|\hat U_3^{(1)} - U_3R_3\right\|\notag\\
	\leq& C_2\frac{p^{3/2}r^{1/2}}{\lambda_{\submin}}.\label{ineq:diff_J_3}
	\end{align}
Combining \eqref{ineq:diff_J_3} and \eqref{ineq:diff_J_1} together, with probability at least $1 - C_1e^{-c_1p}$,
	\begin{align}
	&|\text{\rom{2}}| = |\text{\rom{3}}|\notag\\
	=&\left|\tr\big(U_1\Lambda_1^{-2}U_1^\top\frakJ_1U_{1\perp}U_{1\perp}^\top\frakJ_3U_1\Lambda_1^{-2}U_1^\top\big)\right|\notag\\
	\leq& \left|\tr\big(U_1\Lambda_1^{-2}U_1^\top T_1(\calP_{U_2}\otimes\calP_{U_3})Z_1^{\top}U_{1\perp}U_{1\perp}^\top Z_1(\calP_{U_2}\otimes\calP_{U_3})Z_1^{\top}U_1\Lambda_1^{-2}U_1^\top\big)\right|\notag\\
	&+\left|\tr\big(U_1\Lambda_1^{-2}U_1^\top(\frakJ_1-T_1(\calP_{U_2}\otimes\calP_{U_3})Z_1^{\top})U_{1\perp}U_{1\perp}^\top Z_1(\calP_{U_2}\otimes\calP_{U_3})Z_1^{\top}U_1\Lambda_1^{-2}U_1^\top\big)\right|\notag\\
	&+\left|\tr\big(U_1\Lambda_1^{-2}U_1^\top\frakJ_1U_{1\perp}U_{1\perp}^\top(\frakJ_3-Z_1(\calP_{U_2}\otimes\calP_{U_3})Z_1^{\top})U_1\Lambda_1^{-2}U_1^\top\big)\right|\notag\\
	\leq& \left|\tr\big(U_1\Lambda_1^{-2}U_1^\top T_1(\calP_{U_2}\otimes\calP_{U_3})Z_1^{\top}U_{1\perp}U_{1\perp}^\top Z_1(\calP_{U_2}\otimes\calP_{U_3})Z_1^{\top}U_1\Lambda_1^{-2}U_1^\top\big)\right|\\ 
	&+ r\frac{\|T_1(\calP_{\hat U_2}\otimes\calP_{\hat U_3})Z_1^{\top}\|\|\frakJ_3-Z_1(\calP_{U_2}\otimes\calP_{U_3})Z_1^{\top}\|}{\lambda_{\submin}^4}\notag\\
	&+r\frac{\|\frakJ_1-T_1(\calP_{U_2}\otimes\calP_{U_3})Z_1^{\top}\|\|Z_1(\calP_{U_2}\otimes\calP_{U_3})Z_1^{\top}\|}{\lambda_{\submin}^4}\notag\\
	\leq& \left|\tr\big(U_1\Lambda_1^{-2}U_1^\top T_1(\calP_{U_2}\otimes\calP_{U_3})Z_1^{\top}U_{1\perp}U_{1\perp}^\top Z_1(\calP_{U_2}\otimes\calP_{U_3})Z_1^{\top}U_1\Lambda_1^{-2}U_1^\top\big)\right| + r^{3/2}\kappa_0\frac{p^2}{\lambda_{\submin}^4}\notag\\
	=&\left|\tr\big(U_1\Lambda_1^{-2}G_1W_2^\top W_2W_1^\top\Lambda_1^{-2}U_1^\top\big)\right| + r^{3/2}\kappa_0\frac{p^2}{\lambda_{\submin}^4}\notag\\
	=& \left|\tr\big(\Lambda_1^{-4}G_1W_2^\top W_2W_1^\top\big)\right| + r^{3/2}\kappa_0\frac{p^2}{\lambda_{\submin}^4}\label{ineq1}
	\end{align}
	where $W_1$ and $W_2$ are defined in \eqref{eq:Gaussian}. \\
Observe that
	\begin{equation*}
	\begin{split}
	& \tr\big(\Lambda_1^{-4}G_1W_2^\top W_2W_1^\top\big)\bigg|W_2
	=\big<\Lambda_1^{-4}G_1W_2^{\top}W_2,W_1\big>|W_2
	\sim N\left(0, \left\|\Lambda_1^{-4}G_1W_2^\top W_2\right\|_{\rm F}^2\right)
	\end{split}
	\end{equation*}
By the Gaussian concentration inequality, we have
	\begin{equation}\label{ineq:Gaussian_conditional}
	\begin{split}
	\PP\left(\tr\big(\Lambda_1^{-4}G_1W_2^\top W_2W_1^\top\big) \geq C_2\sqrt{\log(p)}\left\|\Lambda_1^{-4}G_1W_2^\top W_2\right\|_{\rm F}\bigg|W_2\right) \leq p^{-3}.
	\end{split}
	\end{equation}
	Moreover, by \eqref{eq4}, with probability at least $1 - C_1e^{-c_1p}$,
	\begin{align*}
	\left\|\Lambda_1^{-4}G_1W_2^\top W_2\right\|_{\rm F} \leq& \sqrt{r_1}\left\|\Lambda_1^{-4}G_1W_2^\top W_2\right\|
	\leq C_2\sqrt{r_1}\frac{\left\|W_2\right\|^2}{\lambda_{\submin}^3}\\
	\leq& C_2\sqrt{r_1}\frac{\left\|Z_1(U_2 \otimes U_3)\right\|^2}{\lambda_{\submin}^3}
	\leq C_2\sqrt{r_1}\frac{p}{\lambda_{\submin}^3}.
	\end{align*}
	By \eqref{ineq:Gaussian_conditional} and the above inequality, we get with probability at least $1 - C_1p^{-3}$,
	\begin{equation}\label{ineq:Gaussian_trace}
	\tr\big(\Lambda_1^{-4}G_1W_2^\top W_2W_1^\top\big) \leq C_2\sqrt{r_1}\frac{p\sqrt{\log(p)}}{\lambda_{\submin}^3}.
	\end{equation}
	Recall that $\lambda_{\submin} \gg p^{3/4}$, eq. \eqref{ineq1} and \eqref{ineq:Gaussian_trace} together imply that
	\begin{equation}\label{ineq:rom2_3}
	|\text{\rom{2}}| = |\text{\rom{3}}| \leq C_2\left(r^{1/2} p\sqrt{\log(p)}\lambda_{\submin}^{-3} + r^{3/2}\kappa_0p^{2}\lambda_{\submin}^{-4}\right)
	\end{equation}
	with probability $1 - C_1p^{-3}$.
	
\paragraph*{Step 4: treating the leading term of $\langle \calS_{G_1,2}(\frakE_1),U_1U_1^{\top}\rangle$.}
	Now, we consider the leading term \rom{1} $=\tr\big(\frakP_1^{-1}\frakJ_1\frakP_1^{\perp}\frakJ_2\frakP_1^{-1}\big) $. By definition and Algorithm~\ref{algo:am_optimal_PCA}, $\hat{U}_{2}^{(1)}\hat{U}_{2}^{(1)\top}$ is the spectral projector for the top-$r_2$ eigenvectors of 
	\begin{align*}
	A_2(\calP_{\hat U_{1}^{(0)}} \otimes \calP_{\hat U_{3}^{(0)}})A_2^\top
	=& U_2G_2G_2^\top U_2^\top - U_2G_2(U_1^\top\calP^{\perp}_{\hat U_{1}^{(0)}}U_1 \otimes U_3^\top\calP_{\hat U_{3}^{(0)}}U_3)G_2^\top U_2^\top\\
	 &- U_2G_2(I_{r_1} \otimes U_3^\top\calP^{\perp}_{\hat U_{3}^{(0)}}U_3)G_2^\top U_2^\top + T_2(\calP_{\hat U_{1}^{(0)}} \otimes \calP_{\hat U_{3}^{(0)}})Z_2^\top\\ &+ Z_2(\calP_{\hat U_{1}^{(0)}} \otimes \calP_{\hat U_{3}^{(0)}})T_2^\top
	+ Z_2(\calP_{\hat U_{1}^{(0)}} \otimes \calP_{\hat U_{3}^{(0)}})Z_2^\top\\
	 =:& U_2G_2G_2^\top U_2^\top + \hat{\frakE}_2.
	\end{align*}
	Similarly, we can define $\hat{\frakE}_3$. Let $\hat{\Lambda}_2^2$ and $\hat{\Lambda}_3^2$ be the diagonal matrices containing the eigenvalues of $G_2(U_1\calP_{\hat U_{1}^{(0)}}U_1^\top \otimes U_3\calP_{\hat U_{3}^{(0)}}U_3^\top)G_2^\top$ and $G_3(U_1\calP_{\hat U_{1}^{(0)}}U_1^\top \otimes U_2\calP_{\hat U_{2}^{(0)}}U_2^\top)G_3^\top$ with decreasing order, respectively. Let $\hat{\lambda}_{\submin}$ be the smallest eigenvalue among all eigenvalues of $\hat{\Lambda}_2^2$ and $\hat{\Lambda}_3^2$.  
	
Recall that $\frakP_j^k = U_j\Lambda_j^{-2k}U_j^\top$ for positive integer $k$, and $\frakP_j^0 := \frakP_j^{\perp} := \calP_{U_j}^{\perp}$ for $j=2,3$. 
	By Lemma \ref{lem:Ebound} and Lemma \ref{lem:spectral}, with probability at least $1 - C_1e^{-c_1p}$, 
	$$
	\|\hat{\frakE}_j\| \leq C_2\kappa_0\sqrt{p}\lambda_{\submin}
	$$
	and 
	$$\hat{U}_{j}^{(1)}\hat{U}_{j}^{(1)\top} - U_jU_j^\top = \sum_{k \geq 1}\calS_{G_j,k}(\hat{\frakE}_j).$$
	where for positive integer $k$, $$\calS_{G_j,k}(\hat{\frakE}_j)=\sum_{s_1+\cdots+s_{k+1}=k}(-1)^{1+\tau(\bs)}\cdot \frakP_j^{-s_1}\hat{\frakE}_j\frakP_j^{-s_2}\hat{\frakE}_j\frakP_j^{-s_3}\cdots\frakP_j^{-s_k}\hat{\frakE}_j\frakP_j^{-s_{k+1}}.$$
	For $k \geq 2$, similarly to \eqref{ineq:bound_S}, we have
	\begin{equation*}
	\big\|\calS_{G_j,k}(\hat{\frakE}_j)\big\| \leq \bigg(\frac{4\big\|\hat{\frakE}_j\big\|}{\lambda_{\submin}^2}\bigg)^k.
	\end{equation*}
	Then with probability at least $1 - C_1e^{-c_1p}$,
	\begin{equation*}
	\Big\|\sum_{k \geq 2}\calS_{G_j,k}(\hat{\frakE}_j)\Big\| \leq \sum_{k \geq 2}\bigg(\frac{4\big\|\hat{\frakE}_j\big\|}{\lambda_{\submin}^2}\bigg)^k \leq C_2\frac{\kappa_0^2 p}{\lambda_{\submin}^2}.
	\end{equation*}
	Note that 
	$$
	\calP_{U_j}\calS_{G_j,1}(\hat{\frakE}_j) = \calP_{U_j}\left(\frakP_j^{-1}\hat{\frakE}_j\frakP_j^{\perp} + \frakP_j^{\perp}\hat{\frakE}_j\frakP_j^{-1}\right) = U_j\Lambda_j^{-2}U_j^\top \hat{\frakE}_j \calP_{U_j}^{\perp}
	$$
Therefore, with probability at least $1 - C_1e^{-c_1p}$, for $j = 2, 3$, 
	\begin{equation}\label{ineq:projection_product}
	\Big\|U_j^\top\calP_{\hat U_j^{(1)}} - U_j^\top - \Lambda_j^{-2}U_j^\top \hat{\frakE}_j \calP_{U_j}^{\perp}\Big\| = \Big\|\calP_{U_j}\calP_{\hat U_j^{(1)}} - \calP_{U_j} - U_j\Lambda_j^{-2}U_j^\top \hat{\frakE}_j \calP_{U_j}^{\perp}\Big\| \leq C_2\frac{\kappa_0^2 p}{\lambda_{\submin}^2}.
	\end{equation}
For \rom{1} $=\tr\big(\frakP_1^{-1}\frakJ_1\frakP_1^{\perp}\frakJ_2\frakP_1^{-1}\big)$, with probability at least $1 - C_1e^{-c_1p}$,
\begin{align}
	&\left|\text{\rom{1}} - \tr\big(\Lambda_1^{-4}G_1(U_2^\top \otimes U_3^\top)Z_1^\top U_{1\perp}U_{1\perp}^\top Z_1(U_2 \otimes U_3)G_1^\top\big)\right|\notag\\
	=& \left|\tr\big(U_1\Lambda_1^{-2}U_1^\top\frakJ_1U_{1\perp}U_{1\perp}^\top\frakJ_1^\top U_1\Lambda_1^{-2}U_1^\top\big) - \tr\big(\Lambda_1^{-4}G_1(U_2^\top \otimes U_3^\top)Z_1^\top U_{1\perp}U_{1\perp}^\top Z_1(U_2 \otimes U_3)G_1^\top\big)\right|\notag\\
	=& \big|\tr\big(\Lambda_1^{-4}G_1((U_2^\top\calP_{\hat U_2^{(1)}}) \otimes (U_3^\top\calP_{\hat U_{3}^{(1)}}))Z_1^\top U_{1\perp}U_{1\perp}^\top Z_1((\calP_{\hat U_2^{(1)}}U_2) \otimes (\calP_{\hat U_{3}^{(1)}}U_3))G_1^\top\big)\notag\\ 
	&\quad -  \tr\big(\Lambda_1^{-4}G_1(U_2^\top \otimes U_3^\top)Z_1^\top U_{1\perp}U_{1\perp}^\top Z_1(U_2 \otimes U_3)G_1^\top\big)\big|.\label{ineq3}
\end{align}
By \eqref{eq4}, \eqref{ineq:projection_product} and \eqref{ineq:sub_gaussian_fixed}, with probability at least $1-C_1e^{-c_1p}$ that 
\begin{align}
	&\big|\tr\big(\Lambda_1^{-4}G_1((U_2^\top\calP_{\hat U_2^{(1)}}) \otimes (U_3^\top\calP_{\hat U_{3}^{(1)}}))Z_1^\top U_{1\perp}U_{1\perp}^\top Z_1((\calP_{\hat U_2^{(1)}}U_2) \otimes (\calP_{\hat U_{3}^{(1)}}U_3))G_1^\top\big)\notag\\
	&- \tr\big(\Lambda_1^{-4}G_1(U_2^\top \otimes U_3^\top)Z_1^\top U_{1\perp}U_{1\perp}^\top Z_1(U_2 \otimes U_3)G_1^\top\big)\big|\notag\\
	\leq& 2\left|\tr\big(\Lambda_1^{-4}G_1((\Lambda_2^{-2}U_2^\top \hat{\frakE}_2 \calP_{U_2}^{\perp}) \otimes U_3^\top)Z_1^\top U_{1\perp}U_{1\perp}^\top Z_1(U_2 \otimes U_3)G_1^\top\big)\right|\notag\\ 
	&+ 2\left|\tr\big(\Lambda_1^{-4}G_1(U_2^\top \otimes (\Lambda_3^{-2}U_3^\top \hat{\frakE}_3 \calP_{U_3}^{\perp}))Z_1^\top U_{1\perp}U_{1\perp}^\top Z_1(U_2 \otimes U_3)G_1^\top\big)\right|\notag\\
	&+ C_2r_1\max_{j = 2, 3}\left(\|\Lambda_j^{-2}U_j^\top \hat{\frakE}_j \calP_{U_j}^{\perp}\|\right)^2\cdot\left(\sqrt{pr}\right)^2\cdot \lambda_{\submin}^{-2} + C_2r_1\frac{\kappa_0^2p}{\lambda_{\submin}^2}\cdot \sqrt{pr}\cdot \sqrt{p}\cdot \lambda_{\submin}^{-2}\notag\\
	 &+ C_2r_1\frac{\kappa_0^2p}{\lambda_{\submin}^2}\cdot \max_{j = 2, 3}\|\Lambda_j^{-2}U_j^\top \hat{\frakE}_j \calP_{U_j}^{\perp}\|\cdot \left(\sqrt{pr}\right)^2\cdot \lambda_{\submin}^{-2} + C_2r_1\left(\frac{\kappa_0^2p}{\lambda_{\submin}^2}\right)^2\cdot \left(\sqrt{pr}\right)^2\cdot \lambda_{\submin}^{-2}\notag\\
	\leq& 2\left|\tr\big(\Lambda_1^{-4}G_1((\Lambda_2^{-2}U_2^\top \hat{\frakE}_2 \calP_{U_2}^{\perp}) \otimes U_3^\top)Z_1^\top U_{1\perp}U_{1\perp}^\top Z_1(U_2 \otimes U_3)G_1^\top\big)\right|\notag\\ 
	&+ 2\left|\tr\big(\Lambda_1^{-4}G_1(U_2^\top \otimes (\Lambda_3^{-2}U_3^\top \hat{\frakE}_3 \calP_{U_3}^{\perp}))Z_1^\top U_{1\perp}U_{1\perp}^\top Z_1(U_2 \otimes U_3)G_1^\top\big)\right|+ C_2r^2\kappa_0^2p^2\lambda_{\submin}^{-4}.\label{ineq4}
\end{align}
By the definition of $\hat{\frakE}_2$ and recall that $T_2^\top\calP_{U_2}^{\perp}=0$,
	eq. \eqref{eq4} and \eqref{ineq:sub_gaussian_fixed} and Lemma \ref{lem:Ebound} imply that with probability at least $1-C_1e^{-c_p}$, 
\begin{align}
	&\left|\tr\big(\Lambda_1^{-4}G_1((\Lambda_2^{-2}U_2^\top \hat{\frakE}_2 \calP_{U_2}^{\perp}) \otimes U_3^\top)Z_1^\top U_{1\perp}U_{1\perp}^\top Z_1(U_2 \otimes U_3)G_1^\top\big)\right|\notag\\
	\leq& \left|\tr\big(\Lambda_1^{-4}G_1((\Lambda_2^{-2}U_2^\top T_2(\calP_{\hat U_{1}^{(0)}} \otimes \calP_{\hat U_{3}^{(0)}})Z_2^\top \calP_{U_2}^{\perp}) \otimes U_3^\top)Z_1^\top U_{1\perp}U_{1\perp}^\top Z_1(U_2 \otimes U_3)G_1^\top\big)\right|\notag\\
	& + \left|\tr\big(\Lambda_1^{-4}G_1((\Lambda_2^{-2}U_2^\top Z_2(\calP_{\hat U_{1}^{(0)}} \otimes \calP_{\hat U_{3}^{(0)}})Z_2^\top \calP_{U_2}^{\perp}) \otimes U_3^\top)Z_1^\top U_{1\perp}U_{1\perp}^\top Z_1(U_2 \otimes U_3)G_1^\top\big)\right|\notag\\
	\leq& \left|\tr\big(\Lambda_1^{-4}G_1((\Lambda_2^{-2}U_2^\top T_2(\calP_{\hat U_{1}^{(0)}} \otimes \calP_{\hat U_{3}^{(0)}})Z_2^\top \calP_{U_2}^{\perp}) \otimes U_3^\top)Z_1^\top U_{1\perp}U_{1\perp}^\top Z_1(U_2 \otimes U_3)G_1^\top\big)\right|\notag\\
	&+C_2r_1\lambda_{\submin}^{-2}\cdot\sqrt{pr}p\lambda_{\submin}^{-2}\cdot \sqrt{p}\notag\\
	\leq& \left|\tr\big(\Lambda_1^{-4}G_1((\Lambda_2^{-2}U_2^\top T_2(\calP_{U_1} \otimes \calP_{U_3})Z_2^\top \calP_{U_2}^{\perp}) \otimes U_3^\top)Z_1^\top U_{1\perp}U_{1\perp}^\top Z_1(U_2 \otimes U_3)G_1^\top\big)\right|\notag\\
	&+\left|\tr\big(\Lambda_1^{-4}G_1((\Lambda_2^{-2}U_2^\top T_2((\calP_{\hat U_{1}^{(0)}} - \calP_{U_1}) \otimes \calP_{U_3})Z_2^\top \calP_{U_2}^{\perp}) \otimes U_3^\top)Z_1^\top U_{1\perp}U_{1\perp}^\top Z_1(U_2 \otimes U_3)G_1^\top\big)\right|\notag\\
	&+ \left|\tr\big(\Lambda_1^{-4}G_1((\Lambda_2^{-2}U_2^\top T_2(\calP_{\hat U_{1}^{(0)}} \otimes (\calP_{\hat U_{3}^{(0)}} - \calP_{U_3}))Z_2^\top \calP_{U_2}^{\perp}) \otimes U_3^\top)Z_1^\top U_{1\perp}U_{1\perp}^\top Z_1(U_2 \otimes U_3)G_1^\top\big)\right|\notag\\
	&+C_2r^{3/2}p^2\lambda_{\submin}^{-4}\notag\\
	\leq& \left|\tr\big(\Lambda_1^{-4}G_1((\Lambda_2^{-2}G_2(U_1^\top \otimes U_3^\top)Z_2^\top \calP_{U_2}^{\perp}) \otimes U_3^\top)Z_1^\top U_{1\perp}U_{1\perp}^\top Z_1(U_2 \otimes U_3)G_1^\top\big)\right|\notag\\
	&+ C_2r_1\lambda_{\submin}^{-2}\cdot \lambda_{\submin}^{-1}\cdot \sqrt{pr}\frac{\sqrt{p}}{\lambda_{\submin}}\cdot \sqrt{pr}\cdot \sqrt{p}+C_2r^{3/2}p^2\lambda_{\submin}^{-4}\notag\\
	\leq& \left|\tr\big(\Lambda_1^{-4}G_1((\Lambda_2^{-2}G_2W_3^\top \calP_{U_2}^{\perp}) \otimes U_3^\top)W_4^\top W_4(U_2 \otimes U_3)G_1^\top\big)\right| + C_2r^2p^2\lambda_{\submin}^{-4}\notag\\
	=& \left|\tr\big(\Lambda_1^{-4}G_1((\Lambda_2^{-2}G_2W_5^\top) \otimes I_{r_3})W_6^\top W_7G_1^\top\big)\right| + C_2r^2p^2\lambda_{\submin}^{-4}. \label{ineq5}
\end{align}
	where $W_3 = Z_2(U_1 \otimes U_3) \in \RR^{p_2 \times (r_1r_3)}$, $W_4 = U_{1\perp}^\top Z_1  \in \RR^{(p_1 - r_1) \times (p_2p_3)}, W_5 = U_{2\perp}^\top W_3 \in \RR^{(p_2-r_2) \times (r_1r_3)}, W_6 = W_4(U_{2\perp} \otimes U_3) \in \RR^{(p_1 - r_1) \times ((p_2 - r_2)r_3)}, W_7 = W_4(U_{2} \otimes U_3) \in \RR^{(p_1 - r_1) \times (r_2r_3)}$.

By definition, $W_3 \stackrel{i.i.d.}{\sim} N(0, 1)$, $W_4 \stackrel{i.i.d.}{\sim} N(0, 1)$, $W_5 \stackrel{i.i.d.}{\sim} N(0, 1)$, and $W_3$ and $W_4$ are independent. Furthermore, since $W_4([U_2\ U_{2\perp}] \otimes U_3) \stackrel{i.i.d.}{\sim} N(0, 1)$ and $W_6, W_7$ are two disjoint submatrices of $W_4([U_2\ U_{2\perp}] \otimes U_3)$. Therefore, $W_6 \stackrel{iid}{\sim} N(0, 1)$, $W_7 \stackrel{i.i.d.}{\sim} N(0, 1)$,  $W_5, W_6$, and $W_7$ are jointly independent. Then,
\begin{align*}
	&\tr\big(\Lambda_1^{-4}G_1((\Lambda_2^{-2}G_2W_5^\top) \otimes I_{r_3})W_6^\top W_7G_1^\top\big)\bigg|W_5, W_7\\
	=& \tr\big(W_7G_1^\top\Lambda_1^{-4}G_1((\Lambda_2^{-2}G_2W_5^\top) \otimes I_{r_3})W_6^\top\big)\bigg|W_5, W_7\\
	\sim& N(0, \|W_7G_1^\top\Lambda_1^{-4}G_1((\Lambda_2^{-2}G_2W_5^\top) \otimes I_{r_3})\|_{\rm F}^2).
\end{align*}
By the Gaussian concentration inequality, we have
	\begin{equation*}
	\begin{split}
	\PP\Big(\Big|\tr\big(&\Lambda_1^{-4}G_1((\Lambda_2^{-2}G_2W_5^\top) \otimes I_{r_3})W_6^\top W_7G_1^\top\big)\Big| \\
	& > C_2\sqrt{\log(p)}\|W_7G_1^\top\Lambda_1^{-4}G_1((\Lambda_2^{-2}G_2W_5^\top) \otimes I_{r_3})\|_{\rm F}\bigg|W_5, W_7\Big) \leq p^{-3}.
	\end{split}
	\end{equation*}
In addition, with probability at least $1 - C_1e^{-c_1p}$, we have $\|W_5\|, \|W_7\| \leq C_2\sqrt{p_1}$ since $r=O(\sqrt{p})$. By \eqref{eq4}, we obtain
	\begin{equation*}
	\begin{split}
	\|W_7G_1^\top\Lambda_1^{-4}G_1((\Lambda_2^{-2}G_2W_5^\top) \otimes I_{r_3})\|_{\rm F} \leq& \sqrt{r_1}\|W_7G_1^\top\Lambda_1^{-4}G_1((\Lambda_2^{-2}G_2W_5^\top) \otimes I_{r_3})\|\\
	\leq& C_2\sqrt{r_1}\sqrt{p}\lambda_{\submin}^{-2}\cdot \lambda_{\submin}^{-1}\sqrt{p}
	\leq C_2\sqrt{r_1}\frac{p}{\lambda_{\submin}^3}.
	\end{split}
	\end{equation*}
Therefore, with probability at least $1 - C_1p^{-3}$,
	\begin{equation*}
	\left|\tr\big(\Lambda_1^{-4}G_1((\Lambda_2^{-2}G_2W_5^\top) \otimes I_{r_3})W_6^\top W_7G_1^\top\big)\right| \leq C_2\sqrt{r_1}\frac{p\sqrt{\log(p)}}{\lambda_{\submin}^3}.
	\end{equation*}
Combining \eqref{ineq5} and the above inequality, we get with probability at least $1 - C_1p^{-3}$ that
	\begin{equation*}
	\begin{split}
	&\left|\tr\big(\Lambda_1^{-4}G_1((\Lambda_2^{-2}U_2^\top \hat{\frakE}_2 \calP_{U_2}^{\perp}) \otimes U_3^\top)Z_1^\top U_{1\perp}U_{1\perp}^\top Z_1(U_2 \otimes U_3)G_1^\top\big)\right|
	\leq  C_2(\sqrt{r_1}p\sqrt{\log(p)}\lambda_{\submin}^{-3} + r^2p^{2}\lambda_{\submin}^{-4}).
	\end{split}
	\end{equation*}
Similarly, with probability at least $1 - C_1p^{-3}$, 
	\begin{equation*}
	\begin{split}
	&\left|\tr\big(\Lambda_1^{-4}G_1(U_2^\top \otimes (\Lambda_3^{-2}U_3^\top \hat{\frakE}_3 \calP_{U_3}^{\perp}))Z_1^\top U_{1\perp}U_{1\perp}^\top Z_1(U_2 \otimes U_3)G_1^\top\big)\right|
	\leq C_2(\sqrt{r_1}p\sqrt{\log(p)}\lambda_{\submin}^{-3} + r^2p^{2}\lambda_{\submin}^{-4}).
	\end{split}
	\end{equation*}
By \eqref{ineq3}, \eqref{ineq4} and the above two inequalities, we get with probability at least $1 - C_1p^{-3}$ that 
	\begin{equation}\label{ineq6}
	\begin{split}
	&\left|\text{\rom{1}} - \tr\big(\Lambda_1^{-4}G_1(U_2^\top \otimes U_3^\top)Z_1^\top U_{1\perp}U_{1\perp}^\top Z_1(U_2 \otimes U_3)G_1^\top\big)\right|
	\leq C_2(\sqrt{r_1}p\sqrt{\log(p)}\lambda_{\submin}^{-3} + r^2\kappa_0^2p^{2}\lambda_{\submin}^{-4}).
	\end{split}
	\end{equation}
Combining \eqref{eq:main_term}, \eqref{ineq:rom4}, \eqref{ineq:rom2_3} and the above inequality, we get with probability at least $1 - C_1p^{-3}$ that
\begin{align*}
	\Big|\tr\big(\frakP_1^{-1}\frakE_1&\frakP_1^{\perp}\frakE_1\frakP_1^{-1}\big) - \tr\big(\Lambda_1^{-4}G_1(U_2^\top \otimes U_3^\top)Z_1^\top U_{1\perp}U_{1\perp}^\top Z_1(U_2 \otimes U_3)G_1^\top\big)\Big|\\ 
	\leq& C_2(\sqrt{r_1}p\sqrt{\log(p)}\lambda_{\submin}^{-3} + r^2\kappa_0^2p^{2}\lambda_{\submin}^{-4}).
\end{align*}
By \eqref{ineq:second_order}, \eqref{ineq:second_term} and the above inequality, with probability at least $1 - C_1p^{-3}$,
\begin{align}
	\Big|\|\hat U_1\hat U_1^{\top}-U_1U_1^{\top}\|_{\rm F}^2 &- 2\tr\big(\Lambda_1^{-4}G_1(U_2^\top \otimes U_3^\top)Z_1^\top U_{1\perp}U_{1\perp}^\top Z_1(U_2 \otimes U_3)G_1^\top\big)\Big|\notag\\
	 \leq& C_2(\sqrt{r_1}p\sqrt{\log(p)}\lambda_{\submin}^{-3} + r^2\kappa_0^4p^{2}\lambda_{\submin}^{-4}).\label{ineq:diff_main}
\end{align}

\paragraph*{Final step: characterizing the distribution.} By eq. (\ref{ineq:diff_main}), it suffices to prove the distribution of $\tr\big(\Lambda_1^{-4}G_1(U_2^\top \otimes U_3^\top)Z_1^\top U_{1\perp}U_{1\perp}^\top Z_1(U_2 \otimes U_3)G_1^\top\big)$.  We write
\begin{align}
	&\tr\big(\Lambda_1^{-4}G_1(U_2^\top \otimes U_3^\top)Z_1^\top U_{1\perp}U_{1\perp}^\top Z_1(U_2 \otimes U_3)G_1^\top\big)\notag\\
	=& \|\Lambda_1^{-2}G_1(U_2^\top \otimes U_3^\top)Z_1^\top U_{1\perp}\|_{\rm F}^2
	= \sum_{j=r_1+1}^{p_1}\|\Lambda_1^{-2}G_1(U_2^\top \otimes U_3^\top)Z_1^\top u_j\|_{2}^2,\label{ineq:trace_main}
\end{align}
	where $\{u_j\}_{j = r_1+1}^{p_1}$ are the columns of $U_{1\perp}$. For any $r_1+1 \leq j \leq p_1$, $Z_1^\top u_{j} \in N(0, I_{p_2 \times p_3})$, and
	\begin{equation*}
	\EE (Z_1^\top u_{j_1})(Z_1^\top u_{j_2})^\top = 0, \quad \forall r_1+1 \leq j_1 \neq j_2 \leq p_1.
	\end{equation*}
	Therefore, $\{Z_1^\top u_j\}_{j = r_1+1}^{p_1}$ are standard Gaussian random vectors. Recall that $G_1G_1^\top = \Lambda_1^2$ and
	\begin{equation*}
	\begin{split}
	\Lambda_1^{-2}G_1(U_2^\top \otimes U_3^\top)Z_1^\top u_j \sim& N(0, \Lambda_1^{-2}G_1(U_2^\top \otimes U_3^\top)[\Lambda_1^{-2}G_1(U_2^\top \otimes U_3^\top)]^\top)\\
	=& N(0, \Lambda_1^{-2}G_1G_1^\top\Lambda_1) = N(0, \Lambda_1^{-2}).
	\end{split}	
	\end{equation*}
Therefore,
	\begin{equation*}
	\|\Lambda_1^{-2}G_1(U_2^\top \otimes U_3^\top)Z_1^\top U_{1\perp}\|_{\rm F}^2 \stackrel{\rm d.}{=} \sum_{i = 1}^{p_1 - r_1}\|\Lambda_1^{-1}z_i\|_2^2,
	\end{equation*}
	where $z_i \stackrel{i.i.d.}{\sim} N(0, I_r)$. The RHS in above equation is a sum of independent random variables. 
	
Clearly, 
	$
	\EE \|\Lambda_1^{-1}z_i\|_2^2 = \|\Lambda_1^{-1}\|_{\rm F}^2,
	$
	$
	\Var\left(\|\Lambda_1^{-1}z_i\|_2^2\right) = 2\|\Lambda_1^{-2}\|_{\rm F}^2,
	$
	and 
	\begin{equation*}
	\EE \|\Lambda_1^{-1}z_i\|_2^6 \leq C_3\sum_{j_1, j_2, j_3 = 1}^{r_1}\frac{1}{\lambda_{j_1}^{(1)2}\lambda_{j_2}^{(1)2}\lambda_{j_3}^{(1)2}} = C_3\|\Lambda_1^{-1}\|_{\rm F}^6
	\end{equation*}
where we denote $\Lambda_1={\rm diag}\big(\lambda_1^{(1)},\lambda_2^{(1)},\cdots,\lambda_{r_1}^{(1)}\big)$. By Berry-Esseen theorem~\citep{berry1941accuracy, esseen1942liapunov}, we get
	\begin{equation*}
	\sup_{x \in \RR}\left|\PP\left(\frac{\|\Lambda_1^{-2}G_1(U_2^\top \otimes U_3^\top)Z_1^\top U_{1\perp}\|_{\rm F}^2 - 2(p_1 - r_1)\left\|\Lambda_1^{-1}\right\|_{\rm F}^2}{\sqrt{8(p_1 - r_1)}\left\|\Lambda_1^{-2}\right\|_{\rm F}} \leq x\right) - \Phi(x)\right| \leq C_3\left(\frac{\|\Lambda_1^{-1}\|_{\rm F}^4}{\|\Lambda_1^{-2}\|_{\rm F}^2}\right)^{3/2}\cdot \frac{1}{\sqrt{p_1-r_1}}.
	\end{equation*}
Note that $\sqrt{8(p_1 - r_1)\left\|\Lambda_1^{-2}\right\|_{\rm F}} \geq \sqrt{2p_1r_1}\kappa_0^{-2}\lambda_{\submin}^{-2}$. By eq. \eqref{ineq:diff_main} and Lipschitz property of $\Phi(\cdot)$, 
\begin{align}
	&\PP\left(\frac{\|\hat U_1\hat U_1^{\top}-U_1U_1^{\top}\|_{\rm F}^2 - 2(p_1 - r_1)\left\|\Lambda_1^{-1}\right\|_{\rm F}^2}{\sqrt{8(p_1 - r_1)}\left\|\Lambda_1^{-2}\right\|_{\rm F}} \leq x\right)\notag\\
	\leq& \PP\left(\frac{\|\Lambda_1^{-2}G_1(U_2^\top \otimes U_3^\top)Z_1^\top U_{1\perp}\|_{\rm F}^2 - 2(p_1 - r_1)\left\|\Lambda_1^{-1}\right\|_{\rm F}^2}{\sqrt{8(p_1 - r_1)}\left\|\Lambda_1^{-2}\right\|_{\rm F}} \leq x + C_2\Big(\frac{r^{3/2}\kappa_0^{6}p^{3/2}}{\lambda_{\submin}^{2}} + \frac{\kappa_0^2\sqrt{p\log p}}{\lambda_{\submin}}\Big)\right) + C_1p^{-3}\notag\\
	\leq& \Phi\left(x + C_2\Big(\frac{r^{3/2}\kappa_0^{6}p^{3/2}}{\lambda_{\submin}^{2}} + \frac{\kappa_0^2\sqrt{p\log p}}{\lambda_{\submin}}\Big)\right) + C_1p^{-3} + C_3\frac{r_1^{3/2}}{\sqrt{p_1 - r_1}}\notag\\
	\leq& \Phi(x) + C_2\Big(\frac{r^{3/2}\kappa_0^{6}p^{3/2}}{\lambda_{\submin}^{2}} + \frac{\kappa_0^2\sqrt{p\log p}}{\lambda_{\submin}}\Big) + C_3\frac{r^{3/2}}{\sqrt{p}}.\label{ineq:upper_normal_CDF}
\end{align}
	Similarly, we can show that 
	\begin{equation*}
	\begin{split}
	&\PP\bigg(\frac{\|\hat U_1\hat U_1^{\top}-U_1U_1^{\top}\|_{\rm F}^2 - 2(p_1 - r_1)\left\|\Lambda_1^{-1}\right\|_{\rm F}^2}{\sqrt{8(p_1 - r_1)}\left\|\Lambda_1^{-2}\right\|_{\rm F}} \leq x\bigg)\\
	&\hspace{3cm}\geq \Phi(x) - C_2\Big(\frac{r^{3/2}\kappa_0^{6}p^{3/2}}{\lambda_{\submin}^{2}} + \frac{\kappa_0^2\sqrt{p\log p}}{\lambda_{\submin}}\Big)-C_3\frac{r^{3/2}}{\sqrt{p}}.
	\end{split}
	\end{equation*}
    Combining two inequalities above, we know that
    \begin{align*}
    	\sup_{x\in\RR}&\left|\PP\bigg(\frac{\|\hat U_1\hat U_1^{\top}-U_1U_1^{\top}\|_{\rm F}^2 - 2(p_1 - r_1)\left\|\Lambda_1^{-1}\right\|_{\rm F}^2}{\sqrt{8(p_1 - r_1)}\left\|\Lambda_1^{-2}\right\|_{\rm F}} \leq x\bigg) - \Phi(x)\right|\\
    	&\hspace{3cm}\leq C_2\Big(\frac{r^{3/2}\kappa_0^{6}p^{3/2}}{\lambda_{\submin}^{2}} + \frac{\kappa_0^2\sqrt{p\log p}}{\lambda_{\submin}}\Big) + C_3\frac{r^{3/2}}{\sqrt{p}}.
    \end{align*}
    Moreover, by the previous inequality and the Lipschitz property of $\Phi(\cdot)$ and $|x|e^{-x^2/2} < 1$ for all $x \in \RR$, for any $x \in \RR$,
    \begin{align*}
    	&\left|\PP\bigg(\frac{\|\hat U_1\hat U_1^{\top}-U_1U_1^{\top}\|_{\rm F}^2 - 2p_1\left\|\Lambda_1^{-1}\right\|_{\rm F}^2}{\sqrt{8p_1}\left\|\Lambda_1^{-2}\right\|_{\rm F}} \leq x\bigg) - \Phi(x)\right|\\
    	=& \left|\PP\bigg(\frac{\|\hat U_1\hat U_1^{\top}-U_1U_1^{\top}\|_{\rm F}^2 - 2(p_1-r_1)\left\|\Lambda_1^{-1}\right\|_{\rm F}^2}{\sqrt{8(p_1-r_1)}\left\|\Lambda_1^{-2}\right\|_{\rm F}} \leq \sqrt{\frac{p_1}{p_1 - r_1}}x + \frac{r_1\left\|\Lambda_1^{-1}\right\|_{\rm F}^2}{\sqrt{2(p_1-r_1)}\left\|\Lambda_1^{-2}\right\|_{\rm F}}\bigg) - \Phi\left(x\right)\right|\\
    	\leq& C_2\Big(\frac{r^{3/2}\kappa_0^{6}p^{3/2}}{\lambda_{\submin}^{2}} + \frac{\kappa_0^2\sqrt{p\log p}}{\lambda_{\submin}}\Big) + C_3\frac{r^{3/2}}{\sqrt{p}} + \left|\Phi\left(\sqrt{\frac{p_1}{p_1 - r_1}}x + \frac{r_1\left\|\Lambda_1^{-1}\right\|_{\rm F}^2}{\sqrt{2(p_1-r_1)}\left\|\Lambda_1^{-2}\right\|_{\rm F}}\right) - \Phi(x)\right|\\
    	\leq& C_2\Big(\frac{r^{3/2}\kappa_0^{6}p^{3/2}}{\lambda_{\submin}^{2}} + \frac{\kappa_0^2\sqrt{p\log p}}{\lambda_{\submin}}\Big) + C_3\frac{r^{3/2}}{\sqrt{p}} + \left|\Phi\left(\sqrt{\frac{p_1}{p_1 - r_1}}x\right) - \Phi(x)\right|  + C_3\frac{r_1\left\|\Lambda_1^{-1}\right\|_{\rm F}^2}{\sqrt{2(p_1-r_1)}\left\|\Lambda_1^{-2}\right\|_{\rm F}}\\
    	\leq& C_2\Big(\frac{r^{3/2}\kappa_0^{6}p^{3/2}}{\lambda_{\submin}^{2}} + \frac{\kappa_0^2\sqrt{p\log p}}{\lambda_{\submin}}\Big) + C_3\frac{r^{3/2}}{\sqrt{p}} + \left(\sqrt{\frac{p_1}{p_1 - r_1}}-1\right)|x|e^{-x^2/2}\\
    	\leq& C_2\Big(\frac{r^{3/2}\kappa_0^{6}p^{3/2}}{\lambda_{\submin}^{2}} + \frac{\kappa_0^2\sqrt{p\log p}}{\lambda_{\submin}}\Big) + C_3\frac{r^{3/2}}{\sqrt{p}}.
    \end{align*}
	Therefore, we conclude the proof of Theorem~\ref{thm:na_tsvd}.

\subsection{Proof of Theorem \ref{thm:adaptive_tsvd}}
    First, we show that \begin{align}
    &\PP\left(\frac{\|\hat U_1\hat U_1^{\top}-U_1U_1^{\top}\|_{\rm F}^2 - 2p_1\sigma^2\big\|\hat{\Lambda}_1^{-1}\big\|_{\rm F}^2}{\sqrt{8p_1}\sigma^2\big\|\hat{\Lambda}_1^{-2}\big\|_{\rm F}} \leq x\right)\notag\\
    \leq& \Phi(x) + C_2\bigg(\frac{r^{3/2}\kappa_0^{6}p^{3/2}}{(\lambda_{\submin}/\sigma)^2} +  \frac{\kappa_0^3\sqrt{pr(r^2+\log(p))}}{\lambda_{\submin}/\sigma}\bigg)  + C_3\frac{r^{3/2}}{\sqrt{p}}.\label{ineq135}
    \end{align}
	Without loss of generality, we assume $\sigma = 1$ and only prove the result for $\|\hat U_1\hat U_1^{\top}-U_1U_1^{\top}\|_{\rm F}^2$. We denote $\tilde U_1\in\OO_{p_1,r_1}$ the top-$r_1$ left singular vectors of $\calM_1(\calA\times_2 \hat U_2^{\top}\times_3 \hat U_3^{\top})$. By Theorem~\ref{thm:na_tsvd} and \cite{zhang2018tensor}, it is easy to show that $\|\tilde U_1\tilde U_1^{\top}-U_1U_1^{\top}\|\leq C_2\sqrt{p}\lambda_{\submin}^{-1}$ with probability at least $1-C_1p^{-3}$. 
By definition, we know that $\hat{\Lambda}_1^2={\rm diag}(\hat \lambda_1^2,\cdots,\hat\lambda_{r_1}^2)$ contains the eigenvalues of $\tilde{U}_{1}^{\top}A_1(\calP_{\hat U_{2}} \otimes \calP_{\hat U_{3}})A_1^\top\tilde{U}_{1}$. We denote $\Lambda_1={\rm diag}(\lambda_1,\cdots,\lambda_{r_1})$. 
Then,
\begin{align*}
	&\sup_{1 \leq k \leq r_1}\left|\lambda_k^2 - \hat{\lambda}_k^2\right|\\ 
	\leq& \inf_{R \in \OO_{r_1}}\left\|\tilde{U}_{1}^{\top}A_1(\calP_{\hat U_{2}} \otimes \calP_{\hat U_{3}})A_1^\top\tilde{U}_{1} - RG_1G_1^\top R\right\|\\
	\leq& \inf_{R \in \OO_{r_1}}\left\|\tilde{U}_{1}^{\top}T_1(\calP_{\hat U_{2}} \otimes \calP_{\hat U_{3}})T_1^\top\tilde{U}_{1} - R\Lambda_1^2 R\right\|\\
	&+ 2\left\|\tilde{U}_{1}^{\top}T_1(\calP_{\hat U_{2}} \otimes \calP_{\hat U_{3}})Z_1^\top\tilde{U}_{1}\right\|
	+ \left\|\tilde{U}_{1}^{\top}Z_1(\calP_{\hat U_{2}} \otimes \calP_{\hat U_{3}})Z_1^\top\tilde{U}_{1}\right\|\\
	\leq& \left\|\tilde{U}_{1}^{\top}U_1G_1((U_2^\top\calP_{\hat U_{2}}U_2) \otimes (U_3^\top\calP_{\hat U_{3}}U_3))G_1^\top U_1^\top\tilde{U}_{1} - \tilde{U}_{1}^{\top}U_1G_1G_1^\top U_1^\top\tilde{U}_{1}\right\|\\
	&+ \inf_{R \in \OO_{r_1}}\left\|\tilde{U}_{1}^{\top}U_1\Lambda_1^2 U_1^\top\tilde{U}_{1} - R\Lambda_1^2 R^\top\right\|
	+ 2\kappa_0\lambda_{\submin}\left\|\tilde{U}_{1}^{\top}Z_1(\hat U_{2}\otimes \hat U_{3})\right\|
	+ \left\|\tilde U_1^{\top}Z_1(\hat U_{2} \otimes \hat U_{3})\right\|^2.
\end{align*}
	By \eqref{ineq:projection} and the Gaussian concentration inequality, with probability at least $1 - C_1p^{-3}$, 
	\begin{equation*}
	\Big\|\tilde U_1^{\top}Z_1(\hat U_{2} \otimes \hat U_{3})\Big\|^2 \leq C_2p,
	\end{equation*}
	and
\begin{align*}
	\kappa_0\lambda_{\submin}\big\|&\tilde{U}_{1}^{\top}Z_1(\hat U_{2}\otimes \hat U_{3})\big\|
	\leq \kappa_0\lambda_{\submin}\big\|\tilde{U}_{1}^{\top}(\calP_{U_1} + \calP_{U_1}^\perp)Z_1(\hat U_{2}\otimes \hat U_{3})\big\|\\
	\leq& \kappa_0\lambda_{\submin}\Big(\big\|U_1^\top Z_1(\hat U_{2}\otimes \hat U_{3})\big\| + \big\|\tilde{U}_{1}^{\top}U_{1\perp}\big\|\big\|Z_1(\hat U_{2}\otimes \hat U_{3})\big\|\Big)\\
	\leq& \kappa_0\lambda_{\submin}\left(\big\|U_1^\top Z_1(\hat U_{2}\otimes \hat U_{3})\big\| + \big\|\tilde{U}_{1}^{\top}U_{1\perp}\big\|\big\|Z_1(\hat U_{2}\otimes \hat U_{3})\big\|\right)\\
	\leq& \kappa_0\lambda_{\submin}\left(\big\|U_1^\top Z_1(U_2 \otimes U_3)\big\| + C_2\sqrt{pr}\frac{\sqrt{p}}{\lambda_{\submin}} + C_2\frac{\sqrt{p}}{\lambda_{\submin}}\sqrt{p}\right)\\
	\leq& C_2\kappa_0\lambda_{\submin}\left(\sqrt{r^2+\log(p)} + \sqrt{pr}\frac{\sqrt{p}}{\lambda_{\submin}} + \frac{\sqrt{p}}{\lambda_{\submin}}\sqrt{p}\right)\\
	\leq& C_2\kappa_0(\sqrt{r^2+\log(p)}\lambda_{\submin} + p\sqrt{r} ).
\end{align*}
Moreover, with probability at least $1 - C_1e^{-c_1p}$, 
	\begin{align*}
	\Big\|\tilde{U}_{1}^{\top}U_1G_1(&(U_2^\top\calP_{\hat U_{2}}U_2) \otimes (U_3^\top\calP_{\hat U_{3}}U_3))G_1^\top U_1^\top\tilde{U}_{1} - \tilde{U}_{1}^{\top}U_1G_1G_1^\top U_1^\top\tilde{U}_{1}\Big\|\\
	\leq& \Big\|G_1((U_2^\top\calP_{\hat U_{2}}U_2) \otimes (U_3^\top\calP_{\hat U_{3}}U_3))G_1^\top - G_1G_1^\top\Big\|\\
	\leq& \Big\|G_1((U_2^\top\calP_{\hat U_{2}}^\perp U_2) \otimes (U_3^{\top}\calP_{\hat U_{3}}U_3))G_1^\top\Big\| + \Big\|G_1(I_{r_2} \otimes (U_3^\top\calP_{\hat U_{3}}^\perp U_3))G_1^\top\Big\|\\
	\leq& \kappa_0^2\lambda_{\submin}^2\Big(\big\|U_2^\top\calP_{\hat U_{2}}^\perp U_2\big\| + \big\|U_3^\top\calP_{\hat U_{3}}^\top U_3\big\|\Big)
	\leq \kappa_0^2\lambda_{\submin}^2\Big(\big\|U_2^\top\hat U_{2\perp}\big\|^2 + \big\|U_3^\top\hat U_{3\perp}\big\|^2\Big)\\
	\leq& C_2\kappa_0^2p.
	\end{align*}
	%With probability at least $1 - C_1e^{-cp}$, there exists a (random) matrix $R_1 \in \OO_{r_1}$ such that $\left\|U_1 - \widehat{U}_1^{(\tmax)}R_1\right\| \leq C_2\sqrt{p}/\lambda_{\submin}$. 
To deal with $\inf_{R \in \OO_{r_1}}\left\|\tilde{U}_{1}^{\top}U_1\Lambda_1^2 U_1^\top\tilde{U}_{1} - R\Lambda_1^2 R^\top\right\|$, we need the following lemma.
	\begin{lemma}\label{lm:best_rotation}
		For $U, \hat{U} \in \OO_{p, r}$,
		\begin{equation*}
		\inf_{R \in \OO_{r}}\big\|\hat{U}^\top U - R\big\| \leq \|U_{\perp}^\top \hat{U}\|^2.
		\end{equation*}
		\begin{equation*}
			\inf_{R \in \OO_{r}}\big\|\hat{U}^\top U - R\big\|_{\rm F} \leq \min\left\{\|U_{\perp}^\top \hat{U}\|_{\rm F}^2, \sqrt{r}\|U_{\perp}^\top \hat{U}\|^2\right\}.
		\end{equation*}
	\end{lemma}
By Lemma \ref{lm:best_rotation} and \eqref{perturbation_equivalence}, with probability at least $1 - C_1e^{-cp}$, 
	\begin{equation*}
	\begin{split}
	\inf_{R \in \OO_{r_1}}\big\|\tilde{U}_{1}^{\top}&U_1\Lambda_1^2 U_1^\top\tilde{U}_{1} - R\Lambda_1^2 R^\top\big\|\\
	\leq& \inf_{R \in \OO_{r_1}}\left\{\big\|(\tilde{U}_{1}^{\top}U_1 - R) \Lambda_1^2 U_1^\top\tilde{U}_{1}\big\| + \big\|R\Lambda_1^2 (\tilde{U}_{1}^{\top}U_1 - R)^\top\big\|\right\}\\
	\leq& 2\inf_{R \in \OO_{r_1}}\big\|\tilde{U}_{1}^{\top}U_1 - R\big\|\|\Lambda_1^2\|
	\leq C_2(\sqrt{p}\lambda_{\submin}^{-1})^2\cdot\kappa_0^2\lambda_{\submin}^2\\
	\leq& C_2\kappa_0^2p.
	\end{split}
	\end{equation*}
	Combining together the inequalities above, we get with probability at least $1 - C_1p^{-3}$,
	\begin{equation}\label{ineq108}
	\sup_{1 \leq k \leq r_1}\big|\lambda_k^2 - \hat{\lambda}_k^2\big| \leq C_2\left(\kappa_0^2\sqrt{r}p + \kappa_0r\sqrt{\log(p)}\lambda_{\submin}\right).
	\end{equation}
Therefore, with probability at least $1 - C_1p^{-3}$,
	\begin{equation*}
	\begin{split}
	\Big|\left\|\Lambda_1^{-1}\right\|_{\rm F}^2 - \big\|\hat\Lambda_1^{-1}\big\|_{\rm F}^2\Big| \leq r_1\sup_{1 \leq k \leq r_1}\frac{\left|\lambda_k^2 - \hat{\lambda}_k^2\right|}{\lambda_k^2\hat{\lambda}_k^2} \leq C_2\left(\kappa_0^2r^{3/2}p\lambda_{\submin}^{-4} + \kappa_0r\sqrt{r^2+\log(p)}\lambda_{\submin}^{-3}\right)
	\end{split}
	\end{equation*}
and as a result
	\begin{equation*}
	\begin{split}
	& \frac{\Big|\left\|\Lambda_1^{-2}\right\|_{\rm F} - \big\|\hat\Lambda_1^{-2}\big\|_{\rm F}\Big|}{\left\|\Lambda_1^{-2}\right\|_{\rm F}} \leq \frac{\big\|\Lambda_1^{-2} - \hat\Lambda_1^{-2}\big\|_{\rm F}}{\left\|\Lambda_1^{-2}\right\|_{\rm F}} \\
	\leq & \frac{\sup_{1 \leq k \leq r_1}\frac{\left|\lambda_k^2 - \hat{\lambda}_k^2\right|}{\lambda_k^2\hat{\lambda}_k^2}}{\kappa_0^{-2}\lambda_{\submin}^{-2}}\leq C_2\left(\kappa_0^4r^{1/2}p\lambda_{\submin}^{-2} + \kappa_0^3\sqrt{r^2+\log(p)}\lambda_{\submin}^{-1}\right).
	\end{split}
	\end{equation*}
Note that $\sqrt{8(p_1 - r_1)\left\|\Lambda_1^{-2}\right\|_{\rm F}} \geq \sqrt{2p_1r_1}\kappa_0^{-2}\lambda_{\submin}^{-2}$. By eq. \eqref{ineq:upper_normal_CDF}, we have
	\begin{equation*}
	\begin{split}
	&\PP\left(\frac{\|\hat U_1\hat U_1^{\top}-U_1U_1^{\top}\|_{\rm F}^2 - 2p_1\big\|\hat{\Lambda}_1^{-1}\big\|_{\rm F}^2}{\sqrt{8p_1}\left\|\Lambda_1^{-2}\right\|_{\rm F}} \leq x\right)\\
	\leq& \PP\left(\frac{\|\hat U_1\hat U_1^{\top}-U_1U_1^{\top}\|_{\rm F}^2 - 2p_1\left\|\Lambda_1^{-1}\right\|_{\rm F}^2}{\sqrt{8p_1}\left\|\Lambda_1^{-2}\right\|_{\rm F}} \leq x + C_2\left(\kappa_0^4r\frac{p^{3/2}}{\lambda_{\submin}^2} + \kappa_0^3\frac{\sqrt{pr(r^2+\log(p))}}{\lambda_{\submin}}\right)\right) + C_1p^{-3}\\
	\leq& \Phi(x) + C_2\bigg(\frac{r^{3/2}\kappa_0^{6}p^{3/2}}{\lambda_{\submin}^2} +  \frac{\kappa_0^3\sqrt{pr(r^2+\log(p))}}{\lambda_{\submin}}\bigg)  + C_3\frac{r^{3/2}}{\sqrt{p}}.
	\end{split}
	\end{equation*}
Furthermore, we have
\begin{align}
	&\PP\left(\frac{\|\hat U_1\hat U_1^{\top}-U_1U_1^{\top}\|_{\rm F}^2 - 2p_1\big\|\hat{\Lambda}_1^{-1}\big\|_{\rm F}^2}{\sqrt{8p_1}\big\|\hat{\Lambda}_1^{-2}\big\|_{\rm F}} \leq x\right)\notag\\
	=& \PP\left(\frac{\|\hat U_1\hat U_1^{\top}-U_1U_1^{\top}\|_{\rm F}^2 - 2p_1\big\|\hat\Lambda_1^{-1}\big\|_{\rm F}^2}{\sqrt{8p_1}\left\|\Lambda_1^{-2}\right\|_{\rm F}} \leq x\left(1 + \frac{\big\|\hat{\Lambda}_1^{-2}\big\|_{\rm F}-\left\|\Lambda_1^{-2}\right\|_{\rm F}}{\left\|\Lambda_1^{-2}\right\|_{\rm F}}\right)\right)\notag\\
	\leq& \PP\left(\frac{\|\hat U_1\hat U_1^{\top}-U_1U_1^{\top}\|_{\rm F}^2 - 2p_1\big\|\hat\Lambda_1^{-1}\big\|_{\rm F}^2}{\sqrt{8p_1}\left\|\Lambda_1^{-2}\right\|_{\rm F}} \leq x\bigg(1 + C_2\Big(\frac{\kappa_0^4r^{1/2}p}{\lambda_{\submin}^2} + \frac{\kappa_0^3\sqrt{r^2+\log(p)}}{\lambda_{\submin}}\Big)\sgn(x)\bigg)\right)+ \frac{C_1}{p^3}\notag\\
	\leq& \Phi\left( x\bigg(1 + C_2\Big(\frac{\kappa_0^4r^{1/2}p}{\lambda_{\submin}^2} + \frac{\kappa_0^3\sqrt{r^2+\log(p)}}{\lambda_{\submin}}\Big)\sgn(x)\bigg)\right) + \frac{C_1}{p^3}\notag\\ 
	&+ C_2\bigg(\frac{r^{3/2}\kappa_0^{6}p^{3/2}}{\lambda_{\submin}^2} +  \frac{\kappa_0^3\sqrt{pr(r^2+\log(p))}}{\lambda_{\submin}}\bigg)  + C_3\frac{r^{3/2}}{\sqrt{p}}\notag\\
	\leq& \Phi(x) + C_2\left(\Big(\frac{\kappa_0^4r^{1/2}p}{\lambda_{\submin}^2} + \frac{\kappa_0^3\sqrt{r^2+\log(p)}}{\lambda_{\submin}}\Big)|x|\cdot e^{-x^2/2}\right)\notag\\
	&+ C_2\bigg(\frac{r^{3/2}\kappa_0^{6}p^{3/2}}{\lambda_{\submin}^2} +  \frac{\kappa_0^3\sqrt{pr(r^2+\log(p))}}{\lambda_{\submin}}\bigg)  + C_3\frac{r^{3/2}}{\sqrt{p}}\notag\\
	\leq& \Phi(x) + C_2\Big(\frac{\kappa_0^4r^{1/2}p}{\lambda_{\submin}^2} + \frac{\kappa_0^3\sqrt{r^2+\log(p)}}{\lambda_{\submin}}\Big) +C_2\bigg(\frac{r^{3/2}\kappa_0^{6}p^{3/2}}{\lambda_{\submin}^2} +  \frac{\kappa_0^3\sqrt{pr(r^2+\log(p))}}{\lambda_{\submin}}\bigg)  + C_3\frac{r^{3/2}}{\sqrt{p}}\notag\\
	\leq& \Phi(x) + C_2\bigg(\frac{r^{3/2}\kappa_0^{6}p^{3/2}}{\lambda_{\submin}^2} +  \frac{\kappa_0^3\sqrt{pr(r^2+\log(p))}}{\lambda_{\submin}}\bigg)  + C_3\frac{r^{3/2}}{\sqrt{p}},\label{ineq9}
\end{align}
which has proved \eqref{ineq135}.

Then, by Lemma \ref{lm:variance_estimation} and the similar argument for proving \eqref{ineq135}, we further have
\begin{align*}
	&\PP\left(\frac{\|\hat U_1\hat U_1^{\top}-U_1U_1^{\top}\|_{\rm F}^2 - 2p_1\hat\sigma^2\|\hat{\Lambda}_1^{-1}\|_{\rm F}^2}{\sqrt{8p_1}\hat\sigma^2\|\hat{\Lambda}_1^{-2}\|_{\rm F}} \leq x\right)\\
	\leq&\Phi(x)+C_2\bigg(\frac{r^{3/2}\kappa_0^{6}p^{3/2}}{(\lambda_{\submin}/\sigma)^2} +  \frac{\kappa_0^3\sqrt{pr(r^2+\log p)}}{\lambda_{\submin}/\sigma} + \frac{\sqrt{\log(p)}}{p^{1/4}} + \frac{\kappa_0\sqrt{r}}{\sqrt{p}}\bigg) + C_3\frac{r^{3/2}}{\sqrt{p}}.
\end{align*}
Similarly, we have
\begin{align*}
	&\PP\left(\frac{\|\hat U_1\hat U_1^{\top}-U_1U_1^{\top}\|_{\rm F}^2 - 2p_1\hat\sigma^2\|\hat{\Lambda}_1^{-1}\|_{\rm F}^2}{\sqrt{8p_1}\hat\sigma^2\|\hat{\Lambda}_1^{-2}\|_{\rm F}} \leq x\right)\\
	\geq& \Phi(x) - C_2\bigg(\frac{r^{3/2}\kappa_0^{6}p^{3/2}}{(\lambda_{\submin}/\sigma)^2} +  \frac{\kappa_0^3\sqrt{pr(r^2+\log p)}}{\lambda_{\submin}/\sigma} + \frac{\sqrt{\log(p)}}{p^{1/4}} + \frac{\kappa_0\sqrt{r}}{\sqrt{p}}\bigg)  - C_3\frac{r^{3/2}}{\sqrt{p}}.
\end{align*}
Therefore, we conclude the proof of Theorem \ref{thm:adaptive_tsvd}.

\subsection{Proof of Theorem \ref{thm:PCA-sub-Gaussian}}
The outline of proof is similar to the proof of Theorem \ref{thm:na_tsvd}. We note that Lemma \ref{lem:Ebound} still holds for sub-Gaussian noise. Since rank $r_j=1$, we write $u_1$ ($u_2, u_3$ resp.) in short for $U_1$ ($U_2, U_3$ resp.), and $\lambda$ in short for $\lambda_{\submin}$. Similarly, we set $\sigma=1$ without loss of generality. By Lemma \ref{lm:AZB} and the sub-Gaussian concentration inequality, we know that with probability at least $1 - C_1p^{-3}$,
\begin{equation}\label{ineq:sub-Gaussian1}
	\|(u_2\otimes u_3)^{\top}Z_1^{\top}U_{1\perp}\| \leq \|(u_2\otimes u_3)^{\top}Z_1^{\top}\| \leq C\sqrt{p}
\end{equation}
and
\begin{equation}\label{ineq:sub-Gaussian2}
	\left|(u_2\otimes u_3)^{\top}Z_1^{\top}u_1\right| = \left|(u_1 \otimes u_2 \otimes u_3)\text{vec}(\calZ)\right| \leq C\sqrt{\log p}.
\end{equation}

\paragraph*{{\it Step 1: representation of spectral projector $\hat u_1\hat u_1^{\top}$}} Indeed, following the same treatment, we get 
\begin{align}\label{eq:rank1-subGaussian-step1}
\|\hat u_1\hat u_1^{\top}-u_1u_1^{\top}\|_{\rm F}^2=-2\big<\calS_{G_1,2}(\frakE_1), u_1u_1^{\top}\big> -2\big<\calS_{G_1,3}(\frakE_1), u_1u_1^{\top}\big>+O\Big(\frac{p^2}{\lambda^4}\Big)
\end{align}
under event $\calE_0$. 

\paragraph*{{\it Step 2: bounding $\langle \calS_{G_1,3}(\frakE_1), u_1u_1^{\top}\rangle$}} In the Step 2 of proof of Theorem \ref{thm:na_tsvd}, we only use Gaussian assumption to prove the upper bound for 
$$
\left|\tr\big(U_1\Lambda_1^{-2}U_1^\top T_1(\calP_{U_2}\otimes\calP_{U_3})Z_1^{\top}U_{1\perp}U_{1\perp}^\top Z_1(\calP_{U_2}\otimes\calP_{U_3})T_1^{\top}U_1\Lambda_1^{-2}U_1^\top T_1(\calP_{U_2}\otimes\calP_{U_3})Z_1^{\top}U_1\Lambda_1^{-2}U_1^\top\big)\right|
$$
in eq. (\ref{ineq7}). Now we prove an upper bound for this term under the sub-Gaussian noise setting. In the rank one case, this term can be simplified to
\begin{align*}
&\lambda^{-3}\big|\tr\big((u_2\otimes u_3)^{\top}Z_1^{\top}U_{1\perp}U_{1\perp}^{\top}Z_1(u_2\otimes u_3)(u_2\otimes u_3)^{\top}Z_1^{\top}u_1\big)\big|\\
\leq& \lambda^{-3} \|(u_2\otimes u_3)^{\top}Z_1^{\top}U_{1\perp}\|^2\cdot |(u_2\otimes u_3)^{\top}Z_1^{\top}u_1|\\
\leq& C_2 \frac{p\sqrt{\log p}}{\lambda^3},
\end{align*}
which holds with probability at least $1-C_1p^{-3}$. In the last line, we used \eqref{ineq:sub-Gaussian1} and \eqref{ineq:sub-Gaussian2}.

Therefore, we get with probability at least $1-C_1p^{-3}$ that
$$
\big|\tr\big(\frakP_1^{-1}\frakJ_1\frakP_1^{\perp}\frakJ_2\frakP_1^{-1}(\frakJ_1+\frakJ_2)\frakP_1^{-1}\big) \big|\leq C_3\big(\lambda^{-3}p\sqrt{\log p}+\lambda^{-4}p^2\big).
$$
The rest of proof is identical to the {\it Step 2} of proving Theorem \ref{thm:na_tsvd} and we conclude that with the same probability, 
$$
\big|\langle \calS_{G_1,3}(\frakE_1), u_1u_1^{\top}\rangle\big|\leq C_3\big(\lambda^{-3}p\sqrt{\log p}+\lambda^{-4}p^2\big).
$$ 

\paragraph*{{\it Step 3: bounding smaller terms of $\langle \calS_{G_1,2}(\frakE_1), u_1u_1^{\top}\rangle$}} In the Step 3 of proof of Theorem \ref{thm:na_tsvd}, we only use Gaussian assumption to prove the upper bound for 
$$
\left|\tr\big(U_1\Lambda_1^{-2}U_1^\top T_1(\calP_{U_2}\otimes\calP_{U_3})Z_1^{\top}U_{1\perp}U_{1\perp}^\top Z_1(\calP_{U_2}\otimes\calP_{U_3})Z_1^{\top}U_1\Lambda_1^{-2}U_1^\top\big)\right|
$$
in eq. (\ref{ineq1}). We now prove its upper bound under sub-Gaussian assumption. In the rank one case, this term can be written as
\begin{align*}
=&\lambda^{-3}\big| (u_2\otimes u_3)^{\top}Z_1^{\top}U_{1\perp}U_{1\perp}^{\top}Z_1(u_2\otimes u_3)(u_2\otimes u_3)^{\top}Z_1^{\top}u_1\big|\\
&\leq \lambda^{-3}\|(u_2\otimes u_3)^{\top}Z_1^{\top}U_{1\perp}\|^2\cdot \big| (u_2\otimes u_3)Z_1^{\top}u_1\big|\leq C_2\lambda^{-3}p\sqrt{\log p},
\end{align*}
which holds with probability at least $1-C_1p^{-3}$. Therefor, with the same probability, the smaller terms of $\langle \calS_{G_1,2}(\frakE_1), u_1u_1^{\top}\rangle$ can be upper bounded by 
$C_3(\lambda^{-3}p\log p+\lambda^{-4}p^2)$.

\paragraph*{{\it Step 4: treating the leading terms of $\langle \calS_{G_1,2}(\frakE_1), u_1u_1^{\top}\rangle$}} In the step 4 of proof of Theorem \ref{thm:na_tsvd}, we rely Gaussian assumption to prove an upper bound for 
$$
\left|\tr\big(\Lambda_1^{-4}G_1((\Lambda_2^{-2}G_2(U_1^\top \otimes U_3^\top)Z_2^\top \calP_{U_2}^{\perp}) \otimes U_3^\top)Z_1^\top U_{1\perp}U_{1\perp}^\top Z_1(U_2 \otimes U_3)G_1^\top\big)\right|
$$
in eq. (\ref{ineq5}). We now prove its upper bound under sub-Gaussian assumption, which is more involved than the previous steps. In the rank one case, this term can be simplified to 
\begin{align*}
&\lambda^{-3} \big|\tr\big([((u_1^{\top}\otimes u_3^{\top})Z_2^{\top}\calP_{u_2}^{\perp})\otimes u_3^{\top}]Z_1^{\top}\calP_{u_1}^{\perp}Z_1(u_2\otimes u_3)\big) \big|\\
\leq& \lambda^{-3} \big|\tr\big([((u_1^{\top}\otimes u_3^{\top})Z_2^{\top})\otimes u_3^{\top}]Z_1^{\top}\calP_{u_1}^{\perp}Z_1(u_2\otimes u_3)\big) \big|\\
&+ \lambda^{-3} |u_1^{\top}\otimes u_3^{\top})Z_2^{\top}u_2|\big|\tr\big((u_2^\top\otimes u_3^{\top})Z_1^{\top}\calP_{u_1}^{\perp}Z_1(u_2\otimes u_3)\big) \big|\\
=& \lambda^{-3} \big|\tr\big([((u_1^{\top}\otimes u_3^{\top})Z_2^{\top})\otimes u_3^{\top}]Z_1^{\top}\calP_{u_1}^{\perp}Z_1(u_2\otimes u_3)\big) \big|\\
&+ \lambda^{-3} |u_1^{\top}\otimes u_3^{\top})Z_2^{\top}u_2|\left\|(u_2^\top\otimes u_3^{\top})Z_1^{\top}\right\|_2^2\\
\leq& \lambda^{-3} \big|\tr\big([((u_1^{\top}\otimes u_3^{\top})Z_2^{\top})\otimes u_3^{\top}]Z_1^{\top}\calP_{u_1}^{\perp}Z_1(u_2\otimes u_3)\big) \big|+C_3\lambda^{-3}p\sqrt{\log p}\\
\leq& \lambda^{-3} \big|\tr\big([((u_1^{\top}\otimes u_3^{\top})Z_2^{\top})\otimes u_3^{\top}]Z_1^{\top}\calP_{u_1}Z_1(u_2\otimes u_3)\big) \big|\\
&+\lambda^{-3} \big|\tr\big([((u_1^{\top}\otimes u_3^{\top})Z_2^{\top})\otimes u_3^{\top}]Z_1^{\top}Z_1(u_2\otimes u_3)\big) \big|+C_3\lambda^{-3}p\sqrt{\log p},
\end{align*}
which holds with probability at least $1-C_1p^{-3}$ and we used the upper bound of $|(u_1\otimes u_3)^{\top}Z_2^{\top}u_2|$ similarly to \eqref{ineq:sub-Gaussian2}. 

Note that the $j$-th entry of $(u_1^{\top}\otimes u_3^{\top})Z_2^{\top}$ is $\sum_{i,k}Z(i,j,k)u_1(i)u_3(k)$. We know that the $i'$-th entry of $[((u_1^{\top}\otimes u_3^{\top})Z_2^{\top})\otimes u_3^{\top}]Z_1^{\top}$ is 
\begin{align*}
	\sum_{j',k'}Z(i',j',k')\bigg(\sum_{i,k}Z(i,j',k)u_1(i)u_3(k)\bigg)u_3(k').
\end{align*}
Therefore,
\begin{align*}
&\lambda^{-3} \big|\tr\big([((u_1^{\top}\otimes u_3^{\top})Z_2^{\top})\otimes u_3^{\top}]Z_1^{\top}\calP_{u_1}Z_1(u_2\otimes u_3)\big) \big|\\
\leq& \lambda^{-3}\big| [((u_1\otimes u_3)^{\top}Z_2^{\top})\otimes u_3^{\top}]Z_1^{\top}u_1\big|\cdot |u_1^{\top}Z_1(u_2\otimes u_3)|\\
\leq& \frac{|u_1^{\top}Z_1(u_2\otimes u_3)|}{\lambda^3}\cdot \Big|\sum_{i', j', k'}u_1(i')u_3(k')Z(i',j',k') \sum_{i,k}Z(i,j',k)u_1(i)u_3(k)\Big|\\
= & \frac{|u_1^{\top}Z_1(u_2\otimes u_3)|}{\lambda^3}\cdot \Big| \sum_{j'}\Big[\sum_{i',k'}u_1(i')u_3(k')Z(i',j',k')\Big]\Big[\sum_{i,k}u_1(i)u_3(k)Z(i,j',k) \Big]\Big|\\
\leq & \frac{|u_1^{\top}Z_1(u_2\otimes u_3)|}{\lambda^3}\cdot \sum_{j'}\Big[\sum_{i',k'}u_1(i')u_3(k')Z(i',j',k')\Big]^2\\
\leq &C_3\frac{\sqrt{\log p}}{\lambda^3}\cdot p=C_3\lambda^{-3}p\sqrt{\log p},
\end{align*}
which holds with probability at least $1-C_1p^{-3}$ and we used the concentration inequalities of the sum of independent sub-Gaussian and sub-exponential random variables. 

Therefore, it suffices to prove the upper bound for $\lambda^{-3}\big|[((u_1^{\top}\otimes u_3^{\top})Z_2^{\top})\otimes u_3^{\top}]Z_1^{\top}Z_1(u_2\otimes u_3)\big|$. Recall that the $i'$-th entry of $[((u_1\otimes u_3)^{\top}Z_2^{\top})\otimes u_3]Z_1^{\top}$ is 
$$
\sum_{j', k'}Z(i',j',k')\Big(\sum_{i,k}Z(i,j',k)u_1(i)u_3(k)\Big)u_3(k')
$$
and the $i'$-th entry of $Z_1(u_2\otimes u_3)$ is 
$$
\sum_{s,t}Z(i', s,t)u_2(s)u_3(t).
$$
Then, we write
\begin{align*}
&\lambda^{-3}[((u_1^{\top}\otimes u_3^{\top})Z_2^{\top})\otimes u_3^{\top}]Z_1^{\top}Z_1(u_2\otimes u_3)\\
=&\lambda^{-3}\sum_{i'}\Big[\sum_{j',k'}Z(i',j',k')u_3(k')\sum_{i,k}Z(i,j',k)u_1(i)u_3(k)\Big]\Big(\sum_{s,t}Z(i', s,t)u_2(s)u_3(t)\Big)\\
=&\lambda^{-3}\sum_{i'}\sum_{j',k'}\sum_{i,k}\sum_{s,t}Z(i',j',k')Z(i,j',k)Z(i',s,t)u_1(i)u_2(s)u_3(k')u_3(k)u_3(t)\\
=&\lambda^{-3}\sum_{i'}\sum_{i\neq i'}\sum_{j',k'}\sum_{k}\sum_{s,t}Z(i',j',k')Z(i,j',k)Z(i',s,t)u_1(i)u_2(s)u_3(k')u_3(k)u_3(t)\\
& + \lambda^{-3}\sum_{i'}\sum_{j',k'}\sum_{k}\sum_{s,t} Z(i',j',k')Z(i',j',k)Z(i',s,t)u_1(i')u_2(s)u_3(k')u_3(k)u_3(t)\\
=&:\frakJ_1+\frakJ_2.
\end{align*}
To bound $\frakJ_2$, we note
\begin{align*}
\frakJ_2=&\lambda^{-3}\sum_{i'}u_1(i')\Big[\sum_{j',k'}\sum_k Z(i',j',k')Z(i',j',k)u_3(k)u_3(k')\Big]\Big[\sum_{s,t}Z(i',s,t)u_2(s)u_3(t)\Big]\\
=&\lambda^{-3}\sum_{i'}u_1(i')\big\|Z(i',:,:) u_3\big\|_2^2\cdot (u_2^{\top}Z(i',:,:)u_3)
\end{align*}
which is the sum of independent random variables. 

Notice that $\EE \big\|Z(i',:,:) u_3\big\|_2^2 = p_2$ and for any $j \in [p_2]$, we have
\begin{align*}
	\left\|(Z(i',j,:) u_3)^2 - 1\right\|_{\psi_1} \leq C\left\|(Z(i',j,:) u_3)^2\right\|_{\psi_1} \leq C\left\|Z(i',j,:) u_3\right\|_{\psi_2}^2 \leq C\|u_2\|_2^2 = C.
\end{align*}
By \cite[Proposition 5.16]{vershynin2010introduction}, we know that for any $t \geq 0$,
\begin{align*}
	\PP\left(\left| \big\|Z(i',:,:) u_3\big\|_2^2 - p_2\right| \geq t\right) \leq 2\exp\left[-c\min(t^2/p, t)\right] \leq 2\exp(-ct/\sqrt{p}).
\end{align*}
Therefore, $\left\| \big\|Z(i',:,:) u_3\big\|_2^2 - p_2\right\|_{\psi_1} \leq C\sqrt{p}$.
Also note that $\|u_2^{\top}Z(i',:,:)u_3\|_{\psi_2} \leq C$ and
\begin{align*}
	\EE\left(\big\|Z(i',:,:) u_3\big\|_2^2 - p_2\right)(u_2^{\top}Z(i',:,:)u_3) = 0,
\end{align*}
By Lemma \ref{lm:orlicz} and \cite[Lemma 7]{hao2020sparse}, with probability at least $1 - p^{-3}$,
\begin{align*}
	\left|\sum_{i'}\left(\big\|Z(i',:,:) u_3\big\|_2^2 - p_2\right)(u_2^{\top}Z(i',:,:)u_3) \right| \leq C\sqrt{p}\left(\sqrt{p}\sqrt{\log p} + (\log p)^{3/2}\right) \leq Cp\sqrt{\log p}.
\end{align*}
In addition, by the sub-Gaussian concentration inequality, with probability at least $1 - p^{-3}$,
\begin{align*}
	\left|p_2\sum_{i'}u_2^{\top}Z(i',:,:)u_3\right| \leq C_3p\sqrt{\log p}.
\end{align*}
Combining the previous inequalities together, we get 
$$
|\frakJ_2|\leq C_4\lambda^{-3}p\sqrt{\log p}
$$
with probability at least $1-C_1p^{-3}$.

Now, we deal with the term $\frakJ_1$ and write 
\begin{align*}
\frakJ_1=&\lambda^{-3}\sum_{1\leq i\neq i'\leq p_1}\Big[\sum_{j',k'}\sum_{k}\sum_{s,t} Z(i',j',k')Z(i,j',k)Z(i',s,t)u_1(i)u_2(s)u_3(k')u_3(k)u_3(t)\Big]\\
=&\lambda^{-3}\sum_{1\leq i\neq i'\leq p_1}u_1(i)\Big[u_3^{\top}Z(i',:,:)^{\top}Z(i,:,:)u_3\Big]\Big[u_2^{\top}Z(i',:,:)u_3\Big].
\end{align*}
We now apply the decoupling trick to $\{Z(i,:,:)\}_{i=1}^{p_1}$. Let $\tilde Z$ be an independent copy of $Z$. By \cite[Theorem 3.1.1]{de2012decoupling}, it suffices to upper bound 
\begin{align*}
\tilde\frakJ_1=&\lambda^{-3}\sum_{1\leq i\neq i'\leq p_1}u_1(i)\Big[u_3^{\top}Z(i',:,:)^{\top}\tilde Z(i,:,:)u_3\Big]\Big[u_2^{\top}Z(i',:,:)u_3\Big]\\
&=-\lambda^{-3}\sum_{1\leq i\leq p_1}u_1(i)\Big[u_3^{\top}Z(i,:,:)^{\top}\tilde Z(i,:,:)u_3\Big]\Big[u_2^{\top}Z(i,:,:)u_3\Big]\\
&+\lambda^{-3}\Big<\Big[\sum_{1\leq i'\leq p_1}\big(u_2^{\top}Z(i',:,:)u_3\big)Z(i',:,:)u_3\Big],  \Big[\sum_{1\leq i\leq p_1}u_1(i)\tilde Z(i,:,:)u_3\Big]\Big>\\
&=:\tilde\frakJ_{11}+\tilde\frakJ_{12}
\end{align*}
Again, $\tilde\frakJ_{11}$ is the sum of independent random variables and can be treated similarly to $\frakJ_2$. We can show that with probability at least $1-C_1p^{-3}$,
$$
|\tilde\frakJ_{11}|\leq C_4\lambda^{-3}p\sqrt{\log p}.
$$
Note that $\tilde\frakJ_{12}$ is the inner product of two independent random variables. Conditioned on $\sum_{i'}u_2^{\top}Z(i',:,:)u_3Z(i',:,:)u_3$, by Hoeffding-type inequality, we get
$$
|\tilde\frakJ_{12}|\leq C_4\lambda^{-3}\Big\| \sum_{1\leq i'\leq p_1}\big(u_2^{\top}Z(i',:,:)u_3\big)Z(i',:,:)u_3\Big\|\cdot \sqrt{\log p},
$$
which holds with probability at least $1-C_1p^{-3}$. Note that $\sum_{i'}u_2^{\top}Z(i',:,:)u_3Z(i',:,:)u_3$ is the sum of independent random vectors, for any fixed $x \in \SS^{p_2-1}$, we have
\begin{align*}
	&\left|x^\top\sum_{i'}u_2^{\top}Z(i',:,:)u_3Z(i',:,:)u_3\right|\\
	\leq& \frac{1}{2}\sum_{i'}(u_2^{\top}Z(i',:,:)u_3)^2 + \frac{1}{2}\sum_{i'}(x^{\top}Z(i',:,:)u_3)^2.
\end{align*}
Since $\EE(u_2^{\top}Z(i',:,:)u_3)^2 = 1$ and $$\|(u_2^{\top}Z(i',:,:)u_3)^2\|_{\psi_1} \leq C\|u_2^{\top}Z(i',:,:)u_3\|_{\psi_2}^2 \leq C,$$
the Bernstein-type inequality implies that with probability at least $1 - C_1e^{-c_1p}$,
\begin{align*}
	\sum_{i'}(u_2^{\top}Z(i',:,:)u_3)^2 \leq Cp.
\end{align*}
Similarly, with probability at least $1 - C_1e^{-c_1p}$,
\begin{align*}
	\sum_{i'}(x^{\top}Z(i',:,:)u_3)^2 \leq Cp
\end{align*}
and consequently
\begin{align*}
	&\left|x^\top\sum_{i'}u_2^{\top}Z(i',:,:)u_3Z(i',:,:)u_3\right|\\
	\leq& Cp.
\end{align*}
Since $$\|\sum_{i'}u_2^{\top}Z(i',:,:)u_3Z(i',:,:)u_3\| = \sup_{x}\left|x^\top\sum_{i'}u_2^{\top}Z(i',:,:)u_3Z(i',:,:)u_3\right|,$$
the standard $\varepsilon$-net technique shows that with probability at least $1 - C_1e^{-c_1p}$,
$$\|\sum_{i'}u_2^{\top}Z(i',:,:)u_3Z(i',:,:)u_3\| \leq Cp.$$
%whose norm can be upper bounded by, for instance, matrix Bernstein inequality \citep{tropp2012user} together with a truncation trick on $\{u_2^{\top}Z(i',:,:)u_3\}$. We can show that with probability at least $1-C_1p^{-3}$ that $\big\|\sum_{i'}u_2^{\top}Z(i',:,:)u_3Z(i',:,:)u_3\big\|\leq C_4p\log p$. 
Finally, we get with probability at least $1-C_1p^{-3}$ that
$$
|\tilde\frakJ_{12}|\leq C_4 p\log^{1/2} p. 
$$
Therefore, by combining {\it Step 1}-{\it 3}, we conclude that with probability at least $1-C_1p^{-3}$, 
\begin{align*}
\Big|\|\hat u_1\hat u_1^{\top}-u_1u_1^{\top}\|_{\rm F}^2-2\lambda^{-2}\tr\big((u_2\otimes u_3)^{\top}Z_1^{\top}\calP_{u_1}^{\perp}Z_1(u_2\otimes u_3)\big)\Big|\leq C_3(\lambda^{-3}p\sqrt{\log p}+\lambda^{-4}p^2).
\end{align*}

\paragraph*{{\it Final step: characterizing the distribution}} Now it suffices to investigate the distribution of $\tr\big((u_2\otimes u_3)^{\top}Z_1^{\top}\calP_{u_1}^{\perp}Z_1(u_2\otimes u_3)\big)$. We write
\begin{align*}
&\tr\big((u_2\otimes u_3)^{\top}Z_1^{\top}\calP_{u_1}^{\perp}Z_1(u_2\otimes u_3)\big)\\=&-\tr\big((u_2\otimes u_3)^{\top}Z_1^{\top}\calP_{u_1}Z_1(u_2\otimes u_3)\big)
+ \tr\big((u_2\otimes u_3)^{\top}Z_1^{\top}Z_1(u_2\otimes u_3)\big),
\end{align*}
where the first term can be bounded by $C_3\log p$ with probability at least $1-C_1p^{-3}$. We characterize the distribution of the second term. Denote $z_j$ the $j$-th column of $Z_1$ for $j\in[p_1]$. We write
\begin{align}\label{eq:PCA-subGaussian-3}
 \tr\big((u_2\otimes u_3)^{\top}Z_1^{\top}Z_1(u_2\otimes u_3)\big)=\|Z_1(u_2\otimes u_3)\|^2=\sum_{j=1}^{p_1} \langle z_j, u_2\otimes u_3\rangle^2.
\end{align}
The RHS of (\ref{eq:PCA-subGaussian-3}) is the sum of independent random variables. Clearly, $\EE\langle z, u_2\otimes u_3\rangle^2=1$. By Berry-Esseen theorem, we get 
\begin{align*}
\sup_{x\in\RR}\left|\PP\left(\frac{2\lambda^{-2}\|Z_1(u_2\otimes u_3)\|^2-2p_1\lambda^{-2}}{2\lambda^{-2}\sqrt{p_1\cdot{\rm Var}(\langle z, u_2\otimes u_3\rangle^2)}}\leq x\right) -\Phi(x)\right|\leq \frac{C_4}{\sqrt{p}}\cdot \frac{\EE |\langle z, u_2\otimes u_3\rangle|^6}{{\rm Var}^{3/2}(\langle z, u_2\otimes u_3\rangle^2)}.
\end{align*}
We now calculate its variance. Denote $v=u_2\otimes u_3$. Write 
\begin{align*}
\langle z, v\rangle^4=\Big(\sum_{\omega\in [p_2]\times [p_3]}z_{\omega}^2 v_{\omega}^2+\sum_{\omega_1\neq \omega_2} z_{\omega_1}z_{\omega_2}v_{\omega_1}v_{\omega_2}\Big)^2
\end{align*}
and due to independence of the entries of $z$, we get 
\begin{align*}
\EE \langle z, v\rangle^4=&\EE\Big(\sum_{\omega}z_{\omega}^2v_{\omega}^2\Big)^2+2\EE\sum_{\omega_1\neq \omega_2} z_{\omega_1}^2z_{\omega_2}^2v_{\omega_1}^2v_{\omega_2}^2\\
=&\sum_{\omega}\EE z_{\omega}^4v_{\omega}^4 + 3\EE\sum_{\omega_1\neq \omega_2} z_{\omega_1}^2z_{\omega_2}^2v_{\omega_1}^2v_{\omega_2}^2\\
=&\sum_{\omega}\EE z_{\omega}^4v_{\omega}^4 + 3\sum_{\omega_1\neq \omega_2}v_{\omega_1}^2v_{\omega_2}^2\\
=&\sum_{\omega}(\EE z_{\omega}^4-3) v_{\omega}^4+3\|v\|_2^2\|v\|_2^2=3+(\nu-3)\|v\|_{4}^{4},
\end{align*}
where $\nu=\EE z_{\omega}^4$. Therefore, 
$$
{\rm Var}(\langle z, v\rangle^2)=2+(\nu-3)\|v\|_{4}^4.
$$
Since $z$ has i.i.d. sub-Gaussian entries, we have $\EE |\langle z, u_2\otimes u_3\rangle^6|\leq C_5$ for an absolute positive constant $C_5$. Therefore,
\begin{align*}
\sup_{x\in\RR}\left|\PP\left(\frac{2\lambda^{-2}\|Z_1(u_2\otimes u_3)\|^2-2p_1\lambda^{-2}}{2\lambda^{-2}\sqrt{p_1\big(2+(\nu-3)\|u_2\otimes u_3\|_{4}^4\big)}}\leq x\right) -\Phi(x)\right|\leq \frac{C_4}{\sqrt{p}}\cdot \frac{1}{\big(2+(\nu-3)\|u_2\otimes u_3\|_{4}^4\big)^{3/2}}.
\end{align*}
Collecting all the terms in previous steps, we conclude that 
\begin{align*}
\sup_{x\in\RR}&\left|\PP\left(\frac{\|\hat u_1\hat u_1^{\top}-u_1u_1^{\top}\|_{\rm F}^2-2p_1\lambda^{-2}}{2\lambda^{-2}\sqrt{p_1\big(2+(\nu-3)\|u_2\otimes u_3\|_{4}^4\big)}}\leq x\right) -\Phi(x)\right|\leq \frac{C_4}{\sqrt{p}}\cdot \frac{1}{\big(2+(\nu-3)\|u_2\otimes u_3\|_{4}^4\big)^{3/2}}\\
&\quad+C_5\Big(\frac{\sqrt{p\log p}}{\lambda}+\frac{p^{3/2}}{\lambda^2}+\frac{\log p}{\sqrt{p}}\Big)\cdot \frac{1}{\sqrt{2+(\nu-3)\|u_2\otimes u_3\|_{4}^4}}+C_1e^{-c_1p}.
\end{align*}
Finally, recall that $\|\hat u_1\hat u_1^{\top}-u_1u_1^{\top}\|_{\rm F}^2=2\|\sin\Theta(\hat u_1, u_1)\|_{\rm F}^2$, we have finished the proof.

\subsection{Proof of Theorem~\ref{thm:regression_alternating}}
Note that
$$Y_i/\sigma = \langle (\calT/\sigma), \calX_i\rangle + (\xi_i/\sigma).$$ 
We can replace $Y_i, \calT, \xi_i$ by $Y_i/\sigma$, $\calT/\sigma$, and $\xi_i/\sigma$ without changing this problem essentially. Therefore, we assume that $\sigma = 1$ without loss of generality. We only need to focus on the non-trivial case $p \geq r^{1/3}$. Since the proof is technical challenging and long, we divide the proof into several steps. Consider the SVD decomposition $\hat{U}_j^{(t)\top} U_j = L_j^{(t)}S_j^{(t)}D_j^{(t)\top}$ for $t = 0, 1$ and $j=1,2,3$, where $L_j^{(t)}, D_j^{(t)} \in \OO_{r_j}$, and $S_j^{(t)}$ is the diagonal matrix with all singular values of $\hat{U}_j^{(t)\top} U_j$ in decreasing order. Denote $R_j^{(t)} = L_j^{(t)}D_j^{(t)\top} \in \OO_{r_j}$.

For $t = 0, 1$, denote 
\begin{equation}\label{eq:deltaTt+0.5}
	\Delta \calT_{1}^{(t + 0.5)} = \hat \calG^{(t)} \times_1 \hat U_1^{(t+0.5)} \times_2 \hat U_2^{(t)}  \times_3 \hat U_3^{(t)} - \calT
\end{equation}
and denote $\hat G_j^{(t)}=\calM_j(\hat \calG^{(t)})$.

\paragraph*{Step 0: preliminary bounds on $\hat\calG^{(t)}$ and $\hat U_j^{(t+0.5)}$.} Before dealing with $\|\hat U_1^{(2)}\hat U_1^{(2)\top}-U_1U_1^{\top}\|_{\rm F}^2$, we first prove some preliminary results on  $\hat\calG^{(t)}$ and $\hat U_j^{(t+0.5)}$ which shall be used later.

\subparagraph*{Step 0.1: the error of $\hat\calG^{(t)}$.} Without loss of generality, we only prove the bound for $t=0$ and we write $\hat\calG=\hat\calG^{(0)}$ for brevity. 
Consider the SVD decomposition $\hat{U}_i^{(0)\top} U_i = L_iS_iD_i^\top$, where $L_i, D_i \in \OO_{r_i}$, and $S_i$ is the diagonal matrix with all singular values of $\hat{U}_i^{(0)\top} U_i$ in decreasing order. Note that we omitted the superscripts of $L_i, D_i, S_i$ for brevity.  Let $R_i = L_iD_i^\top \in \OO_{r_i}$ and
$$
\frakJ_1 = \frac{1}{n}R_1\bigg[U_1^\top\Big(\sum_{j = 1}^n \xi_j\calM_1(\calX_j)\Big)\left(U_2 \otimes U_3\right)\bigg](R_2^\top \otimes R_3^\top).
$$ 
We aim to show that with probability at least $1 - C_1e^{-c_1p} - p^{-3}$, 
\begin{equation}\label{ineq:error_G}
		\left\|\hat{G}_1 - R_1G_1(R_2^\top \otimes R_3^\top) - \frakJ_1\right\| \leq C\left(\kappa_0\frac{pr}{n} + \kappa_0\frac{p\sqrt{r}}{n\lambda_{\submin}}\right)
\end{equation}
where $\hat G_1=\calM_1(\hat\calG)$.
		
Since $\frac{\partial}{\partial \calG}\ell_n\left(\hat{\calG}\times_1 \hat{U}_1^{(0)} \times_2 \hat{U}_2^{(0)} \times_3 \hat{U}_3^{(0)}\right) = 0$, we have
\begin{equation}\label{ineq44}
			2\sum_{i=1}^{n}\hat{U}_1^{(0)\top} \calM_1(\calX_i)(\hat{U}_2^{(0)} \otimes \hat{U}_3^{(0)})\langle \hat{G}_1, \hat{U}_1^{(0)\top} \calM_1(\calX_i)(\hat{U}_2^{(0)} \otimes \hat{U}_3^{(0)})\rangle - 2\sum_{i=1}^{n}Y_i\hat{U}_1^{(0)\top} \calM_1(\calX_i)(\hat{U}_2^{(0)} \otimes \hat{U}_3^{(0)}) = 0.
\end{equation}
Denote
$$
\Delta\calG = \hat\calG - \calG \times_1 R_1 \times_2 R_2 \times_3 R_3, \quad \Delta G_1 = \calM_1(\Delta \calG) = \hat{G}_1 - R_1G_1(R_2^\top \otimes R_3^\top)
$$ 
and 
\begin{align*}
			\calA_{\calT}^{(0)} =& \calG \times_1 (\hat U_1^{(0)}R_1) \times_2 (\hat U_2^{(0)}R_2) \times_3 (\hat U_2^{(0)}R_2) - \calT,\\ \quad A_{T_1}^{(0)} =& \calM_1(\calA_{\calT}^{(0)}) = (\hat U_1^{(0)}R_1)G_1\left((\hat U_2^{(0)}R_2)^\top \otimes (\hat U_3^{(0)}R_3)^\top\right) - T_1.
\end{align*}
By \eqref{ineq44}, we have 
\begin{align}
			&\Delta G_1 - \frakJ_1\notag\\ 
			=& \Big(\calM_1(\Delta \calG) - \frac{1}{n}\sum_{i=1}^{n}\hat{U}_1^{(0)\top} \calM_1(\calX_i)(\hat{U}_2^{(0)} \otimes \hat{U}_3^{(0)})\big\langle \Delta \calG \times_1 \hat{U}_1^{(0)} \times_2 \hat{U}_2^{(0)} \times_3 \hat{U}_3^{(0)}, \calX_i\big\rangle\Big)\notag\\
			&- \Big(\frac{1}{n}\sum_{i=1}^{n}\hat{U}_1^{(0)\top} \calM_1(\calX_i)(\hat{U}_2^{(0)} \otimes \hat{U}_3^{(0)})\big\langle \calA_{\calT}^{(0)}, \calX_i\big\rangle - \hat{U}_1^{(0)\top} \calM_1(\calA_{\calT}^{(0)})(\hat{U}_2^{(0)} \otimes \hat{U}_3^{(0)})\Big)\notag\\
			&+ \Big(\frac{1}{n}\sum_{i=1}^{n}\xi_i\hat{U}_1^{(0)\top} \calM_1(\calX_i)(\hat{U}_2^{(0)} \otimes \hat{U}_3^{(0)}) - \frakJ_1\Big) - \hat{U}_1^{(0)\top} A_{T_1}^{(0)}(\hat{U}_2^{(0)} \otimes \hat{U}_3^{(0)}).\label{ineq45}
\end{align}
Notice that $\rank(\calM_i(\Delta \calG)) \leq 2r_i$ and $\rank(\calM_i(\calA_{\calT}^{(0)})) \leq 2r_i$ for $i \in [3]$, by Lemma \ref{lm:Gaussian_ensemble}, with probability at least $1 - e^{-C_1pr}$, 
\begin{align}\label{ineq50}
			\Big\|\calM_1(\Delta \calG) - \frac{1}{n}\sum_{i=1}^{n}\hat{U}_1^{(0)\top} \calM_1(\calX_i)(\hat{U}_2^{(0)} \otimes \hat{U}_3^{(0)})&\big\langle \Delta \calG \times_1 \hat{U}_1^{(0)} \times_2 \hat{U}_2^{(0)} \times_3 \hat{U}_3^{(0)}, \calX_i\big\rangle\Big\|\notag\\
			 &\leq C_2\sqrt{\frac{pr}{n}}\|\Delta \calG\|_{\rm F},
\end{align}
and
\begin{align}
			\Big\|\frac{1}{n}\sum_{i=1}^{n}\hat{U}_1^{(0)\top} \calM_1(\calX_i)(\hat{U}_2^{(0)} \otimes \hat{U}_3^{(0)})&\big\langle \calA_{\calT}^{(0)}, \calX_i\big\rangle - \hat{U}_1^{(0)\top} \calM_1(\calA_{\calT}^{(0)})(\hat{U}_2^{(0)} \otimes \hat{U}_3^{(0)})\Big\|\notag\\
			 &\leq C_2\sqrt{\frac{pr}{n}}\left\|\calA_{\calT}^{(0)}\right\|_{\rm F}.\label{ineq51}
\end{align}
By the definition of $R_i$ and Lemma \ref{lm:best_rotation}, with probability at least $1 - C_1e^{-c_1p}$, 
\begin{equation}\label{ineq14}
		\big\|\hat{U}_i^{(0)\top} U_i - R_i\big\| \leq \big\|U_{i\perp}^\top \hat{U}_i^{(0)}\big\|^2 \leq C\frac{p}{n\lambda_{\submin}^2},
\end{equation}
and
\begin{align}\label{ineq15}
		\big\|U_i - \hat{U}_i^{(0)}R_i\big\| \leq&\Big\|\calP_{\hat{U}_i^{(0)}}(U_i - \hat{U}_i^{(0)}R_i)\Big\| + \Big\|\calP_{\hat{U}_i^{(0)}}^\perp(U_i - \hat{U}_i^{(0)}R_i)\Big\|\notag\\ 
		 =& \big\|\hat{U}_i^{(0)\top} U_i - R_i\big\| + \big\|\big(\hat{U}_i^{(0)}\big)_{\perp}^{\top} U_i\big\|
		 \leq C\frac{\sqrt{p/n}}{\lambda_{\submin}}.
\end{align}
Thus with probability at least $1 - C_1e^{-c_1p}$,
\begin{align*}
			\big\|\calA_{\calT}^{(0)}\big\|_{\rm F}=& \big\|\calG \times_1 (\hat U_1^{(0)}R_1) \times_2 (\hat U_2^{(0)}R_2) \times_3 (\hat U_3^{(0)}R_3) - \calG \times_1 U_1 \times_2 U_2 \times_3 U_3\big\|_{\rm F}\\
			\leq& \Big\|\calG \times_1 (\hat U_1^{(0)}R_1 - U_1) \times_2 (\hat U_2^{(0)}R_2) \times_3 (\hat U_3^{(0)}R_3)\Big\|_{\rm F} \\
			& + \big\|\calG \times_1 U_1 \times_2 (\hat U_2^{(0)}R_2 - U_2) \times_3 (\hat U_3^{(0)}R_3)\big\|_{\rm F} + \big\|\calG \times_1 U_1 \times_2 U_2 \times_3 (\hat U_3^{(0)}R_3 - U_3)\big\|_{\rm F}\\
			\leq& \big\|\hat U_1^{(0)}R_1 - U_1\big\|_{\rm F}\|G_1\| + \big\|\hat U_2^{(0)}R_2 - U_2\big\|_{\rm F}\|G_2\| + \big\|\hat U_3^{(0)}R_3 - U_3\big\|_{\rm F}\|G_3\|\\
			\leq& C\frac{\sqrt{pr/n}}{\lambda_{\submin}}\cdot \kappa_0\lambda_{\submin} = C\kappa_0\sqrt{pr/n}.
\end{align*}
The previous inequality and \eqref{ineq51}, with probability at least $1 - e^{-C_1pr} - C_1e^{-c_1p}$,
\begin{equation}\label{ineq52}
			\bigg\|\frac{1}{n}\sum_{i=1}^{n}\hat{U}_1^{(0)\top} \calM_1(\calX_i)(\hat{U}_2^{(0)} \otimes \hat{U}_3^{(0)})\big\langle \calA_{\calT}^{(0)}, \calX_i\big\rangle - \hat{U}_1^{(0)\top} \calM_1(\calA_{\calT}^{(0)})(\hat{U}_2^{(0)} \otimes \hat{U}_3^{(0)})\bigg\| \leq C_2\kappa_0\frac{pr}{n}.
\end{equation}
By Lemma \ref{lm:sub_Gaussian} and \eqref{ineq15}, with probability at least $1 - e^{-C_1pr} - C_1e^{-c_1p}$,
\begin{align}\label{ineq53}
			&\Big\|\frac{1}{n}\sum_{i=1}^{n}\xi_i\hat{U}_1^{(0)\top} \calM_1(\calX_i)\left(\hat{U}_2^{(0)} \otimes \hat{U}_3^{(0)}\right) - \frakJ_1\Big\|\notag\\
			\leq& \Big\|\frac{1}{n}\sum_{i=1}^{n}\xi_i\big(\hat{U}_1^{(0)} - U_1R_1^\top\big)^\top\calM_1(\calX_i)\left(\hat{U}_2^{(0)} \otimes \hat{U}_3^{(0)}\right)\Big\|\notag\\ 
			&+ \Big\|\frac{1}{n}\sum_{i=1}^{n}\xi_i\hat{U}_1^{(0)\top} \calM_1(\calX_i)\left(\big(\hat{U}_2^{(0)} - U_2R_2^\top\big) \otimes \hat{U}_3^{(0)}\right)\Big\|\notag\\
			& + \Big\|\frac{1}{n}\sum_{i=1}^{n}\xi_i\hat{U}_1^{(0)\top} \calM_1(\calX_i)\left(\hat{U}_2^{(0)} \otimes \big(\hat{U}_3^{(0)} - U_3R_3^\top\big)\right)\Big\|\notag\\
			\leq& C\sqrt{\frac{pr}{n}}\left(\big\|\hat{U}_1^{(0)} - U_1R_1^\top\big\| + \big\|\hat{U}_2^{(0)} - U_2R_2^\top\big\| + \big\|\hat{U}_3^{(0)} - U_3R_3^\top\big\|\right)
			\leq C_2\frac{p\sqrt{r}}{n\lambda_{\submin}}.
\end{align}
In addition, by \eqref{ineq14}, with probability at least $1 - C_1e^{-c_1p}$,
\begin{align}\label{ineq54}
			&\big\|\hat{U}_1^{(0)\top} A_{T_1}^{(0)}(\hat{U}_2^{(0)} \otimes \hat{U}_3^{(0)})\big\|\notag\\ 
			= & \Big\|R_1G_1\big(R_2^\top \otimes R_3^\top\big) - \big(\hat U_1^{(0)\top}U_1\big)G_1\Big(\big(\hat U_2^{(0)\top}U_2\big)^\top \otimes \big(\hat U_3^{(0)\top}U_3\big)^\top\Big)\Big\|\notag\\
			\leq& \Big\|\big(R_1 - \hat U_1^{(0)\top}U_1\big)G_1\big(R_2^\top \otimes R_3^\top\big)\Big\| + \Big\|\big(\hat U_1^{(0)\top}U_1\big)G_1\Big(\big(R_2 - \hat U_2^{(0)\top}U_2\big)^\top \otimes R_3^\top\Big)\Big\|\notag\\
			&+ \Big\|\big(\hat U_1^{(0)\top}U_1\big)G_1\Big(\big(\hat U_2^{(0)\top}U_2\big)^\top \otimes \big(R_3 - \hat U_3^{(0)\top}U_3\big)^\top\Big)\Big\|\notag\\
			\leq& \left(\big\|R_1 - \hat U_1^{(0)\top}U_1\big\| + \big\|R_2 - \hat U_2^{(0)\top}U_2\big\| + \big\|R_3 - \hat U_3^{(0)\top}U_3\big\|\right)\kappa_0\lambda_{\submin}\notag\\
			\leq& C_2\frac{p}{n\lambda_{\submin}^2}\kappa_0\lambda_{\submin} = C_2\kappa_0\frac{p}{n\lambda_{\submin}}.
\end{align}
Putting \eqref{ineq45}, \eqref{ineq50}, \eqref{ineq52}, \eqref{ineq53} and \eqref{ineq54} and Lemma \ref{lm:Gaussian_ensemble} together, we get with probability $1 - p^{-3} - C_1e^{-c_1p}$ that
\begin{align*}
			\|\Delta G_1\| \leq& \|\frakJ_1\| + C_2\left(\sqrt{\frac{pr}{n}}\|\Delta \calG\|_{\rm F} + \kappa_0\frac{pr}{n} + \frac{p\sqrt{r}}{n\lambda_{\submin}} + \kappa_0\frac{p}{n\lambda_{\submin}}\right)\\
			 \leq& C_2\sqrt{\frac{pr^2}{n}}\|\Delta G_1\| + C_2\left(\sqrt{\frac{r^2+\log(p)}{n}} + \kappa_0\frac{pr}{n} + \kappa_0\frac{p\sqrt{r}}{n\lambda_{\submin}}\right)\\
			\leq& \frac{1}{2}\|\Delta G_1\| + C_2\left(\sqrt{\frac{r^2+\log(p)}{n}} + \kappa_0\frac{pr}{n} + \kappa_0\frac{p\sqrt{r}}{n\lambda_{\submin}}\right)
\end{align*}
and as a result
\begin{equation}\label{ineq57}
			\|\Delta G_1\| \leq C_2\left(\sqrt{\frac{r^2+\log(p)}{n}} + \kappa_0\frac{pr}{n} + \kappa_0\frac{p\sqrt{r}}{n\lambda_{\submin}}\right)
\end{equation}
Therefore, with probability at least $1 - p^{-3} - C_1e^{-c_1p}$,
\begin{equation*}
			\|\Delta G_1 - \frakJ_1\| \leq C_2\left(\sqrt{\frac{pr^2}{n}}\|\Delta G_1\| + \kappa_0\frac{pr}{n} + \kappa_0\frac{p\sqrt{r}}{n\lambda_{\submin}}\right) \leq C_2\left(\kappa_0\frac{pr}{n} + \kappa_0\frac{p\sqrt{r}}{n\lambda_{\submin}}\right)
\end{equation*}
which proves (\ref{ineq:error_G}).

\subparagraph*{Step 0.2: the error of $\hat U_j^{(t+0.5)}$.} Without loss of generality, we only prove the bound for $t=0$ and $j=1$. Again, we denote $\hat\calG=\hat\calG^{(0)}$ and $\hat G_1=\calM_1(\hat\calG)$ for brevity.

We aim to show that with probability at least $1 - C_1p^{-3} - C_1e^{-c_1p}$,
\begin{equation}\label{eq:regression_step_2}
	\hat U_1^{(0.5)} = U_1R_1^\top + \frac{1}{n}\calP_{U_1}^\perp\Big(\sum_{j=1}^{n}\xi_j\calM_1(\calX_j)\Big)(U_2 \otimes U_3)G_1^\top(G_1G_1^\top)^{-1}R_1^\top + \frakE,
\end{equation}
where $\|\frakE\| \leq C_2\left(\kappa_0\frac{pr}{n\lambda_{\submin}} + \kappa_0^2\frac{p\sqrt{r}}{n\lambda_{\submin}^2}\right)$.

Since $\frac{\partial}{\partial U_1}\ell_n(\hat{\calG}\times_1 \hat{U}_1^{(0.5)} \times_2 \hat{U}_2^{(0)} \times_3 \hat{U}_3^{(0)}) = 0$, we have
\begin{equation*}
	2\sum_{i=1}^{n}\calM_1(\calX_i)(\hat U_2^{(0)} \otimes \hat U_3^{(0)})\hat G_1^\top\langle \hat U_1^{(0.5)}, \calM_1(\calX_i)(\hat U_2^{(0)} \otimes \hat U_3^{(0)})\hat G_1^\top\rangle - 2\sum_{i=1}^{n}\calM_1(\calX_i)(\hat U_2^{(0)} \otimes \hat U_3^{(0)})\hat G_1^\top Y_i = 0.
\end{equation*}
Denote
\begin{equation*}
	\Delta \calT_{1} = \hat \calG \times_1 \hat U_1^{(0.5)} \times_2 \hat U_2^{(0)}  \times_3 \hat U_3^{(0)} - \calT.
\end{equation*}
For brevity, denote $\hat U_1=\hat U_1^{(0.5)}$ for simplicity. Then, we write
\begin{align*}
	&\hat{U}_1(\hat G_1\hat G_1^\top) - U_1R_1^\top(\hat G_1\hat G_1^\top)\\
	=& \bigg(\calM_1(\Delta \calT_{1})(\hat U_2^{(0)} \otimes \hat U_3^{(0)}) - \frac{1}{n}\sum_{i=1}^{n}\langle \Delta \calT_{1}, \calX_i\rangle\calM_1(\calX_i)(\hat U_2^{(0)} \otimes \hat U_3^{(0)})\bigg)\hat G_1^\top\\
	&+ \Big(U_1G_1\big((\hat U_2^{(0)\top}U_2)^\top \otimes (\hat U_3^{(0)\top}U_3)^\top\big)\hat G_1^\top - U_1R_1^\top\big(\hat G_1\hat G_1^\top\big) \Big)\\
	&+ \frac{1}{n}\sum_{i=1}^{n}\xi_i\calM_1(\calX_i)(\hat U_2^{(0)} \otimes \hat U_3^{(0)})\hat G_1^\top,
\end{align*}
which is equivalent to 
\begin{align}\label{ineq55}
	&\hat{U}_1 - U_1R_1^\top\notag\\
	=& \bigg(\calM_1(\Delta \calT_{1})(\hat U_2^{(0)} \otimes \hat U_3^{(0)}) - \frac{1}{n}\sum_{i=1}^{n}\langle \Delta \calT_{1}, \calX_i\rangle\calM_1(\calX_i)(\hat U_2^{(0)} \otimes \hat U_3^{(0)})\bigg)\hat G_1^\top\big(\hat G_1\hat G_1^\top\big)^{-1}\notag\\
	&+\Big(U_1G_1\big((\hat U_2^{(0)\top}U_2)^\top \otimes (\hat U_3^{(0)\top}U_3)^\top\big)\hat G_1^\top - U_1R_1^\top\big(\hat G_1\hat G_1^\top\big) \Big)\hat G_1^\top\big(\hat G_1\hat G_1^\top\big)^{-1}\notag\\
	&+ \frac{1}{n}\sum_{i=1}^{n}\xi_i\calM_1(\calX_i)(\hat U_2^{(0)} \otimes \hat U_3^{(0)})\hat G_1^\top\big(\hat G_1\hat G_1^\top\big)^{-1}\notag\\
	=:& \text{\rom{1}} + \text{\rom{2}} + \text{\rom{3}}.
\end{align}
By \eqref{ineq:error_G} and Lemma \ref{lm:Gaussian_ensemble}, with probability at least $1 - C_1p^{-3} - C_1e^{-c_1p}$, 
\begin{align*}
	&\bigg\| \bigg(\calM_1(\Delta \calT_{1})(\hat U_2^{(0)} \otimes \hat U_3^{(0)}) - \frac{1}{n}\sum_{i=1}^{n}\langle \Delta \calT_{1}, \calX_i\rangle\calM_1(\calX_i)(\hat U_2^{(0)} \otimes \hat U_3^{(0)})\bigg)\hat G_1^\top\big(\hat G_1\hat G_1^\top\big)^{-1}\bigg\|\\
	\leq& C\sqrt{\frac{pr}{n}}\|\Delta \calT_{1}\|_{\rm F}\Big\|\hat G_1^\top\big(\hat G_1\hat G_1^\top\big)^{-1}\Big\|\\
	\leq& C\sqrt{\frac{pr}{n}}\left\|\hat \calG \times_1 \hat U_1 \times_2 \hat U_2^{(0)} \times_3 \hat U_3^{(0)} - \calT\right\|_{\rm F}\bigg(\lambda_{\submin} - C_2\Big(\sqrt{\frac{r^2+\log(p)}{n}} + \kappa_0\frac{pr}{n} + \kappa_0\frac{p\sqrt{r}}{n\lambda_{\submin}}\Big)\bigg)^{-1}\\
	\leq& C\sqrt{\frac{pr}{n}}\lambda_{\submin}^{-1}\left\|\hat \calG \times_1 \hat U_1 \times_2 \hat U_2^{(0)} \times_3 \hat U_3^{(0)} - \calT\right\|_{\rm F}.
\end{align*}
Moreover, with probability at least $1 - C_1p^{-3} - C_1e^{-c_1p}$, 
\begin{align}
	&\big\|\hat \calG \times_1 \hat U_1 \times_2 \hat U_2^{(0)} \times_3 \hat U_3^{(0)} - \calT\big\|_{\rm F}\notag\\
	\leq& \big\|\big(\hat \calG - \calG \times_1 R_1 \times_2 R_2 \times_3 R_3\big) \times_1 \hat U_1 \times_2 \hat U_2^{(0)} \times_3 \hat U_3^{(0)}\big\|_{\rm F}\notag\\ 
	&+ \big\|\calG \times_1 \big(\hat U_1R_1\big) \times_2 \big(\hat U_2^{(0)}R_2\big) \times_3 \big(\hat U_3^{(0)}R_3\big) - \calG \times_1 U_1 \times_2 U_2 \times_3 U_3\big\|_{\rm F}\notag\\
	\leq& \sqrt{r}\big\|\hat \calG - \calG \times_1 R_1 \times_2 R_2 \times_3 R_3\big\| + \sqrt{r}\big\|U_1 - \hat U_1R_1\big\|\|G_1\| + \sqrt{r}\big\|U_2 - \hat U_2^{(0)}R_2\big\|\|G_2\|\notag\\ 
	&+ \sqrt{r}\big\|U_3 - \hat U_3^{(0)}R_3\big\|\|G_3\|\notag\\
	\leq& C_2\sqrt{r}\bigg(\sqrt{\frac{r^2+\log(p)}{n}} + \kappa_0\frac{pr}{n} + \kappa_0\frac{p\sqrt{r}}{n\lambda_{\submin}}\bigg) + \sqrt{r}\kappa_0\lambda_{\submin}\big\|U_1 - \hat U_1R_1\big\| + C_2\kappa_0\sqrt{\frac{pr}{n}}\notag\\
	\leq& C_2\bigg(\kappa_0\sqrt{\frac{pr}{n}} + \kappa_0\frac{pr^{3/2}}{n} + \kappa_0\frac{pr}{n\lambda_{\submin}}\bigg) + \sqrt{r}\kappa_0\lambda_{\submin}\big\|U_1 - \hat U_1R_1\big\|.\label{ineq13}
\end{align}
The two previous inequalities together imply that with probability at least $1 - C_1p^{-3} - C_1e^{-c_1p}$, 
\begin{align}
	\left\|\text{\rom{1}}\right\|=&\bigg\|\bigg(\calM_1(\Delta \calT_{1})(\hat U_2^{(0)} \otimes \hat U_3^{(0)}) - \frac{1}{n}\sum_{i=1}^{n}\langle \Delta \calT_{1}, \calX_i\rangle\calM_1(\calX_i)(\hat U_2^{(0)} \otimes \hat U_3^{(0)})\bigg)\hat G_1^\top\big(\hat G_1\hat G_1^\top\big)^{-1}\bigg\|\notag\\
	\leq& C_2\left(\frac{\kappa_0pr}{n\lambda_{\submin}} + \kappa_0\frac{p^{3/2}r^2}{n^{3/2}\lambda_{\submin}} + \kappa_0\frac{p^{3/2}r^{3/2}}{n^{3/2}\lambda_{\submin}^2}\right) + C_2\kappa_0r\sqrt{\frac{p}{n}}\big\|U_1 - \hat U_1R_1\big\|.\label{ineq58}
\end{align}
By \eqref{ineq57}, with probability at least $1 - C_1p^{-3} - C_1e^{-c_1p}$, we have
\begin{align}\label{ineq59}
	&\Big\|\hat G_1^\top\big(\hat G_1\hat G_1^\top\big)^{-1} - (R_2 \otimes R_3)G_1^\top(G_1G_1^\top)^{-1}R_1^\top\Big\|\notag\\
	\leq& \Big\|\Big(\hat G_1 - R_1G_1(R_2^\top \otimes R_3^\top)\Big)^\top\big(\hat G_1\hat G_1^\top\big)^{-1}\Big\|\notag\\ 
	&+ \Big\|(R_2 \otimes R_3)G_1^\top R_1^\top\left(R_1(G_1G_1^\top)^{-1}R_1^\top - \big(\hat G_1\hat G_1^\top\big)^{-1}\right)\Big\|\notag\\
	\leq& C_2\bigg(\sqrt{\frac{r^2+\log(p)}{n}} + \kappa_0\frac{pr}{n} + \kappa_0\frac{p\sqrt{r}}{n\lambda_{\submin}}\bigg)\lambda_{\submin}^{-2}\notag\\ 
	&+ \|G_1\|\big\|(G_1G_1^\top)^{-1}\big\|\big\|\big(\hat G_1\hat G_1^\top\big)^{-1}\big\|\big\|R_1G_1G_1^\top R_1^\top - \hat G_1\hat G_1^\top\big\|\notag\\
	\leq& C_2\bigg(\sqrt{\frac{r^2+\log(p)}{n}} + \kappa_0\frac{pr}{n} + \kappa_0\frac{p\sqrt{r}}{n\lambda_{\submin}}\bigg)\lambda_{\submin}^{-2}\notag\\ 
	&+ C_2\kappa_0\lambda_{\submin}^{-3}\Big(\left\|\left(R_1G_1(R_2^\top \otimes R_3^\top) - \hat G_1\right)(R_2 \otimes R_3)G_1^\top R_1^\top\right\| + 
	\Big\|\hat{G}_1\left(R_1G_1(R_2^\top \otimes R_3^\top) - \hat G_1\right)^\top\Big\|\Big)\notag\\
	\leq& C_2\kappa_0^2\lambda_{\submin}^{-2}\bigg(\sqrt{\frac{r^2+\log(p)}{n}} + \kappa_0\frac{pr}{n} + \kappa_0\frac{p\sqrt{r}}{n\lambda_{\submin}}\bigg),
\end{align}
and 
\begin{align*}
	&\Big\|U_1G_1\Big((\hat U_2^{(0)\top}U_2)^\top \otimes (\hat U_3^{(0)\top}U_3)^\top\Big) - U_1R_1^\top\hat G_1\Big\|\\
	\leq& \left\|G_1\big((\hat U_2^{(0)\top}U_2 - R_2)^\top \otimes (\hat U_3^{(0)\top}U_3)^\top\big)\right\| + \left\|G_1\big(R_2^\top \otimes (\hat U_3^{(0)\top}U_3 -R_3)^\top\big)\right\|\\ &
	+ \big\|R_1G_1\big(R_2^\top \otimes R_3^\top\big) - \hat G_1\big\|\\
	\leq& C_2\kappa_0\lambda_{\submin}\bigg(\frac{\sqrt{p/n}}{\lambda_{\submin}}\bigg)^2 + C_2\bigg(\sqrt{\frac{r^2+\log(p)}{n}} + \kappa_0\frac{pr}{n} + \kappa_0\frac{p\sqrt{r}}{n\lambda_{\submin}}\bigg)\\
	\leq& C_2\bigg(\sqrt{\frac{r^2+\log(p)}{n}} + \kappa_0\frac{pr}{n} + \kappa_0\frac{p\sqrt{r}}{n\lambda_{\submin}}\bigg).
\end{align*}
By the two previous inequality and \eqref{ineq:error_G}, with probability at least $1 - C_1p^{-3} - C_1e^{-c_1p}$, 
\begin{align}\label{ineq60}
	&\bigg\|\text{\rom{2}} + \frac{1}{n}\calP_{U_1}\Big(\sum_{j=1}^{n}\xi_j\calM_1(\calX_j)\Big)(U_2 \otimes U_3)G_1^\top(G_1G_1^\top)^{-1}R_1^\top\bigg\|\notag\\
	\leq& \left\|U_1G_1\left((\hat U_2^{(0)\top}U_2)^\top \otimes (\hat U_3^{(0)\top}U_3)^\top\right) - U_1R_1^\top\hat G_1\right\|\Big\|\hat G_1^\top\big(\hat G_1\hat G_1^\top\big)^{-1} - (R_2 \otimes R_3)G_1^\top(G_1G_1^\top)^{-1}R_1^\top\Big\|\notag\\
	&+ \bigg\|U_1G_1\left(R_2^\top \otimes R_3^\top\right) - U_1R_1^\top\hat G_1 + \frac{1}{n}\calP_{U_1}\Big(\sum_{j=1}^{n}\xi_j\calM_1(\calX_j)\Big)\big((U_2R_2^\top) \otimes (U_3R_3^\top)\big)\bigg\|\notag\\
	&\cdot\big\|(R_2 \otimes R_3)G_1^\top(G_1G_1^\top)^{-1}R_1^\top\big\|\notag\\
	&+\left(\left\|G_1\big((\hat U_2^{(0)\top}U_2 - R_2)^\top \otimes (\hat U_3^{(0)\top}U_3)^\top\big)\right\| + \left\|G_1\big(R_2^\top \otimes (\hat U_3^{(0)\top}U_3 -R_3)^\top\big)\right\|\right)\notag\\
	&\cdot\big\|(R_2 \otimes R_3)G_1^\top(G_1G_1^\top)^{-1}R_1^\top\big\|\notag\\
	\leq& C_2\kappa_0^2\lambda_{\submin}^{-2}\left(\frac{r^2+\log(p)}{n} + \kappa_0^2\frac{p^2r^2}{n^2} + \kappa_0^2\frac{p^2r}{n^2\lambda_{\submin}^2}\right) + C_2\kappa_0\left(\frac{pr}{n} + \frac{p\sqrt{r}}{n\lambda_{\submin}}\right)\lambda_{\submin}^{-1} + C_2\kappa_0\lambda_{\submin}\left(\frac{\sqrt{p/n}}{\lambda_{\submin}}\right)^2\lambda_{\submin}^{-1}\notag\\
	\leq& C_2\left(\kappa_0\frac{pr}{n\lambda_{\submin}} + \kappa_0^2\frac{p\sqrt{r}}{n\lambda_{\submin}^2}\right).
\end{align}
For term \rom{3}, by \eqref{ineq59} and Lemma \ref{lm:sub_Gaussian} Part 3, with probability $1 - C_1p^{-3} - C_1e^{-c_1p}$ that
\begin{align}\label{ineq61}
	&\bigg\|\text{\rom{3}} - \frac{1}{n}\Big(\sum_{j=1}^{n}\xi_j\calM_1(\calX_j)\Big)(U_2 \otimes U_3)G_1^\top(G_1G_1^\top)^{-1}R_1^\top\bigg\|\notag\\
	\leq& \bigg\|\frac{1}{n}\sum_{i=1}^{n}\xi_i\calM_1(\calX_i)(\hat U_2^{(0)} \otimes \hat U_3^{(0)})\Big(\hat G_1^\top\big(\hat G_1\hat G_1^\top\big)^{-1} - (R_2 \otimes R_3)G_1^\top(G_1G_1^\top)^{-1}R_1^\top\Big)\bigg\|\notag\\
	&+ \bigg\|\frac{1}{n}\Big(\sum_{j=1}^{n}\xi_j\calM_1(\calX_j)\Big)\big((\hat U_2^{(0)}R_2 - U_2) \otimes (\hat U_3^{(0)}R_3)\big)G_1^\top(G_1G_1^\top)^{-1}R_1^\top\bigg\|\notag\\
	&+ \bigg\|\frac{1}{n}\Big(\sum_{j=1}^{n}\xi_j\calM_1(\calX_j)\Big)\big(U_2 \otimes (\hat U_3^{(0)}R_3 - U_3)\big)G_1^\top(G_1G_1^\top)^{-1}R_1^\top\bigg\|\notag\\
	\leq& C_2\sqrt{\frac{pr}{n}}\Big\|\hat G_1^\top\big(\hat G_1\hat G_1^\top\big)^{-1} - (R_2 \otimes R_3)G_1^\top(G_1G_1^\top)^{-1}R_1^\top\Big\|\notag\\
	&+ C_2\sqrt{\frac{pr}{n}}\big\|\hat U_2^{(0)}R_2 - U_2\big\|\lambda_{\submin}^{-1} + C_2\sqrt{\frac{pr}{n}}\big\|\hat U_3^{(0)}R_3 - U_3\big\|\lambda_{\submin}^{-1}
	\leq C_2\kappa_0^2\frac{p\sqrt{r}}{n\lambda_{\submin}^2}.
\end{align}
Putting \eqref{ineq55}, \eqref{ineq58}, \eqref{ineq60} and \eqref{ineq61} and Lemma \ref{lm:sub_Gaussian} Part 3 together, we get with probability at least $1 - C_1p^{-3} - C_1e^{-c_1p}$ that
\begin{align*}
	\big\|\hat U_1 - U_1R_1^\top\big\| \leq& \bigg\|\frac{1}{n}\calP_{U_1}^\perp\Big(\sum_{j=1}^{n}\xi_j\calM_1(\calX_j)\Big)(U_2 \otimes U_3)G_1^\top(G_1G_1^\top)^{-1}R_1^\top\bigg\|\\
	 &+  C_2\kappa_0r\sqrt{\frac{p}{n}}\big\|\hat U_1 - U_1R_1^\top\big\| + C_2\left(\kappa_0\frac{pr}{n\lambda_{\submin}} + \kappa_0^2\frac{p\sqrt{r}}{n\lambda_{\submin}^2}\right)\\
	\leq& C_2\sqrt{\frac{p}{n}}\lambda_{\submin}^{-1} + \frac{1}{2}\big\|\hat U_1 - U_1R_1^\top\big\| + C_2\left(\kappa_0\frac{pr}{n\lambda_{\submin}} + \kappa_0^2\frac{p\sqrt{r}}{n\lambda_{\submin}^2}\right)\\
	\leq& \frac{1}{2}\big\|\hat U_1 - U_1R_1^\top\big\| + C_2\sqrt{\frac{p}{n}}\lambda_{\submin}^{-1}.
\end{align*}
Therefore, with probability at least $1 - C_1p^{-3} - C_1e^{-c_1p}$,
\begin{equation}\label{ineq67}
	\big\|\hat U_1 - U_1R_1^\top\big\| \leq C_2\sqrt{\frac{p}{n}}\lambda_{\submin}^{-1}.
\end{equation}
Thus with probability at least $1 - C_1p^{-3} - C_1e^{-c_1p}$,
\begin{align*}
	&\bigg\|\hat U_1^{(0.5)} - U_1R_1^\top - \frac{1}{n}\calP_{U_1}^\perp\Big(\sum_{j=1}^{n}\xi_j\calM_1(\calX_j)\Big)(U_2 \otimes U_3)G_1^\top(G_1G_1^\top)^{-1}R_1^\top\bigg\|\\
	\leq& C_2\left(\kappa_0\frac{pr}{n\lambda_{\submin}} + \kappa_0^2\frac{p\sqrt{r}}{n\lambda_{\submin}^2}\right) + C_2\kappa_0r\sqrt{\frac{p}{n}}\cdot\sqrt{\frac{p}{n}}\lambda_{\submin}^{-1}\\
	\leq& C_2\left(\kappa_0\frac{pr}{n\lambda_{\submin}} + \kappa_0^2\frac{p\sqrt{r}}{n\lambda_{\submin}^2}\right).
\end{align*}

Now, we continue from eq. (\ref{eq:deltaTt+0.5}) and prove the distribution of $\|\hat U_1^{(2)}\hat U_1^{(2)\top}-U_1U_1^{\top}\|_{\rm F}^2$.
\paragraph*{Step 1: bounding $\|\hat U_j^{(1)\top}U_j-R_j^{(1)}\|_{\rm F}$ and $\|U_j-\hat U_j^{(1)} R_j^{(1)}\|$.} Without loss of generality, we only prove the bound for $j=1$. By definition of $\hat U_1^{(t+0.5)}$ in Algorithm~\ref{algo:am_optimal_regression}, we write
\begin{align}
		&\hat{U}_1^{(t+0.5)} - U_1R_1^{(t)\top}\notag\\
		=& \Big(\calM_1(\Delta \calT_{1}^{(t + 0.5)})(\hat U_2^{(t)} \otimes \hat U_3^{(t)}) - \frac{1}{n}\sum_{i=1}^{n}\langle \Delta \calT_{1}^{(t + 0.5)}, \calX_i\rangle\calM_1(\calX_i)(\hat U_2^{(t)} \otimes \hat U_3^{(t)})\Big)\hat G_1^{(t)\top}\big(\hat G_1^{(t)}\hat G_1^{(t)\top}\big)^{-1}\notag\\
		&+ \left(U_1G_1\big((\hat U_2^{(t)\top}U_2)^\top \otimes (\hat U_3^{(t)\top}U_3)^\top\big) - U_1R_1^{(t)\top}\hat G_1^{(t)}\right)\hat G_1^{(t)\top}\big(\hat G_1^{(t)}\hat G_1^{(t)\top}\big)^{-1}\notag\\
		&+ \frac{1}{n}\sum_{i=1}^{n}\xi_i\calM_1(\calX_i)(\hat U_2^{(t)} \otimes \hat U_3^{(t)})\hat G_1^{(t)\top}\big(\hat G_1^{(t)}\hat G_1^{(t)\top}\big)^{-1}\notag\\
		=:& \frakJ_{U_1, 1}^{(t)} + \frakJ_{U_1, 2}^{(t)} + \frakJ_{U_1, 3}^{(t)}.\label{ineq62}
\end{align}
Denote $\frakE_1^{(t)} = \frakJ_{U_1, 1}^{(t)} + \frakJ_{U_1, 2}^{(t)} + \frakJ_{U_1, 3}^{(t)}$.  Recall that $\hat U_1^{(t+1)}$ are the left singular vectors of $\hat U_1^{(t+0.5)}$. We can also apply the spectral representation formula (Lemma~\ref{lem:spectral}) to investigate $\hat U^{(t+1)}_1$. Toward that end, we define
\begin{equation*}
\begin{pmatrix}
		0 & \hat{U}_1^{(t+0.5)}\\
		\hat{U}_1^{(t+0.5)\top} & 0
		\end{pmatrix} = 
		\begin{pmatrix}
		0 & U_1R_1^{(t)\top}\\
		R_1^{(t)}U_1^\top & 0
		\end{pmatrix} +
		\begin{pmatrix}
		0 & \frakE_1^{(t)}\\
		\frakE_1^{(t)\top} & 0
\end{pmatrix}.
\end{equation*}
Note that the non-zero eigenvalues of the symmetric matrix 
$$
\begin{pmatrix}
		0 & U_1R_1^{(t)\top}\\
		R_1^{(t)}U_1^\top & 0
\end{pmatrix}
$$
are $\mu_1^{(t)} = \dots = \mu_{r_1}^{(t)} = 1$ and $\mu_{r_1 +1}^{(t)} = \cdots = \mu_{2r_1}^{(t)} = -1$, and for $1 \leq i \leq r_1$, the corresponding eigenvectors of $\mu_i^{(t)}$ and $\mu_{r_1 + i}^{(t)}$ are
\begin{equation*}
		\theta_i^{(t)} = \frac{1}{\sqrt{2}}\begin{pmatrix}
		\bar{u}_i^{(t)}\\
		e_i
		\end{pmatrix} \quad \text{and} \quad \theta_{r_1 + i}^{(t)} = \frac{1}{\sqrt{2}}\begin{pmatrix}
		\bar{u}_i^{(t)}\\
		-e_i
		\end{pmatrix},
\end{equation*}
where $\bar{u}_i^{(t)}$ is the $i$-th column of $U_1R_1^{(t)\top}$ and $e_i$ is the $i$-th canonical basis of $\RR^{r_1}$.

Denote a $(p_1+r_1)\times (2r_1)$ matrix
$$
\Theta^{(t)} = \left(\theta_1^{(t)} \ \dots \theta_{2r_1}^{(t)}\right)
$$
 and $\Theta_{\perp}^{(t)} \in \OO_{p_1 + r_1, p_1 - r_1}$ such that $\left(\Theta^{(t)} \ \Theta_{\perp}^{(t)}\right) \in \OO_{p_1 + r_1}$. Then, we write
$$
\Theta^{(t)}\Theta^{(t)\top} = \sum_{1 \leq j \leq 2r_1}\theta_j^{(t)}\theta_j^{(t)\top} = 
\begin{pmatrix}
		U_1U_1^\top & 0\\
		0 & I_{r_1}
\end{pmatrix}
$$ 
For $k \geq 1$, denote
$$
\left(\frakP_1^{(t)}\right)^{-k} = \sum_{1 \leq j \leq 2r_1}\frac{1}{(\mu_j^{(t)})^k}\theta_j^{(t)}\theta_j^{(t)\top} = \begin{cases}
		\begin{pmatrix}
		0 & U_1R_1^{(t)\top}\\
		R_1^{(t)}U_1^\top & 0
		\end{pmatrix}, & \text{ if } k \text{ is old},\\
		\begin{pmatrix}
		UU^\top & 0\\
		0 & I_{r_1}
		\end{pmatrix}, &\text{ if } k \text{ is even},
		\end{cases}
$$
and
$$
\left(\frakP_1^{(t)}\right)^{0} = \Theta_{\perp}^{(t)}\Theta_{\perp}^{(t)\top} = \begin{pmatrix}
		U_{1\perp}U_{1\perp}^\top & 0\\
		0 & 0
		\end{pmatrix}.
$$
Let 
$$
E^{(t)} = \begin{pmatrix}
		0 & \frakE_1^{(t)}\\
		\frakE_1^{(t)\top} & 0
\end{pmatrix}.
$$
By Lemma \ref{lm:sub_Gaussian} and together with \eqref{eq:regression_step_2}, with probability at least $1 - C_1p^{-3} - C_1e^{-c_1p}$,
\begin{equation}\label{ineq2}
		\left\|E^{(0)}\right\| = \left\|\frakE_1^{(0)}\right\| \leq C_2\frac{\sqrt{p/n}}{\lambda_{\submin}} + C_2\left(\kappa_0\frac{pr}{n\lambda_{\submin}} + \kappa_0^2\frac{p\sqrt{r}}{n\lambda_{\submin}^2}\right) \leq C_2\frac{\sqrt{p/n}}{\lambda_{\submin}} < \frac{1}{8}.
\end{equation}
By Lemma \ref{lem:spectral}, with probability at least $1 - C_1p^{-3} - C_1e^{-c_1p}$,
\begin{equation}\label{ineq:representation_regression}
		\begin{pmatrix}
		\hat U_1^{(1)}\hat U_1^{(1)\top} & 0\\
		0 & I_{r_1}
		\end{pmatrix} - \begin{pmatrix}
		U_1U_1^{\top} & 0\\
		0 & I_{r_1}
		\end{pmatrix} = \sum_{k \geq 1}\calS_{U_1^{(0)}, k}\big(E^{(0)}\big)
\end{equation}
where 
$$
\calS_{U_1^{(t)}, k}(X) = \sum_{s_1+\cdots+s_{k+1}=k}(-1)^{1+\tau(\bs)}\cdot \big(\frakP_1^{(t)}\big)^{-s_1}X\big(\frakP_1^{(t)}\big)^{-s_2}X\big(\frakP_1^{(t)}\big)^{-s_3}\cdots\big(\frakP_1^{(t)}\big)^{-s_k}X\big(\frakP_1^{(t)}\big)^{-s_{k+1}}
$$
where $s_1,\cdots,s_{k+1}$ are non-negative integers and $\tau(\bs)=\sum_{j=1}^{k+1}\II(s_j>0)$.

Clearly, we have 
\begin{equation}\label{ineq63}
		\left\|\calS_{U_1^{(0)}, k}\big(E^{(0)}\big)\right\| \leq \binom{2k}{k}\|E^{(0)}\|^{k-1}\|E^{(0)}\|_{\rm F} \leq \big(4\|E^{(0)}\|\big)^{k}.
\end{equation}
By \eqref{ineq2}, \eqref{ineq:representation_regression} and \eqref{ineq63}, with probability at least $1 - C_1p^{-3} -C_1e^{-c_1p}$,
\begin{align*}
	\big\|\hat U_1^{(1)}\hat U_1^{(1)\top} - U_1U_1^\top\big\| \leq \sum_{k \geq 1}\big\|\calS_{U_1^{(0)}, k}\big(E^{(0)}\big)\big\| \leq C_2\frac{\sqrt{p/n}}{\lambda_{\submin}}.
\end{align*}

Thus with probability at least $1 - C_1p^{-3} - C_1e^{-c_1p}$,
\begin{equation}\label{ineq70}
		\left\|\hat U_1^{(1)\top}U_1 - R_1^{(1)}\right\| \leq C_2\frac{p}{n\lambda_{\submin}^2}, \quad \left\|\hat U_1^{(1)\top}U_1 - R_1^{(1)}\right\|_{\rm F} \leq  C_2\frac{p\sqrt{r}}{n\lambda_{\submin}^2}
\end{equation}
and
\begin{equation}\label{ineq71}
		\left\|U_1 - \hat U_1^{(1)}R_1^{(1)}\right\| \leq C_2\frac{\sqrt{p/n}}{\lambda_{\submin}}, \quad \left\|U_1 - \hat U_1^{(1)}R_1^{(1)}\right\|_{\rm F} \leq C_2\frac{\sqrt{pr/n}}{\lambda_{\submin}}.
\end{equation}

\paragraph*{Step 2: representation of $\|\hat U_1^{(2)}\hat U_1^{(2)\top}-U_1U_1^{\top}\|_{\rm F}^2$ and its first order approximation. } For convenience, we denote 
$$
E = E^{(1)},\quad \frakE_1 = \frakE_1^{(1)}\quad{\rm and}\quad \frakP_1^{-k} = \big(\frakP_1^{(1)}\big)^{-k}.
$$
Similarly to {\it Step 1}, we apply Lemma \ref{lem:spectral} to $\hat U_1^{(2)}\hat U_1^{(2)\top}$ and get with probability at least $1 - C_1p^{-3} - C_1e^{-c_1p}$ that
\begin{align}
&\big\|\hat U_1^{(2)}\hat U_1^{(2)\top} - U_1U_1^\top\big\|_{\rm F}^2\notag\\ 
=& -2\left\langle \begin{pmatrix}
		U_1U_1^\top & 0\\
		0 & I_{r_1}
		\end{pmatrix},
		\begin{pmatrix}
		\hat U_1^{(2)}\hat U_1^{(2)\top} - U_1U_1^\top & 0\\
		0 & 0
		\end{pmatrix}
		\right\rangle\notag\\		
=& -2\left\langle \begin{pmatrix}
		U_1U_1^\top & 0\\
		0 & I_{r_1}
		\end{pmatrix}, 
		\sum_{k \geq 1}\calS_{U_1^{(1)}, k}\left(E\right)
		\right\rangle\notag\\
=& -2\left\langle \begin{pmatrix}
		U_1U_1^\top & 0\\
		0 & I_{r_1}
		\end{pmatrix}, 
		\sum_{k \geq 4}\calS_{U_1^{(1)}, k}\left(E\right)
		\right\rangle - 2\left\langle \begin{pmatrix}
		U_1U_1^\top & 0\\
		0 & I_{r_1}
		\end{pmatrix}, 
		\calS_{U_1^{(1)}, 2}\left(E\right)
		\right\rangle\notag\\
& - 2\left\langle \begin{pmatrix}
		U_1U_1^\top & 0\\
		0 & I_{r_1}
		\end{pmatrix}, 
		\calS_{U_1^{(1)}, 3}\left(E\right)
		\right\rangle \label{eq3}
\end{align}
where we use the fact $\left\langle \begin{pmatrix}
U_1U_1^\top & 0\\
0 & I_{r_1}
\end{pmatrix}, 
\calS_{U_1^{(1)}, 1}\left(E\right)
\right\rangle = 0.$\\
Similarly to \eqref{ineq63}, we know that with probability at least $1 - C_1p^{-3} - C_1e^{-c_1p}$,
\begin{align}
		\left|2\left\langle 
		\begin{pmatrix}
		U_1U_1^\top & 0\\
		0 & I_{r_1}
		\end{pmatrix}, 
		\sum_{k \geq 4}\calS_{U_1^{(1)}, k}\left(E\right)
		\right\rangle\right| \leq 4r\sum_{k \geq 4}\left\|\calS_{U_1^{(1)}, k}\left(E\right)\right\| \leq &C_2r\left(\frac{\sqrt{p/n}}{\lambda_{\submin}}\right)^4 \notag\\
		&\leq C_2r\frac{p^2}{n^2\lambda_{\submin}^4}.\label{ineq64}
\end{align}
We now bound the third order term. 	Notice that
\begin{align*}
		\frakP_1^{0}E\frakP_1^{0} = \begin{pmatrix}
		U_{1\perp}U_{1\perp}^\top & 0\\
		0 & 0
		\end{pmatrix}
		\begin{pmatrix}
		0 & \frakE_1\\
		\frakE_1^\top & 0
		\end{pmatrix}
		\begin{pmatrix}
		U_{1\perp}U_{1\perp}^\top & 0\\
		0 & 0
		\end{pmatrix} = 0.
\end{align*}
Therefore, we have
\begin{align}
		&\left\langle \begin{pmatrix}
		U_1U_1^\top & 0\\
		0 & I_{r_1}
		\end{pmatrix}, 
		\calS_{U_1^{(1)}, 3}\left(E\right)
		\right\rangle\notag\\
		=& -2\left\langle \begin{pmatrix}
		U_1U_1^\top & 0\\
		0 & I_{r_1}
		\end{pmatrix}, 
		\frakP_1^{-1}E\frakP_1^{0}E\frakP_1^{0}E\frakP_1^{-2}
		\right\rangle
		+ 2\left\langle \begin{pmatrix}
		U_1U_1^\top & 0\\
		0 & I_{r_1}
		\end{pmatrix}, 
		\frakP_1^{-1}E\frakP_1^{-1}E\frakP_1^{0}E\frakP_1^{-1}
		\right\rangle\notag\\
		=& 2\left\langle \begin{pmatrix}
		U_1U_1^\top & 0\\
		0 & I_{r_1}
		\end{pmatrix}, 
		\frakP_1^{-1}E\frakP_1^{-1}E\frakP_1^{0}E\frakP_1^{-1}
		\right\rangle. \label{ineq65}
\end{align}
By simple calculation, we have
		\begin{equation*}
		\frakP_1^{-1}E\frakP_1^{-1}E\frakP_1^{0}E\frakP_1^{-1} = \begin{pmatrix}
		U_1R_1^{(1)}\frakE_1^\top U_1R_1^{(1)}\frakE_1^\top \calP_{U_1}^\perp \frakE_1R_1^{(1)\top}U_1^\top & 0\\
		0 & 0 
		\end{pmatrix}.
		\end{equation*}
Therefore, with probability at least $1 - C_1p^{-3} - C_1e^{-c_1p}$,
\begin{align*}
		&\left|\left\langle \begin{pmatrix}
		U_1U_1^\top & 0\\
		0 & I_{r_1}
		\end{pmatrix}, 
		\frakP_1^{-1}E\frakP_1^{-1}E\frakP_1^{0}E\frakP_1^{-1}
		\right\rangle\right|\\
		=& \left|\left\langle \begin{pmatrix}
		U_1U_1^\top & 0\\
		0 & I_{r_1}
		\end{pmatrix}, 
		\begin{pmatrix}
		U_1R_1^{(1)\top}\frakE_1^\top U_1R_1^{(1)\top}\frakE_1^\top \calP_{U_1}^\perp \frakE_1R_1^{(1)}U_1^\top & 0\\
		0 & 0 
		\end{pmatrix}
		\right\rangle\right|\\
		=& \left|\tr\left(U_1R_1^{(1)\top}\frakE_1^\top U_1R_1^{(1)\top}\frakE_1^\top \calP_{U_1}^\perp \frakE_1R_1^{(1)}U_1^\top\right)\right|
		= \left|\tr\left(\frakE_1^\top U_1R_1^{(1)\top}\frakE_1^\top \calP_{U_1}^\perp \frakE_1\right)\right|\\
		\leq& r\|\frakE_1^\top U_1\|\|\frakE_1\|^2
		\leq C_2\frac{pr}{n\lambda_{\submin}^2}\|\frakE_1^\top U_1\|.
\end{align*}
Similarly to {\it Step 0.2} and by Lemma~\ref{lm:sub_Gaussian} and $\calP_{U_1}^\perp U_1 = 0$, with probability at least $1 - C_1p^{-3} - C_1e^{-c_1p}$,
\begin{align*}
		\big\|\frakE_1^\top U_1\big\| \leq& \bigg\|\bigg(\frac{1}{n}\calP_{U_1}\Big(\sum_{j=1}^{n}\xi_j\calM_1(\calX_j)\Big)(U_2 \otimes U_3)G_1^{\top}(G_1 G_1^{\top})^{-1}R_1^{(1)\top}\bigg)^\top U_1\bigg\|\\
		& + C_2\left(\kappa_0\frac{pr}{n\lambda_{\submin}} + \kappa_0^2\frac{p\sqrt{r}}{n\lambda_{\submin}^2}\right)
		\leq  C_2\left(\sqrt{\frac{r^2+\log(p)}{n\lambda_{\submin}^2}} + \kappa_0\frac{pr}{n\lambda_{\submin}} + \kappa_0^2\frac{p\sqrt{r}}{n\lambda_{\submin}^2}\right).
\end{align*}
Combining \eqref{ineq65} and the above two inequalities together, we get with probability at least $1 - C_1p^{-3} - C_1e^{-c_1p}$,
\begin{align}
		\left|\left\langle \begin{pmatrix}
		U_1U_1^\top & 0\\
		0 & I_{r_1}
		\end{pmatrix}, 
		\calS_{U_1^{(1)}, 3}\left(E\right)
		\right\rangle\right|
		\leq& C_2\frac{pr}{n\lambda_{\submin}^2}\left(\sqrt{\frac{r^2+\log(p)}{n\lambda_{\submin}^2}} + \kappa_0\frac{pr}{n\lambda_{\submin}} + \kappa_0^2\frac{p\sqrt{r}}{n\lambda_{\submin}^2}\right)\notag\\
		\leq& C_2\left(\frac{pr(r + \sqrt{\log(p)})}{n^{3/2}\lambda_{\submin}^3} + \kappa_0\frac{p^2r^2}{n^2\lambda_{\submin}^3} + \kappa_0^2\frac{p^2r^{3/2}}{n^2\lambda_{\submin}^4}\right).\label{ineq:S_{3}}
\end{align}
Therefore, we conclude that with probability at least $1 - C_1p^{-3} - C_1e^{-c_1p}$, 
		\begin{align}\label{ineq129}
			\begin{split}
			\left|\left\|\hat U_1^{(2)}\hat U_1^{(2)\top} - U_1U_1^\top\right\|_{\rm F}^2 + 2\left\langle \begin{pmatrix}
			U_1U_1^\top & 0\\
			0 & I_{r_1}
			\end{pmatrix}, 
			\calS_{U_1^{(1)}, 2}\left(E\right)
			\right\rangle\right|\\ \leq C_2\left(\frac{pr(r + \sqrt{\log(p)})}{n^{3/2}\lambda_{\submin}^3} + \kappa_0\frac{p^2r^2}{n^2\lambda_{\submin}^3} + \kappa_0^2\frac{p^2r^{3/2}}{n^2\lambda_{\submin}^4}\right).
			\end{split}
		\end{align}
		
\paragraph*{Step 3: representing the leading term of $\|\hat U_1^{(2)}\hat U_1^{(2)\top}-U_1U_1^{\top}\|_{\rm F}^2$.} Recall from {\it Step 2}, the leading term of $\|\hat U_1^{(2)}\hat U_1^{(2)\top} - U_1U_1^\top\|_{\rm F}^2$ is 
$$
-2\left\langle \begin{pmatrix}
			U_1U_1^\top & 0\\
			0 & I_{r_1}
			\end{pmatrix}, 
			\calS_{U_1^{(1)}, 2}\left(E\right)
\right\rangle
$$
In {\it Step 3}, we aim to approximate this leading term by a sum of independent random variables. By definition of $\calS_{U_1^{(1)}, 2}(E)$, we have 
\begin{align}\label{ineq66}
		\left\langle \begin{pmatrix}
		U_1U_1^\top & 0\\
		0 & I_{r_1}
		\end{pmatrix}, 
		\calS_{U_1^{(1)}, 2}\left(E\right)
		\right\rangle
		=& -\left\langle \begin{pmatrix}
		U_1U_1^\top & 0\\
		0 & I_{r_1}
		\end{pmatrix}, 
		\frakP_1^{-1}E\frakP_1^{0}E\frakP_1^{-1}\right\rangle\notag\\
		=& -\left\langle \begin{pmatrix}
		U_1U_1^\top & 0\\
		0 & I_{r_1}
		\end{pmatrix}, 
		\begin{pmatrix}
		U_1R_1^{(1)\top}\frakE_1^\top \calP_{U_1}^\perp \frakE_1R_1^{(1)}U_1^\top & 0\\
		0 & 0
		\end{pmatrix}\right\rangle\notag\\
		=& -\tr\left(U_1R_1^{(1)\top}\frakE_1^\top \calP_{U_1}^\perp \frakE_1R_1^{(1)}U_1^\top\right)
		= -\tr\left(\frakE_1^\top \calP_{U_1}^\perp \frakE_1\right)\notag\\
		=& -\tr\left(\big(\frakJ_{U_1, 1}^{(1)} + \frakJ_{U_1, 3}^{(1)}\big)^\top \calP_{U_1}^\perp \big(\frakJ_{U_1, 1}^{(1)} + \frakJ_{U_1, 3}^{(1)}\big)\right).
\end{align}
The last equation holds since $\frakJ_{U_1, 2}^{(1)\top}U_{1\perp} = 0$. Therefore, we only need to bound $\tr\big(\frakJ_{U_1, 1}^{(1)\top} \calP_{U_1}^\perp\frakJ_{U_1, 1}^{(1)}\big)$, $\tr\big(\frakJ_{U_1, 1}^{(1)\top} \calP_{U_1}^\perp\frakJ_{U_1, 3}^{(1)}\big)$ and $\tr\big(\frakJ_{U_1, 3}^{(1)\top} \calP_{U_1}^\perp\frakJ_{U_1, 3}^{(1)}\big)$, respectively.

\subparagraph*{Step 3.1: bounding $\tr\big(\frakJ_{U_1, 1}^{(1)\top} \calP_{U_1}^\perp\frakJ_{U_1, 1}^{(1)}\big)$.}
Recall the definitions of $\frakJ_{U_1, 1}^{(1)}$ in {\it Step 1}. Similarly to \eqref{ineq58} and by \eqref{ineq67}, we get with probability at least $1 - C_1p^{-3} - C_1e^{-c_1p}$ that,
\begin{align}
		\left\|\frakJ_{U_1, 1}^{(1)}\right\| \leq& C_2\left(\frac{\kappa_0pr}{n\lambda_{\submin}} + \kappa_0\frac{p^{3/2}r^2}{n^{3/2}\lambda_{\submin}} + \kappa_0\frac{p^{3/2}r^{3/2}}{n^{3/2}\lambda_{\submin}^2}\right) + C_2\kappa_0r\sqrt{\frac{p}{n}}\cdot\sqrt{\frac{p}{n}}\lambda_{\submin}^{-1}\notag\\
		\leq& C_2\kappa_0\frac{pr}{n\lambda_{\submin}}.\label{ineq:J_{U_1,1}}
\end{align}
Therefore, with probability at least $1 - C_1p^{-3} - C_1e^{-c_1p}$,
\begin{equation}\label{ineq68}
		\left|\tr\big(\frakJ_{U_1, 1}^{(1)\top} \calP_{U_1}^\perp\frakJ_{U_1, 1}^{(1)}\big)\right| \leq r\big\|\frakJ_{U_1, 1}^{(1)}\big\|^2 \leq C_2 \kappa_0^2\frac{p^2r^2}{n^2\lambda_{\submin}^2}.
\end{equation}

\subparagraph*{Step 3.2: bounding $\tr\big(\frakJ_{U_1, 1}^{(1)\top} \calP_{U_1}^\perp\frakJ_{U_1, 3}^{(1)}\big)$.} Denote 
$$
K_{U_1}^{(1)} = \Big(\calM_1(\Delta \calT_{1}^{(1.5)}) - \frac{1}{n}\sum_{i=1}^{n}\langle \Delta \calT_{1}^{(1.5)}, \calX_i\rangle\calM_1(\calX_i)\Big)\left(U_2 \otimes U_3\right)
$$
and
$$
L_1 = \frac{1}{n}\Big(\sum_{j=1}^{n}\xi_j\calM_1(\calX_j)\Big)(U_2 \otimes U_3).
$$ 
Lemma \ref{lm:sub_Gaussian} immediately implies that 
\begin{equation}\label{ineq:L_1}
		\PP\left(\left\|L_1\right\| \geq C_2\sqrt{\frac{p}{n}}\right) \leq 1 - p^{-3}.
\end{equation}
By \eqref{ineq61}, with probability at least $1 - C_1p^{-3} - C_1e^{-c_1p}$,
\begin{equation}\label{ineq78}
		\big\|\frakJ_{U_1, 3}^{(1)} - L_1G_1^\top(G_1G_1^\top)^{-1}R_1^{(1)\top}\big\| \leq C_2\kappa_0^2\frac{p\sqrt{r}}{n\lambda_{\submin}^2}.
\end{equation}
By \eqref{ineq13}, \eqref{ineq59}, \eqref{ineq67} and Lemma \ref{lm:Gaussian_ensemble}, with probability at least $1 - C_1p^{-3} - C_1e^{-c_1p}$,
\begin{align*}
		&\big\|\frakJ_{U_1, 1}^{(1)} - K_{U_1}^{(1)}G_1^\top(G_1G_1^\top)^{-1}R_1^{(1)\top}\big\|\\
		\leq& \Big\|\Big(\calM_1(\Delta \calT_{1}^{(1.5)}) - \frac{1}{n}\sum_{i=1}^{n}\langle \Delta \calT_{1}^{(1.5)}, \calX_i\rangle\calM_1(\calX_i)\Big)\Big(\big(\hat U_2^{(1)}R_2^{(1)} - U_2\big) \otimes U_3\Big)G_1^\top(G_1G_1^\top)^{-1}R_1^{(1)\top}\Big\|\\
		& + \Big\|\Big(\calM_1(\Delta \calT_{1}^{(1.5)}) - \frac{1}{n}\sum_{i=1}^{n}\langle \Delta \calT_{1}^{(1.5)}, \calX_i\rangle\calM_1(\calX_i)\Big)\Big(\big(\hat U_2^{(1)}R_2^{(1)}\big) \otimes \big(\hat U_3^{(1)}R_3^{(1)} - U_3\big)\Big)G_1^\top(G_1G_1^\top)^{-1}R_1^{(1)\top}\Big\|\\
		& + \Big\|\big(\calM_1(\Delta \calT_{1}^{(1.5)}) - \frac{1}{n}\sum_{i=1}^{n}\langle \Delta \calT_{1}^{(1.5)}, \calX_i\rangle\calM_1(\calX_i)\big)\big(\hat U_2^{(1)} \otimes \hat U_3^{(1)}\big)\Big\|\\ 
		&\hspace{3cm}\times\Big\|\big(R_2^{(1)} \otimes R_3^{(1)}\big)G_1^\top(G_1G_1^\top)^{-1}R_1^{(1)\top} - \hat G_1^\top(\hat G_1\hat G_1^\top)^{-1}\Big\|\\
		\leq& C_2\sqrt{\frac{pr}{n}}\big\|\calM_1(\Delta \calT_{1}^{(1.5)})\big\|_{\rm F}\big(\big\|\hat U_2^{(1)}R_2^{(1)} - U_2\big\| + \big\|\hat U_3^{(1)}R_3^{(1)} - U_3\big\|\big)\lambda_{\submin}^{-1}\\ 
		&+ C_2\sqrt{\frac{pr}{n}}\big\|\calM_1(\Delta \calT_{1}^{(1.5)})\big\|_{\rm F}\kappa_0^2\lambda_{\submin}^{-2}\Big(\sqrt{\frac{r^2+\log(p)}{n}} + \kappa_0\frac{pr}{n} + \kappa_0\frac{p\sqrt{r}}{n\lambda_{\submin}}\Big)\\
		\leq& C_2\sqrt{\frac{pr}{n}}\big\|\calM_1(\Delta \calT_{1}^{(1.5)})\big\|_{\rm F}\kappa_0^2\sqrt{\frac{p}{n}}\lambda_{\submin}^{-2}\\
		\leq& C_2\kappa_0^2\frac{p\sqrt{r}}{n\lambda_{\submin}^2}C_2\bigg(\kappa_0\sqrt{\frac{pr}{n}} + \kappa_0\frac{pr^{3/2}}{n} + \kappa_0\frac{pr}{n\lambda_{\submin}}\bigg)\\
		\leq& C_2\kappa_0^3\frac{p^{3/2}r}{n^{3/2}\lambda_{\submin}^2}.
\end{align*}
\eqref{ineq:J_{U_1,1}}, \eqref{ineq68}, \eqref{ineq:L_1}, \eqref{ineq78} and the previous inequality together imply with probability at least $1 - C_1p^{-3} - C_1e^{-c_1p}$,
\begin{align}
		&\left|\tr\big(\frakJ_{U_1, 1}^{(1)\top} \calP_{U_1}^\perp\frakJ_{U_1, 3}^{(1)}\big)\right|\notag\\
		\leq& \Big|\tr\Big(\big(\frakJ_{U_1, 1}^{(1)} - K_{U_1}^{(1)}G_1^\top(G_1G_1^\top)^{-1}R_1^{(1)\top}\big)^\top \calP_{U_1}^\perp\frakJ_{U_1, 3}^{(1)}\Big)\Big|\notag\\
		& + \Big|\tr\Big(\frakJ_{U_1, 1}^{(1)\top} \calP_{U_1}^\perp\big(\frakJ_{U_1, 3}^{(1)} - L_1G_1^\top(G_1G_1^\top)^{-1}R_1^{(1)\top}\big)\Big)\Big|\notag\\
		& + \Big|\tr\Big(R_1^{(1)}(G_1G_1^\top)^{-1}G_1K_{U_1}^{(1)\top} \calP_{U_1}^\perp L_1G_1^\top(G_1G_1^\top)^{-1}R_1^{(1)\top}\Big)\Big|\notag\\
		\leq& r\big\|\frakJ_{U_1, 1}^{(1)} - K_{U_1}^{(1)}G_1^\top(G_1G_1^\top)^{-1}R_1^{(1)\top}\big\|\big\|\frakJ_{U_1, 3}^{(1)}\big\| + r\big\|\frakJ_{U_1, 1}^{(1)}\big\|\big\|\frakJ_{U_1, 3}^{(1)} - L_1G_1^\top(G_1G_1^\top)^{-1}R_1^{(1)\top}\big\|\notag\\
		&+ \big|\tr\big(G_1^\top(G_1G_1^\top)^{-2}G_1K_{U_1}^{(1)\top} \calP_{U_1}^\perp L_1\big)\big|\notag\\
		\leq& C_2r\cdot \kappa_0^3\frac{p^{3/2}r}{n^{3/2}\lambda_{\submin}^2}\cdot \frac{\sqrt{p/n}}{\lambda_{\submin}} + C_2r\cdot \kappa_0\frac{pr}{n\lambda_{\submin}}\cdot \kappa_0^2\frac{p\sqrt{r}}{n\lambda_{\submin}^2} + \big|\tr\big(G_1^\top(G_1G_1^\top)^{-2}G_1K_{U_1}^{(1)\top} \calP_{U_1}^\perp L_1\big)\big|\notag\\
		\leq& C_2\kappa_0^3\frac{p^2r^{5/2}}{n^2\lambda_{\submin}^3} + \big|\tr\big(G_1^\top(G_1G_1^\top)^{-2}G_1K_{U_1}^{(1)\top} \calP_{U_1}^\perp L_1\big)\big|.\label{ineq69}
\end{align}
By \eqref{ineq:error_G} and Lemma \ref{lm:sub_Gaussian}, with probability at least $1 - C_1p^{-3} - C_1e^{-c_1p}$,
\begin{align}
		&\Big\|\calM_1\big(\Delta \calT_{1}^{(1.5)}\big) - \Big((\hat U_1^{(1.5)}R_1^{(1)})G_1\big((\hat U_2^{(1)}R_1^{(1)})^\top \otimes (\hat U_3^{(1)}R_1^{(1)})^\top\big) - U_1G_1(U_2 \otimes U_3)\Big)\Big\|_{\rm F}\notag\\ 
		=& \Big\|\hat U_1^{(1.5)}\big(\hat G_1^{(1)} - R_1^{(1)}G_1(R_2^{(1)} \otimes R_3^{(1)})\big)(\hat U_2^{(1)\top} \otimes \hat U_3^{(1)\top})\Big\|_{\rm F}\notag\\
		\leq& C_2\left(\kappa_0\frac{pr^{3/2}}{n} + \kappa_0\frac{pr}{n\lambda_{\submin}} + \sqrt{\frac{r^3 + \log(p)}{n}}\right).\label{ineq72}
\end{align}
Define
\begin{align*}
		&\Delta T_{U_1,1}^{(1.5)} = (\hat{U}_1^{(1.5)}R_1^{(1)} - U_1)G_1(U_2^\top \otimes U_3^\top), \quad \Delta T_{U_2,1}^{(1)} = U_1G_1\left((\hat{U}_2^{(1)}R_2^{(1)} - U_2)^\top \otimes U_3^\top\right),\\
		&\Delta T_{U_3,1}^{(1)} = U_1G_1\left(U_2^\top \otimes (\hat{U}_3^{(1)}R_3^{(1)} - U_3)^\top\right).
\end{align*}
By \eqref{eq:regression_step_2}, \eqref{ineq70}, \eqref{ineq71} and Lemma \ref{lm:sub_Gaussian}, with probability at least $1 - C_1p^{-3} - C_1e^{-c_1p}$, we get
\begin{align*}
		\max\left\{\|U_1 - \hat U_1^{(1.5)}R_1^{(1)}\|, \|U_2 - \hat U_2^{(1)}R_2^{(1)}\|, \|U_3 - \hat U_3^{(1)}R_3^{(1)}\|\right\} \leq C_2\frac{\sqrt{p/n}}{\lambda_{\submin}}.
\end{align*}
Therefore, with probability at least $1 - C_1p^{-3} - C_1e^{-c_1p}$, 
\begin{align*}
		&\left\|(\hat U_1^{(1.5)}R_1^{(1)})G_1\left((\hat U_2^{(1)}R_1^{(1)})^\top \otimes (\hat U_3^{(1)}R_1^{(1)})^\top\right) - U_1G_1(U_2 \otimes U_3) - \left(\Delta T_{U_1,1}^{(1.5)} + \Delta T_{U_2,1}^{(1)} + \Delta T_{U_3,1}^{(1)}\right)\right\|_{\rm F}\\
		\leq& C_2\sqrt{r}\Big(\frac{\sqrt{p/n}}{\lambda_{\submin}}\Big)^2\cdot \kappa_0\lambda_{\submin} = C_2\kappa_0\frac{p\sqrt{r}}{n\lambda_{\submin}}.
\end{align*}
Combining \eqref{ineq72} and the above inequality, with probability at least $1 - C_1p^{-3} - C_1e^{-c_1p}$, 
\begin{equation}\label{ineq73}
		\left\|\calM_1\big(\Delta \calT_{1}^{(1.5)}\big) - \left(\Delta T_{U_1,1}^{(1.5)} + \Delta T_{U_2,1}^{(1)} + \Delta T_{U_3,1}^{(1)}\right)\right\|_{\rm F} \leq C_2\Big(\kappa_0\frac{pr^{3/2}}{n} + \kappa_0\frac{pr}{n\lambda_{\submin}} + \sqrt{\frac{r^3 + \log(p)}{n}}\Big).
\end{equation}
Define
\begin{align*}
		S_1^{(1.5)} =& \tr\bigg(G_1^\top(G_1G_1^\top)^{-2}G_1\bigg\{\Big[\Delta T_{U_1,1}^{(1.5)} - \frac{1}{n}\sum_{i=1}^{n}\langle \Delta T_{U_1,1}^{(1.5)}, \calM_1(\calX_i)\rangle\calM_1(\calX_i)\Big](U_2 \otimes U_3)\bigg\}^\top\calP_{U_1}^\perp L_1\bigg),\\
		S_2^{(1)} =& \tr\bigg(G_1^\top(G_1G_1^\top)^{-2}G_1\bigg\{\Big[\Delta T_{U_2,1}^{(1)} - \frac{1}{n}\sum_{i=1}^{n}\langle \Delta T_{U_2,1}^{(1)}, \calM_1(\calX_i)\rangle\calM_1(\calX_i)\Big](U_2 \otimes U_3)\bigg\}^\top\calP_{U_1}^\perp L_1\bigg),\\
		S_3^{(1)} =& \tr\bigg(G_1^\top(G_1G_1^\top)^{-2}G_1\bigg\{\Big[\Delta T_{U_3,1}^{(1)} - \frac{1}{n}\sum_{i=1}^{n}\langle \Delta T_{U_3,1}^{(1)}, \calM_1(\calX_i)\rangle\calM_1(\calX_i)\Big](U_2 \otimes U_3)\bigg\}^\top\calP_{U_1}^\perp L_1\bigg).
\end{align*}
By \eqref{ineq73}, \eqref{ineq:L_1} and Lemma \ref{lm:Gaussian_ensemble}, with probability at least $1 - C_1p^{-3} - C_1e^{-c_1p}$, 
\begin{align}
		&\left|\tr\Big(G_1^\top(G_1G_1^\top)^{-2}G_1K_{U_1}^{(1)\top} \calP_{U_1}^\perp L_1\Big) - S_1^{(1.5)} - S_2^{(1)} - S_3^{(1)}\right|\notag\\
		\leq& C_2r\big\|G_1^\top(G_1G_1^\top)^{-2}G_1\big\|\left\|L_1\right\|\left(\sqrt{\frac{pr}{n}}\left\|\calM_1\left(\Delta \calT_{1}^{(1.5)}\right) - \left(\Delta T_{U_1,1}^{(1.5)} + \Delta T_{U_2,1}^{(1)} + \Delta T_{U_3,1}^{(1)}\right)\right\|_{\rm F}\right)\notag\\
		\leq& C_2r\lambda_{\submin}^{-2}\sqrt{\frac{p}{n}}\cdot \sqrt{\frac{pr}{n}}\Big(\kappa_0\frac{pr^{3/2}}{n} + \kappa_0\frac{pr}{n\lambda_{\submin}} + \sqrt{\frac{r^3 + \log(p)}{n}}\Big)\notag\\
		\leq& C_2\bigg(\kappa_0\frac{p^2r^{3}}{n^2\lambda_{\submin}^2} + \kappa_0\frac{p^2r^{5/2}}{n^2\lambda_{\submin}^3} + \frac{p\big(r^{3}+r^{3/2}\sqrt{\log p}\big)}{n^{3/2}\lambda_{\submin}^2}\bigg).\label{ineq74}
\end{align}
Therefore, it suffices to bound $\big|S_1^{(1.5)}\big|$, $\big|S_2^{(1)}\big|$, and $\big|S_3^{(1)}\big|$, respectively.

\subparagraph*{-Step 3.2.1: bounding $\big|S_1^{(1.5)}\big|$.} 		
We consider $\big|S_1^{(1.5)}\big|$ first. The proof of this part is involved and highly non-trivial, and some decoupling techniques (e.g., \cite{de2012decoupling}) are needed. Let
		\begin{equation*}
		\frakE_{U_1, i}^{(1)} = \frakJ_{U_1, i}^{(1)}R_1^{(1)}G_1(U_2^\top \otimes U_3^\top), \quad \forall i \in [3].
		\end{equation*}
By \eqref{ineq62}, 
\begin{align}
		&S_1^{(1.5)}\notag\\ 
		=& \tr\bigg(G_1^\top(G_1G_1^\top)^{-2}G_1\Big\{\Big[\frakE_{U_1, 1}^{(1)} - \frac{1}{n}\sum_{i=1}^{n}\langle \frakE_{U_1, 1}^{(1)}, \calM_1(\calX_i)\rangle\calM_1(\calX_i)\Big](U_2 \otimes U_3)\Big\}^\top\calP_{U_1}^\perp L_1\bigg)\notag\\
		& + \tr\bigg(G_1^\top(G_1G_1^\top)^{-2}G_1\Big\{\Big[\frakE_{U_1, 2}^{(1)} - \frac{1}{n}\sum_{i=1}^{n}\langle \frakE_{U_1, 2}^{(1)}, \calM_1(\calX_i)\rangle\calM_1(\calX_i)\Big](U_2 \otimes U_3)\Big\}^\top\calP_{U_1}^\perp L_1\bigg)\notag\\
		& + \tr\bigg(G_1^\top(G_1G_1^\top)^{-2}G_1\Big\{\Big[\frakE_{U_1, 3}^{(1)} - \frac{1}{n}\sum_{i=1}^{n}\langle \frakE_{U_1, 3}^{(1)}, \calM_1(\calX_i)\rangle\calM_1(\calX_i)\Big](U_2 \otimes U_3)\Big\}^\top\calP_{U_1}^\perp L_1\bigg).\label{ineq98}
\end{align}
By \eqref{ineq:J_{U_1,1}}, \eqref{ineq:L_1} and Lemma \ref{lm:Gaussian_ensemble}, we have
\begin{align}
		&\bigg|\tr\bigg(G_1^\top(G_1G_1^\top)^{-2}G_1\Big\{\Big[\frakE_{U_1, 1}^{(1)} - \frac{1}{n}\sum_{i=1}^{n}\langle \frakE_{U_1, 1}^{(1)}, \calM_1(\calX_i)\rangle\calM_1(\calX_i)\Big](U_2 \otimes U_3)\Big\}^\top\calP_{U_1}^\perp L_1\bigg)\bigg|\notag\\
		\leq& r\big\|G_1^\top(G_1G_1^\top)^{-2}G_1\big\|\|L_1\|\Big(\sqrt{\frac{pr}{n}}\big\|\frakE_{U_1, 1}^{(1)}\big\|_{\rm F}\Big)\notag\\
		\leq& C_2r\lambda_{\submin}^{-2}\sqrt{\frac{p}{n}}\cdot \sqrt{\frac{pr}{n}}\sqrt{r}\big\|\frakJ_{U_1, 1}^{(1)}\big\|\|G_1\|
		\leq C_2\kappa_0^2\frac{p^{2}r^{3}}{n^2\lambda_{\submin}^2}.\label{ineq75}
\end{align}
In addition, by \eqref{ineq:error_G}, \eqref{ineq70} and Lemma \ref{lm:sub_Gaussian},  with probability at least $1 - C_1p^{-3} - C_1e^{-c_1p}$,
\begin{align*}
		\big\|\frakJ_{U_1, 2}^{(1)}\big\|_{\rm F} =& \Big\|\left(U_1G_1\big((\hat U_2^{(1)\top}U_2)^\top \otimes (\hat U_3^{(1)\top}U_3)^\top\big) - U_1R_1^{(1)\top}\hat G_1^{(1)}\right)\hat G_1^{(1)\top}\big(\hat G_1^{(1)}\hat G_1^{(1)\top}\big)^{-1}\Big\|_{\rm F}\\
		\leq& C_2\lambda_{\submin}^{-1}\left\|G_1\left((\hat U_2^{(1)\top}U_2)^\top \otimes (\hat U_3^{(1)\top}U_3)^\top\right) - R_1^{(1)\top}\hat G_1^{(1)}\right\|_{\rm F}\\
		\leq& C_2\lambda_{\submin}^{-1}\left\|G_1\left((\hat U_2^{(1)\top}U_2)^\top \otimes (\hat U_3^{(1)\top}U_3)^\top\right) - G_1\left(R_2^{(1)\top} \otimes R_3^{(1)\top}\right)\right\|_{\rm F}\\
		&+ C_2\lambda_{\submin}^{-1}\left\|R_1^{(1)\top}\left(\hat G_1^{(1)} - R_1^{(1)}G_1\left(R_2^{(1)\top} \otimes R_3^{(1)\top}\right)\right)\right\|_{\rm F}\\
		\leq& C_2\lambda_{\submin}^{-1}\cdot\kappa_0\lambda_{\submin}\frac{p\sqrt{r}}{n\lambda_{\submin}^2} + C_2\lambda_{\submin}^{-1}\Big(\kappa_0\frac{pr^{3/2}}{n} + \kappa_0\frac{pr}{n\lambda_{\submin}} + \sqrt{\frac{r^3 + \log p}{n}}\Big)\\
		=& C_2\Big(\kappa_0\frac{pr}{n\lambda_{\submin}^2} + \kappa_0\frac{pr^{3/2}}{n\lambda_{\submin}} + \sqrt{\frac{r^3+\log p}{n\lambda_{\submin}^2}}\Big).
\end{align*}
Therefore, \eqref{ineq:L_1}, Lemma \ref{lm:Gaussian_ensemble} and the previous inequality together imply that with probability at least $1 - C_1p^{-3} - C_1e^{-c_1p}$,
\begin{align}
		&\bigg|\tr\bigg(G_1^\top(G_1G_1^\top)^{-2}G_1\Big\{\Big[\frakE_{U_1, 2}^{(1)} - \frac{1}{n}\sum_{i=1}^{n}\langle \frakE_{U_1, 2}^{(1)}, \calM_1(\calX_i)\rangle\calM_1(\calX_i)\Big](U_2 \otimes U_3)\Big\}^\top\calP_{U_1}^\perp L_1\bigg)\bigg|\notag\\
		\leq& r\left\|G_1^\top(G_1G_1^\top)^{-2}G_1\right\|\|L_1\|\Big(\sqrt{\frac{pr}{n}}\big\|\frakE_{U_1, 2}^{(1)}\big\|_{\rm F}\Big)\notag\\
		\leq& C_2r\lambda_{\submin}^{-2}\sqrt{\frac{p}{n}}\cdot \sqrt{\frac{pr}{n}}\big\|\frakJ_{U_1, 2}^{(1)}\big\|_{\rm F}\|G_1\|\notag\\
		\leq& C_2\left(\kappa_0^2\frac{p^2r^{5/2}}{n^2\lambda_{\submin}^{3}} + \kappa_0^2\frac{p^2r^{3}}{n^2\lambda_{\submin}^{2}} + \kappa_0\frac{r^{3}p + r^{3/2}p\sqrt{\log(p)}}{n^{3/2}\lambda_{\submin}^2}\right).\label{ineq76}
\end{align}
Let $$J_1 = L_1G_1^\top(G_1G_1^\top)^{-1}G_1(U_2^\top \otimes U_3^\top).$$ Consider the SVD decomposition $G_1 = U_{G_1}\Lambda_{1}V_{G_1}^\top$, where $U_{G_1} \in \OO_{r_1}, V_{G_1} \in \OO_{r_2r_3, r_1}$ and $\Lambda_1 \in \RR^{r_1 \times r_1}$ is a diagonal matrix containing all singular values of $G_1$. Then 
$$
G_1^\top(G_1G_1^\top)^{-1}G_1 = V_{G_1}V_{G_1}^\top.
$$
By \eqref{ineq78}, with probability at least $1 - C_1p^{-3} - C_1e^{-c_1p}$,
\begin{align}
		&\bigg|\tr\bigg(G_1^\top(G_1G_1^\top)^{-2}G_1\Big\{\Big[\frakE_{U_1, 3}^{(1)} - \frac{1}{n}\sum_{i=1}^{n}\langle \frakE_{U_1, 3}^{(1)}, \calM_1(\calX_i)\rangle\calM_1(\calX_i)\Big](U_2 \otimes U_3)\Big\}^\top\calP_{U_1}^\perp L_1\bigg)\bigg|\notag\\
		\leq& \bigg|\tr\bigg(G_1^\top(G_1G_1^\top)^{-2}G_1\Big\{\Big[J_1 - \frac{1}{n}\sum_{i=1}^{n}\langle J_1, \calM_1(\calX_i)\rangle\calM_1(\calX_i)\Big](U_2 \otimes U_3)\Big\}^\top\calP_{U_1}^\perp L_1\bigg)\bigg|\notag\\
		&+ r\left\|G_1^\top(G_1G_1^\top)^{-2}G_1\right\|\|L_1\|\Big(\sqrt{\frac{pr}{n}}\big\|\frakE_{U_1, 3}^{(1)} - J_1\big\|_{\rm F}\Big)\notag\\
		\leq& \bigg|\tr\bigg(G_1^\top(G_1G_1^\top)^{-2}G_1\Big\{\Big[J_1 - \frac{1}{n}\sum_{i=1}^{n}\langle J_1, \calM_1(\calX_i)\rangle\calM_1(\calX_i)\Big](U_2 \otimes U_3)\Big\}^\top\calP_{U_1}^\perp L_1\bigg)\bigg|\notag\\
		&+ C_2r\lambda_{\submin}^{-2}\sqrt{\frac{p}{n}}\cdot \sqrt{\frac{pr}{n}}\sqrt{r}\left\|\frakJ_{U_1, 3}^{(1)} - L_1G_1^\top(G_1G_1^\top)^{-1}R_1^{(1)\top}\right\|\|G_1\|\notag\\
		\leq& \bigg|\tr\bigg(G_1^\top(G_1G_1^\top)^{-2}G_1\Big\{\Big[J_1 - \frac{1}{n}\sum_{i=1}^{n}\langle J_1, \calM_1(\calX_i)\rangle\calM_1(\calX_i)\Big](U_2 \otimes U_3)\Big\}^\top\calP_{U_1}^\perp L_1\bigg)\bigg|\notag\\
		&+ C_2r\lambda_{\submin}^{-2}\sqrt{\frac{p}{n}}\cdot \sqrt{\frac{pr}{n}}\cdot \sqrt{r}\kappa_0^2\frac{p\sqrt{r}}{n\lambda_{\submin}^2}\cdot \kappa_0\lambda_{\submin}\notag\\
		\leq& \bigg|\tr\bigg(G_1^\top(G_1G_1^\top)^{-2}G_1\Big\{\Big[J_1 - \frac{1}{n}\sum_{i=1}^{n}\langle J_1, \calM_1(\calX_i)\rangle\calM_1(\calX_i)\Big](U_2 \otimes U_3)\Big\}^\top L_1\bigg)\bigg|\notag\\
		&+\bigg|\tr\bigg(G_1^\top(G_1G_1^\top)^{-2}G_1\Big\{\Big[J_1 - \frac{1}{n}\sum_{i=1}^{n}\langle J_1, \calM_1(\calX_i)\rangle\calM_1(\calX_i)\Big](U_2 \otimes U_3)\Big\}^\top\calP_{U_1}L_1\bigg)\bigg|\notag\\
		&+ C_2\kappa_0^3\frac{p^2r^{5/2}}{n^2\lambda_{\submin}^3}.\label{ineq77}
\end{align}
For $i \in [n]$, let 
$$
Z_i = \calM_1(\calX_i)(U_2 \otimes U_3)V_{G_1} \in \RR^{p_1 \times r_1},
$$
and
$$
Z_{\xi} = L_1V_{G_1} = \sum_{i=1}^{n}\xi_i\calM_1(\calX_i)(U_2 \otimes U_3)V_{G_1} \in \RR^{p_1 \times r_1}.
$$
Then
$$
Z_i \stackrel{i.i.d.}{\sim} N(0, 1) \quad \text{and} \quad Z_{\xi} = \sum_{j=1}^{n}\xi_jZ_j.
$$
Thus
\begin{align}
		&\bigg|\tr\bigg(G_1^\top(G_1G_1^\top)^{-2}G_1\Big\{\Big[J_1 - \frac{1}{n}\sum_{i=1}^{n}\langle J_1, \calM_1(\calX_i)\rangle\calM_1(\calX_i)\Big](U_2 \otimes U_3)\Big\}^\top L_1\bigg)\bigg|\notag\\
		=& \bigg|\frac{1}{n^3}\sum_{i=1}^{n}\sum_{j=1}^{n}\sum_{k=1}^{n}\xi_j\xi_k\langle Z_i, Z_j\rangle\langle Z_i\Lambda_1^{-1}, Z_k\Lambda_1^{-1}\rangle - \frac{1}{n^2}\sum_{j=1}^{n}\sum_{k=1}^{n}\xi_j\xi_k\langle Z_j\Lambda_1^{-1}, Z_k\Lambda_1^{-1}\rangle\bigg|\notag\\
		\leq& \bigg|\frac{1}{n^3}\sum_{i=1}^{n}\xi_i^2\|Z_i\|_{\rm F}^2\|Z_i\Lambda_1^{-1}\|_{\rm F}^2 - \frac{1}{n^3}\sum_{i=1}^{n}\xi_i^2\|Z_i\Lambda_1^{-1}\|_{\rm F}^2\bigg|\notag\\
		&+ \bigg|\frac{1}{n^3}\sum_{i = 1}^n\bigg[\sum_{j \neq i}\xi_j^2\left(\langle Z_i, Z_j\rangle\langle Z_i, Z_j\Lambda_1^{-2}\rangle - \langle Z_j, Z_j\Lambda_1^{-2}\rangle\right)\bigg]\bigg|\notag\\
		&+ \bigg|\frac{1}{n^3}\sum_{i=1}^{n}\xi_i\|Z_i\|_{\rm F}^2\sum_{k \neq i}\xi_k\langle Z_i, Z_k\Lambda_1^{-2}\rangle\bigg|
		+\bigg|\frac{1}{n^3}\sum_{i=1}^{n}\xi_i\left\|Z_i\Lambda_1^{-1}\right\|_{\rm F}^2\sum_{j \neq i}\xi_j\langle Z_i, Z_j\rangle\bigg|\notag\\
		&+ \bigg|\frac{2}{n^3}\sum_{i=1}^{n}\xi_i\sum_{k \neq i}\xi_k\langle Z_i, Z_k\Lambda_1^{-2}\rangle\bigg|\notag\\
		&+ \bigg|\frac{1}{n^3}\sum_{i=1}^{n}\sum_{j \neq k \neq i}\xi_j\xi_k\left[\langle Z_j, Z_i\rangle\langle Z_i, Z_k\Lambda_1^{-2}\rangle - \langle Z_j, Z_k\Lambda_1^{-2}\rangle\right]\bigg|.\label{ineq90}
\end{align}
By \cite[Corollary 5.35]{vershynin2010introduction}, for any $i \in [n]$, with probability at least $1 - e^{-c_1(p + \log(n))}$, 
\begin{equation}\label{ineq:Z_i}
		\|Z_i\| \leq C_2\sqrt{p + \log(n)} \quad \text{and} \quad \|Z_i\|_{\rm F} \leq C_2\sqrt{r(p+\log(n))}
\end{equation}
and 
\begin{equation}\label{ineq79}
		\|Z_i\Lambda_1^{-1}\|_{\rm F} \leq \|Z_i\|_{\rm F}\|\Lambda_1^{-1}\| \leq C_2\frac{\sqrt{r(p + \log(n))}}{\lambda_{\submin}}.
\end{equation}
By the union bound and Bernstein-type inequality, with probability at least $1 - e^{-C_1(pr + \log(n))}$,
\begin{align*}
		\frac{1}{n^3}\sum_{i=1}^{n}\xi_i^2\|Z_i\|_{\rm F}^2\|Z_i\Lambda_1^{-1}\|_{\rm F}^2 \leq \frac{1}{n^3}\sum_{i = 1}^n\xi_i^2\cdot C_2\frac{r^2(p + \log(n))^2}{\lambda_{\submin}^2} \leq& C_2\frac{r^2(p + \log(n))^2}{n^3\lambda_{\submin}^2}\cdot Cn\\
		 \leq& C_2\frac{p^2r^2 + r^2\log^2(n)}{n^2\lambda_{\submin}^2}
\end{align*}
and 
\begin{equation}\label{ineq86}
		\frac{1}{n^3}\sum_{i=1}^{n}\xi_i^2\|Z_i\Lambda_1^{-1}\|_{\rm F}^2 \leq \frac{1}{n^3}\sum_{i=1}^{n}\xi_i^2\cdot\frac{r(p+\log n)}{\lambda_{\submin}} \leq C_2\frac{r(p+\log n)}{n^2\lambda_{\submin}^2}.
\end{equation}
Therefore, with probability at least $1 - e^{-c_1pr}$,
\begin{equation}\label{ineq80}
		\bigg|\frac{1}{n^3}\sum_{i=1}^{n}\xi_i^2\|Z_i\|_{\rm F}^2\|Z_i\Lambda_1^{-1}\|_{\rm F}^2 - \frac{1}{n^3}\sum_{i=1}^{n}\xi_i^2\|Z_i\Lambda_1^{-1}\|_{\rm F}^2\bigg| \leq C_2\frac{p^2r^2 + r^2\log^2n}{n^2\lambda_{\submin}^2}.
\end{equation}
		Since $Z_i$ and $Z_j$ are independent for all $1 \leq i \neq j \leq n$, we have $\langle Z_i, Z_j\rangle\big|Z_i \sim N(0, \|Z_i\|_{\rm F}^2)$ and $\langle Z_i, Z_j\Lambda_1^{_2}\rangle\big|Z_i \sim N(0, \|Z_i\Lambda_1^{-2}\|_{\rm F}^2)$, which imply that
\begin{equation*}
		\left\|\langle Z_i, Z_j\rangle\right\|_{\psi_2}\bigg|Z_i \leq C\|Z_i\|_{\rm F} \quad \text{and} \quad \left\|\langle Z_i, Z_j\Lambda_1^{_2}\rangle\right\|_{\psi_2}\bigg|Z_i \leq C\|Z_i\Lambda_1^{-2}\|_{\rm F}.
\end{equation*}
Since $\EE \left[\langle Z_i, Z_j\rangle\langle Z_i, Z_j\Lambda_1^{-2}\rangle\right]\big|Z_i = \langle Z_i, Z_i\Lambda_1^{-2}\rangle$, by \cite[Remark 5.18]{vershynin2010introduction} and \eqref{ineq:psi_1},
\begin{align}
		\big\|\langle Z_i, Z_j\rangle&\langle Z_i, Z_j\Lambda_1^{-2}\rangle - \langle Z_j, Z_j\Lambda_1^{-2}\rangle\big\|_{\psi_1}\big| Z_i
		\leq C\left\|\langle Z_i, Z_j\rangle\langle Z_i, Z_j\Lambda_1^{-2}\rangle\right\|_{\psi_1}\big| Z_i\notag\\
		\leq& C\|\langle Z_i, Z_j\rangle\|_{\psi_2}\|\langle Z_i, Z_j\Lambda_1^{_2}\rangle\|_{\psi_2}\big|Z_i
		\leq C\|Z_i\|_{\rm F}\|Z_i\Lambda_1^{-2}\|_{\rm F}.\label{ineq93}
\end{align}
By Bernstein-type inequality, we have
\begin{equation*}
		\begin{split}
		&\PP\bigg(\bigg|\sum_{j \neq i}\xi_j^2\left(\langle Z_i, Z_j\rangle\langle Z_i, Z_j\Lambda_1^{-2}\rangle - \langle Z_j, Z_j\Lambda_1^{-2}\rangle\right)\bigg| \geq C_2\|Z_i\|_{\rm F}\|Z_i\Lambda_1^{-2}\|_{\rm F}\big(\sum_{j \neq i}\xi_j^4\big)^{1/2}\log(n)\big|Z_i, \xi_1, \dots, \xi_n\bigg)\\ \leq& n^{-3}.
		\end{split}
\end{equation*}
The union bound and \eqref{ineq:Z_i} together imply that
\begin{equation*}
		\begin{split}
		&\PP\bigg(\bigg|\frac{1}{n^3}\sum_{i = 1}^n\Big[\sum_{j \neq i}\xi_j^2\left(\langle Z_i, Z_j\rangle\langle Z_i, Z_j\Lambda_1^{-2}\rangle - \langle Z_j, Z_j\Lambda_1^{-2}\rangle\right)\Big]\bigg| \geq C_2\frac{r(p + \log(n))\log(n)}{n^2\lambda_{\submin}^2}\big(\sum_{j=1}^n\xi_j^4\big)^{1/2}\bigg|\xi_1, \dots, \xi_n\bigg)\\ \leq &n^{-3}.
		\end{split}
\end{equation*}
Notice that 
$$
\EE\xi_i^4 \leq \bigg(2\sup_{q \geq 1}q^{-1/2}\left(\EE \left|\xi_i\right|^q\right)^{1/q}\bigg)^{4} \leq C\|\xi_i\|_{\psi_2}^4 \leq C,
$$
by \cite[Lemmas 7 and 8]{hao2020sparse}, 
\begin{equation}\label{ineq:Gaussian_fourth}
		\PP\bigg(\sum_{j=1}^n\xi_j^4 - Cn \geq C_2\left(\sqrt{n\log(p)} + \log^2(p)\right)\bigg) \leq p^{-3}.
\end{equation}
By combining the above two inequalities together, we know that with probability at least $1 - C_1p^{-3}$,
\begin{equation}\label{ineq81}
		\begin{split}
		\bigg|\frac{1}{n^3}\sum_{i = 1}^n\bigg[\sum_{j \neq i}\xi_j^2\left(\langle Z_i, Z_j\rangle\langle Z_i, Z_j\Lambda_1^{-2}\rangle - \langle Z_j, Z_j\Lambda_1^{-2}\rangle\right)\bigg]\bigg| \leq C_2\frac{r(p + \log(n))\log(n)}{n^{3/2}\lambda_{\submin}^2}.
		\end{split}
\end{equation}
Note that 
\begin{align*}
		&\sum_{i=1}^{n}\xi_i\|Z_i\|_{\rm F}^2\sum_{k \neq i}\xi_k\langle Z_i, Z_k\Lambda_1^{-2}\rangle
		= \sum_{i=1}^n\|Z_i\|_{\rm F}^2\langle Z_i, \xi_i\sum_{k \neq i}\xi_kZ_k\Lambda_1^{-2}\rangle
\end{align*}
By \cite[Theorem 3.4.1]{de2012decoupling}, there exists a constant $C > 0$, for any $t > 0$, we have
\begin{align}
		&\PP\bigg(\bigg|\sum_{i=1}^{n}\xi_i\|Z_i\|_{\rm F}^2\sum_{k \neq i}\xi_k\langle Z_i, Z_k\Lambda_1^{-2}\rangle\bigg| \geq t\bigg) \leq C\PP\bigg(\bigg|\sum_{i=1}^n\|Z_i^{(1)}\|_{\rm F}^2\langle Z_i^{(1)}, \xi_i\sum_{k \neq i}\xi_kZ_k^{(2)}\Lambda_1^{-2}\rangle\bigg| \geq t/C\bigg)\notag\\
		\leq& C\PP\bigg(\bigg|\sum_{i=1}^n\|Z_i^{(1)}\|_{\rm F}^2\langle Z_i^{(1)}, \xi_i\sum_{k = 1}^n\xi_kZ_k^{(2)}\Lambda_1^{-2}\rangle\bigg| \geq \frac{t}{2C}\bigg)
		+C\PP\bigg(\bigg|\sum_{i=1}^n\|Z_i^{(1)}\|_{\rm F}^2\langle Z_i^{(1)}, \xi_i^2Z_i^{(2)}\Lambda_1^{-2}\rangle\bigg| \geq \frac{t}{2C}\bigg),\label{ineq85}
\end{align}
where $\{Z_1^{(1)}, \dots, Z_n^{(1)}\}$ and $\{Z_1^{(2)}, \dots, Z_n^{(2)}\}$ are two independent copies of $\{Z_1, \dots, Z_n\}$. By Lemma \ref{lm:third_moment_gaussian}, with probability at least $1 - p^{-3}$, 
\begin{align*}
		&\Big|\sum_{i=1}^n\|Z_i^{(1)}\|_{\rm F}^2\langle Z_i^{(1)}, \xi_i\sum_{k = 1}^n\xi_kZ_k^{(2)}\Lambda_1^{-2}\rangle\Big|\bigg|\left\{\xi_k, Z_k^{(2)}\right\}_{k=1}^n\\
		\leq& Cpr\Big\|\sum_{k = 1}^n\xi_kZ_k^{(2)}\Lambda_1^{-2}\Big\|_{\rm F}\|\vec{\xi}\|_2\sqrt{\log(p)}
		\leq Cpr\Big\|\sum_{k = 1}^n\xi_kZ_k\Big\|_{\rm F}\lambda_{\submin}^{-2}\|\vec{\xi}\|_2\sqrt{\log(p)}.
\end{align*}
By Lemma \ref{lm:sub_Gaussian}, we have
\begin{equation*}
		\PP\Big(\Big\|\sum_{k = 1}^n\xi_kZ_k\Big\|_{\rm F} \geq C_2\sqrt{npr}\Big) \leq e^{-C_1pr},
\end{equation*}
The previous two inequalities and \eqref{ineq37} together imply that
\begin{equation}\label{ineq84}
		\begin{split}
		\PP\Big(\Big|\sum_{i=1}^n\|Z_i^{(1)}\|_{\rm F}^2\langle Z_i^{(1)}, \xi_i\sum_{k = 1}^n\xi_kZ_k^{(2)}\Lambda_1^{-2}\rangle\Big| \geq C_2np^{3/2}r^{3/2}\sqrt{\log(p)}\lambda_{\submin}^{-2}\Big) \leq p^{-3}.
		\end{split}
\end{equation}
Since
$$
\sum_{i=1}^n\|Z_i^{(1)}\|_{\rm F}^2\langle Z_i^{(1)}, \xi_i^2Z_i^{(2)}\Lambda_1^{-2}\rangle\bigg|\left\{\xi_k, Z_k^{(1)}\right\}_{k=1}^n \sim N(0, \sum_{i=1}^{n}\xi_i^4\|Z_i^{(1)}\|_{\rm F}^4\|Z_i^{(1)}\Lambda_1^{-2}\|_{\rm F}^2),
$$
we know that
\begin{equation*}
		\PP\bigg(\bigg|\sum_{i=1}^n\|Z_i^{(1)}\|_{\rm F}^2\langle Z_i^{(1)}, \xi_i^2Z_i^{(2)}\Lambda_1^{-2}\rangle\bigg| \geq C\sqrt{\sum_{i=1}^{n}\xi_i^4\|Z_i^{(1)}\|_{\rm F}^4\|Z_i^{(1)}\Lambda_1^{-2}\|_{\rm F}^2\log(p)}\bigg|\left\{\xi_k, Z_k^{(1)}\right\}_{k=1}^n\bigg) \leq p^{-3}.
\end{equation*}
By Cauchy-Schwarz inequality, 
\begin{align*}
		\sum_{i=1}^{n}\xi_i^4\|Z_i^{(1)}\|_{\rm F}^4\|Z_i^{(1)}\Lambda_1^{-2}\|_{\rm F}^2 \leq& \sum_{i=1}^{n}\xi_i^4\|Z_i^{(1)}\|_{\rm F}^4\|Z_i^{(1)}\Lambda_1^{-2}\|_{\rm F}^2 \leq \sum_{i=1}^{n}\xi_i^4\|Z_i^{(1)}\|_{\rm F}^6\|\Lambda_1^{-2}\|^2\\ 
		\leq& \lambda_{\submin}^{-4}\Big(\sum_{i=1}^{n}\xi_i^8\Big)^{1/2}\Big(\sum_{i=1}^{n}\|Z_i^{(1)}\|_{\rm F}^{12}\Big)^{1/2}.
\end{align*}
Similarly to \eqref{ineq:Gaussian_fourth},
\begin{equation*}
		\PP\left(\sum_{i=1}^{n}\xi_i^8 \geq Cn\right) \leq n^{-3}.
\end{equation*}
By \eqref{ineq:Z_i} and the union bound, with probability at least $1 - e^{-C_1(p+\log(n))}$,
\begin{equation*}
		\sum_{i=1}^{n}\|Z_i^{(1)}\|_{\rm F}^{12} \leq C_2n\left(\sqrt{r(p + \log(n))}\right)^{12} \leq Cnr^6\left(p + \log(n)\right)^6.
\end{equation*}
By combining the previous four inequalities together, we have
\begin{equation*}
		\PP\bigg(\bigg|\sum_{i=1}^n\|Z_i^{(1)}\|_{\rm F}^2\langle Z_i^{(1)}, \xi_i^2Z_i^{(2)}\Lambda_1^{-2}\rangle\bigg| \geq C_2\sqrt{nr^3\left(p + \log(n)\right)^3\log(p)}\lambda_{\submin}^{-2}\bigg) \leq C_1p^{-3}.
\end{equation*}
By \eqref{ineq85}, \eqref{ineq84} and the previous inequality, 
\begin{equation}\label{ineq88}
		\PP\bigg(\bigg|\frac{1}{n^3}\sum_{i=1}^{n}\xi_i\|Z_i\|_{\rm F}^2\sum_{k \neq i}\xi_k\langle Z_i, Z_k\Lambda_1^{-2}\rangle\bigg| \geq C_2\frac{p^{3/2}r^{3/2}\sqrt{\log(p)}}{n^2\lambda_{\submin}^2}\bigg) \leq C_1p^{-3}.
\end{equation}
Similarly, we have
\begin{equation}\label{ineq89}
		\PP\bigg(\bigg|\frac{1}{n^3}\sum_{i=1}^{n}\xi_i\|Z_i\Lambda_1^{-1}\|_{\rm F}^2\sum_{k \neq i}\xi_k\langle Z_i, Z_k\rangle\bigg| \geq C_2\frac{p^{3/2}r^{3/2}\sqrt{\log(p)}}{n^2\lambda_{\submin}^2}\bigg) \leq C_1p^{-3}.
\end{equation}
By Lemma \ref{lm:sub_Gaussian}, with probability at least $1 - e^{-C_1pr}$,
\begin{align*}
		\bigg|\sum_{i=1}^{n}\xi_i\sum_{k = 1}^n\xi_k\langle Z_i, Z_k\Lambda_1^{-2}\rangle\bigg|
		=& \bigg|\langle \sum_{i=1}^{n}\xi_iZ_i, \sum_{k=1}^{n}\xi_kZ_k\Lambda_1^{-2}\rangle\bigg|
		\leq \Big\|\sum_{i=1}^{n}\xi_iZ_i\Big\|_{\rm F}\Big\|\sum_{k=1}^{n}\xi_kZ_k\Lambda_1^{-2}\Big\|_{\rm F}\\
		\leq& \Big\|\sum_{i=1}^{n}\xi_iZ_i\bigg\|_{\rm F}^2\left\|\Lambda_1^{-2}\right\|
		\leq C_2npr\lambda_{\submin}^{-2}.
\end{align*}
The previous inequality and \eqref{ineq86} together show that
\begin{equation}\label{ineq87}
		\bigg|\frac{1}{n^3}\sum_{i=1}^{n}\xi_i\sum_{k \neq i}\xi_k\langle Z_i, Z_k\Lambda_1^{-2}\rangle\bigg| \leq C_2\frac{r(p + \log(n))}{n^2\lambda_{\submin}^2}.
\end{equation}
Now, we consider $\big|\frac{1}{n^3}\sum_{i=1}^{n}\sum_{j \neq k \neq i}\xi_j\xi_k\left[\langle Z_j, Z_i\rangle\langle Z_i, Z_k\Lambda_1^{-2}\rangle - \langle Z_j, Z_k\Lambda_1^{-2}\rangle\right]\big|$. By \cite[Theorem 3.4.1]{de2012decoupling}, for any $t \geq 0$,
\begin{align}
		&\PP\bigg(\bigg|\frac{1}{n^3}\sum_{i=1}^{n}\sum_{j \neq k \neq i}\xi_j\xi_k\left[\langle Z_j, Z_i\rangle\langle Z_i, Z_k\Lambda_1^{-2}\rangle - \langle Z_j, Z_k\Lambda_1^{-2}\rangle\right]\bigg| \geq t\bigg)\notag\\
		\leq& C\PP\bigg(\bigg|\frac{1}{n^3}\sum_{i=1}^{n}\sum_{j \neq k \neq i}\xi_j\xi_k\left[\langle Z_j^{(2)}, Z_i^{(1)}\rangle\langle Z_i^{(1)}, Z_k^{(3)}\Lambda_1^{-2}\rangle - \langle Z_j^{(2)}, Z_k^{(3)}\Lambda_1^{-2}\rangle\right]\bigg| \geq \frac{t}{C}\bigg)\notag\\
		\leq& C\PP\bigg(\bigg|\frac{1}{n^3}\sum_{i=1}^{n}\bigg[\langle \sum_{j \neq i}\xi_jZ_j^{(2)}, Z_i^{(1)}\rangle\langle Z_i^{(1)}, \sum_{k \neq i}\xi_kZ_k^{(3)}\Lambda_1^{-2}\rangle - \langle \sum_{j \neq i}\xi_jZ_j^{(2)}, \sum_{k \neq i}\xi_kZ_k^{(3)}\Lambda_1^{-2}\rangle\bigg]\bigg| \geq \frac{t}{2C}\bigg)\notag\\
		&+C\PP\bigg(\bigg|\frac{1}{n^3}\sum_{i=1}^{n}\sum_{j \neq i}\xi_j^2\left[\langle Z_j^{(2)}, Z_i^{(1)}\rangle\langle Z_i^{(1)}, Z_j^{(3)}\Lambda_1^{-2}\rangle - \langle Z_j^{(2)}, Z_j^{(3)}\Lambda_1^{-2}\rangle\right]\bigg| \geq \frac{t}{2C}\bigg).\label{ineq91}
\end{align}
Here, $\{Z_i^{(1)}\}_{i=1}^n, \{Z_i^{(2)}\}_{i=1}^n$ and $\{Z_i^{(3)}\}_{i=1}^n$ are independent copies of $\{Z_i\}_{i=1}^n$. Conditioning on $\{\xi_i\}_{i=1}^n, \{Z_i^{(2)}\}_{i=1}^n$ and $\{Z_i^{(3)}\}_{i=1}^n$, we know that
$$
\bigg\{\langle \sum_{j \neq i}\xi_jZ_j^{(2)}, Z_i^{(1)}\rangle\langle Z_i^{(1)}, \sum_{k \neq i}\xi_kZ_k^{(3)}\Lambda_1^{-2}\rangle - \langle \sum_{j \neq i}\xi_jZ_j^{(2)}, \sum_{k \neq i}\xi_kZ_k^{(3)}\Lambda_1^{-2}\rangle\bigg\}_{i=1}^n
$$
are independent. In addition, 
$$
\big< \sum\nolimits_{j \neq i}\xi_jZ_j^{(2)}, Z_i^{(1)}\big>\bigg|\left\{\xi_i, Z_i^{(2)}, Z_i^{(3)}\right\}_{i=1}^n \sim N\big(0, \|\sum\nolimits_{j \neq i}\xi_jZ_j^{(2)}\|_{\rm F}^2\big)
$$
and 
$$
\big< Z_i^{(1)}, \sum\nolimits_{k \neq i}\xi_kZ_k^{(3)}\Lambda_1^{-2}\big> \bigg|\left\{\xi_i, Z_i^{(2)}, Z_i^{(3)}\right\}_{i=1}^n \sim N\big(0, \|\sum\nolimits_{k \neq i}\xi_kZ_k^{(3)}\Lambda_1^{-2}\|_{\rm F}^2\big).
$$
Note that 
$$
\EE\bigg(\langle \sum_{j \neq i}\xi_jZ_j^{(2)}, Z_i^{(1)}\rangle\langle Z_i^{(1)}, \sum_{k \neq i}\xi_kZ_k^{(3)}\Lambda_1^{-2}\rangle\bigg|\left\{\xi_i, Z_i^{(2)}, Z_i^{(3)}\right\}_{i=1}^n\bigg) = \Big< \sum_{j \neq i}\xi_jZ_j^{(2)}, \sum_{k \neq i}\xi_kZ_k^{(3)}\Lambda_1^{-2}\Big>
$$
and
\begin{align*}
		&\bigg\|\langle \sum_{j \neq i}\xi_jZ_j^{(2)}, Z_i^{(1)}\rangle\langle Z_i^{(1)}, \sum_{k \neq i}\xi_kZ_k^{(3)}\Lambda_1^{-2}\rangle - \langle \sum_{j \neq i}\xi_jZ_j^{(2)}, \sum_{k \neq i}\xi_kZ_k^{(3)}\Lambda_1^{-2}\rangle\bigg\|_{\psi_1}\bigg|\left\{\xi_i, Z_i^{(2)}, Z_i^{(3)}\right\}_{i=1}^n\\
		\leq& C\Big\|\langle \sum_{j \neq i}\xi_jZ_j^{(2)}, Z_i^{(1)}\rangle\langle Z_i^{(1)}, \sum_{k \neq i}\xi_kZ_k^{(3)}\Lambda_1^{-2}\rangle\Big\|_{\psi_1}\bigg|\left\{\xi_i, Z_i^{(2)}, Z_i^{(3)}\right\}_{i=1}^n\\
		\leq& C\Big\|\langle \sum_{j \neq i}\xi_jZ_j^{(2)}, Z_i^{(1)}\rangle\Big\|_{\psi_2}\Big\|\langle Z_i^{(1)}, \sum_{k \neq i}\xi_kZ_k^{(3)}\Lambda_1^{-2}\rangle\Big\|_{\psi_2}\bigg|\left\{\xi_i, Z_i^{(2)}, Z_i^{(3)}\right\}_{i=1}^n\\
		\leq& C\|\sum_{j \neq i}\xi_jZ_j^{(2)}\|_{\rm F}\|\sum_{k \neq i}\xi_kZ_k^{(3)}\Lambda_1^{-2}\|_{\rm F}
		\leq C\|\sum_{j \neq i}\xi_jZ_j^{(2)}\|_{\rm F}\|\sum_{k \neq i}\xi_kZ_k^{(3)}\|_{\rm F}\lambda_{\submin}^{-2}.
\end{align*}
By Bernstein-type inequality, for any $t \geq 0$,
\begin{align*}
		&\PP\bigg(\bigg|\sum_{i=1}^{n}\Big[\langle \sum_{j \neq i}\xi_jZ_j^{(2)}, Z_i^{(1)}\rangle\langle Z_i^{(1)}, \sum_{k \neq i}\xi_kZ_k^{(3)}\Lambda_1^{-2}\rangle - \langle \sum_{j \neq i}\xi_jZ_j^{(2)}, \sum_{k \neq i}\xi_kZ_k^{(3)}\Lambda_1^{-2}\rangle\Big]\bigg| \\
		& \qquad \geq \frac{t}{\lambda_{\submin}^2}\bigg|\left\{\xi_i, Z_i^{(2)}, Z_i^{(3)}\right\}_{i=1}^n\bigg)\\
		\leq& 2\exp\bigg(-C_1\min\bigg\{\frac{t^2}{\sum_{i=1}^{n}\|\sum_{j \neq i}\xi_jZ_j^{(2)}\|_{\rm F}^2\|\sum_{k \neq i}\xi_kZ_k^{(3)}\|_{\rm F}^2}, \\
		& \qquad\qquad\qquad \frac{t}{\max_{1 \leq i \leq n}\|\sum_{j \neq i}\xi_jZ_j^{(2)}\|_{\rm F}\|\sum_{k \neq i}\xi_kZ_k^{(3)}\|_{\rm F}}\bigg\}\bigg).
\end{align*}
By Lemma \ref{lm:sub_Gaussian}, for any $i \in [n]$, with probability at least $1 - e^{-C_1(pr + \log(n))}$,
\begin{equation*}
		\max\Big\{\|\sum_{j \neq i}\xi_jZ_j^{(2)}\|_{\rm F}, \|\sum_{k \neq i}\xi_kZ_k^{(3)}\|_{\rm F}\Big\} \leq C_2\sqrt{n(pr + \log(n))}.
\end{equation*}
The union bound shows that with probability at least $1 - e^{-C_1(pr + \log(n))}$,
\begin{equation*}
		\max\bigg\{\|\sum_{j \neq i}\xi_jZ_j^{(2)}\|_{\rm F}, \|\sum_{k \neq i}\xi_kZ_k^{(3)}\|_{\rm F}\bigg\} \leq C_2\sqrt{n(pr + \log(n))}, \quad \forall i \in [n].
\end{equation*}
Therefore, with probability at least $1 - e^{-C_1(pr + \log(n))}$,
\begin{equation*}
		\sum_{i=1}^{n}\|\sum_{j \neq i}\xi_jZ_j^{(2)}\|_{\rm F}^2\|\sum_{k \neq i}\xi_kZ_k^{(3)}\|_{\rm F}^2 \leq C_2n^3\left(pr + \log(n)\right)^2,
\end{equation*}
and
\begin{equation*}
		\max_{1 \leq i \leq n}\|\sum_{j \neq i}\xi_jZ_j^{(2)}\|_{\rm F}\|\sum_{k \neq i}\xi_kZ_k^{(3)}\|_{\rm F} \leq C_2n(pr + \log(n)).
\end{equation*}
Thus
\begin{align}
		\PP\bigg(\bigg|\sum_{i=1}^{n}\Big[\langle \sum_{j \neq i}\xi_jZ_j^{(2)}, Z_i^{(1)}\rangle\langle Z_i^{(1)}, \sum_{k \neq i}\xi_kZ_k^{(3)}\Lambda_1^{-2}\rangle - &\langle \sum_{j \neq i}\xi_jZ_j^{(2)}, \sum_{k \neq i}\xi_kZ_k^{(3)}\Lambda_1^{-2}\rangle\Big]\bigg| \notag\\
		\geq& C_2\frac{n^{3/2}(pr + \log(n))\sqrt{\log(p)}}{\lambda_{\submin}^2}\bigg)
		\leq C_1p^{-3}.
\label{ineq92}
\end{align}
Similarly to \eqref{ineq93}, for any $i \in [n]$, we have
\begin{align*}
		\Big\|\langle Z_j^{(2)}, Z_i^{(1)}\rangle\langle Z_i^{(1)}, Z_j^{(3)}\Lambda_1^{-2}\rangle &- \langle Z_j^{(2)}, Z_j^{(3)}\Lambda_1^{-2}\rangle\Big\|_{\psi_1}\bigg|Z_i^{(1)}, \xi_j\\
		\leq& C\|Z_i^{(1)}\|_{\rm F}\|Z_i^{(1)}\Lambda_1^{-2}\|_{\rm F} \leq C\|Z_i^{(1)}\|_{\rm F}^2\lambda_{\submin}^{-2}.
\end{align*}
By Bernstein-type inequality,
\begin{align*}
		\PP\bigg(\bigg|\sum_{j \neq i}\xi_j^2\Big[\langle Z_j^{(2)}, Z_i^{(1)}\rangle\langle Z_i^{(1)}, Z_j^{(3)}\Lambda_1^{-2}\rangle - &\langle Z_j^{(2)}, Z_j^{(3)}\Lambda_1^{-2}\rangle\Big]\bigg| \notag\\
		\geq& C_2\|Z_i^{(1)}\|_{\rm F}^2\lambda_{\submin}^{-2}\Big(\sum_{j \neq i}\xi_j^4\Big)^{1/2}\log(n)\bigg|\{Z_i, \xi_i\}_{i = 1}^n\bigg)
		\leq n^{-3}.
\end{align*}
By \eqref{ineq:Z_i}, \eqref{ineq:Gaussian_fourth}, the previous inequality and the union bound together show that
\begin{align*}
		&\PP\bigg(\bigg|\sum_{i=1}^{n}\sum_{j \neq i}\xi_j^2\Big[\langle Z_j^{(2)}, Z_i^{(1)}\rangle\langle Z_i^{(1)}, Z_j^{(3)}\Lambda_1^{-2}\rangle - \langle Z_j^{(2)}, Z_j^{(3)}\Lambda_1^{-2}\rangle\Big]\bigg| \geq C_2\frac{n^{3/2}r(p + \log(n))\log(n)}{\lambda_{\submin}^2}\bigg)
		\leq n^{-3}.
\end{align*}
By \eqref{ineq91}, \eqref{ineq92} and the previous inequality, we have
\begin{align}\label{ineq94}
		&\PP\bigg(\bigg|\frac{1}{n^3}\sum_{i=1}^{n}\sum_{j \neq k \neq i}\xi_j\xi_k\Big[\langle Z_j, Z_i\rangle\langle Z_i, Z_k\Lambda_1^{-2}\rangle - \langle Z_j, Z_k\Lambda_1^{-2}\rangle\Big]\bigg| \geq C_2\frac{r(p + \log(n))\log(n)}{n^{3/2}\lambda_{\submin}^2}\bigg)
		\leq C_1p^{-3}.
\end{align}
By combining \eqref{ineq90}, \eqref{ineq80}, \eqref{ineq81}, \eqref{ineq88}, \eqref{ineq89}, \eqref{ineq87} and the previous inequality, we conclude that with probability at least $1 - C_1p^{-3}$,
\begin{align}
		\bigg|\tr\bigg(G_1^\top(G_1G_1^\top)^{-2}G_1&\Big\{\Big[J_1 - \frac{1}{n}\sum_{i=1}^{n}\langle J_1, \calM_1(\calX_i)\rangle\calM_1(\calX_i)\Big](U_2 \otimes U_3)\Big\}^\top L_1\bigg)\bigg|\notag\\
		\leq& C_2\left(\frac{p^2r^2}{n^2\lambda_{\submin}^2} + \frac{rp\log(n)}{n^{3/2}\lambda_{\submin}^2} + \frac{r \log^2(n)}{n^{3/2}\lambda_{\submin}^2}\right).\label{ineq95}
\end{align}
By Lemma \ref{lm:sub_Gaussian}, with probability at least $1 - p^{-3}$,
\begin{equation*}
		\|U_1^\top L_1\| = \bigg\|\frac{1}{n}U_1^\top \bigg(\sum_{j=1}^n\xi_j\calM_1(\calX_j) \bigg)(U_2 \otimes U_3) \bigg\| \leq C_2\sqrt{\frac{r^2 + \log(p)}{n}}
\end{equation*}
and
\begin{equation*}
		\| J_1\|_{\rm F} \leq \|L_1\|_{\rm F}\|G_1^\top(G_1G_1^\top)^{-1}G_1\| = \|L_1\|_{\rm F} \leq C_2\sqrt{\frac{p}{n}}.
\end{equation*}
Therefore, by Lemma \ref{lm:Gaussian_ensemble},
\begin{align}
		&\bigg|\tr\bigg(G_1^\top(G_1G_1^\top)^{-2}G_1\Big\{\Big[J_1 - \frac{1}{n}\sum_{i=1}^{n}\langle J_1, \calM_1(\calX_i)\rangle\calM_1(\calX_i)\Big](U_2 \otimes U_3)\Big\}^\top\calP_{U_1}L_1\bigg)\bigg|\notag\\
		\leq& r\|G_1^\top(G_1G_1^\top)^{-2}G_1\|\bigg\|U_1^\top\Big[J_1 - \frac{1}{n}\sum_{i=1}^{n}\langle J_1, \calM_1(\calX_i)\rangle\calM_1(\calX_i)\Big](U_2 \otimes U_3)\bigg\|\|U_1^\top L_1\|\notag\\
		\leq& C_2\lambda_{\submin}^{-2}\cdot r\sqrt{\frac{r^2 + \log(p)}{n}}\cdot \sqrt{\frac{pr}{n}}\| J_1\|_{\rm F}
		\leq C_2\left(\frac{r^{5/2}p}{n^{3/2}\lambda_{\submin}^2} + \frac{r^{3/2}p\sqrt{\log(p)}}{n^{3/2}\lambda_{\submin}^2}\right).\label{ineq96}
\end{align}
By \eqref{ineq77}, \eqref{ineq95} and \eqref{ineq96}, with probability at least $1 - C_1p^{-3}+C_2e^{-c_0p}$, 
\begin{align}
		&\bigg|\tr\bigg(G_1^\top(G_1G_1^\top)^{-2}G_1\Big\{\Big[\frakE_{U_1, 3}^{(1)} - \frac{1}{n}\sum_{i=1}^{n}\langle \frakE_{U_1, 3}^{(1)}, \calM_1(\calX_i)\rangle\calM_1(\calX_i)\Big](U_2 \otimes U_3)\Big\}^\top\calP_{U_1}^\perp L_1\bigg)\bigg|\notag\\
		\leq& C_2\left(\kappa_0^3\frac{p^2r^{5/2}}{n^2\lambda_{\submin}^3} + \frac{p^2r^2}{n^2\lambda_{\submin}^2} + \frac{r^{5/2}p\log(n)}{n^{3/2}\lambda_{\submin}^2} + \frac{r\log^2(n)}{n^{3/2}\lambda_{\submin}^2}\right).\label{ineq97}
\end{align}
Combining \eqref{ineq98}, \eqref{ineq75}, \eqref{ineq76} and \eqref{ineq97} together, we have
\begin{equation}\label{ineq115}
		\big|S_1^{(1.5)}\big| \leq C_2\left(\kappa_0^3\frac{p^2r^{5/2}}{n^2\lambda_{\submin}^3} + \kappa_0^2\frac{p^2r^{3}}{n^2\lambda_{\submin}^2} + \kappa_0\frac{r^3p\log(n)}{n^{3/2}\lambda_{\submin}^2} + \frac{r\log^2(n)}{n^{3/2}\lambda_{\submin}^2}\right).
\end{equation}

\subparagraph*{-Step 3.2.2: bounding $|S_2^{(1)}|$ and $|S_3^{(1)}|$.}
Let $\frakE_2^{(0)} = \hat{U}_2^{(0.5)} - U_2R_2^{(0)\top}.$ For $k \geq 1$, denote
$$
\left(\frakP_2^{(0)}\right)^{-k} = \begin{cases}
		\begin{pmatrix}
		0 & U_2R_2^{(0)\top}\\
		R_2^{(0)}U_2^\top & 0
		\end{pmatrix}, & \text{ if } k \text{ is old},\\
		\begin{pmatrix}
		U_2U_2^\top & 0\\
		0 & I_{r_2}
		\end{pmatrix}, &\text{ if } k \text{ is even},
		\end{cases}
$$
and 
$$
\left(\frakP_2^{(0)}\right)^{0} = \begin{pmatrix}
		U_{2\perp}U_{2\perp}^\top & 0\\
		0 & 0
		\end{pmatrix}.
$$
Let 
$$
E_2^{(0)} = \begin{pmatrix}
		0 & \frakE_2^{(0)}\\
		\frakE_2^{(0)\top} & 0
		\end{pmatrix}.
$$
Similarly to \eqref{ineq2}, with probability at least $1 - C_1p^{-3} - C_1e^{-c_1p}$, 
\begin{equation*}
		\big\|E_2^{(0)}\big\| \leq C\frac{\sqrt{p/n}}{\lambda_{\submin}}.
\end{equation*}
By Lemma \ref{lem:spectral}, with probability at least $1 - C_1p^{-3} - C_1e^{-c_1p}$, 
\begin{equation}
		\begin{pmatrix}
		\calP_{\hat U_2^{(1)}} - \calP_{U_2} & 0\\
		0 & 0
		\end{pmatrix}
		=
		\sum_{k \geq 1}\calS_{U_2^{(0)}, k}\big(E_2^{(0)}\big)
\end{equation}
where 
$$
\calS_{U_2^{(0)}, k}(X) = \sum_{s_1+\cdots+s_{k+1}=k}(-1)^{1+\tau(\bs)}\cdot \big(\frakP_2^{(0)}\big)^{-s_1}X\big(\frakP_2^{(0)}\big)^{-s_2}X\big(\frakP_2^{(0)}\big)^{-s_3}\cdots\big(\frakP_2^{(0)}\big)^{-s_k}X\big(\frakP_2^{(0)}\big)^{-s_{k+1}}.
$$
By \eqref{ineq63}, with probability at least $1 - C_1p^{-3} - C_1e^{-c_1p}$, 
\begin{equation}\label{ineq100}
		\sum_{k \geq 2}\Big\|\calS_{U_2^{(0)}, k}\big(E_2^{(0)}\big)\Big\| \leq \sum_{k \geq 2}\big(4\|E_2^{(0)}\|\big)^{k} = \frac{16\|E_2^{(0)}\|^2}{1 - 4\|E_2^{(0)}\|} \leq C\frac{p}{n\lambda_{\submin}^2}.
\end{equation}
Note that
\begin{align*}
		\calS_{U_2^{(0)}, 1}\big(E_2^{(0)}\big) =&\big(\frakP_2^{(0)}\big)^{-1}E_2^{(0)}\big(\frakP_2^{(0)}\big)^{0} + \big(\frakP_2^{(0)}\big)^{0}E_2^{(0)}\big(\frakP_2^{(0)}\big)^{-1}\\
		=& \begin{pmatrix}
		0 & U_2R_2^{(0)\top}\\
		R_2^{(0)}U_2^\top & 0
		\end{pmatrix}
		\begin{pmatrix}
		0 & \frakE_2^{(0)}\\
		\frakE_2^{(0)\top} & 0
		\end{pmatrix}
		\begin{pmatrix}
		U_{2\perp}U_{2\perp}^\top & 0\\
		0 & 0
		\end{pmatrix}\\
		&+ \begin{pmatrix}
		U_{2\perp}U_{2\perp}^\top & 0\\
		0 & 0
		\end{pmatrix}
		\begin{pmatrix}
		0 & \frakE_2^{(0)}\\
		\frakE_2^{(0)\top} & 0
		\end{pmatrix}
		\begin{pmatrix}
		0 & U_2R_2^{(0)\top}\\
		R_2^{(0)}U_2^\top & 0
		\end{pmatrix}\\
		=& \begin{pmatrix}
		U_2R_2^{(0)\top}\frakE_2^{(0)\top}U_{2\perp}U_{2\perp}^\top & 0\\
		0 & 0
		\end{pmatrix}
		+\begin{pmatrix}
		U_{2\perp}U_{2\perp}^\top\frakE_2^{(0)}R_2^{(0)}U_2^\top & 0\\
		0 & 0
		\end{pmatrix}\\
		=& \begin{pmatrix}
		U_2R_2^{(0)\top}\frakE_2^{(0)\top}\calP_{U_{2}}^\perp + \calP_{U_{2}}^\perp\frakE_2^{(0)}R_2^{(0)}U_2^\top & 0\\
		0 & 0
		\end{pmatrix},
\end{align*}
with probability at least $1 - C_1p^{-3} - C_1e^{-c_1p}$,
\begin{align*}
		&\Big\|\Big(\calP_{\hat U_2^{(1)}} - \calP_{U_2}\Big)U_2 - \calP_{U_{2}}^\perp\frakE_2^{(0)}R_2^{(0)}\Big\|\\
		=& \Big\|\Big(\calP_{\hat U_2^{(1)}} - \calP_{U_2} - \big(U_2R_2^{(0)\top}\frakE_2^{(0)\top}\calP_{U_{2}}^\perp + \calP_{U_{2}}^\perp\frakE_2^{(0)}R_2^{(0)}U_2^\top\big)\Big)U_2\Big\|\\
		=& \left\|\begin{pmatrix}
		\calP_{\hat U_2^{(1)}} - \calP_{U_2} & 0\\
		0 & 0
		\end{pmatrix} - \calS_{U_2^{(0)}, 1}\left(E_2^{(0)}\right)\right\|
		\leq \sum_{k \geq 2}\left\|\calS_{U_2^{(0)}, k}\left(E_2^{(0)}\right)\right\|
		\leq C\frac{p}{n\lambda_{\submin}^2}.
\end{align*}
In addition,  with probability at least $1 - C_1p^{-3} - C_1e^{-c_1p}$,
\begin{align*}
		&\Big\|\hat{U}_2^{(1)}R_2^{(1)} - U_2 - \big(\calP_{\hat U_2^{(1)}} - \calP_{U_2}\big)U_2\Big\| = \Big\|\hat{U}_2^{(1)}\big(R_2^{(1)} - \hat U_2^{(1)\top}U_2\big)\Big\| \leq \big\|R_2^{(1)} - \hat U_2^{(1)\top}U_2\big\| \leq C\frac{p}{n\lambda_{\submin}^2},
\end{align*}
we know that with probability $1 - C_1p^{-3} - C_1e^{-c_1p}$,
\begin{equation}\label{ineq119}
		\left\|\hat{U}_2^{(1)}R_2^{(1)} - U_2 -  \calP_{U_{2}}^\perp\frakE_2^{(0)}R_2^{(0)}\right\| \leq C\frac{p}{n\lambda_{\submin}^2},
\end{equation}
and
\begin{align*}
		&\Big\|\Delta T_{U_2,1}^{(1)} - U_1G_1\Big(\big(\calP_{U_{2}}^\perp\frakE_2^{(0)}R_2^{(0)}\big)^\top \otimes U_3^\top\Big)\Big\|_{\rm F}\\
		=&\Big\|U_1G_1\Big(\big(\hat{U}_2^{(1)}R_2^{(1)} - U_2 - \calP_{U_{2}}^\perp\frakE_2^{(0)}R_2^{(0)}\big)^\top \otimes U_3^\top\Big)\Big\|_{\rm F}\\
		=&\Big\|\calG \times_1 U_1 \times_2 \Big(\hat{U}_2^{(1)}R_2^{(1)} - U_2 - \calP_{U_{2}}^\perp\frakE_2^{(0)}R_2^{(0)}\Big) \times U_3\Big\|_{\rm F}\\
		\leq& \|G_2\|\Big\|\hat{U}_2^{(1)}R_2^{(1)} - U_2 - \calP_{U_{2}}^\perp\frakE_2^{(0)}R_2^{(0)}\Big\|_{\rm F}
		\leq \kappa_0\lambda_{\submin}\cdot C\frac{p\sqrt{r}}{n\lambda_{\submin}^2}
		\leq C\kappa_0\frac{p\sqrt{r}}{n\lambda_{\submin}}.
\end{align*}
Let 
$$
\widetilde\Delta T_{U_2,1}^{(1)} = U_1G_1\Big(\big(\calP_{U_{2}}^\perp\frakE_2^{(0)}R_2^{(0)}\big)^\top \otimes U_3^\top\Big)
$$
and 
$$
\widetilde S_2^{(1)} = \tr\bigg(G_1^\top(G_1G_1^\top)^{-2}G_1\Big\{\Big[\widetilde\Delta T_{U_2,1}^{(1)} - \frac{1}{n}\sum_{i=1}^{n}\langle \widetilde\Delta T_{U_2,1}^{(1)}, \calM_1(\calX_i)\rangle\calM_1(\calX_i)\Big](U_2 \otimes U_3)\Big\}^\top\calP_{U_1}^\perp L_1\bigg).
$$
By Lemmas \ref{lm:Gaussian_ensemble} and \ref{lm:sub_Gaussian}, with probability at least $1 - C_1p^{-3} - C_1e^{-c_1p}$,
\begin{align}
		\big|S_2^{(1)} - \widetilde S_2^{(1)}\big| \leq& r\|G_1^\top(G_1G_1^\top)^{-2}G_1\|\cdot\sqrt{\frac{p}{n}}\Big\|\Delta T_{U_2,1}^{(1)} - \widetilde\Delta T_{U_2,1}^{(1)}\Big\|_{\rm F}\|L_1\|\notag\\
		\leq& C_2r\lambda_{\submin}^{-2}\sqrt{\frac{p}{n}}\kappa_0\frac{p\sqrt{r}}{n\lambda_{\submin}}\cdot\sqrt{\frac{p}{n}}
		\leq C_2\kappa_0\frac{p^2r^{3/2}}{n^2\lambda_{\submin}^3}.\label{ineq99}
\end{align}
By \eqref{ineq55}, we have
\begin{align}
		\frakE_2^{(0)}
		=& \bigg(\calM_2(\Delta \calT_{2}^{(0.5)})(\hat U_1^{(0)} \otimes \hat U_3^{(0)}) - \frac{1}{n}\sum_{i=1}^{n}\langle \Delta \calT_{2}^{(0.5)}, \calX_i\rangle\calM_2(\calX_i)(\hat U_1^{(0)} \otimes \hat U_3^{(0)}) \bigg)\hat G_2^{(0)\top}\big(\hat G_2^{(0)}\hat G_2^{(0)\top}\big)^{-1}\notag\\
		&+ \Big(U_2G_2\big((\hat U_1^{(0)\top}U_1)^\top \otimes (\hat U_3^{(0)\top}U_3)^\top\big) - U_2R_2^{(0)\top}\hat G_2^{(0)}\Big)\hat G_2^{(0)\top}\big(\hat G_2^{(0)}\hat G_2^{(0)\top}\big)^{-1}\notag\\
		&+ \frac{1}{n}\sum_{i=1}^{n}\xi_i\calM_2(\calX_i)(\hat U_1^{(0)} \otimes \hat U_3^{(0)})\hat G_2^{(0)\top}\big(\hat G_2^{(0)}\hat G_2^{(0)\top}\big)^{-1}\notag\\
		=:& \frakJ_{U_2, 1}^{(0)} + \frakJ_{U_2, 2}^{(0)} + \frakJ_{U_2, 3}^{(0)}.\label{ineq101}
\end{align}
Here, $\Delta\calT_2^{(0.5)} = \hat\calG^{(0)} \times_1 \hat{U}_1^{(0)} \times_2 \hat{U}_2^{(0.5)} \times_3 \hat{U}_3^{(0)} - \calT$.
Let 
$$
\frakE_{U_2, i}^{(0)} = U_1G_1\left(\left(\calP_{U_{2}}^\perp\frakJ_{U_2, i}^{(0)}R_2^{(0)}\right)^\top \otimes U_3^\top\right), \quad \forall i \in [3].
$$
Then 
$$
\widetilde\Delta T_{U_2,1}^{(1)} = \frakE_{U_2, 1}^{(0)} + \frakE_{U_2, 2}^{(0)} + \frakE_{U_2, 3}^{(0)}.
$$
Similarly to \eqref{ineq75} and \eqref{ineq76}, with probability at least $1 - C_1p^{-3} - C_1e^{-c_1p}$,
\begin{align}\label{ineq102}
		\bigg|\tr\bigg(G_1^\top(G_1G_1^\top)^{-2}&G_1\Big\{\Big[\frakE_{U_2, 1}^{(0)} - \frac{1}{n}\sum_{i=1}^{n}\langle \frakE_{U_2, 1}^{(0)}, \calM_1(\calX_i)\rangle\calM_1(\calX_i)\Big](U_2 \otimes U_3)\Big\}^\top\calP_{U_1}^\perp L_1\bigg)\bigg|\notag\\
		& \leq C_2\kappa_0^2\frac{p^{2}r^{3}}{n^2\lambda_{\submin}^2},
\end{align}
and
\begin{align}
		\bigg|\tr\bigg(G_1^\top(G_1G_1^\top)^{-2}G_1&\Big\{\Big[\frakE_{U_2, 2}^{(0)} - \frac{1}{n}\sum_{i=1}^{n}\langle \frakE_{U_2, 2}^{(0)}, \calM_1(\calX_i)\rangle\calM_1(\calX_i)\Big](U_2 \otimes U_3)\Big\}^\top\calP_{U_1}^\perp L_1\bigg)\bigg|\nonumber\\ 
		\leq& C_2\bigg(\kappa_0^2\frac{p^2r^{5/2}}{n^2\lambda_{\submin}^{3}} + \kappa_0^2\frac{p^2r^{3}}{n^2\lambda_{\submin}^{2}} + \kappa_0\frac{r^{3}p + r^{3/2}p\sqrt{\log(p)}}{n^{3/2}\lambda_{\submin}^2}\bigg).\label{ineq103}
\end{align}
Similarly to \eqref{ineq78}, with probability at least $1 - C_1p^{-3} - C_1e^{-c_1p}$,
\begin{equation}\label{ineq105}
		\Big\|\frakJ_{U_2, 3}^{(0)} - \frac{1}{n}\sum_{i=1}^{n}\xi_i\calM_2(\calX_i)(U_1 \otimes U_3)G_2^\top\left(G_2G_2^{\top}\right)^{-1}R_2^{(0)\top}\Big\| \leq C_2\kappa_0^2\frac{p\sqrt{r}}{n\lambda_{\submin}^2}
\end{equation}
and
\begin{align}
		&\bigg\|\frakE_{U_2, 3}^{(0)} - U_1G_1\bigg(\Big(\calP_{U_{2}}^\perp\big(\frac{1}{n}\sum_{i=1}^{n}\xi_i\calM_2(\calX_i)(U_1 \otimes U_3)G_2^\top(G_2G_2^{\top})^{-1}\big)\Big)^\top \otimes U_3^\top\bigg)\bigg\|\notag\\
		=&\bigg\|U_1G_1\bigg(\Big(\calP_{U_{2}}^\perp\big(\frakJ_{U_2, 3}^{(0)} - \frac{1}{n}\sum_{i=1}^{n}\xi_i\calM_2(\calX_i)(U_1 \otimes U_3)G_2^\top(G_2G_2^{\top})^{-1}R_2^{(0)\top}\big)R_2^{(0)}\Big)^\top \otimes U_3^\top\bigg)\bigg\|\notag\\
		\leq& \|G_1\|\Big\|\frakJ_{U_2, 3}^{(0)} - \frac{1}{n}\sum_{i=1}^{n}\xi_i\calM_2(\calX_i)(U_1 \otimes U_3)G_2^\top\big(G_2G_2^{\top}\big)^{-1}R_2^{(0)\top}\Big\|
		\leq C_2\kappa_0^3\frac{p\sqrt{r}}{n\lambda_{\submin}}.\label{ineq104}
\end{align}
Let 
$$
\widetilde\frakE_{U_2, 3}^{(0)} = U_1G_1\bigg(\Big(\calP_{U_{2}}^\perp\big(\frac{1}{n}\sum_{i=1}^{n}\xi_i\calM_2(\calX_i)(U_1 \otimes U_3)G_2^\top(G_2G_2^{\top})^{-1}\big)\Big)^\top \otimes U_3^\top\bigg).
$$
The same argument for proving \eqref{ineq77} shows that with probability at least $1 - C_1p^{-3} - C_1e^{-c_1p}$,
\begin{align}
		&\bigg|\tr\bigg(G_1^\top(G_1G_1^\top)^{-2}G_1\Big\{\Big[\frakE_{U_2, 3}^{(0)} - \frac{1}{n}\sum_{i=1}^{n}\langle \frakE_{U_2, 3}^{(0)}, \calM_1(\calX_i)\rangle\calM_1(\calX_i)\Big](U_2 \otimes U_3)\Big\}^\top\calP_{U_1}^\perp L_1\bigg)\bigg|\notag\\
		\leq& \bigg|\tr\bigg(G_1^\top(G_1G_1^\top)^{-2}G_1\Big\{\Big[\widetilde\frakE_{U_2, 3}^{(0)} - \frac{1}{n}\sum_{i=1}^{n}\langle \widetilde\frakE_{U_2, 3}^{(0)}, \calM_1(\calX_i)\rangle\calM_1(\calX_i)\Big](U_2 \otimes U_3)\Big\}^\top\calP_{U_1}^\perp L_1\bigg)\bigg|\notag\\
		&+C_2\kappa_0^3\frac{p^2r^{5/2}}{n^2\lambda_{\submin}^3}.\label{ineq106}
\end{align}
Let 
$$
G_2 = U_{G_2}\Lambda_2V_{G_2}^\top
$$
be the SVD decomposition of $G_2$. Let $W_i = U_{1\perp}^\top\calM_1(\calX_i)(U_2 \otimes U_3)V_{G_1} \in \RR^{(p_1-r_1) \times r_1}$ and $\widetilde W_i = U_{2\perp}^\top\calM_2(\calX_i)(U_1 \otimes U_3)V_{G_2} \in \RR^{(p_2 - r_2) \times r_2}$. Then $W_i \stackrel{i.i.d.}{\sim} N(0,1)$ and $\widetilde{W}_i \stackrel{i.i.d.}{\sim} N(0,1)$. In addition, since $\calX_i \times_1 [U_1 \ U_{1\perp}] \times_2 [U_2 \ U_{2\perp}] \times_3 [U_3 \ U_{3\perp}] \stackrel{i.i.d.}{\sim} N(0, 1)$, $U_{1\perp}^\top\calM_1(\calX_i)(U_2 \otimes U_3)$ and $U_{2\perp}^\top\calM_2(\calX_i)(U_1 \otimes U_3)$ are independent. Therefore, $W_i$ and $\widetilde{W}_i$ are independent.
		
Note that $\widetilde\frakE_{U_2, 3}^{(0)\top}\calP_{U_1}^\perp = 0$ and
\begin{align*}
		&\langle \widetilde\frakE_{U_2, 3}^{(0)}, \calM_1(\calX_i)\rangle\\ 
		=& \left\langle U_1G_1\bigg(\Big(\calP_{U_{2}}^\perp\big(\frac{1}{n}\sum_{j=1}^{n}\xi_j\calM_2(\calX_j)(U_1 \otimes U_3)G_2^\top(G_2G_2^{\top})^{-1}\big)\Big)^\top \otimes U_3^\top\bigg), \calM_1(\calX_i)\right\rangle\\
		=& \left\langle \calG \times_1 U_1 \times_2 \bigg(\calP_{U_{2}}^\perp\Big(\frac{1}{n}\sum_{j=1}^{n}\xi_j\calM_2(\calX_j)(U_1 \otimes U_3)G_2^\top(G_2G_2^{\top})^{-1}\Big)\bigg) \times_3 U_3, \calX_i\right\rangle\\
		=& \left\langle \calP_{U_{2}}^\perp\frac{1}{n}\sum_{j=1}^{n}\xi_j\calM_2(\calX_j)(U_1 \otimes U_3)G_2^\top\big(G_2G_2^{\top}\big)^{-1}G_2\big(U_1^\top \otimes U_3^\top\big), \calM_2(\calX_i)\right\rangle\\
		=& \left\langle\frac{1}{n}\sum_{j=1}^{n}\xi_jU_{2\perp}^\top\calM_2(\calX_j)(U_1 \otimes U_3)V_{G_2}, U_{2\perp}^\top\calM_2(\calX_i)(U_1 \otimes U_3)V_{G_2}\right\rangle\\
		=& \left\langle \frac{1}{n}\sum_{j=1}^{n}\xi_j\widetilde W_j, \widetilde W_i\right\rangle,
\end{align*}
we have
\begin{align*}
		&\bigg|\tr\bigg(G_1^\top(G_1G_1^\top)^{-2}G_1\Big\{\Big[\widetilde\frakE_{U_2, 3}^{(0)} - \frac{1}{n}\sum_{i=1}^{n}\langle \widetilde\frakE_{U_2, 3}^{(0)}, \calM_1(\calX_i)\rangle\calM_1(\calX_i)\Big](U_2 \otimes U_3)\Big\}^\top\calP_{U_1}^\perp L_1\bigg)\bigg|\\
		=&\bigg|\tr\bigg(\bigg(\frac{1}{n}\sum_{i=1}^{n}\Big< \frac{1}{n}\sum_{j=1}^{n}\xi_j\widetilde W_j, \widetilde W_i\Big> U_{1\perp}\calM_1(\calX_i)(U_2 \otimes U_3)V_{G_1}\bigg)^\top U_{1\perp}L_1V_{G_1}\Lambda_1^{-2}\bigg)\bigg|\\
		=& \bigg|\frac{1}{n^3}\sum_{i=1}^{n}\sum_{j=1}^{n}\sum_{k=1}^{n}\xi_j\xi_k\langle \widetilde{W}_i, \widetilde{W}_j\rangle\langle W_i, W_k\Lambda_1^{-2}\rangle\bigg|\\
		\leq& \bigg|\frac{1}{n^3}\sum_{i=1}^{n}\xi_i^2\big\|\widetilde{W}_i\big\|_{\rm F}^2\left\|W_i\Lambda_1^{-2}\right\|_{\rm F}^2\bigg| + \bigg|\frac{1}{n^3}\sum_{i=1}^{n}\sum_{k \neq i}\xi_i\xi_k\big\|\widetilde{W}_i\big\|_{\rm F}^2\langle W_i, W_k\Lambda_1^{-2}\rangle\bigg|\\
		&+\bigg|\frac{1}{n^3}\sum_{i=1}^{n}\sum_{j \neq i}\xi_i\xi_j\langle \widetilde{W}_i, \widetilde{W}_j\rangle\left\|W_i\Lambda_1^{-1}\right\|_{\rm F}^2\bigg| + \bigg|\frac{1}{n^3}\sum_{i=1}^{n}\sum_{j \neq i}\xi_j^2\langle \widetilde{W}_i, \widetilde{W}_j\rangle\langle W_i, W_j\Lambda_1^{-2}\rangle\bigg|\\ 
		&+ \bigg|\frac{1}{n^3}\sum_{i \neq j \neq k}\xi_j\xi_k\langle \widetilde{W}_i, \widetilde{W}_j\rangle\langle W_i, W_k\Lambda_1^{-2}\rangle\bigg|.
\end{align*}
Similarly to \eqref{ineq95}, with probability at least $1 - C_1p^{-3}$, 
\begin{align*}
		\bigg|\tr\bigg(G_1^\top(G_1G_1^\top)^{-2}G_1&\Big\{\Big[\widetilde\frakE_{U_2, 3}^{(0)} - \frac{1}{n}\sum_{i=1}^{n}\langle \widetilde\frakE_{U_2, 3}^{(0)}, \calM_1(\calX_i)\rangle\calM_1(\calX_i)\Big](U_2 \otimes U_3)\Big\}^\top\calP_{U_1}^\perp L_1\bigg)\bigg|\\
		 \leq&C_2\left(\frac{p^2r^2}{n^2\lambda_{\submin}^2} + \frac{rp\log(n)}{n^{3/2}\lambda_{\submin}^2} + \frac{r \log^2(n)}{n^{3/2}\lambda_{\submin}^2}\right).
\end{align*}
By \eqref{ineq99}, \eqref{ineq102}, \eqref{ineq103}, \eqref{ineq106} and the above inequality, with probability at least $1 - C_1p^{-3} - C_1e^{-c_1p}$,
\begin{equation}\label{ineq116}
		\big|S_2^{(1)}\big| \leq C_2\left(\kappa_0^3\frac{p^2r^{5/2}}{n^2\lambda_{\submin}^3} + \kappa_0^2\frac{p^2r^{3}}{n^2\lambda_{\submin}^2} + \kappa_0\frac{r^3p\log(n)}{n^{3/2}\lambda_{\submin}^2} + \frac{r \log^2(n)}{n^{3/2}\lambda_{\submin}^2}\right).
\end{equation}
Similarly, with probability at least $1 - C_1p^{-3} - C_1e^{-c_1p}$,
\begin{equation}\label{ineq117}
		\big|S_3^{(1)}\big| \leq C_2\left(\kappa_0^3\frac{p^2r^{5/2}}{n^2\lambda_{\submin}^3} + \kappa_0^2\frac{p^2r^{3}}{n^2\lambda_{\submin}^2} + \kappa_0\frac{r^3p\log(n)}{n^{3/2}\lambda_{\submin}^2} + \frac{r \log^2(n)}{n^{3/2}\lambda_{\submin}^2}\right).
\end{equation}
Putting \eqref{ineq69}, \eqref{ineq74}, \eqref{ineq115}, \eqref{ineq116} and \eqref{ineq117} together, we get with probability at least $1 - C_1p^{-3} - C_1e^{-c_1p}$ that
\begin{equation}\label{ineq118}
		\big|\tr\big(\frakJ_{U_1, 1}^{(1)\top} \calP_{U_1}^\perp\frakJ_{U_1, 3}^{(1)}\big)\big| \leq C_2\left(\kappa_0^3\frac{p^2r^{5/2}}{n^2\lambda_{\submin}^3} + \kappa_0^2\frac{p^2r^{3}}{n^2\lambda_{\submin}^2} + \kappa_0\frac{r^3p\log(n)}{n^{3/2}\lambda_{\submin}^2} + \frac{r \log^2(n)}{n^{3/2}\lambda_{\submin}^2}\right).
\end{equation}

\paragraph*{Final step: characterizing the distribution of $\big|\tr\big(\frakJ_{U_1, 3}^{(1)\top} \calP_{U_1}^\perp\frakJ_{U_1, 3}^{(1)}\big)\big|$.}
By \eqref{ineq59}, with probability at least $1 - p^{3} - C_1e^{-c_1p}$,
\begin{align*}
		&\Big\|\frakJ_{U_1, 3}^{(1)} - \frac{1}{n}\sum_{i=1}^{n}\xi_i\calM_1(\calX_i)\Big((\hat U_2^{(1)}R_2^{(1)}) \otimes (\hat U_3^{(1)}R_3^{(1)})\Big)G_1^\top(G_1G_1^\top)^{-1}R_1^{(1)\top}\Big\|\\
		=&\Big\|\frac{1}{n}\sum_{i=1}^{n}\xi_i\calM_1(\calX_i)(\hat U_2^{(1)} \otimes \hat U_3^{(1)})\Big(\hat G_1^{(1)\top}\left(\hat G_1^{(1)}\hat G_1^{(1)\top}\right)^{-1} - (R_2^{(1)} \otimes R_3^{(1)})G_1^\top(G_1G_1^\top)^{-1}R_1^{(1)\top}\Big)\Big\|\\
		\leq& C_2\kappa_0^2\lambda_{\submin}^{-2}\left(\sqrt{\frac{r^2+\log(p)}{n}} + \kappa_0\frac{pr}{n} + \kappa_0\frac{p\sqrt{r}}{n\lambda_{\submin}}\right)\Big\|\frac{1}{n}\sum_{i=1}^{n}\xi_i\calM_1(\calX_i)(\hat U_2^{(1)} \otimes \hat U_3^{(1)})\Big\|\\
		\leq&  C_2\kappa_0^2\lambda_{\submin}^{-2}\left(\sqrt{\frac{r^2+\log(p)}{n}} + \kappa_0\frac{pr}{n} + \kappa_0\frac{p\sqrt{r}}{n\lambda_{\submin}}\right)\\
		&\cdot\left(\sqrt{\frac{p}{n}}+ \sqrt{\frac{pr^2}{n}}\Big\|\hat U_2^{(1)} - U_2R_2^{(1)\top}\Big\| + \sqrt{\frac{pr^2}{n}}\Big\|\hat U_3^{(1)} - U_3R_3^{(1)\top}\Big\|\right)\\
		\leq& C_2\kappa_0^2\lambda_{\submin}^{-2}\left(\sqrt{\frac{r^2+\log(p)}{n}} + \kappa_0\frac{pr}{n} + \kappa_0\frac{p\sqrt{r}}{n\lambda_{\submin}}\right)\sqrt{\frac{p}{n}}\\
		\leq& C_2\left(\kappa_0^2\frac{\sqrt{pr^2+p\log(p)}}{n\lambda_{\submin}^2} + \kappa_0^3\frac{p^{3/2}r}{n^{3/2}\lambda_{\submin}^2} + \kappa_0^3\frac{p^{3/2}r^{1/2}}{n^{3/2}\lambda_{\submin}^3}\right).
\end{align*}
Moreover, by \eqref{ineq119} and Lemma \ref{lm:sub_Gaussian}, with probability at least $1 - C_1p^{-3} - C_1e^{-c_1p}$,
\begin{align*}
		&\Big\|\frac{1}{n}\sum_{i=1}^{n}\xi_i\calM_1(\calX_i)\left((\hat U_2^{(1)}R_2^{(1)}) \otimes (\hat U_3^{(1)}R_3^{(1)}) - (U_2 +  \calP_{U_{2}}^\perp\frakE_2^{(0)}R_2^{(0)}) \otimes (U_3 +  \calP_{U_{3}}^\perp\frakE_3^{(0)}R_3^{(0)})\right)G_1^\top(G_1G_1^\top)^{-1}\Big\|\\
		\leq& C_2\sqrt{\frac{pr}{n}}\big\|\hat U_2^{(1)}R_2^{(1)}-U_2-\calP_{U_{2}}^\perp\frakE_2^{(0)}R_2^{(0)}\big\|\big\|G_1^\top(G_1G_1^\top)^{-1}\big\|\\
		&+C_2\sqrt{\frac{pr}{n}}\big\|\hat U_3^{(1)}R_3^{(1)}-U_3-\calP_{U_{3}}^\perp\frakE_3^{(0)}R_3^{(0)}\big\|\big\|G_1^\top(G_1G_1^\top)^{-1}\big\|\\
		&+C_2\sqrt{\frac{pr}{n}}\big\|\hat U_2^{(1)}R_2^{(1)}-U_2-\calP_{U_{2}}^\perp\frakE_2^{(0)}R_2^{(0)}\big\|\big\|\hat U_3^{(1)}R_3^{(1)}-U_3-\calP_{U_{3}}^\perp\frakE_3^{(0)}R_3^{(0)}\big\|_{\rm F}\big\|G_1^\top(G_1G_1^\top)^{-1}\big\|\\
		\leq& C_2\frac{p^{3/2}r^{1/2}}{n^{3/2}\lambda_{\submin}^3}.
\end{align*}
Similarly to \eqref{ineq:J_{U_1,1}}, with probability at least $1 - C_1p^{-3} - C_1e^{-c_1p}$,
\begin{equation}\label{ineq120}
		\left\|\frakJ_{U_2, 1}^{(0)}\right\| \leq C_2\kappa_0\frac{pr}{n\lambda_{\submin}}.
\end{equation}
Note that $\calP_{U_{2}}^\perp\frakJ_{U_2, 2}^{(0)} = 0$, by \eqref{ineq105} and \eqref{ineq120}, with probability at least $1 - C_1p^{-3} - C_1e^{-c_1p}$,
\begin{align}
		&\Big\|\calP_{U_{2}}^\perp\frakE_2^{(0)}R_2^{(0)}-\frac{1}{n}\sum_{i=1}^{n}\xi_i\calP_{U_{2}}^\perp\calM_2(\calX_i)(U_1 \otimes U_3)G_2^\top(G_2G_2^{\top})^{-1}\Big\|\notag\\
		 \leq& \left\|\frakJ_{U_2, 1}^{(0)}\right\| + \Big\|\frakJ_{U_2, 3}^{(0)} - \frac{1}{n}\sum_{i=1}^{n}\xi_i\calM_2(\calX_i)(U_1 \otimes U_3)G_2^\top(G_2G_2^{\top})^{-1}R_2^{(0)\top}\Big\|\notag\\
		\leq& C_2\left(\kappa_0\frac{pr}{n\lambda_{\submin}} + \kappa_0^2\frac{p\sqrt{r}}{n\lambda_{\submin}^2}\right).\label{ineq121}
\end{align}
Similarly, with probability at least $1 - C_1p^{-3} - C_1e^{-c_1p}$,
\begin{equation}\label{ineq122}
		\Big\|\calP_{U_{3}}^\perp\frakE_3^{(0)}R_3^{(0)}-\frac{1}{n}\sum_{i=1}^{n}\xi_i\calP_{U_{3}}^\perp\calM_3(\calX_i)(U_1 \otimes U_2)G_3^\top(G_3G_3^{\top})^{-1}\Big\| \leq C_2\left(\kappa_0\frac{pr}{n\lambda_{\submin}} + \kappa_0^2\frac{p\sqrt{r}}{n\lambda_{\submin}^2}\right).
\end{equation}
In addition, with probability at least $1 - C_1p^{-3} - C_1e^{-c_1p}$,
\begin{equation*}
		\big\|\calP_{U_{2}}^\perp\frakE_2^{(0)}R_2^{(0)}\big\|\big\|\calP_{U_{3}}^\perp\frakE_3^{(0)}R_3^{(0)}\big\| \leq \big\|\frakE_2^{(0)}\big\|\big\|\frakE_3^{(0)}\big\| \leq C_2\frac{p}{n\lambda_{\submin}^2}.
\end{equation*}
Let
$$
Q_2 = \frac{1}{n}\sum_{i=1}^{n}\xi_i\calP_{U_{2}}^\perp\calM_2(\calX_i)(U_1 \otimes U_3)G_2^\top\big(G_2G_2^{\top}\big)^{-1}
$$
and
$$
Q_3 = \frac{1}{n}\sum_{i=1}^{n}\xi_i\calP_{U_{3}}^\perp\calM_3(\calX_i)(U_1 \otimes U_2)G_3^\top\big(G_3G_3^{\top}\big)^{-1}.
$$
Combining the previous six inequalities together, with probability at least $1 - C_1p^{-3} - C_1e^{-c_1p}$, 
\begin{align}
		&\Big\|\frakJ_{U_1, 3}^{(1)} - \frac{1}{n}\sum_{i=1}^{n}\xi_i\calM_1(\calX_i)\left(U_2 \otimes U_3 + Q_2 \otimes U_3 + U_2 \otimes Q_3\right)G_1^\top(G_1G_1^\top)^{-1}R_1^{(1)\top}\Big\|\notag\\ 
		\leq& C_2\left(\kappa_0^2\frac{\sqrt{pr^2+p\log(p)}}{n\lambda_{\submin}^2} + \kappa_0^3\frac{p^{3/2}r^{3/2}}{n^{3/2}\lambda_{\submin}^2} + \kappa_0^3\frac{p^{3/2}r}{n^{3/2}\lambda_{\submin}^3}\right).\label{ineq123}
\end{align}
Let  $\bar\frakJ_{U_1, 3}^{(1)} = \frac{1}{n}\sum_{i=1}^{n}\xi_i\calM_1(\calX_i)\left(U_2 \otimes U_3 + Q_2 \otimes U_3 + U_2 \otimes Q_3\right)G_1^\top(G_1G_1^\top)^{-1}$. By Lemma \ref{lm:sub_Gaussian}, with probability at least $1 - e^{-C_1p}$,
\begin{equation*}
		\left\|\bar\frakJ_{U_1, 3}^{(1)}\right\| \leq C_2\sqrt{\frac{p}{n}}\lambda_{\submin}^{-1} + C_2\sqrt{\frac{pr}{n}}\sqrt{\frac{p}{n}}\lambda_{\submin}^{-2} \leq C_2\sqrt{\frac{p}{n}}\lambda_{\submin}^{-1}.
\end{equation*}
Therefore, with probability at least $1 - C_1p^{-3} - C_1e^{-c_1p}$, 
\begin{align}
		&\left|\tr\big(\frakJ_{U_1, 3}^{(1)\top} \calP_{U_1}^\perp\frakJ_{U_1, 3}^{(1)}\big) - \tr\big(\bar\frakJ_{U_1, 3}^{(1)\top} \calP_{U_1}^\perp\bar\frakJ_{U_1, 3}^{(1)}\big)\right|\notag\\
		\leq& 2\Big|\tr\Big(\big(\frakJ_{U_1, 3}^{(1)} - \bar\frakJ_{U_1, 3}^{(1)}R_1^{(1)\top}\big)^\top \calP_{U_1}^\perp\bar\frakJ_{U_1, 3}^{(1)}R_1^{(1)\top}\Big)\Big| + \Big|\tr\Big(\big(\frakJ_{U_1, 3}^{(1)} - \bar\frakJ_{U_1, 3}^{(1)}R_1^{(1)\top}\big)^\top \calP_{U_1}^\perp\big(\frakJ_{U_1, 3}^{(1)} - \bar\frakJ_{U_1, 3}^{(1)}R_1^{(1)\top}\big)\Big)\Big|\notag\\
		\leq& r\big\|\frakJ_{U_1, 3}^{(1)} - \bar\frakJ_{U_1, 3}^{(1)}R_1^{(1)\top}\big\|\big\|\bar\frakJ_{U_1, 3}^{(1)}\big\| + r\big\|\frakJ_{U_1, 3}^{(1)} - \bar\frakJ_{U_1, 3}^{(1)}R_1^{(1)\top}\big\|^2\notag\\
		\leq& C_2\left(\kappa_0^2\frac{pr^2+pr\sqrt{\log(p)}}{n^{3/2}\lambda_{\submin}^3} + \kappa_0^3\frac{p^{2}r^{5/2}}{n^{2}\lambda_{\submin}^3} + \kappa_0^3\frac{p^{2}r^2}{n^{2}\lambda_{\submin}^4}\right).\label{ineq124}
\end{align}
Let
$$
\widetilde Z_j = U_{1\perp}^\top\calM_1(\calX_j)(U_2 \otimes U_3)V_{G_1} \in \RR^{(p_1 - r_1) \times r_1}.
$$
Then $\widetilde Z_j \stackrel{i.i.d.}{\sim} N(0, 1)$. With probability at least $1 - e^{-C_1p}$, 
\begin{align}
		&\bigg|\tr\left(\bar\frakJ_{U_1, 3}^{(1)\top} \calP_{U_1}^\perp\bar\frakJ_{U_1, 3}^{(1)}\right) - \frac{1}{n^2}\Big\|\sum_{i=1}^{n}\xi_i\widetilde Z_i\Lambda_1^{-1}\Big\|_{\rm F}^2\bigg|\notag\\
		\leq&2\left|\left\langle\calP_{U_1}^\perp\frac{1}{n}\sum_{j=1}^{n}\xi_j\calM_1(\calX_j)(U_2 \otimes U_3)V_{G_1}\Lambda_1^{-2}, \frac{1}{n}\sum_{i=1}^{n}\xi_i\calM_1(\calX_i)\left(Q_2 \otimes U_3\right)V_{G_1}\right\rangle\right|\notag\\
		&+2\left|\left\langle\calP_{U_1}^\perp\frac{1}{n}\sum_{j=1}^{n}\xi_j\calM_1(\calX_j)(U_2 \otimes U_3)V_{G_1}\Lambda_1^{-2}, \frac{1}{n}\sum_{i=1}^{n}\xi_i\calM_1(\calX_i)\left(U_2 \otimes Q_3\right)V_{G_1}\right\rangle\right|\notag\\
		&+C_2r\left(\sqrt{\frac{pr}{n}}\sqrt{\frac{p}{n}}\lambda_{\submin}^{-2}\right)^2.\label{ineq126}
\end{align}
Note that 
\begin{align*}
		&\left|\left\langle\calP_{U_1}^\perp\frac{1}{n}\sum_{j=1}^{n}\xi_j\calM_1(\calX_j)(U_2 \otimes U_3)V_{G_1}\Lambda_1^{-2}, \frac{1}{n}\sum_{i=1}^{n}\xi_i\calM_1(\calX_i)\left(Q_2 \otimes U_3\right)V_{G_1}\right\rangle\right|\\
		=&\bigg|\bigg\langle\left(U_{1\perp}^\top\frac{1}{n}\sum_{k=1}^{n}\xi_k\calM_1(\calX_k)(U_{2\perp} \otimes U_3)\right)^\top \left(U_{1\perp}^\top\frac{1}{n}\sum_{j=1}^{n}\xi_j\calM_1(\calX_j)(U_2 \otimes U_3)V_{G_1}\right)\Lambda_1^{-2},\\ 
		&\left(\left(\frac{1}{n}\sum_{i=1}^{n}\xi_iU_{2\perp}^\top\calM_2(\calX_i)(U_1 \otimes U_3)G_2^\top\left(G_2G_2^{\top}\right)^{-1}\right) \otimes I_{r_3}\right)V_{G_1}\bigg\rangle\bigg|\\
		=&\left|\left\langle\left(\frac{1}{n}\sum_{j=1}^{n}\xi_jW_jV_{G_1}\Lambda_1^{-2}\right)\left[\left(\left(\frac{1}{n}\sum_{i=1}^{n}\xi_i\widetilde{W}_iV_{G_2}\Lambda_2^{-1}\right) \otimes I_{r_3}\right)V_{G_1}\right]^\top, \frac{1}{n}\sum_{k=1}^{n}\xi_k\bar{W}_k\right\rangle\right|.
\end{align*}
Here, $\bar{W}_i = U_{1\perp}^\top\calM_1(\calX_i)(U_{2\perp} \otimes U_3) \in \RR^{(p_1 - r_1) \times ((p_2 - r_2)r_3)}, W_i = U_{1\perp}^\top\calM_1(\calX_i)(U_2 \otimes U_3) \in \RR^{(p_1 - r_1) \times r_2r_3}$ and $\widetilde{W}_i = U_{2\perp}^\top\calM_2(\calX_i)(U_1 \otimes U_3) \in \RR^{(p_2 - r_2) \times r_1r_3}$. Since $\calX_i \times_1 [U_1\ U_1^\top] \times_2 [U_2\ U_2^\top] \times_3 [U_3\ U_3^\top] \stackrel{i.i.d.}{\sim} N(0, 1)$, we know that $\bar{W}_i \stackrel{i.i.d.}{\sim} N(0, 1), W_i \stackrel{i.i.d.}{\sim} N(0, 1), \widetilde{W}_i\stackrel{i.i.d.}{\sim} N(0, 1)$, and $\bar{W}_i, W_i$ and $\widetilde{W}_i$ are independent. Therefore,
\begin{align*}
		&\left\langle\calP_{U_1}^\perp\frac{1}{n}\sum_{j=1}^{n}\xi_j\calM_1(\calX_j)(U_2 \otimes U_3)V_{G_1}\Lambda_1^{-2}, \frac{1}{n}\sum_{i=1}^{n}\xi_i\calM_1(\calX_i)\left(Q_2 \otimes U_3\right)V_{G_1}\right\rangle\bigg|\left\{W_k, \widetilde{W}_k, \xi_k\right\}_{k=1}^n\\
		\sim& N\bigg(0, \bigg\|\bigg(\frac{1}{n}\sum_{j=1}^{n}\xi_jW_jV_{G_1}\Lambda_1^{-2}\bigg)\bigg[\Big(\big(\frac{1}{n}\sum_{i=1}^{n}\xi_i\widetilde{W}_iV_{G_2}\Lambda_2^{-1}\big) \otimes I_{r_3}\Big)V_{G_1}\bigg]^\top\bigg\|_{\rm F}^2\frac{\|\vec{\xi}\|_2^2}{n^2}\bigg).
\end{align*}
Note that with probability at least $1 - e^{-C_1p}$,
\begin{align*}
		& \bigg\|\bigg(\frac{1}{n}\sum_{j=1}^{n}\xi_jW_jV_{G_1}\Lambda_1^{-2}\bigg)\bigg[\Big(\big(\frac{1}{n}\sum_{i=1}^{n}\xi_i\widetilde{W}_iV_{G_2}\Lambda_2^{-1}\big) \otimes I_{r_3}\Big)V_{G_1}\bigg]^\top\bigg\|_{\rm F}\\
		\leq& \sqrt{r}\Big\|\frac{1}{n}\sum_{j=1}^{n}\xi_jW_jV_{G_1}\Lambda_1^{-2}\Big\|\Big\|\frac{1}{n}\sum_{i=1}^{n}\xi_i\widetilde{W}_iV_{G_2}\Lambda_2^{-1}\Big\|
		\leq C_2\sqrt{r}\left(\sqrt{\frac{p}{n}}\lambda_{\submin}^{-2}\right)\left(\sqrt{\frac{p}{n}}\lambda_{\submin}^{-1}\right)\\
		\leq& C_2\frac{pr^{1/2}}{n\lambda_{\submin}^3}, 
\end{align*}
and as a result
\begin{align*}
		&\PP\left(\bigg|\bigg<\calP_{U_1}^\perp\frac{1}{n}\sum_{j=1}^{n}\xi_j\calM_1(\calX_j)(U_2 \otimes U_3)V_{G_1}\Lambda_1^{-2}, \frac{1}{n}\sum_{i=1}^{n}\xi_i\calM_1(\calX_i)\left(Q_2 \otimes U_3\right)V_{G_1}\bigg>\bigg| \geq C_2\frac{\sqrt{r}p\sqrt{\log(p)}}{n^{3/2}\lambda_{\submin}^3}\right) \leq \frac{C_1}{p^3}.
\end{align*}
Similarly, 
\begin{align*}
		&\PP\left(\bigg|\bigg<\calP_{U_1}^\perp\frac{1}{n}\sum_{j=1}^{n}\xi_j\calM_1(\calX_j)(U_2 \otimes U_3)V_{G_1}\Lambda_1^{-2}, \frac{1}{n}\sum_{i=1}^{n}\xi_i\calM_1(\calX_i)\left(U_2 \otimes Q_3\right)V_{G_1}\bigg>\bigg| \geq C_2\frac{\sqrt{r}p\sqrt{\log(p)}}{n^{3/2}\lambda_{\submin}^3}\right)
		\leq \frac{C_1}{p^3}.
\end{align*}
By \eqref{ineq124}, \eqref{ineq126} and the two inequalities above, with probability at least $1 - C_1p^{-3} - C_1e^{-c_1p}$,
\begin{equation*}
		\bigg|\tr\big(\frakJ_{U_1, 3}^{(1)\top} \calP_{U_1}^\perp\frakJ_{U_1, 3}^{(1)}\big) - \frac{1}{n^2}\Big\|\sum_{i=1}^{n}\xi_i\widetilde Z_i\Lambda_1^{-1}\Big\|_{\rm F}^2\bigg| \leq C_2\left(\kappa_0^2\frac{pr^2+pr\sqrt{\log(p)}}{n^{3/2}\lambda_{\submin}^3} + \kappa_0^3\frac{p^{2}r^{5/2}}{n^{2}\lambda_{\submin}^3} + \kappa_0^3\frac{p^{2}r^2}{n^{2}\lambda_{\submin}^4}\right).
\end{equation*}
By \eqref{ineq66}, \eqref{ineq68}, \eqref{ineq118} and the previous inequality, with probability at least $1 - C_1p^{-3} - C_1e^{-c_1p}$, we have proved
\begin{align}\label{ineq130}
			&\left|\left\langle \begin{pmatrix}
			U_1U_1^\top & 0\\
			0 & I_{r_1}
			\end{pmatrix}, 
			\calS_{U_1^{(1)}, 2}\left(E\right)
			\right\rangle + \frac{1}{n^2}\left\|\sum_{i=1}^{n}\xi_i\widetilde Z_i\Lambda_1^{-1}\right\|_{\rm F}^2\right|\notag\\ 
			\leq& C_2\left(\kappa_0^2\frac{pr^2+pr\sqrt{\log(p)}}{n^{3/2}\lambda_{\submin}^3} + \kappa_0^3\frac{p^{2}r^{5/2}}{n^{2}\lambda_{\submin}^3} + \kappa_0^3\frac{p^{2}r^2}{n^{2}\lambda_{\submin}^4} + \kappa_0^2\frac{p^2r^{3}}{n^2\lambda_{\submin}^2} + \kappa_0\frac{r^3p\log(n)}{n^{3/2}\lambda_{\submin}^2} + \frac{r \log^2(n)}{n^{3/2}\lambda_{\submin}^2}\right).
\end{align}
By \eqref{eq3}, \eqref{ineq129} and \eqref{ineq130}, with probability at least $1 - C_1p^{-3} - C_1e^{-c_1p}$,
\begin{align}\label{ineq125}
&\bigg|\left\|\hat U_1^{(2)}\hat U_1^{(2)\top} - U_1U_1^\top\right\|_{\rm F}^2 - \frac{2}{n^2}\Big\|\sum_{i=1}^{n}\xi_i\widetilde Z_i\Lambda_1^{-1}\Big\|_{\rm F}^2\bigg|\notag\\ 
\leq& C_2\left(\kappa_0^2\frac{pr^2+pr\sqrt{\log(p)}}{n^{3/2}\lambda_{\submin}^3} + \kappa_0^3\frac{p^{2}r^{5/2}}{n^{2}\lambda_{\submin}^3} + \kappa_0^3\frac{p^{2}r^2}{n^{2}\lambda_{\submin}^4} + \kappa_0^2\frac{p^2r^{3}}{n^2\lambda_{\submin}^2} + \kappa_0\frac{r^3p\log(n)}{n^{3/2}\lambda_{\submin}^2} + \frac{r \log^2(n)}{n^{3/2}\lambda_{\submin}^2}\right).
\end{align}
%where we assumed $n\leq e^p$. 

Recall that $\widetilde Z_i = U_{1\perp}^\top\calM_1(\calX_i)(U_2 \otimes U_3)V_{G_1} \in \RR^{(p_1 - r_1) \times r_1}$ so that $\widetilde Z_i \stackrel{i.i.d.}{\sim} N(0,1)$. For fixed $a = (a_1, \dots, a_n) \in \RR^n$, the rows of $(\sum_{i=1}^{n}a_i\widetilde Z_i)\Lambda_1^{-1} \stackrel{i.i.d.}{\sim} N\left(0, \|a\|_2^2\Lambda_1^{-2}\right)$. Therefore, for any $a \in \RR^n$,
\begin{equation*}
	\frac{1}{\|a\|_2^2}\Big\|\Big(\sum_{i=1}^{n}a_i\widetilde Z_i\Big)\Lambda_1^{-1}\Big\|_{\rm F}^2 \stackrel{\rm d.}{=} \frac{1}{n}\Big\|\sum_{i=1}^{n}\widetilde Z_i\Lambda_1^{-1}\Big\|_{\rm F}^2
\end{equation*}
which means that
\begin{equation*}
\frac{1}{\|\vec{\xi}\|_2^2}\Big\|\Big(\sum_{i=1}^{n}\xi_i\widetilde Z_i\Big)\Lambda_1^{-1}\Big\|_{\rm F}^2 \stackrel{\rm d.}{=} \frac{1}{n}\Big\|\sum_{i=1}^{n}\widetilde Z_i\Lambda_1^{-1}\Big\|_{\rm F}^2.
\end{equation*}
For any $1 \leq i \leq p_1 - r_1$,
\begin{equation*}
	\EE\Big\|\Big(\sum_{j=1}^{n}(\widetilde Z_j)_{[i, :]}\Big)\Lambda_1^{-1}\Big\|_2^2 = n\|\Lambda_1^{-1}\|_{\rm F}^2,
\end{equation*}
\begin{equation*}
\Var\bigg(\Big\|\Big(\sum_{j=1}^{n}(\widetilde Z_j)_{[i, :]}\Big)\Lambda_1^{-1}\Big\|_2^2\bigg) = 2n^2\|\Lambda_1^{-2}\|_{\rm F}^2.
\end{equation*}
and
\begin{equation*}
	\EE\bigg\|\Big(\sum_{j=1}^{n}(\widetilde Z_j)_{[i, :]}\Big)\Lambda_1^{-1}\bigg\|_2^6 \leq C_3n^3\sum_{j_1, j_2, j_3 = 1}^{r_1}\frac{1}{\lambda_{j_1}^{(1)2}\lambda_{j_2}^{(1)2}\lambda_{j_3}^{(1)2}} \leq C_3n^3\|\Lambda_1^{-1}\|_{\rm F}^6.
\end{equation*}
By Berry-Esseen theorem, we have
\begin{equation}\label{ineq127}
	\sup_{x \in \RR}\bigg|\PP\bigg(\frac{\frac{2}{n\|\vec{\xi}\|_2^2}\big\|(\sum_{j=1}^{n}\xi_j\widetilde Z_j)\Lambda_1^{-1}\big\|_{\rm F}^2 - 2(p_1-r_1)\|\Lambda_1^{-1}\|_{\rm F}^2/n}{\sqrt{8(p_1-r_1)}\|\Lambda_1^{-2}\|_{\rm F}/n} \leq x\bigg) - \Phi(x)\bigg| \leq C\Big(\frac{\|\Lambda_1^{-1}\|_{\rm F}^4}{\|\Lambda_1^{-2}\|_{\rm F}^2}\Big)^{3/2}\frac{1}{\sqrt{p}}.
\end{equation}
By Lemma \ref{lm:sub_Gaussian}, with probability $1 - e^{-C_1p}$,
\begin{equation*}
	\big\|\sum_{j=1}^{n}\xi_j\widetilde Z_j\big\| \leq C_2\sqrt{np}.
\end{equation*}
Therefore, with probability $1 - e^{-C_1p}$,
\begin{equation*}
	\Big\|(\sum_{j=1}^{n}\xi_j\widetilde Z_j)\Lambda_1^{-1}\Big\|_{\rm F}^2 \leq r_1\Big\|\sum_{j=1}^{n}\xi_j\widetilde Z_j\Big\|^2\|\Lambda_1^{-1}\|^2 \leq C_2rnp\lambda_{\submin}^{-2}.
\end{equation*}
By Bernstein-type inequality (\cite[Proposition 5.16]{vershynin2010introduction}),
\begin{equation}\label{ineq:noise}
	\begin{split}
	\PP\left(|\|\vec{\xi}\|_2^2 - n| \geq C_2\sqrt{n\log(p)}\right) \leq 2\exp\left[-C_1\min\left\{\frac{n\log(p)}{n}, \sqrt{n\log(p)}\right\}\right] \leq p^{-C_1}.
	\end{split}
\end{equation}
Therefore, with probability at least $1 - C_1p^{-3}$,
\begin{equation}\label{ineq128}
	\left|\frac{2}{n\|\vec{\xi}\|_2^2} - \frac{2}{n^2}\right| \leq \frac{2|\|\vec{\xi}\|_2^2 - n|}{n^2\|\vec{\xi}\|_2^2} \leq C_2\frac{\sqrt{n\log(p)}}{n^3}.
\end{equation}
By \eqref{ineq127}, \eqref{ineq128} and the previous inequality and the similar proof in Theorem \ref{thm:na_tsvd}, we have
\begin{align*}
&\sup_{x \in \RR}\bigg|\PP\bigg(\frac{\big\|\hat U_1^{(2)}\hat U_1^{(2)\top} - U_1U_1^\top\big\|_{\rm F}^2 - \frac{2p_1}{n}\|\Lambda_1^{-1}\|_{\rm F}^2}{\frac{\sqrt{8p_1}}{n}\|\Lambda_1^{-2}\|_{\rm F}} \leq x\bigg) - \Phi(x)\bigg|\\ 
\leq& C_2\left(\kappa_0^4\sqrt{\frac{pr^{3}+pr\log(p)}{n\lambda_{\submin}^2}} + \kappa_0^5\frac{p^{3/2}r^{2}}{n\lambda_{\submin}} + \kappa_0^5\frac{p^{3/2}r^{3/2}}{n\lambda_{\submin}^2} + \kappa_0^4\frac{p^{3/2}r^{5/2}}{n} + \kappa_0^3\sqrt{\frac{r^5p\log^2(n)}{n}}\right)\\ &+ C_3\frac{r^{3/2}}{\sqrt{p}} + C_1e^{-c_1p}
\end{align*}
where we use the fact that $C_2\frac{\kappa_0^2\sqrt{r}\log^2(n)}{\sqrt{pn}} \leq C_3\frac{r^{3/2}}{\sqrt{p}}$.

\subsection{Proof of Theorem \ref{thm:regression_adaptive}}
Without loss of generality, we assume $\sigma = 1$. By \eqref{ineq:error_G} and Lemma \ref{lm:sub_Gaussian}, with probability at least $1 - C_1p^{-3} - C_1e^{-c_1p}$,
\begin{align*}
		\max_{1 \leq i \leq r_1}\big|\lambda_i - \hat\lambda_i^{(1)}\big| \leq \big\|\hat{G}_1 - R_1G_1(R_2^\top \otimes R_3^\top)\big\| \leq C_2\left(\kappa_0\frac{pr}{n} + \kappa_0\frac{p\sqrt{r}}{n\lambda_{\submin}} + \sqrt{\frac{r^2 + \log(p)}{n}}\right).
\end{align*}
Therefore, with probability at least $1 - C_1p^{-3} - C_1e^{-c_1p}$,
\begin{equation*}
		\Big|\left\|\Lambda_1^{-1}\right\|_{\rm F}^2 - \big\|\big(\hat\Lambda_1^{(1)}\big)^{-1}\big\|_{\rm F}^2\Big| \leq r\max_{1 \leq i \leq r_1}\frac{\big|\lambda_i^2 - \hat\lambda_i^2\big|}{\lambda_i^2\hat\lambda_i^2} \leq C_2r\left(\kappa_0\frac{pr}{n} + \kappa_0\frac{p\sqrt{r}}{n\lambda_{\submin}} + \sqrt{\frac{r^2 + \log(p)}{n}}\right)\lambda_{\submin}^{-3},
\end{equation*}
	and
\begin{align*}
	\frac{\Big|\left\|\Lambda_1^{-1}\right\|_{\rm F}^2 - \big\|\big(\hat\Lambda_1^{(1)}\big)^{-1}\big\|_{\rm F}^2\Big|}{\left\|\Lambda_1^{-2}\right\|_{\rm F}} &\leq \frac{\max_{1 \leq i \leq r_1}\frac{\left|\lambda_i^2 - \hat\lambda_i^2\right|}{\lambda_i^2\hat\lambda_i^2}}{\kappa_0^{-2}\lambda_{\submin}^{-2}} \leq C_2\kappa_0^2\left(\kappa_0\frac{pr}{n} + \kappa_0\frac{p\sqrt{r}}{n\lambda_{\submin}} + \sqrt{\frac{r^2 + \log(p)}{n}}\right)\lambda_{\submin}^{-1}.
\end{align*}
The rest of the proof is essentially the same as the proof of Theorem \ref{thm:adaptive_tsvd}.

\subsection{Proof of Theorem~\ref{thm:T-orth}}
Without loss of generality, we assume $\sigma=1$, $\pi(j) = j$ and only prove the normal approximation for $\langle \hat u_1, u_1\rangle^2$. We denote $\tilde U_1=(u_2,\cdots,u_r)\in \OO_{p_1,r-1}$, $\tilde V_1=(v_2,\cdots,v_r)\in\OO_{p_2,r-1}$ and $\tilde W_1=(w_2,\cdots,w_r)\in\OO_{p_3,r-1}$. Denote the $(r-1)\times (r-1)\times (r-1)$ diagonal tensor $\tilde \bLambda={\rm diag}(\lambda_2,\cdots,\lambda_r)$, and $\tilde\Lambda_1=\calM_1(\tilde{\bLambda}), \tilde\Lambda_2=\calM_2(\tilde{\bLambda}), \tilde\Lambda_3=\calM_3(\tilde{\bLambda})$. 

By definition, $\hat u_1$ is the left singular vector of $\calA\times_2 \hat v_1^{(1)\top}\times_3 \hat w_1^{(1)\top}$, for which we write
\begin{align*}
\calA\times_2 \hat v_1^{(1)\top}\times_3 \hat w_1^{(1)\top}=&\lambda_1\langle \hat v_1^{(1)}, v_1\rangle \langle \hat w_1^{(1)}, w_1\rangle u_1+\tilde U_1\tilde\Lambda_1(\tilde V_1\otimes \tilde W_1)^{\top}(\hat v_1^{(1)}\otimes\hat w_1^{(1)})+\underbrace{\calZ\times_2 \hat v_1^{(1)\top}\times_3 \hat w_1^{(1)\top}}_{\hat z_1^{(1)}}\\
=:&\tilde\lambda_1 u_1+\hat E_1^{(1)}
\end{align*}
where we define $\tilde\lambda_1=\lambda_1\langle \hat v_1^{(1)}, v_1\rangle \langle \hat w_1^{(1)}, w_1\rangle$ and $\hat E_1^{(1)}=\tilde U_1\tilde\Lambda_1(\tilde V_1\otimes \tilde W_1)^{\top}(\hat v_1^{(1)}\otimes\hat w_1^{(1)})+\hat z_1^{(1)}$. 

Similarly to the proof of Theorem~\ref{thm:na_tsvd}, the following bounds hold with probability at least $1-C_1e^{-c_1p}$,
$$
\max\big\{\|\tilde V_1^{\top}\hat v_1^{(1)}\|_2, \|\tilde W_1^{\top}\hat w_1^{(1)}\|_2\big\}\leq C_2\sqrt{p}/\lambda_{\submin}\quad {\rm and}\quad \|\hat z_1^{(1)}\|_2\leq C_3(\sqrt{p}+p\lambda_{\submin}^{-1}).
$$
As a result, with the same probability, $\|\hat E_1^{(1)}\|\leq C_3(\sqrt{p}+\kappa_0p\lambda_{\submin}^{-1})$. Then, we write
$$
\left( \begin{array}{cc}0& \tilde\lambda_1 u_1+\hat E_1^{(1)}\\\tilde\lambda_1 u_1^{\top}+\hat E_1^{(1)\top} &0\end{array}\right)=\tilde\lambda_1\left( \begin{array}{cc}0&  u_1\\ u_1^{\top} &0\end{array}\right)+\left( \begin{array}{cc}0&  \hat E_1^{(1)}\\ \hat E_1^{(1)\top} &0\end{array}\right).
$$
We now apply Lemma~\ref{lem:spectral} and represent $\langle u_1u_1^{\top}, \hat u_1 \hat u_1^{\top}-u_1u_1^{\top}\rangle$. Similarly to the proof of Theorem~\ref{thm:na_tsvd}, the $1$st-order term does not matter, and the $4$-th and higher-order terms can be simply bounded. Therefore, we obtain,
\begin{align*}
    \langle u_1u_1^{\top}, \hat u_1 \hat u_1^{\top}-u_1u_1^{\top}\rangle=-\frac{1}{\tilde{\lambda}_1^2}\tr\big(\hat E_1^{(1)\top}U_{1\perp}U_{1\perp}^{\top}\hat E_1^{(1)}\big)+\frac{2}{\tilde{\lambda}_1^3}\tr\big(\hat E_1^{(1)\top}u_1\hat E_1^{(1)\top}U_{1\perp}U_{1\perp}^{\top}\hat E_1^{(1)}\big)+\tilde{R}_1^{(1)}
\end{align*}
where $(u_1,U_{1\perp})\in\OO_{p_1}$, and $\|\tilde{R}_1^{(1)}\|\leq C_3\kappa_0^4p^2/\lambda_{\submin}^4$ with probability at least $1-C_1e^{-c_1p}$. 

Similarly to the proof of Theorem~\ref{thm:na_tsvd}, we have $|\hat E_1^{(1)\top}u_1|\leq C_1(p\lambda_{\submin}^{-1}+\sqrt{\log p})$ which holds with probability at least $1-p^{-3}$. Therefore, with probability at least $1-2p^{-3}$, $|\tr(\hat E_1^{(1)\top}u_1\hat E_1^{(1)\top}U_{1\perp}U_{1\perp}^{\top}\hat E_1^{(1)})|\leq C_2\kappa_0^2p(p\lambda_{\submin}^{-1}+\sqrt{\log p})$. Therefore, we conclude with probability at least $1-3p^{-3}$ that 
\begin{align*}
    \bigg| \langle u_1u_1^{\top}, \hat u_1 \hat u_1^{\top}-u_1u_1^{\top}\rangle+\frac{1}{\tilde{\lambda}_1^2}\tr\big(\hat E_1^{(1)\top}U_{1\perp}U_{1\perp}^{\top}\hat E_1^{(1)}\big) \bigg|\leq C_2\Big(\frac{\kappa_0^2p\sqrt{\log p}}{\lambda_{\submin}^3}+\frac{\kappa_0^4p^2}{\lambda_{\submin}^4}\Big).
\end{align*}

Similarly to the proof of Theorem~\ref{thm:na_tsvd}, we have $|\tilde{\lambda}_1^2-\lambda_1^2|\leq C_2\kappa_0^2p$ with probability at least $1-C_1e^{-c_1p}$. 

Therefore, with probability at least $1-4p^{-3}$ that 
\begin{align*}
\bigg| \langle u_1u_1^{\top}, \hat u_1 \hat u_1^{\top}-u_1u_1^{\top}\rangle+\frac{1}{\lambda_1^2}\tr\big(\hat E_1^{(1)\top}U_{1\perp}U_{1\perp}^{\top}\hat E_1^{(1)}\big) \bigg|\leq C_2\Big(\frac{\kappa_0^2p\sqrt{\log p}}{\lambda_{\submin}^3}+\frac{\kappa_0^4p^2}{\lambda_{\submin}^4}\Big).
\end{align*}
It then suffices to prove the normal approximation of $\tr(\hat E_1^{(1)\top}U_{1\perp}U_{1\perp}^{\top}\hat E_1^{(1)})$. Recall that $\hat E_1^{(1)}=\tilde U_1\tilde\Lambda_1(\tilde V_1\otimes \tilde W_1)^{\top}(\hat v_1^{(1)}\otimes\hat w_1^{(1)})+\hat z_1^{(1)}$. Note that 
$$
\|\tilde U_1\tilde\Lambda_1(\tilde V_1\otimes \tilde W_1)^{\top}(\hat v_1^{(1)}\otimes\hat w_1^{(1)})\|_2\leq C_2\frac{\kappa_0p}{\lambda_{\submin}}
$$
implying that, with probability at least $1-C_1e^{-c_1p}$, 
$$
\Big| \tr\Big(\big(\tilde U_1\tilde\Lambda_1(\tilde V_1\otimes \tilde W_1)^{\top}(\hat v_1^{(1)}\otimes\hat w_1^{(1)})\big)^{\top}U_{1\perp}U_{1\perp}^{\top}\big(\tilde U_1\tilde\Lambda_1(\tilde V_1\otimes \tilde W_1)^{\top}(\hat v_1^{(1)}\otimes\hat w_1^{(1)})\big)\Big)\Big|\leq C_3\frac{\kappa_0^2p^2}{\lambda_{\submin}^2}.
$$
Now, we consider the cross term (recall $z_1=\calZ\times_2 v_1^{\top}\times_3 w_1^{\top}$) and conclude with probability at least $1-p^{-3}$,
\begin{align*}
    \Big|\tr\Big(&\hat z_1^{(1)\top}U_{1\perp}U_{1\perp}^{\top}\big(\tilde U_1\tilde\Lambda_1(\tilde V_1\otimes \tilde W_1)^{\top}(\hat v_1^{(1)}\otimes\hat w_1^{(1)})\big)\Big)\Big|=\Big|\tr\Big(\hat z_1^{(1)\top}\big(\tilde U_1\tilde\Lambda_1(\tilde V_1\otimes \tilde W_1)^{\top}(\hat v_1^{(1)}\otimes\hat w_1^{(1)})\big)\Big)\Big|\\
    \leq&\Big|\tr\Big(z_1^{\top}\big(\tilde U_1\tilde\Lambda_1(\tilde V_1\otimes \tilde W_1)^{\top}(\hat v_1^{(1)}\otimes\hat w_1^{(1)})\big)\Big)\Big|+\Big|\tr\Big((\hat z_1^{(1)}-z_1)^{\top}\big(\tilde U_1\tilde\Lambda_1(\tilde V_1\otimes \tilde W_1)^{\top}(\hat v_1^{(1)}\otimes\hat w_1^{(1)})\big)\Big)\Big|\\
\leq &C_2\sqrt{r\log p}\cdot \kappa_0p/\lambda_{\submin}+C_3\kappa_0p^2/\lambda_{\submin}^2
\end{align*}
where the first term is due to $\|z_1^{\top}\tilde U_1\|=O(\sqrt{r\log p})$ with probability at least $1-p^{-3}/2$, and the second term is similar as the proof of Lemma~2. 

To this end, we obtain with probability at least $1-5p^{-3}$ that
\begin{align}\label{eq:cor3ineq1}
\bigg| \langle u_1u_1^{\top}, \hat u_1 \hat u_1^{\top}-u_1u_1^{\top}\rangle+\frac{1}{\lambda_1^2}\tr\big(\hat z_1^{(1)\top}U_{1\perp}U_{1\perp}^{\top}\hat z_1^{(1)}\big) \bigg|\leq C_2\Big(\frac{\kappa_0^2p\sqrt{\log p}+\kappa_0p\sqrt{r\log p}}{\lambda_{\submin}^3}+\frac{\kappa_0^4p^2}{\lambda_{\submin}^4}\Big).
\end{align}
Now, we investigate the main term $\tr\big(\hat z_1^{(1)\top}U_{1\perp}U_{1\perp}^{\top}\hat z_1^{(1)}\big)$ for which we write
\begin{align*}
&\tr\big(\hat z_1^{(1)\top}U_{1\perp}U_{1\perp}^{\top}\hat z_1^{(1)}\big)=\tr\big(U_{1\perp}U_{1\perp}^{\top}\hat z_1^{(1)}\hat z_1^{(1)\top}\big)=\tr\big(U_{1\perp}U_{1\perp}^{\top}Z_1\big((\hat v_1^{(1)}\hat v_1^{(1)\top})\otimes (\hat w_1^{(1)}\hat w_1^{(1)\top})\big)Z_1^{\top}\big)\\
=&\tr\big(U_{1\perp}U_{1\perp}^{\top}Z_1\big((v_1 v_1^{\top})\otimes (w_1 w_1^{\top})\big)Z_1^{\top}\big)+\tr\big(U_{1\perp}U_{1\perp}^{\top}Z_1\big((\hat v_1^{(1)}\hat v_1^{(1)\top}-v_1 v_1^{\top})\otimes (w_1 w_1^{\top})\big)Z_1^{\top}\big)\\
&+\tr\big(U_{1\perp}U_{1\perp}^{\top}Z_1\big((v_1 v_1^{\top})\otimes (\hat w_1^{(1)}\hat w_1^{(1)\top}-w_1 w_1^{\top})\big)Z_1^{\top}\big)\\&+\tr\big(U_{1\perp}U_{1\perp}^{\top}Z_1\big((\hat v_1^{(1)}\hat v_1^{(1)\top}-v_1 v_1^{\top})\otimes (\hat w_1^{(1)}\hat w_1^{(1)\top}-w_1 w_1^{\top})\big)Z_1^{\top}\big)
\end{align*}
where we denote $Z_1=\calM_1(\calZ)$, and the last term can be simply bounded by $C_4p^2/\lambda_{\submin}^2$ with probability at least $1-C_1e^{-c_1p}$.

The idea of bounding the term $\tr\big(U_{1\perp}U_{1\perp}^{\top},Z_1\big((\hat v_1^{(1)}\hat v_1^{(1)\top}-v_1 v_1^{\top})\otimes (w_1 w_1^{\top})\big)Z_1^{\top}\big)$ is the same as the Step 4 in the proof of Theorem~\ref{thm:na_tsvd}. Indeed, we shall recall that $\hat v_1^{(1)}$ is the singular vector of 
\begin{align*}
\hat v_1^{(0.5)}=\calA\times_1 \hat u_1^{(0)}\times_3 \hat w_1^{(0)}=&\lambda_1\langle u_1, \hat u_1^{(0)}\rangle \langle w_1, \hat w_1^{(0)}\rangle v_1+\tilde V_1\tilde\Lambda_2(\tilde U_1\otimes \tilde W_1)^{\top}(\hat u_1^{(0)}\otimes \hat w_1^{(0)})+\calZ\times_1 \hat u_1^{(0)}\times_3\hat w_1^{(0)}\\
=&\tilde{\lambda}_1^{(0)}v_1+\hat E_2^{(0)}.
\end{align*}
Similarly to the proof of Theorem~\ref{thm:na_tsvd}, it suffices to consider the $1$st-order term in $\hat v_1^{(1)}\hat v_1^{(1)\top}-v_1v_1^{\top}$. It is then easy to show that $\big|\tr\big(U_{1\perp}U_{1\perp}^{\top},Z_1\big((\hat v_1^{(1)}\hat v_1^{(1)\top}-v_1 v_1^{\top})\otimes (w_1 w_1^{\top})\big)Z_1^{\top}\big)\big|\leq C_4\kappa_0p\sqrt{\log p}/\lambda_{\submin}+C_5\kappa_0^4p^2/\lambda_{\submin}^2$ with probability at least $1-p^{-3}$. 
Together with (\ref{eq:cor3ineq1}), we conclude with probability at least $1-6p^{-3}$ that 
$$
\bigg| \langle u_1u_1^{\top}, \hat u_1 \hat u_1^{\top}-u_1u_1^{\top}\rangle+\frac{1}{\lambda_1^2}\tr\big(z_1^{\top}U_{1\perp}U_{1\perp}^{\top} z_1\big) \bigg|\leq C_2\Big(\frac{\kappa_0^2p\sqrt{\log p}+k_0p\sqrt{r\log p}}{\lambda_{\submin}^3}+\frac{\kappa_0^4p^2}{\lambda_{\submin}^4}\Big).
$$
The rest of the proof is identical to the final step in the proof of Theorem~\ref{thm:na_tsvd}.

\subsection{Proof of Theorem~\ref{thm:entry_inference}}
	Without loss of generality, we assume $\sigma = 1$. For random variables (or vectors) $A$ and $B$, we use $A \stackrel{\rm d.}{\to} (\text{resp.} \stackrel{\rm p.}{\to}, \stackrel{\rm a.s.}{\to}) B$ as a shorthand for $A \stackrel{\rm d.}{\to} (\text{resp.} \stackrel{\rm p.}{\to}, \stackrel{\rm a.s.}{\to}) B$ as $p \to \infty$. 
	By \cite[Theorem 1]{zhang2018tensor}, with probability at least $1 - C_1e^{-c_1p}$, 
	$$\|\hat u - u\|_2, \|\hat v - v\|_2, \|\hat w - w\|_2 \leq C_2\frac{\sqrt{p}}{\lambda}.$$ 
	Let $U_{\perp} \in \OO_{p_1, p_1-1}, V_{\perp} \in \OO_{p_2, p_2-1}, W_{\perp} \in \OO_{p_3, p_3-1}$ be orthogonal complements of $u, v$ and $w$, respectively. For any $O_i \in \OO_{p_i-1}, i \in [3]$, let $$\tilde{O}_1 = uu^\top + U_{\perp}O_1U_{\perp}^\top \in \OO_{p_1}, \quad \tilde{O}_2 = vv^\top + V_{\perp}O_2V_{\perp}^\top \in \OO_{p_2}, \quad \tilde{O}_3 = ww^\top + W_{\perp}O_3W_{\perp}^\top \in \OO_{p_3}.$$ 
	Let $\tilde{\calA} = \calA \times_1 \tilde{O}_1^\top \times_2 \tilde{O}_2^\top \times_3 \tilde{O}_3^\top$. Notice that $\tilde{O}_1u = u, \tilde{O}_2v = v, \tilde{O}_3w = w$, we have
	\begin{align*}
	\tilde{\calA} =& \calA \times_1 \tilde{O}_1^\top \times_2 \tilde{O}_2^\top \times_3 \tilde{O}_3^\top = \calT \times_1 \tilde{O}_1^\top \times_2 \tilde{O}_2^\top \times_3 \tilde{O}_3^\top + \calZ \times_1 \tilde{O}_1^\top \times_2 \tilde{O}_2^\top \times_3 \tilde{O}_3^\top\\
	=& \calT + \calZ \times_1 \tilde{O}_1^\top \times_2 \tilde{O}_2^\top \times_3 \tilde{O}_3^\top.
	\end{align*}
	Since $O_i \in \OO_{p_{i-1}}$, the entries of $\calZ \times_1 \tilde{O}_1^\top \times_2 \tilde{O}_2^\top \times_3 \tilde{O}_3^\top \stackrel{i.i.d.}{\sim} N(0,1)$. Consequently, 
	\begin{equation}\label{eq5}
	\tilde{\calA} \stackrel{\rm d.}{=} \calA.
	\end{equation}
	Let $\tilde{u}, \tilde{v}, \tilde{w}$ be the outputs of Algorithm \ref{algo:hooi} after $t_{\submax}$ iterations. Then we have 
	$$\tilde{u} = \tilde{O}_1^\top\hat{u}, \quad \tilde{v} = \tilde{O}_2^\top\hat{v}, \quad \tilde{w} = \tilde{O}_3^\top\hat{w}.$$ In addition, we have
	$$\langle u, \hat{u}\rangle = \langle \tilde{O}_1u, \hat{u}\rangle = \langle u, \tilde{O}_1^\top\hat{u}\rangle = \langle u, \tilde{u}\rangle.$$ Similarly, we have $\langle v, \hat{v}\rangle = \langle v, \tilde{v}\rangle$ and $\langle w, \hat{w}\rangle = \langle w, \tilde{w}\rangle$.
	By \eqref{eq5}, for any $O_i \in \OO_{p_i-1}, i \in [3]$,
	\begin{align*}
	& (\hat{u}^\top U_{\perp}, \hat{v}^\top V_{\perp}, \hat{w}^\top W_{\perp})\bigg|(\langle u, \hat{u}\rangle, \langle v, \hat{v}\rangle, \langle w, \hat{w}\rangle) \\
	\stackrel{\rm d.}{=}& (\tilde{u}^\top U_{\perp}, \tilde{v}^\top V_{\perp}, \tilde{w}^\top W_{\perp})\bigg|(\langle u, \tilde{u}\rangle, \langle v, \tilde{v}\rangle, \langle w, \tilde{w}\rangle)\\ =& (\hat{u}^\top\tilde O_1U_{\perp}, \hat{v}^\top\tilde O_2V_{\perp}, \hat{w}^\top\tilde O_3W_{\perp})\bigg|(\langle u, \hat{u}\rangle, \langle v, \hat{v}\rangle, \langle w, \hat{w}\rangle)\\ =& (\hat{u}^\top U_{\perp}O_1, \hat{v}^\top V_{\perp}O_2, \hat{w}^\top W_{\perp}O_3)\bigg|(\langle u, \hat{u}\rangle, \langle v, \hat{v}\rangle, \langle w, \hat{w}\rangle).
	\end{align*}
	Therefore, for any $O_i \in \OO_{p_i-1}, i \in [3]$,
	\begin{align*}
	&\bigg(\frac{\hat{u}^\top U_{\perp}}{\|U_{\perp}^\top\hat{u}\|_2}, \frac{\hat{v}^\top V_{\perp}}{\|V_{\perp}^\top\hat{v}\|_2}, \frac{\hat{w}^\top W_{\perp}}{\|W_{\perp}^\top\hat{w}\|_2}\bigg)\bigg|(\langle u, \hat{u}\rangle, \langle v, \hat{v}\rangle, \langle w, \hat{w}\rangle)\\ \stackrel{\rm d.}{=}& \bigg(\frac{\hat{u}^\top U_{\perp}O_1}{\|(U_{\perp}O_1)^\top\hat{u}\|_2}, \frac{\hat{v}^\top V_{\perp}O_2}{\|(V_{\perp}O_2)^\top\hat{v}\|_2}, \frac{\hat{w}^\top W_{\perp}O_3}{\|(W_{\perp}O_3)^\top\hat{w}\|_2}\bigg)\bigg|(\langle u, \hat{u}\rangle, \langle v, \hat{v}\rangle, \langle w, \hat{w}\rangle)\\
	=& \bigg(\frac{\hat{u}^\top U_{\perp}}{\|U_{\perp}^\top\hat{u}\|_2}O_1, \frac{\hat{v}^\top V_{\perp}}{\|V_{\perp}^\top\hat{v}\|_2}O_2, \frac{\hat{w}^\top W_{\perp}}{\|W_{\perp}^\top\hat{w}\|_2}O_3\bigg)\bigg|(\langle u, \hat{u}\rangle, \langle v, \hat{v}\rangle, \langle w, \hat{w}\rangle).
	\end{align*}
	Set $O_1 = I_{p_1-1}$, for any $v_1, v_2 \in \SS^{p_2-2} = \{x \in \RR^{p_2-1}: \|x\|_2 = 1\}$ and $w_1, w_2 \in \SS^{p_3-2}$, we know that there exists $O_2 \in \OO_{p_2-1}$ and $O_3 \in \OO_{p_3-1}$ such that $v_2 = O_2v_1, w_2 = O_3v_2$. Then for any Borel set $A \subseteq \SS^{p_1-2}$,
	\begin{align*}
	&\PP\left(\frac{U_{\perp}^\top\hat{u}}{\|U_{\perp}^\top\hat{u}\|_2} \in A\bigg|\frac{V_{\perp}^\top\hat{v}}{\|V_{\perp}^\top\hat{v}\|_2} = v_1, \frac{W_{\perp}^\top\hat{w}}{\|W_{\perp}^\top\hat{w}\|_2} = w_1\right)\\ =& \PP\left(\frac{U_{\perp}^\top\hat{u}}{\|U_{\perp}^\top\hat{u}\|_2} \in A\bigg|O_2^\top\frac{V_{\perp}^\top\hat{v}}{\|V_{\perp}^\top\hat{v}\|_2} = v_1, O_3^\top\frac{W_{\perp}^\top\hat{w}}{\|W_{\perp}^\top\hat{w}\|_2} = w_1\right)\\
	=& \PP\left(\frac{U_{\perp}^\top\hat{u}}{\|U_{\perp}^\top\hat{u}\|_2} \in A\bigg|\frac{V_{\perp}^\top\hat{v}}{\|V_{\perp}^\top\hat{v}\|_2} = v_2, \frac{W_{\perp}^\top\hat{w}}{\|W_{\perp}^\top\hat{w}\|_2} = w_2\right).
	\end{align*}
	Therefore, $\frac{U_{\perp}^\top\hat{u}}{\|U_{\perp}^\top\hat{u}\|_2}$ and $\big(\frac{V_{\perp}^\top\hat{v}}{\|V_{\perp}^\top\hat{v}\|_2}, \frac{W_{\perp}^\top\hat{w}}{\|W_{\perp}^\top\hat{w}\|_2}\big)$ are independent. Similarly, we know that $\frac{U_{\perp}^\top\hat{u}}{\|U_{\perp}^\top\hat{u}\|_2}, \frac{V_{\perp}^\top\hat{v}}{\|V_{\perp}^\top\hat{v}\|_2}$ and $\frac{W_{\perp}^\top\hat{w}}{\|W_{\perp}^\top\hat{w}\|_2}$ are independent. In addition, by using the property of Haar measure, we know that $\big(\frac{\hat{u}^\top U_{\perp}}{\|U_{\perp}^\top\hat{u}\|_2},  \frac{\hat{v}^\top V_{\perp}}{\|V_{\perp}^\top\hat{v}\|_2}, \frac{\hat{w}^\top W_{\perp}}{\|W_{\perp}^\top\hat{w}\|_2}\big)^\top$ and $(\langle u, \hat{u}\rangle, \langle v, \hat{v}\rangle, \langle w, \hat{w}\rangle)^\top$ are independent, and
	\begin{align*}
	\frac{U_{\perp}^\top\hat{u}}{\|U_{\perp}^\top\hat{u}\|_2}\stackrel{\rm d.}{=}&\bigg(\frac{g_1^{(1)}}{\sqrt{\sum_{i=1}^{p_1-1}g_i^{(1)2}}}, \dots, \frac{g_{p_1-1}^{(1)}}{\sqrt{\sum_{i=1}^{p_1-1}g_i^{(1)2}}}\bigg)^\top,\\
	\frac{V_{\perp}^\top\hat{v}}{\|V_{\perp}^\top\hat{v}\|_2}\stackrel{\rm d.}{=}&\bigg(\frac{g_1^{(2)}}{\sqrt{\sum_{i=1}^{p_2-1}g_i^{(2)2}}}, \dots, \frac{g_{p_2-1}^{(2)}}{\sqrt{\sum_{i=1}^{p_2-1}g_i^{(2)2}}}\bigg)^\top,\\
	\frac{W_{\perp}^\top\hat{w}}{\|W_{\perp}^\top\hat{w}\|_2}\stackrel{\rm d.}{=}&\bigg(\frac{g_1^{(3)}}{\sqrt{\sum_{i=1}^{p_3-1}g_i^{(3)2}}}, \dots, \frac{g_{p_3-1}^{(3)}}{\sqrt{\sum_{i=1}^{p_3-1}g_i^{(3)2}}}\bigg)^\top
	\end{align*}
	where $g^{(1)} = (g_1^{(1)}, \dots, g_{p_1-1}^{(1)})^\top, g^{(2)} = (g_1^{(2)}, \dots, g_{p_2-1}^{(2)})^\top, g^{(3)} = (g_{1}^{(3)}, \dots, g_{p_3-1}^{(3)})^\top$ are independent standard Gaussian random vectors. Moreover, by SLLN, 
	\begin{align*}
	\frac{1}{p_1}\sum_{i=1}^{p_1-1}g_i^{(1)2} \stackrel{\rm a.s.}{\to} 1, \quad \frac{1}{p_2}\sum_{i=1}^{p_2-1}g_i^{(2)2} \stackrel{\rm a.s.}{\to} 1, \quad \frac{1}{p_3}\sum_{i=1}^{p_3-1}g_i^{(3)2} \stackrel{\rm a.s.}{\to} 1.
	\end{align*}
	For any fixed $f_1 \in \SS^{p_1-2}, f_2 \in \SS^{p_2-2}, f_3 \in \SS^{p_3-2}$, notice that $(f_1^\top g_1, f_2^\top g_2, f_3^\top g_3)^\top \sim N(0, I_3)$, we have
	\begin{align}\label{ineq109}
	   &\bigg(\sqrt{p_1}\frac{\hat{u}^\top U_{\perp}}{\|U_{\perp}^\top\hat{u}\|_2}f_1, \sqrt{p_2}\frac{\hat{v}^\top V_{\perp}}{\|V_{\perp}^\top\hat{v}\|_2}f_2, \sqrt{p_3}\frac{\hat{w}^\top W_{\perp}}{\|W_{\perp}^\top\hat{w}\|_2}f_3\bigg)^\top\notag\\
		\stackrel{\rm d.}{=}& \diag\bigg(\sqrt{\frac{p_1}{\sum_{i=1}^{p_1-1}g_i^{(1)2}}}, \sqrt{\frac{p_2}{\sum_{i=1}^{p_2-1}g_i^{(2)2}}}, \sqrt{\frac{p_3}{\sum_{i=1}^{p_3-1}g_i^{(3)2}}}\bigg)\cdot(f_1^\top g_1, f_2^\top g_2, f_3^\top g_3)^\top\notag\\
		\stackrel{\rm d.}{\to}& N(0, I_3).
	\end{align}
	By Theorem \ref{thm:na_tsvd}, 
	\begin{equation}\label{ineq131}
		\bigg(\frac{\langle u, \hat u\rangle^2 - (1 - p_1\lambda^{-2})}{\sqrt{2p_1}\lambda^{-2}},
		\frac{\langle v, \hat v\rangle^2 - (1 - p_2\lambda^{-2})}{\sqrt{2p_2}\lambda^{-2}},
		\frac{\langle w, \hat w\rangle^2 - (1 - p_3\lambda^{-2})}{\sqrt{2p_3}\lambda^{-2}}
		\bigg)^\top \stackrel{\rm d.}{\to} N(0, I_3).
	\end{equation}
	The delta method and the fact that $1 - p_i\lambda^{-2} \stackrel{\rm a.s.}{\to} 1$ for $i \in [3]$ together show that
	\begin{equation}\label{ineq132}
	\bigg(\frac{\langle u, \hat u\rangle - \sqrt{1 - p_1\lambda^{-2}}}{\sqrt{p_1/2}\lambda^{-2}},
	\frac{\langle v, \hat v\rangle - \sqrt{1 - p_2\lambda^{-2}}}{\sqrt{p_2/2}\lambda^{-2}},
	\frac{\langle w, \hat w\rangle - \sqrt{1 - p_3\lambda^{-2}}}{\sqrt{p_3/2}\lambda^{-2}}
	\bigg)^\top \stackrel{\rm d.}{\to} N(0, I_3).
	\end{equation}
	Also note that $$1 - p_i\lambda^{-2}/2 - \sqrt{1 - p_i\lambda^{-2}} = \frac{p_i^2\lambda^{-4}}{4(\sqrt{1 - p_i\lambda^{-2}} + 1 - p_i\lambda^{-2}/2)} \asymp \frac{p_i^2}{\lambda^4} \ll \frac{\sqrt{p_i}}{\lambda^2}$$ and $\big(\frac{\hat{u}^\top U_{\perp}}{\|U_{\perp}^\top\hat{u}\|_2},  \frac{\hat{v}^\top V_{\perp}}{\|V_{\perp}^\top\hat{v}\|_2}, \frac{\hat{w}^\top W_{\perp}}{\|W_{\perp}^\top\hat{w}\|_2}\big)^\top$ is independent of $(\langle u, \hat{u}\rangle, \langle v, \hat{v}\rangle, \langle w, \hat{w}\rangle)^\top$, by \eqref{ineq109} and \eqref{ineq132}, for any fixed $f_1 \in \SS^{p_1-2}, f_2 \in \SS^{p_2-2}, f_3 \in \SS^{p_3-2}$, 
	\begin{align}\label{ineq133}
		\begin{pmatrix}
		\sqrt{p_1}\frac{\hat{u}^\top U_{\perp}}{\|U_{\perp}^\top\hat{u}\|_2}f_1\\
		\sqrt{p_2}\frac{\hat{v}^\top V_{\perp}}{\|V_{\perp}^\top\hat{v}\|_2}f_2\\
		\sqrt{p_3}\frac{\hat{w}^\top W_{\perp}}{\|W_{\perp}^\top\hat{w}\|_2}f_3\\
		\frac{\langle u, \hat u\rangle - (1 - p_1\lambda^{-2}/2)}{\sqrt{p_1/2}\lambda^{-2}}\\
		\frac{\langle v, \hat v\rangle - (1 - p_2\lambda^{-2}/2)}{\sqrt{p_2/2}\lambda^{-2}}\\
		\frac{\langle w, \hat w\rangle - (1 - p_3\lambda^{-2}/2)}{\sqrt{p_3/2}\lambda^{-2}}
		\end{pmatrix} \stackrel{\rm d.}{\to} N(0, I_6).
	\end{align}
    By \eqref{ineq131},
	\begin{equation*}
	\frac{\|U_{\perp}^\top\hat{u}\|_2^2}{p_1/\lambda^2} = \frac{1-\langle\hat u, u\rangle^2}{p_1/\lambda^2} \stackrel{\rm p.}{\to} 1.
	\end{equation*}
	By \eqref{ineq133}, for any fixed $f_1 \in \SS^{p_1-2}, f_2 \in \SS^{p_2-2}, f_3 \in \SS^{p_3-2}$, 
	\begin{align}\label{ineq110}
	 &\bigg(\lambda\hat{u}^\top U_{\perp}f_1, \lambda\hat{v}^\top V_{\perp}f_2, \lambda\hat{w}^\top W_{\perp}f_3, \\
	 & \qquad \qquad \frac{\langle u, \hat u\rangle - (1 - p_1\lambda^{-2}/2)}{\sqrt{p_1/2}\lambda^{-2}}, \frac{\langle v, \hat v\rangle - (1 - p_2\lambda^{-2}/2)}{\sqrt{p_2/2}\lambda^{-2}}, \frac{\langle w, \hat w\rangle - (1 - p_3\lambda^{-2}/2)}{\sqrt{p_3/2}\lambda^{-2}}\bigg)^\top\notag\\
	=& \diag\bigg(\frac{\|U_{\perp}^\top\hat{u}\|_2}{\sqrt{p_1}/\lambda}, \frac{\|V_{\perp}^\top\hat{v}\|_2}{\sqrt{p_2}/\lambda}, \frac{\|W_{\perp}^\top\hat{w}\|_2}{\sqrt{p_3}/\lambda}, 1, 1, 1\bigg)\begin{pmatrix}
	\sqrt{p_1}\frac{\hat{u}^\top U_{\perp}}{\|U_{\perp}^\top\hat{u}\|_2}f_1\\
	\sqrt{p_2}\frac{\hat{v}^\top V_{\perp}}{\|V_{\perp}^\top\hat{v}\|_2}f_2\\
	\sqrt{p_3}\frac{\hat{w}^\top W_{\perp}}{\|W_{\perp}^\top\hat{w}\|_2}f_3\\
	\frac{\langle u, \hat u\rangle - (1 - p_1\lambda^{-2}/2)}{\sqrt{p_1/2}\lambda^{-2}}\\
	\frac{\langle v, \hat v\rangle - (1 - p_2\lambda^{-2}/2)}{\sqrt{p_2/2}\lambda^{-2}}\\
	\frac{\langle w, \hat w\rangle^2 - (1 - p_3\lambda^{-2}/2)}{\sqrt{p_3/2}\lambda^{-2}}
	\end{pmatrix}\notag\\ \stackrel{\rm d.}{\to}& N(0, I_{6}).
	\end{align}
	For simplicity, let $q_i = q_i^{(p_i)}$ for $i \in [3]$. Note that
	\begin{equation}\label{eq:u_decomposition}
	\langle \hat u, q_1\rangle = \langle \hat u, \calP_{u}q_1\rangle + \langle\hat u, \calP_{u}^\perp q_1\rangle = (q_1^\top u)\hat u^\top u + (U_{\perp}^\top q_1)^\top U_{\perp}^\top\hat{u}.
	\end{equation}
	If $q_1 \neq \pm u, q_2 \neq \pm v, q_3 \neq \pm w$ for $i \in [3]$, since $U_{\perp}^\top q_1, V_{\perp}^\top q_2, W_{\perp}^\top q_3$ are fixed vectors, by \eqref{ineq110}, we have
	\begin{align}\label{ineq107}
		&\begin{pmatrix}
		\lambda\frac{(U_{\perp}^\top q_1)^\top}{\|U_{\perp}^\top q_1\|_2} U_{\perp}^\top\hat{u},
		\frac{\langle u, \hat u\rangle - (1 - p_1\lambda^{-2}/2)}{\sqrt{p_1/2}\lambda^{-2}},
		\lambda\frac{(V_{\perp}^\top q_2)^\top}{\|V_{\perp}^\top q_2\|_2} V_{\perp}^\top\hat{v},
		\frac{\langle v, \hat v\rangle - (1 - p_2\lambda^{-2}/2)}{\sqrt{p_2/2}\lambda^{-2}},
		\lambda\frac{(W_{\perp}^\top q_3)^\top}{\|W_{\perp}^\top q_3\|_2} W_{\perp}^\top\hat{w},
		\frac{\langle w, \hat w\rangle - (1 - p_3\lambda^{-2}/2)}{\sqrt{p_3/2}\lambda^{-2}}
		\end{pmatrix}^\top\notag\\ &\stackrel{\rm d.}{\to} N(0, I_6).
	\end{align}
	Since 
	\begin{align*}
		\frac{\langle q_1, \hat u - u\rangle + \frac{p_1\langle q_1, u\rangle}{2\lambda^2}}{\sqrt{\frac{p_1\langle q_1, u\rangle^2}{2\lambda^4} + \frac{1 - \langle q_1, u\rangle^2}{\lambda^2}}} =& \bigg(\frac{\frac{\sqrt{1 - \langle u, q_1\rangle^2}}{\lambda}}{\sqrt{\frac{1 - \langle u, q_1\rangle^2}{\lambda^2} + \frac{p_1\langle u, q_1\rangle^2}{2\lambda^2}}}, \frac{\frac{\sqrt{p_1/2}\langle u, q_1\rangle}{\lambda^2}}{\sqrt{\frac{1 - \langle u, q_1\rangle^2}{\lambda^2} + \frac{p_1\langle u, q_1\rangle^2}{2\lambda^2}}}\bigg)\\&\cdot\bigg(\lambda\frac{(U_{\perp}^\top q_1)^\top}{\|U_{\perp}^\top q_1\|_2} U_{\perp}^\top\hat{u},
		\frac{\langle u, \hat u\rangle - (1 - p_1\lambda^{-2}/2)}{\sqrt{p_1/2}\lambda^{-2}}\bigg)^\top
	\end{align*}
	where $\bigg(\frac{\frac{\sqrt{1 - \langle u, q_1\rangle^2}}{\lambda}}{\sqrt{\frac{1 - \langle u, q_1\rangle^2}{\lambda^2} + \frac{p_1\langle u, q_1\rangle^2}{2\lambda^2}}}, \frac{\frac{\sqrt{p_1/2}\langle u, q_1\rangle}{\lambda^2}}{\sqrt{\frac{1 - \langle u, q_1\rangle^2}{\lambda^2} + \frac{p_1\langle u, q_1\rangle^2}{2\lambda^2}}}\bigg)^\top$ is a fixed unit vector, by Lemma \ref{lm:asymptotic}, we have
	\begin{align}\label{ineq134}
		&\bigg(\frac{\langle q_1, \hat u - u\rangle + \frac{p_1\langle q_1, u\rangle}{2\lambda^2}}{\sqrt{\frac{p_1\langle q_1, u\rangle^2}{2\lambda^4} + \frac{1 - \langle q_1, u\rangle^2}{\lambda^2}}}, \frac{\langle q_2, \hat v - v\rangle + \frac{p_2\langle q_2, v\rangle}{2\lambda^2}}{\sqrt{\frac{p_2\langle q_2, v\rangle^2}{2\lambda^4} + \frac{1 - \langle q_2, v\rangle^2}{\lambda^2}}}, \frac{\langle q_3, \hat w - w\rangle + \frac{p_3\langle q_3, w\rangle}{2\lambda^2}}{\sqrt{\frac{p_3\langle q_3, w\rangle^2}{2\lambda^4} + \frac{1 - \langle q_3, w\rangle^2}{\lambda^2}}}\bigg)^\top \stackrel{\rm d.}{\to} N(0, I_3).
	\end{align}
	If $q_1 = \pm u$, $q_2 = \pm v$ or $q_3 = \pm w$, by \eqref{ineq131}, we still have \eqref{ineq134}.
	
	Specifically, if if $|u_i|, |v_j|, |w_k| \ll \min\{\lambda/p, 1\}$ for some $i \in [p_1], j \in [p_2], k \in [p_3]$, by setting $q_1 = e_i, q_2 = e_j, q_3 = e_k$ and noticing that $\frac{p_1u_i^2}{\lambda^4} \ll \frac{p_1^2u_i^2}{\lambda^4} \ll \frac{1}{\lambda^2}$ and $u_i \stackrel{\rm a.s.}{\to} 0$, we know that \eqref{ineq:hat_u_asymp} holds.
	
	Given $\lambda^{-1} \ll |u_i|, |v_j|, |w_k| \ll \min\{\lambda/p, 1/\sqrt{\log(p)}\}$, immediately we have $\frac{\hat u_i}{u_i} \stackrel{\rm p.}{\to} 1, \frac{\hat v_j}{v_j} \stackrel{\rm p.}{\to} 1, \frac{\hat w_k}{w_k} \stackrel{\rm p.}{\to} 1$. Then
	\begin{align*}
		\begin{pmatrix}
		\lambda\frac{\hat u_i\hat v_j\hat w_k - u_i\hat v_j\hat w_k}{v_jw_k}\\
		\lambda\frac{u_i\hat v_j\hat w_k - u_iv_j\hat w_k}{u_iw_k}\\
		\lambda\frac{u_iv_j\hat w_k - u_iv_jw_k}{u_iv_j}
		\end{pmatrix} = \begin{pmatrix}
		\frac{\hat v_j\hat w_k}{v_jw_k} & 0 & 0\\
		0 & \frac{\hat w_k}{w_k} & 0\\
		0 & 0 & 1
		\end{pmatrix}\begin{pmatrix}
		\lambda(\hat u_i-u_i)\\
		\lambda(\hat v_j-v_j)\\
		\lambda(\hat w_k-w_k)
		\end{pmatrix} \stackrel{\rm d.}{\to} N(0, I_3).
	\end{align*}
	By Lemma \ref{lm:asymptotic}, we have
	\begin{align}\label{ineq111}
		\lambda\frac{\hat u_i\hat v_j\hat w_k - u_iv_jw_k}{\sqrt{u_i^2v_j^2 + v_j^2w_k^2 + w_k^2u_i^2}} = \begin{pmatrix}
		\frac{v_jw_k}{\sqrt{u_i^2v_j^2 + v_j^2w_k^2 + w_k^2u_i^2}}\\ \frac{w_ku_i}{\sqrt{u_i^2v_j^2 + v_j^2w_k^2 + w_k^2u_i^2}}\\ \frac{u_iv_j}{\sqrt{u_i^2v_j^2 + v_j^2w_k^2 + w_k^2u_i^2}}
		\end{pmatrix}^\top\begin{pmatrix}
		\lambda\frac{\hat u_i\hat v_j\hat w_k - u_i\hat v_j\hat w_k}{v_jw_k}\\
		\lambda\frac{u_i\hat v_j\hat w_k - u_iv_j\hat w_k}{u_iw_k}\\
		\lambda\frac{u_iv_j\hat w_k - u_iv_jw_k}{u_iv_j}
		\end{pmatrix} \stackrel{\rm d.}{\to} N(0, 1).
	\end{align}
	Notice that $(\hat u_i^2\hat v_j^2 + \hat v_j^2\hat w_k^2 + \hat w_k^2\hat u_i^2)/(u_i^2v_j^2 + v_j^2w_k^2 + w_k^2u_i^2) \stackrel{\rm p.}{\to} 1$, we have
	\begin{equation*}
		\lambda\frac{\hat u_i\hat v_j\hat w_k - u_iv_jw_k}{\sqrt{\hat u_i^2\hat v_j^2 + \hat v_j^2\hat w_k^2 + \hat w_k^2\hat u_i^2}} \stackrel{\rm d.}{\to} N(0, 1).
	\end{equation*}
    Finally, by \eqref{ineq108}, with probability at least $1 - Cp^{-3}$, 
    \begin{equation}\label{ineq112}
    	\big|(\hat\lambda - \lambda)\frac{\hat u_i\hat v_j\hat w_k}{\sqrt{\hat u_i^2\hat v_j^2 + \hat v_j^2\hat w_k^2 + \hat w_k^2\hat u_i^2}}\big| \leq C_2(\frac{p}{\lambda} + \sqrt{\log(p)})|w_k| \ll C_2,
    \end{equation}
    i.e., 
    \begin{equation*}
    	\bigg|(\hat\lambda - \lambda)\frac{\hat u_i\hat v_j\hat w_k}{\sqrt{\hat u_i^2\hat v_j^2 + \hat v_j^2\hat w_k^2 + \hat w_k^2\hat u_i^2}}\bigg| \stackrel{\rm p.}{\to} 0.
    \end{equation*}
    Therefore, we conclude that 
    \begin{equation*}
    	\frac{\hat{\calT}_{ijk} - \calT_{ijk}}{\sqrt{\hat u_i^2\hat v_j^2 + \hat v_j^2\hat w_k^2 + \hat w_k^2\hat u_i^2}} \stackrel{\rm d.}{\to} N(0, 1).
    \end{equation*}
\subsection{Proof of Theorem~\ref{th:CI}}
	Without loss of generality, we assume $\sigma = 1$. We discuss in four scenarios:
	\begin{itemize}[leftmargin=*]
		\item[(1).]{$|u_i|, |v_j|, |w_k|\geq (\log(p))^{1/8}\lambda^{-1}$.} By Theorem \ref{thm:entry_inference}, 
		\begin{equation*}
		\lim_{p\to\infty}\PP\left(|\calT_{ijk} - T_{ijk}| \leq z_{\alpha/2}\sqrt{\hat u_i^2\hat v_j^2 + \hat v_j^2\hat w_k^2 + \hat w_k^2\hat u_i^2}\right) = 1 - \alpha.
		\end{equation*}
		Therefore \eqref{ineq:CI} holds.
		\item[(2).]{Exactly two of $|u_i|, |v_j|, |w_k|\geq (\log(p))^{1/8}\lambda^{-1}$}. Without loss of generality, we assume $|v_j|, |w_k|\geq (\log(p))^{1/8}\lambda^{-1}$. By the essentially same proof of \eqref{ineq111}, we have
		\begin{equation*}
			\lambda\frac{\hat u_i\hat v_j\hat w_k - u_iv_jw_k}{\sqrt{u_i^2v_j^2 + v_j^2w_k^2 + w_k^2u_i^2}} \stackrel{\rm d.}{\to} N(0, 1).
		\end{equation*}
		(If $u_i = 0$, then immediately we have $\lambda\frac{\hat u_i\hat v_j\hat w_k - u_iv_jw_k}{\sqrt{u_i^2v_j^2 + v_j^2w_k^2 + w_k^2u_i^2}} = \frac{\hat v_j\hat w_k}{v_jw_k}\cdot \lambda\hat u_i \stackrel{\rm d.}{\to} N(0, 1).$)\\
		If $|u_i| \geq (\log(p))^{1/16}\lambda^{-1}$, then by \eqref{ineq:hat_u_asymp}, $\hat u_i/u_i, \hat v_j/v_j$, $\hat w_k/v_k \stackrel{\rm p.}{\to} 1$. Therefore, $(\hat u_i^2\hat v_j^2 + \hat v_j^2\hat w_k^2 + \hat w_k^2\hat u_i^2)/(u_i^2v_j^2 + v_j^2w_k^2 + w_k^2u_i^2) \stackrel{\rm p.}{\to} 1$. If $|u_i| < (\log(p))^{1/16}\lambda^{-1}$, then \eqref{ineq:hat_u_asymp} shows that $\hat{u}_i/((\log(p))^{1/16}\lambda^{-1}) \stackrel{\rm p.}{\to} 0$. Thus
		\begin{equation*}
			\frac{\hat u_i^2\hat v_j^2 + \hat v_j^2\hat w_k^2 + \hat w_k^2\hat u_i^2}{u_i^2v_j^2 + v_j^2w_k^2 + w_k^2u_i^2} = \frac{\hat u_i^2\hat v_j^2 + \hat v_j^2\hat w_k^2 + \hat w_k^2\hat u_i^2}{\hat v_j^2\hat w_k^2}\cdot \frac{v_j^2w_k^2}{u_i^2v_j^2 + v_j^2w_k^2 + w_k^2u_i^2}\cdot \frac{\hat v_j^2\hat w_k^2}{v_j^2w_k^2} \stackrel{\rm p.}{\to} 1.
		\end{equation*}
		As a consequence,
		\begin{equation}\label{ineq113}
			\lambda\frac{\hat u_i\hat v_j\hat w_k - u_iv_jw_k}{\sqrt{\hat u_i^2\hat v_j^2 + \hat v_j^2\hat w_k^2 + \hat w_k^2\hat u_i^2}} \stackrel{\rm d.}{\to} N(0, 1).
		\end{equation}
	    By combining \eqref{ineq112} and \eqref{ineq113} together, we have
	    \begin{equation*}
	    	\frac{\hat{\calT}_{ijk} - \calT_{ijk}}{\sqrt{\hat u_i^2\hat v_j^2 + \hat v_j^2\hat w_k^2 + \hat w_k^2\hat u_i^2}} \stackrel{\rm d.}{\to} N(0, 1),
	    \end{equation*}
	    which indicates that \eqref{ineq:CI} holds.
	    \item[(3).] {At least two of $|u_i|, |v_j|, |w_k| < (\log(p))^{1/8}\lambda^{-1}$.} Without loss of generality, we assume $|v_j|, |w_k| < (\log(p))^{1/8}\lambda^{-1}$. By \eqref{ineq:hat_u_asymp}, $$\frac{|\hat v_j|}{(\log(p))^{1/6}\lambda^{-1}}  \stackrel{\rm p.}{\to} 0, \quad \frac{|\hat w_k|}{(\log(p))^{1/6}\lambda^{-1}} \stackrel{\rm p.}{\to} 0.$$
	    By \eqref{ineq108}, we have $\hat\lambda/\lambda \stackrel{\rm p.}{\to} 1$. Then
	    \begin{align*}
	    	& \bigg|\frac{\hat \calT_{ijk}}{\sqrt{s(\hat u_i^2)s(\hat v_j^2)+s(\hat v_j^2)s(\hat w_k^2)+s(\hat w_k^2)s(\hat u_i^2)}}\bigg| = \frac{|\hat\lambda\hat u_i\hat v_j\hat w_k|}{\sqrt{s(\hat u_i^2)s(\hat v_j^2)+s(\hat v_j^2)s(\hat w_k^2)+s(\hat w_k^2)s(\hat u_i^2)}} \\
	    	\leq & \hat\lambda|\hat w_k|\frac{|\hat v_k|}{\sqrt{s(\hat v_k^2)}} \leq (\log(p))^{-1/6}\frac{\hat\lambda}{\lambda}\cdot\frac{|\hat v_k|}{(\log(p))^{1/6}\lambda^{-1}}\frac{|\hat w_k|}{(\log(p))^{1/6}\lambda^{-1}}
	    	\stackrel{\rm p.}{\to} 0,
	    \end{align*}
	    and
	    \begin{align*}
	    	& \bigg|\frac{\calT_{ijk}}{\sqrt{s(\hat u_i^2)s(\hat v_j^2)+s(\hat v_j^2)s(\hat w_k^2)+s(\hat w_k^2)s(\hat u_i^2)}}\bigg| = \frac{|\lambda u_i v_j w_k|}{\sqrt{s(\hat u_i^2)s(\hat v_j^2)+s(\hat v_j^2)s(\hat w_k^2)+s(\hat w_k^2)s(\hat u_i^2)}} \\
	    	\leq & \lambda|w_k|\frac{|v_k|}{\sqrt{s(\hat v_k^2)}} \leq \lambda\frac{((\log(p))^{1/8}\lambda^{-1})^2}{\sqrt{\log(p)}\hat\lambda^{-1}}
	    	\stackrel{\rm p.}{\to} 0.
	    \end{align*}
	    Therefore, 
	    \begin{align*}
	    	\lim_{p \to \infty}\PP(\calT_{ijk} \in \widetilde{CI}_{\alpha}(\hat \calT_{ijk})) = \lim_{p \to \infty}\PP\bigg(\frac{|\hat \calT_{ijk} - \calT_{ijk}|}{\sqrt{s(\hat u_i^2)s(\hat v_j^2)+s(\hat v_j^2)s(\hat w_k^2)+s(\hat w_k^2)s(\hat u_i^2)}} \leq z_{\alpha/2}\bigg) = 1.
	    \end{align*}
	\end{itemize}
In conclusion, we have proved \eqref{ineq:CI}.

\newpage

\subsection{Proof of supporting lemmas}

\begin{lemma}\label{lm:equivalence-different-error}
	For any $0 \leq \delta \leq 1$, if either of the following inequality holds, (1) $\|\calM_j(\hat{\calT}^{(0)} - \calT)\|\leq \delta\lambda_{\submin}/2$; (2) $\|\hat{U}_j^{(0)\top} \hat{U}_j^{(0)\top} - U_j U_j^\top\| \leq \delta$; (3) $\|\hat{U}_j^{(0)\top} \hat{U}_j^{(0)\top} - U_j U_j^\top\|_{\F} \leq \sqrt{2}\delta$; (4) $\|\hat{U}_j^{(0)\top} U_j\|\geq \sqrt{1-\delta^2}$; (5) $\|\hat{U}_j^{(0)\top} U_j\|_\F\geq \sqrt{r_j-\delta^2}$, we have $\|\sin\Theta(\hat{U}_j^{(0)}, U_j)\|\leq \delta.$
\end{lemma}
\begin{proof}[Proof of Lemma~\ref{lm:equivalence-different-error}]
	For simplicity, let $T_j$ and $\hat T_j^{(0)}$ denote $\calM_j(\calT)$ and $\calM_j(\hat\calT^{(0)})$, respectively. Suppose $\|\hat T_j^{(0)} - T_j\| \leq \delta\lambda_{\submin}/2$. By \cite[Lemma 6]{zhang2018tensor}, we have
	$$\|\hat U_{j\perp}^{(0)} T_j\| \leq 2\|\hat T_j^{(0)} - T_j\| \leq \delta\lambda_{\submin}$$
	and consequently,
	\begin{align*}
	\big\|\sin\Theta(\hat U_j^{(0)}, U_j)\big\| = \|\hat U_{j\perp}^{(0)\top}U_j\| \leq \frac{\|\hat U_{j\perp}^{(0)\top}U_jU_j^\top T_j\|}{\sigma_{\min}(U_j^\top T_j)} = \frac{\|\hat U_{j\perp}^{(0)\top} T_j\|}{\lambda_{\submin}} \leq \delta.
\end{align*}
In addition, by \cite[Lemma 1]{cai2018rate}, we have
\begin{equation*}
\big\|\sin\Theta(\hat U_j^{(0)}, U_j)\big\| = \sqrt{1 - \|\hat U_j^{(0)\top}U_j\|^2} \leq \|\hat U_j^{(0)}\hat U_j^{(0)\top} - U_jU_j^\top\|
\end{equation*}
and 
\begin{equation*}
\big\|\sin\Theta(\hat U_j^{(0)}, U_j)\big\| \leq \big\|\sin\Theta(\hat U_j^{(0)}, U_j)\big\|_{\F} = \sqrt{r_j - \|\hat U_j^{(0)\top}U_j\|_{\F}^2} = \|\hat U_j^{(0)}\hat U_j^{(0)\top} - U_jU_j^\top\|_{\F}/\sqrt{2},
\end{equation*}
which have finished the proof of Lemma \ref{lm:equivalence-different-error}.
\end{proof}
\begin{proof}[Proof of Lemma~\ref{lm:variance_estimation}]
	Notice that $$\|\calA - \calA \times_1 \calP_{\hat U_1} \times_2 \calP_{\hat U_2} \times_3 \calP_{\hat U_3}\|_{\F} = \|\calZ - \calZ \times_1 \calP_{\hat U_1} \times_2 \calP_{\hat U_2} \times_3 \calP_{\hat U_3} + \calT - \calT \times_1 \calP_{\hat U_1} \times_2 \calP_{\hat U_2} \times_3 \calP_{\hat U_3}\|_{\F},$$ we have
	\begin{align*}
	&\big|\|\calA - \calA \times_1 \calP_{\hat U_1} \times_2 \calP_{\hat U_2} \times_3 \calP_{\hat U_3}\|_{\F} - \|\calZ\|_{\F}\big|\\ \leq& \|\calZ \times_1 \calP_{\hat U_1} \times_2 \calP_{\hat U_2} \times_3 \calP_{\hat U_3}\|_{\F} + \|\calT - \calT \times_1 \calP_{\hat U_1} \times_2 \calP_{\hat U_2} \times_3 \calP_{\hat U_3}\|_{\F}.
	\end{align*}
	By \eqref{ineq:projection} and \eqref{ineq:sub-Gaussian}, with probability at least $1-C_1e^{-c_1p}$,
	\begin{align*}
	\|\calZ \times_1 \calP_{\hat U_1} \times_2 \calP_{\hat U_2} \times_3 \calP_{\hat U_3}\|_{\F} = \|\hat U_1^\top Z_1(\hat U_2 \otimes \hat U_3)\|_{\F} \leq \sqrt{r_1}\|Z_1(\hat U_2 \otimes \hat U_3)\| \leq C_2\sigma\sqrt{pr}.
	\end{align*}
	In addition, we have
	\begin{align*}
	&\|\calT - \calT \times_1 \calP_{\hat U_1} \times_2 \calP_{\hat U_2} \times_3 \calP_{\hat U_3}\|_{\F}\\ =& \|\calT \times_1 \calP_{U_1} \times_2 \calP_{U_2} \times_3 \calP_{U_3} - \calT \times_1 \calP_{\hat U_1} \times_2 \calP_{\hat U_2} \times_3 \calP_{\hat U_3}\|_{\F}\\
	\leq& \|(\calP_{U_1} - \calP_{\hat U_1})T_1(\calP_{U_{2}} \otimes \calP_{U_{3}})\|_{\F} + \|(\calP_{U_2} - \calP_{\hat U_2})T_2(\calP_{\hat U_{1}} \otimes \calP_{U_{3}})\|_{\F}
	+ \|(\calP_{U_3} - \calP_{\hat U_3})T_3(\calP_{\hat U_{1}} \otimes \calP_{U_{2}})\|_{\F}\\
	\leq& \big(\|\calP_{U_1} - \calP_{\hat U_1}\| + \|\calP_{U_2} - \calP_{\hat U_2}\| + \|\calP_{U_2} - \calP_{\hat U_2}\|\big)\|\calT\|_{\F}\\
	\leq& C_2\frac{\sqrt{p}\sigma}{\lambda_{\submin}}\cdot \sqrt{r}\kappa_0\lambda_{\submin} = C_2\kappa_0\sigma\sqrt{pr}
	\end{align*}
	with probability at least $1-C_1e^{-c_1p}$.
	Therefore, with probability at least $1-C_1e^{-c_1p}$, we have
	\begin{equation}\label{ineq114}
	\big|\|\calA - \calA \times_1 \calP_{\hat U_1} \times_2 \calP_{\hat U_2} \times_3 \calP_{\hat U_3}\|_{\F} - \|\calZ\|_{\F}\big| \leq C_2\kappa_0\sigma\sqrt{pr}.
	\end{equation}
	By \cite[Lemma 1]{laurent2000adaptive}, 
	\begin{equation*}
	\PP\left(\left|\frac{\|\calZ\|_{\F}^2}{\sigma^2} - p_1p_2p_3\right| \geq C_2(\sqrt{p_1p_2p_3}\sqrt{\log(p)} + \log(p))\right) \leq p^{-3}.
	\end{equation*}
	As a consequence, with probability at least $1 - p^{-3}$, 
	\begin{equation*}
	\big|\|\calZ\|_{\F} - \sqrt{p_1p_2p_3}\sigma\big| \leq C_2\sqrt{\log(p)}\sigma.
	\end{equation*}
	Combing \eqref{ineq114} and the previous inequality together, we know that with probability at least $1 - C_1p^{-3}$, 
	\begin{equation*}
	|\hat\sigma/\sigma - 1| \leq C_2(\kappa_0\sqrt{r}p^{-1} + p^{-3/4}\sqrt{\log(p)})
	\end{equation*}
	and 
	\begin{equation*}
	|\hat\sigma^2/\sigma^2 - 1| = |\hat\sigma/\sigma - 1||\hat\sigma/\sigma + 1| \leq 2|\hat\sigma/\sigma - 1| + |\hat\sigma/\sigma - 1|^2 \leq C_2(\kappa_0\sqrt{r}p^{-1} + p^{-3/4}\sqrt{\log(p)}).
	\end{equation*}
\end{proof}

\begin{proof}[Proof of Lemma~\ref{lem:Ebound}]
	By definition, $\|\frakE_1\|\leq \|\frakJ_1\|+\|\frakJ_2\|+\|\frakJ_3\|+\|\frakJ_4\|$. 
	\begin{equation}\label{ineq:E_1}
	\|\frakE_1\|\leq \|\frakJ_1\|+\|\frakJ_2\|+\|\frakJ_3\|+\|\frakJ_4\|.
	\end{equation}
	We first proved the upper bound for $\|\frakJ_1\|$. By the definition of $\|\frakJ_1\|$, \begin{equation}\label{eq:J1init}
	\|\frakJ_1\|\leq \left\|T_1(\calP_{\hat U_2^{(1)}}\otimes\calP_{\hat U_3^{(1)}})Z_1^{\top}\right\| \leq \left\|T_1\|\|(\calP_{\hat U_2^{(1)}}\otimes\calP_{\hat U_3^{(1)}})Z_1^{\top}\right\| \leq \kappa_0\lambda_{\submin}\left\|Z_1(\hat{U}_2^{(1)} \otimes \hat{U}_3^{(1)})\right\|.
	\end{equation}
	For any fixed matrices $X \in \RR^{p_2 \times r_2}, Y \in \RR^{p_3 \times r_3}$ satisfying $\|X\|, \|Y\| \leq 1$, 
	\begin{equation}\label{ineq:sub_gaussian_fixed}
	\PP\left(\|Z_1(X \otimes Y)\| \geq C_2\sqrt{pr}\right) \leq C_1e^{-c_1pr}.
	\end{equation}
	Let $\mathcal{X}_{p_k, r_k} = \{X \in \RR^{p_k \times r_k}: \|X\| \leq 1\}$. By \cite[Lemma 7]{zhang2018tensor}, there exists an $1/4$-net $\bar{\mathcal{X}}_{p_k, r_k}$ with cardinality at most $9^{p_kr_k}$ for $\mathcal{X}_{p_k, r_k}$. That is, for any $X \in \mathcal{X}_{p_k, r_k}$, there exists $X' \in \mathcal{X}_{p_k, r_k}$ such that $\|X' - X\| \leq 1/4$. For any $X \in \mathcal{X}_{p_2, r_2}$ and $Y \in \mathcal{X}_{p_3, r_3}$, let $X' \in \bar{\mathcal{X}}_{p_2, r_2}$ and $Y' \in \bar{\mathcal{X}}_{p_3, r_3}$ satisfying $\|X - X'\| \leq 1/4, \|Y - Y'\| \leq 1/4$.
	Then
	\begin{align*}
	&\|Z_1(X \otimes Y)\|\\ \leq& \|Z_1(X' \otimes Y')\| + \|Z_1((X - X') \otimes Y)\| + \|Z_1(X \otimes (Y - Y'))\| + \|Z_1((X - X') \otimes (Y - Y'))\|\\
	\leq& \|Z_1(X' \otimes Y')\| + \frac{3}{4}\sup_{\stackrel{X \in \RR^{p_2 \times r_2}, Y \in \RR^{p_3 \times r_3}}{\|X\|, \|Y\| \leq 1}}\|Z_1(X \otimes Y)\|.
	\end{align*}
	By taking the supremum over any $X \in \mathcal{X}_{p_2, r_2}$ and $Y \in \mathcal{X}_{p_3, r_3}$, we have
	\begin{equation*}
	\sup_{\stackrel{X \in \RR^{p_2 \times r_2}, Y \in \RR^{p_3 \times r_3}}{\|X\|, \|Y\| \leq 1}}\|Z_1(X \otimes Y)\| \leq 4\sup_{X' \in \bar{\mathcal{X}}_{p_2, r_2}, Y' \in \bar{\mathcal{X}}_{p_3, r_3}}\|Z_1(X' \otimes Y')\|.
	\end{equation*}
	The union bound shows that  
	\begin{align*}
		&\PP\left(\sup_{\stackrel{X \in \RR^{p_2 \times r_2}, Y \in \RR^{p_3 \times r_3}}{\|X\|, \|Y\| \leq 1}}\|Z_1(X \otimes Y)\| \geq C_2\sqrt{pr}\right)\\
		\leq& \PP\left(\sup_{X' \in \bar{\mathcal{X}}_{p_2, r_2}, Y' \in \bar{\mathcal{X}}_{p_3, r_3}}\|Z_1(X' \otimes Y')\| \geq C_2\sqrt{pr}\right)\\
		\leq& \sum_{X' \in \bar{\mathcal{X}}_{p_2, r_2}, Y' \in \bar{\mathcal{X}}_{p_3, r_3}}\PP\left(\|Z_1(X' \otimes Y')\| \geq C_2\sqrt{pr}\right)\\ \leq& C_1e^{-c_1pr}.
	\end{align*}
	By \cite[Lemma 1]{cai2018rate}, with probability at least $1 - C_1e^{-c_1p}$, 
	\begin{equation}\label{perturbation_equivalence}
	\|\calP_{U_{k}}^{\perp} \hat{U}_k^{(1)}\| = \|U_{k\perp}^\top \hat{U}_k^{(1)}\| \leq \|\hat{U}_k^{(1)}\hat{U}_k^{(1)\top} - U_{k}U_{k}^\top\| \leq C_2\sqrt{p}/\lambda_{\submin}, \quad 1 \leq k \leq 3.
	\end{equation}
	Therefore, with probability at least $1 - C_1e^{-c_1p}$,
	\begin{align}\label{ineq:projection}
	&\left\|Z_1(\hat U_2^{(1)} \otimes \hat U_3^{(1)})\right\|\notag\\ =& \left\|Z_1(\calP_{U_2 \otimes U_3} + \calP_{U_{2\perp} \otimes U_3} + \calP_{U_2 \otimes U_{3\perp}} + \calP_{U_{2\perp} \otimes U_{3\perp}})(\hat U_2^{(1)} \otimes \hat U_3^{(1)})\right\|\notag\\
	\leq& \left\|Z_1\calP_{U_2 \otimes U_3}(\hat U_2^{(1)} \otimes \hat U_3^{(1)})\right\| + \left\|Z_1\calP_{U_{2\perp} \otimes U_3}(\hat U_2^{(1)} \otimes \hat U_3^{(1)})\right\|\notag\\
	&+ \left\|Z_1\calP_{U_2 \otimes U_{3\perp}}(\hat U_2^{(1)} \otimes \hat U_3^{(1)})\right\| + \left\|Z_1\calP_{U_{2\perp} \otimes U_{3\perp}}(\hat U_2^{(1)} \otimes \hat U_3^{(1)})\right\|\notag\\ 
	=& \left\|Z_1(U_2 \otimes U_3)(U_2 \otimes U_3)^\top(\hat U_2^{(1)} \otimes \hat U_3^{(1)})\right\| + \left\|Z_1\left((\calP_{U_{2}}^\perp\hat U_2^{(1)}) \otimes (\calP_{U_3}\hat U_3^{(1)})\right)\right\|\notag\\
	&+ \left\|Z_1\left((\calP_{U_{2}}\hat U_2^{(1)}) \otimes (\calP_{U_{3}}^\perp\hat U_3^{(1)})\right)\right\| + \left\|Z_1\left((\calP_{U_{2}}^\perp\hat U_2^{(1)}) \otimes (\calP_{U_{3}}^\perp\hat U_3^{(1)})\right)\right\|\notag\\
	\leq& \left\|Z_1(U_2 \otimes U_3)\right\| + C_2\sqrt{pr}\left\|\calP_{U_{2}}^\perp\hat U_2^{(1)}\right\|\left\|\calP_{U_3}\hat U_3^{(1)}\right\|\notag\\& + C_2\sqrt{pr}\left\|\calP_{U_{2}}\hat U_2^{(1)}\right\|\left\|\calP_{U_{3}}^\perp\hat U_3^{(1)}\right\|+ C_2\sqrt{pr}\left\|\calP_{U_{2}}^\perp\hat U_2^{(1)}\right\|\left\|\calP_{U_{3}}^\perp\hat U_3^{(1)}\right\|\notag\\
	\leq& \left\|Z_1(U_2 \otimes U_3)\right\| + C_2\sqrt{pr}\sqrt{p}\lambda_{\submin}^{-1} + C_2\sqrt{pr}\sqrt{p}\lambda_{\submin}^{-1} + C_2\sqrt{pr}p\lambda_{\submin}^{-2}\notag\\
	\leq& \left\|Z_1(U_2 \otimes U_3)\right\| + C_2\sqrt{pr}\sqrt{p}\lambda_{\submin}^{-1}\notag\\
	\leq& \left\|Z_1(U_2 \otimes U_3)\right\| + C_2\sqrt{p}.
	\end{align}
	By the Gaussian concentration inequality,
	\begin{equation}\label{ineq:sub-Gaussian}
	\PP\left(\left\|Z_1(U_2 \otimes U_3)\right\| \geq C_3\sqrt{p}\right) \leq C_1e^{-c_1p}.
	\end{equation}
	\eqref{eq:J1init}, \eqref{ineq:projection} and \eqref{ineq:sub-Gaussian} together imply that
	\begin{equation}\label{ineq:J_1}
	\PP\left(\|\frakJ_1\| \geq C_2\kappa_0\lambda_{\submin}\sqrt{p}\right) \leq C_1e^{-c_1p}.
	\end{equation}
	Since $\frakJ_2 = \frakJ_1^\top$, we also have
	\begin{equation}\label{ineq:J_2}
	\PP\left(\|\frakJ_2\| \geq C_2\kappa_0\lambda_{\submin}\sqrt{p}\right) \leq C_1e^{-c_1p}.
	\end{equation}
	
	For $\frakJ_3$, by definition,
	\begin{equation}\label{eq:J3init}
	\|\frakJ_3\|= \|Z_1(\hat U_2^{(1)} \otimes \hat U_3^{(1)})\|^2.
	\end{equation}
	Combining \eqref{eq:J3init}, \eqref{ineq:projection} and \eqref{ineq:sub-Gaussian} together, we have
	\begin{equation}\label{ineq:J_3}
	\PP\left(\|\frakJ_3\| \geq C_2p\right) \leq C_1e^{-c_1p}.
	\end{equation}
	
	Then, we consider $\frakJ_4$. By \eqref{perturbation_equivalence}, with probability at least $1 - C_1e^{-c_1p}$,
	\begin{align*}
	\|\frakJ_4 \leq& \left\|T_1((\calP_{\hat U_2^{(1)}} - \calP_{U_{2}})\otimes \calP_{\hat U_3^{(1)}})T_1^{\top}\right\| + \left\|T_1(\calP_{U_2} \otimes (\calP_{\hat U_3^{(1)}} - \calP_{U_3})T_1^{\top}\right\|\\
	=& \left\|U_1G_1((U_2^\top(\calP_{\hat U_2^{(1)}} - \calP_{U_{2}})U_2)\otimes (U_3^\top\calP_{\hat U_3^{(1)}}U_3))G_1^{\top}U_1^\top\right\|\\
	&+ \left\|U_1G_1(\calP_{U_{2}}\otimes (U_3^\top(\calP_{\hat U_3^{(1)}} - \calP_{U_{3}})U_3))G_1^{\top}U_1^\top\right\|\\
	=& \left\|U_1G_1((U_2^\top\calP_{\hat U_2^{(1)}}^\perp U_2)\otimes (U_3^\top\calP_{\hat U_3^{(1)}}U_3))G_1^{\top}U_1^\top\right\|\\
	&+ \left\|U_1G_1(\calP_{U_{2}}\otimes (U_3^\top\calP_{\hat U_3^{(1)}}^{\perp}U_3))G_1^{\top}U_1^\top\right\|\\
	\leq& \|G_1\|^2\left\|U_2^\top \hat U_{2\perp}^{(1)}\right\|^2 + \|G_1\|^2\left\|U_3^\top \hat U_{3\perp}^{(1)}\right\|^2\\
	\leq& C_2\kappa_0^2\lambda_{\submin}^2(\sqrt{p}/\lambda_{\submin})^2\\
	=&C_2\kappa_0^2p.
	\end{align*}

	Therefore, by \eqref{ineq:E_1}, \eqref{ineq:J_1}, \eqref{ineq:J_2} and \eqref{ineq:J_3} and notice that $\lambda_{\submin} \geq C_2\kappa_0\sqrt{p}$, we conclude with 
	$$
	\PP\big(\|\frakE_1\|\geq C_2\kappa_0\lambda_{\submin}\sqrt{p}\big)\leq C_1e^{-c_1p}.
	$$

\end{proof}

\begin{proof}[Proof of Lemma~\ref{lm:best_rotation}]
	Consider the SVD decomposition $\hat{U}^\top U = LSW^\top$, where $L, W \in \OO_r$, and $S \in \RR^{r \times r} = \diag(s_1, \dots, s_r)$ is a diagonal matrix with diagonal entries $1 \geq s_1 \geq \cdots \geq s_r \geq 0$. By setting $R = LW^\top$, we have $\hat{U}^\top U - R = L(S - I_r)W^\top$. Therefore, $\left\|\hat{U}^\top U - R\right\| = \|S -I_r\|$. Since $|x - 1| \leq |x^2 - 1|$ for all $x \geq 0$, we have
	\begin{equation}\label{ineq10}
	\left\|\hat{U}^\top U - R\right\| \leq \left\|S^2 - I_r\right\| = \left\|\hat{U}^\top UU^\top\hat{U} - I_r\right\| =\left\|\hat{U}^\top U_{\perp}U_{\perp}^\top\hat{U}\right\| = \left\|U_{\perp}^\top\hat{U}\right\|^2.
	\end{equation}
	For $\left\|\hat{U}^\top U - R\right\|_{\rm F}$, we have 
	\begin{equation}\label{ineq11}
		\left\|\hat{U}^\top U - R\right\|_{\rm F} \leq \sqrt{r}\left\|\hat{U}^\top U - R\right\| \leq \sqrt{r}\left\|U_{\perp}^\top\hat{U}\right\|^2
	\end{equation}
	and
	\begin{equation}\label{ineq12}
		\left\|\hat{U}^\top U - R\right\|_{\rm F} \leq \left\|S^2 - I_r\right\|_{\rm F} = \left\|\hat{U}^\top UU^\top\hat{U} - I_r\right\|_{\rm F} =\left\|\hat{U}^\top U_{\perp}U_{\perp}^\top\hat{U}\right\|_{\rm F} \leq \left\|U_{\perp}^\top\hat{U}\right\|_{\rm F}^2,
	\end{equation}
	which have finish the proof of Lemma \ref{lm:best_rotation}.
\end{proof}

\begin{lemma}\label{lm:subspace_regression}
	Under tensor regression model (\ref{eq:tr_model}) with $\calX(i_1,i_2,i_3)\stackrel{i.i.d.}{\sim} N(0,1)$, $\Var(\xi_i) = \sigma^2$ and $\|\xi_i\|_{\psi_2} \leq C\sigma$ for some constant $C > 0$, if $\|\tilde\calT-\calT\|_{\rm F}^2\leq C_2pr_{\submax}\sigma^2/n$, $n(\lambda_{\submin}/\sigma)^2\geq C_0(p^{3/2} \vee \kappa_0^4pr_{\submax}^2)$ and $n\geq C_0(p^{3/2} \vee \kappa_0^2pr_{\submax}^3)$ for some constants $C_0, C_2 > 0$, then there exists some constants $C_1, c_1, C_3 > 0$ such that with probability at least $1- C_1e^{-c_1p}$, 
	\begin{equation*}
	\left\|\sin\Theta(\hat U_j^{(1)},U_j)\right\| \leq C_3\sqrt{p/n}\sigma/\lambda_{\submin}, \quad \forall j=1,2,3,
	\end{equation*}
	where $U_j^{(1)}$ is the one-step alternating minimization defined in Algorithm \ref{algo:am_optimal_regression}.
\end{lemma}
\begin{proof}[Proof of Lemma~\ref{lm:subspace_regression}]
	Without loss of generality, we assume $\sigma = 1$. By Assumption~\ref{assump:tr} and \cite[Lemma 6]{zhang2018tensor}, with probability $1 - C_1e^{-c_1p}$,
	\begin{equation*}
	\big\|\hat{U}_{1\perp}^{(0)\top} T_1\big\|_{\rm F} \leq 2\big\|T_1 - \hat T_1^{(0)}\big\|_{\rm F} \leq C\frac{pr}{n}.
	\end{equation*}
	By \cite[Lemma 1]{cai2018rate}, we get with probability $1 - C_1e^{-c_1p}$ that
	\begin{align*}
	\inf_{O \in \OO_{r_1}}\big\|\hat U_1^{(0)} - U_1O\big\|_{\rm F} \leq& \big\|\hat U_1^{(0)}\hat U_1^{(0)\top} - U_1U_1^\top\big\|_{\rm F} =\sqrt{2}\big\|\sin\Theta(\hat U_1^{(0)}, U_1)\big\|_{\rm F} = \sqrt{2}\big\|\hat{U}_{1\perp}^{(0)\top} U_1\big\|_{\rm F} \\\leq& \sqrt{2}\frac{\|\hat{U}_{1\perp}^{(0)\top} T_1\|_{\rm F}}{\|G_1(U_2^\top \otimes U_3^\top)\|}
	\leq C\frac{\sqrt{pr/n}}{\|G_1\|} \leq C\frac{\sqrt{pr/n}}{\lambda_{\submin}}.
	\end{align*}
	Similarly, with the same probability, we get
	\begin{equation*}
	\big\|\hat U_2^{(0)}\hat U_2^{(0)\top} - U_2U_2^\top\big\|_{\rm F}, \big\|\hat U_3^{(0)}\hat U_3^{(0)\top} - U_3U_3^\top\big\|_{\rm F} \leq C\frac{\sqrt{pr/n}}{\lambda_{\submin}}.
	\end{equation*}
	Based on the two equations above, we can prove Lemma~\ref{lm:subspace_regression} by similar proof of \eqref{eq:regression_step_2}.
\end{proof}

\begin{lemma}\label{lm:epsilon_net}
	There exists an $\epsilon$-net $\bar{\OO}_{p, r} = \{U^{(j)} \in \OO_{p, r}, 1 \leq j \leq N\}$ in $\|\cdot\|$ norm with cardinality $N \leq \left((4 + \epsilon)/\epsilon\right)^{pr}$ for $\OO_{p, r}$. That is, for any $U \in \OO_{p, r}$, there exists $j \in [N]$ such that $\|U - U^{(j)}\| \leq \epsilon$.
\end{lemma}
\begin{proof}[Proof of Lemma \ref{lm:epsilon_net}]
		By \cite[Lemma 7]{zhang2018tensor}, for $\calU_{p, r} = \{U \in \RR^{p \times r}, \|U\| \leq 1\}$, there exists an $\epsilon/2$-net $\bar\calU_{p, r} = \{\bar U^{(j)} \in \RR^{p \times r}, \|\bar U^{(j)}\| \leq 1, 1 \leq j \leq N\}$ in $\|\cdot\|$ norm with $N \leq ((4+\epsilon)/\epsilon)^{pr}$ for $\calU_{p, r}$. Let $U^{(j)} \in \argmin_{U \in \OO_{p, r}}\|\bar U^{(j)} - U\|$, $1 \leq j \leq N$. For any $U \in \OO_{p, r}$, there exists $\bar U^{(j)}$ such that $\|\bar U^{(j)} - U\| \leq \epsilon/2$. Then $\|U^{(j)} - U\| \leq \|\bar U^{(j)} - U\| + \|\bar U^{(j)} - U^{(j)}\| \leq 2\|\bar U^{(j)} - U\| \leq \epsilon$. 
\end{proof}

\begin{lemma}\label{lm:AZB}
	Suppose $Z \in \RR^{p \times q}$ is a matrix with independent zero-mean $\sigma$-sub-Gaussian entries. $A \in \RR^{m \times p}, B \in \RR^{q \times n}$ satisfy $\|A\|, \|B\| \leq 1$, $m \leq p, n \leq q$. Then 
	\begin{equation}\label{ineq31}
	\PP\left(\|AZB\| \geq 2\sigma\sqrt{m + t}\right) \leq 2\cdot 5^n\exp\left[-c\min\left(\frac{t^2}{m}, t\right)\right].
	\end{equation}
		\begin{equation}\label{ineq34}
	\PP\left(\|AZB\|_{\rm F} \geq \sigma\sqrt{mn + t}\right) \leq 2\exp\left[-c\min\left(\frac{t^2}{mn}, t\right)\right].
	\end{equation}
\end{lemma}
\begin{proof}[Proof of Lemma \ref{lm:AZB}]
	Without loss of generality, assume $\sigma = 1$. For fixed $x \in \RR^n$ satisfying $\|x\|_2 = 1$, we have $AZBx = \text{vec}(AZBx) = (x^\top B^\top \otimes A)\text{vec}(Z)$. Since $Z_{ij}$ is $1$-sub-Gaussian, we know that $\Var(Z_{ij}) \leq 1$. In addition,
	\begin{equation}\label{ineq32}
	\begin{split}
	\EE \|(x^\top B^\top \otimes A)\text{vec}(Z)\|_{2}^2 =& \EE\left[\text{trace}\left(\text{vec}(Z)^\top(x^\top B^\top \otimes A)^\top (x^\top B^\top \otimes A)\text{vec}(Z)\right)\right]\\
	=& \text{trace}\left[\EE\left((x^\top B^\top \otimes A)^\top (x^\top B^\top \otimes A)\text{vec}(Z)\text{vec}(Z)^\top\right)\right]\\
	=& \text{trace}\left[(x^\top B^\top \otimes A)^\top (x^\top B^\top \otimes A)\EE\left(\text{vec}(Z)\text{vec}(Z)^\top\right)\right]\\
	\leq& \text{trace}\left((x^\top B^\top \otimes A)^\top (x^\top B^\top \otimes A)\right)\\
	=& \left\|x^\top B^\top \otimes A\right\|_{\rm F}^2 = \|Bx\|_2^2\|A\|_{\rm F}^2 \leq \|x\|_2^2\|A\|_{\rm F}^2 \\
	\leq& m.
	\end{split}
	\end{equation}
	The first inequality holds since $\EE\left(\text{vec}(Z)\text{vec}(Z)^\top\right)$ is a diagonal matrix with diagonal entries $\Var(Z_{ij}) \leq 1$; the last inequality is due to $\|A\|_{\rm F} \leq \min\{m, p\}\|A\|_2 \leq m$.\\
	By Hanson-Wright inequality, we have
	\begin{equation*}
	\PP\left(\|AZBx\|_{2}^2 - m \geq t\right) \leq 2\exp\left[-c\min\left(\frac{t^2}{\|(Bxx^\top B^\top) \otimes (A^\top A)\|_{\rm F}^2}, \frac{t}{\|(Bxx^\top B^\top) \otimes (A^\top A)\|}\right)\right].
	\end{equation*}
	Since $\|x\|_2 = 1$ and $\|A\|, \|B\| \leq 1$, 
	\begin{equation*}
	\begin{split}
	\|(Bxx^\top B^\top) \otimes (A^\top A)\|_{\rm F}^2 =& \|Bxx^\top B^\top\|_{\rm F}^2\|A^\top A\|_{\rm F}^2 = (x^\top B^\top Bx)^2\|A^\top A\|_{\rm F}^2\\ \leq& (x^\top x)^2\|A^\top A\|_{\rm F}^2 = \sum_{i=1}^{\min\{m, p\}}\sigma_i^4(A) \leq m,
	\end{split}
	\end{equation*}
	\begin{equation*}
	\|(Bxx^\top B^\top) \otimes (A^\top A)\| \leq \|Bxx^\top B^\top\|\|A^\top A\| \leq \|xx^\top\|\|A^\top A\| \leq 1.
	\end{equation*}
	Thus, for fixed $x$ satisfying $\|x\|_2 = 1$, we have
	\begin{equation}\label{ineq33}
	\PP\left(\|AZBx\|_{2}^2 \geq m + t\right) \leq 2\exp\left[-c\min\left(\frac{t^2}{m}, t\right)\right].
	\end{equation}
	By \cite{vershynin2010introduction}[Lemma 5.2], there exists $\mathcal{N}_{1/2}$, a $1/2$-net of $\{x \in \RR^n: \|x\|_2 = 1\}$, such that $\left|\mathcal{N}_{1/2}\right| \leq 5^n$. The union bound, \cite{vershynin2010introduction}[Lemma 5.2] and \eqref{ineq33} together imply that 
	\begin{equation*}
	\begin{split}
	\PP\left(\|AZB\| \geq 2\sqrt{m + t}\right) \leq \PP\left(\max_{x \in \mathcal{N}_{1/2}}\|AZBx\|_{2} \geq \sqrt{m + t}\right) \leq 2\cdot 5^n\exp\left[-c\min\left(\frac{t^2}{m}, t\right)\right].
	\end{split}
	\end{equation*}
	
	For $\|AZB\|_{\rm F}$, note that $AZB = (B^\top \otimes A)\text{vec}(Z)$, 
	Similarly to \eqref{ineq32}, we have
	\begin{equation*}
	\begin{split}
	\EE\|(B^\top \otimes A)\text{vec}(Z)\|_2^2 =& \EE \left[\text{vec}(Z)^\top (B^\top \otimes A)^\top(B^\top \otimes A)\text{vec}(Z)\right]\\
	=& \EE \text{trace}\left[\text{vec}(Z)^\top (B^\top \otimes A)^\top(B^\top \otimes A)\text{vec}(Z)\right]\\
	=& \text{trace}\EE\left[(B^\top \otimes A)^\top(B^\top \otimes A)\text{vec}(Z)\text{vec}(Z)^\top\right]\\
	=& \text{trace}\left[(B^\top \otimes A)^\top(B^\top \otimes A)\EE\left(\text{vec}(Z)\text{vec}(Z)^\top\right)\right]\\
	\leq& \text{trace}\left[(B^\top \otimes A)^\top(B^\top \otimes A)\right]\\
	=& \|B^\top \otimes A\|_{\rm F}^2 = \|B\|_{\rm F}^2\|A\|_{\rm F}^2\\ \leq& mn.
	\end{split}
	\end{equation*}
	By Hanson-Wright inequality, we have
	\begin{equation*}
	\PP\left(\|AZB\|_{\rm F}^2 - mn \geq t\right) \leq 2\exp\left[-c\min\left(\frac{t^2}{\|(BB^\top) \otimes (A^\top A)\|_{\rm F}^2}, \frac{t}{\|(BB^\top) \otimes (A^\top A)\|}\right)\right].
	\end{equation*}
	Since $\|A\|, \|B\| \leq 1$, we have 
	\begin{equation*}
	\begin{split}
	&\|(BB^\top) \otimes (A^\top A)\|_{\rm F} = \sqrt{\|A^\top A\|_{\rm F}^2\|BB^\top\|_{\rm F}^2} = \sqrt{\sum_{i=1}^{\min\{m, p\}}\sigma^4_i(A)\sum_{i=1}^{\min\{q, n\}}\sigma^4_i(B)} \leq \sqrt{mn},\\
	&\|(BB^\top) \otimes (A^\top A)\| \leq 1.
	\end{split}
	\end{equation*}
	Therefore,
	\begin{equation*}
	\PP\left(\|AZB\|_{\rm F}^2\geq mn + t\right) \leq 2\exp\left[-c\min\left(\frac{t^2}{mn}, t\right)\right].
	\end{equation*}
\end{proof}

\begin{lemma}\label{lm:epsilon_net_Frobenius}
		For the class of low-rank tensors under the Frobenius norm $\calX_{\bp, \br} = \{\calA \in \RR^{p_1 \times p_2 \times p_3}: \rank(\calM_i(\calA)) \leq r_i, i \in [3], \|\calA\|_{\rm F} \leq 1\}$, there exists 
		\begin{equation}\label{ineq48}
			\bar\calX_{\bp, \br} = \{\calA^{(1)}, \dots, \calA^{(N)}\}
		\end{equation}
		with $N \leq ((8 + \epsilon)/\epsilon)^{r_1r_2r_3 + \sum_{i=1}^{3}p_ir_i}$ satisfying $\calA^{(i)} \in \RR^{p_1 \times p_2 \times p_3}: \|\calA^{(i)}\|_{\rm F} \leq 1$, such that for all $\calA \in \calX_{\bp, \br}$, there exists $i \in [N]$ satisfying $\|\calA^{(i)} - \calA\|_{\rm F} \leq \epsilon$.
\end{lemma}
\begin{proof}[Proof of Lemma \ref{lm:epsilon_net_Frobenius}]
	By \cite[Lemma 7]{zhang2018tensor}, there exist $\epsilon/4$-nets $\bar{\calX}_{p_i, r_i}$ for $\calX_{p_i, r_i} = \{U \in \RR^{p_i \times r_i}: \|U\| \leq 1\}$ under the spectral norm with cardinality at most $((8+\epsilon)/\epsilon)^{p_ir_i}$, $i \in [3]$, and $\bar{\calX}_{r_1, r_2r_3}$ for $\calX_{r_1, r_2r_3} = \{B \in \RR^{r_1 \times (r_2r_3)}: \|B\|_{\rm F} \leq 1\}$ under the Frobenius norm with cardinality at most $((8+\epsilon)/\epsilon)^{r_1r_2r_3}$. Let $$\bar\calX_{\bp, \br} = \{\calB \times_1 U_1 \times_2 U_2 \times_3 U_3: U_i \in \bar\calX_{p_i, r_i}, \calM_1(\calB) \in \bar{\calX}_{r_1, r_2r_3}\}.$$ For any $\calA \in \calX_{\bp, \br}$, there exist $U_i \in \calX_{p_i, r_i}$ and $D_1 \in \calX_{r_1, r_2r_3}$ such that $\calM_1(\calA) = V_1D_1(V_2^\top \otimes V_3^\top)$. Then we can find $U_i^* \in \bar{\calX}_{p_i, r_i}$ and $\calB^* \in \RR^{r_1 \times r_2 \times r_3}, B_1^* = \calM_1(\calB^*) \in \bar{\calX}_{r_1, r_2r_3}$, and $\calB^* \times_1 U_1^* \times_2 U_2^* \times_3 U_3^* \in \bar\calX_{\bp, \br}$ satisfying
	\begin{align*}
		&\|\calA - \calB^* \times_1 U_1^* \times_2 U_2^* \times_3 U_3^*\|_{\rm F}\\ =& \|V_1D_1(V_2^\top \otimes V_3^\top) - U_1^*B_1^*(U_2^{*\top} \otimes U_3^{*\top})\|_{\rm F}\\ \leq& \|(V_1 - U_1^*)D_1(V_2^\top \otimes V_3^\top)\|_{\rm F} + \|U_1^*(D_1 - B_1^*)(V_2^\top \otimes V_3^\top)\|_{\rm F} + \|U_1^*B_1^*((V_2 - U_2^*)^\top \otimes V_3^\top)\|_{\rm F}\\
		&+ \|U_1^*B_1^*(U_2^{*\top} \otimes (V_3 - U_3^*)^\top)\|_{\rm F}\\
		\leq& \|V_1 - U_1^*\|\|V_2\|\|V_3\|\|D_1\|_{\rm F} + \|U_1^*\|\|V_2\|\|V_3\|\|D_1 - B_1^*\|_{\rm F} + \|U_1^*\|\|V_2 - U_2^*\|\|V_3\|\|B_1^*\|_{\rm F}\\ &+ \|U_1^*\|\|U_2^*\|\|V_3 - U_3^*\|\|B_1^*\|_{\rm F}\\
		\leq& \frac{\epsilon}{4} + \frac{\epsilon}{4} + \frac{\epsilon}{4} + \frac{\epsilon}{4} = \epsilon.
	\end{align*}
	Notice that $|\bar\calX_{\bp, \br}| \leq |\bar{\calX}_{p_1, r_1}||\bar{\calX}_{p_2, r_2}||\bar{\calX}_{p_3, r_3}||\bar{\calX}_{r_1, r_2r_3}| \leq ((8 + \epsilon)/\epsilon)^{r_1r_2r_3 + \sum_{i=1}^{3}p_ir_i}$, we have finished the proof of Lemma \ref{lm:epsilon_net_Frobenius}.
\end{proof}

\begin{lemma}\label{lm:Gaussian_ensemble}
	\begin{itemize}
		\item[(1)] Suppose $X \in \RR^{p_1 \times p_2}, X(i_1,i_2) \stackrel{i.i.d.}{\sim} N(0, 1)$ and $X_1, \dots, X_n$ are i.i.d. copies of $X$. Then there exist two universal constants $C, C_1 > 0$ such that for any fixed $U \in \OO_{p_1, r_1}, V \in \OO_{p_2, r_2}$ and $\Delta \in \RR^{p_1 \times p_2}$,
		\begin{equation}\label{ineq47}
		\PP\left(\left\|\frac{1}{n}\sum_{i=1}^{n}\langle X_i, \Delta\rangle U^\top X_iV - U^\top\Delta V\right\| \geq C\|\Delta\|_{\rm F}t\right) \leq 2\cdot 7^{r_1 + r_2}e^{-C_1\min\{nt^2, nt\}}.
		\end{equation} 
		\item[(2)] Suppose $\calX \in \RR^{p_1 \times p_2 \times p_3}, \calX(i_1,i_2,i_3) \stackrel{i.i.d.}{\sim} N(0, 1)$ and $\calX_1, \dots, \calX_n$ are i.i.d. copies of $\calX$. Then
		\begin{equation*}
		\begin{split}
		&\PP\left(\sup_{\begin{subarray}{c}
			U_i \in \RR^{p_i \times r_i}, \|U_i\| \leq 1\\\calA \in \RR^{p_1 \times p_2 \times p_3}, \|\calA\|_{\rm F} \leq 1\\ \rank(\calA) \leq (\bar r_1, \bar r_2, \bar r_3)
			\end{subarray}}\left\|\frac{1}{n}\sum_{i=1}^{n}\langle \calX_i, \calA\rangle \calM_1(\calX_i)(U_2 \otimes U_3) - \calM_1(\calA) (U_2 \otimes U_3)\right\| \geq Ct\right)\\ \leq& 2\cdot 7^{p_1 + r_2r_3}9^{p_2r_2 + p_3r_3}33^{\bar r_1\bar r_2\bar r_3 + \sum_{i=1}^{3}p_i\bar r_i}e^{-C_1\min\{nt^2, nt\}}.
		\end{split}
		\end{equation*}
	\end{itemize}
\end{lemma}
\begin{proof}[Proof of Lemma \eqref{lm:Gaussian_ensemble}]
	\begin{itemize}[leftmargin = *]
		\item[(1)] We only need to show that \eqref{ineq47} holds for any fixed $U \in \OO_{p_1, r_1}, V \in \OO_{p_2, r_2}$. For any fixed $a \in \RR^{r_1}, b \in \RR^{r_2}$ satisfying $\|a\|_2 = 1, \|b\|_2 = 1$, notice that $\EE[\langle X_i, \Delta\rangle X_i] = \Delta$ for $i \in [n]$, we have
		$$\EE[\langle X_i, \Delta\rangle a^\top U^\top X_iVb] = a^\top U^\top\Delta Vb, \quad \forall i \in [n].$$
For any random variable $Y_1$ and $Y_2$, by Cauchy-Schwarz inequality, we have
\begin{equation}\label{ineq:psi_1}
	\begin{split}
	\|Y_1Y_2\|_{\psi_1} \leq & C\sup_{q \geq 1}\frac{1}{q}\left(\EE|Y_1Y_2|^q\right)^{1/q} \leq C\left[\sup_{q \geq 1}\frac{1}{\sqrt{2q}}\left(\EE|Y_1|^{2q}\right)^{\frac{1}{2q}}\right]\left[\sup_{q \geq 1}\frac{1}{\sqrt{2q}}\left(\EE|Y_1|^{2q}\right)^{\frac{1}{2q}}\right]\\ 
	\leq & C\|Y_1\|_{\psi_2}\|Y_2\|_{\psi_2}.
	\end{split}
\end{equation} 
Since $\langle X_i, \Delta\rangle \sim N(0, \|\Delta\|_{\rm F}^2)$ and $a^\top U^\top X_iVb \sim N(0,1)$, by \cite[Remark 5.18]{vershynin2010introduction} and the above inequality, we have
\begin{equation*}
	\begin{split}
	&\left\|\langle X_i, \Delta\rangle a^\top U^\top X_iVb - a^\top U^\top\Delta Vb\right\|_{\psi_1}\\ \leq& C\left\|\langle X_i, \Delta\rangle a^\top U^\top X_iVb\right\|_{\psi_1}\\ \leq& C\|\langle X_i, \Delta\rangle\|_{\psi_2}\|a^\top U^\top X_iVb\|_{\psi_2}\\
	\leq& C\|\Delta\|_{\rm F}.
	\end{split}
\end{equation*}
By Bernstein-type inequality, we have
\begin{equation*}
	\PP\left(\left|\sum_{i=1}^{n}\frac{1}{n}\left(\langle X_i, \Delta\rangle a^\top U^\top X_iVb - a^\top U^\top\Delta Vb\right)\right| \geq C\|\Delta\|_{\rm F}t\right) \leq 2\exp\left[-C_1\min\{nt^2, nt\}\right].
\end{equation*}
By \cite[Lemma 5.2]{vershynin2010introduction}, there exist a $1/3$-net $\calN_1$ for $S^{r_1 - 1} = \{x: x \in \RR^{r_1}, \|x\|_2 = 1\}$ with cardinality at most $7^{r_1}$ and a $1/3$-net $\calN_2$ for $S^{r_2 - 1} = \{x: x \in \RR^{r_2}, \|x\|_2 = 1\}$ with cardinality at most $7^{r_2}$. By the union bound, we have
\begin{equation}\label{ineq46}
	\PP\left(\sup_{a \in \calN_1, b \in \calN_2}\left|a^\top\left[\sum_{i=1}^{n}\frac{1}{n}\left(\langle X_i, \Delta\rangle U^\top X_iV -  U^\top\Delta V\right)\right]b\right| \geq C\|\Delta\|_{\rm F}t\right) \leq 2\cdot 7^{r_1 + r_2}e^{-C_1\min\{nt^2, nt\}}.
\end{equation}
Let $a^* \in S^{r_1 - 1}$ and $b^* \in S^{r_2 - 1}$ satisfy 
\begin{equation*}
	\begin{split}
	& \left|a^{*\top}\left[\sum_{i=1}^{n}\frac{1}{n}\left(\langle X_i, \Delta\rangle U^\top X_iV -  U^\top\Delta V\right)\right]b^*\right|\\ = & \sup_{a \in S^{r_1 - 1}, b \in S^{r_2 - 1}}\left|a^\top\left[\sum_{i=1}^{n}\frac{1}{n}\left(\langle X_i, \Delta\rangle U^\top X_iV -  U^\top\Delta V\right)\right]b\right|.
	\end{split}
\end{equation*}
Then there exist $\tilde{a} \in \calN_1$ and $\tilde{b} \in \calN_2$ such that $\left\|\tilde{a} - a^*\right\| \leq \frac{1}{3}$ and $\left\|\tilde{b} - b^*\right\| \leq \frac{1}{3}$. Therefore,
\begin{equation*}
	\begin{split}
	&\left\|\sum_{i=1}^{n}\frac{1}{n}\left(\langle X_i, \Delta\rangle U^\top X_iV -  U^\top\Delta V\right)\right\|\\ =& \left|a^{*\top}\left[\sum_{i=1}^{n}\frac{1}{n}\left(\langle X_i, \Delta\rangle U^\top X_iV -  U^\top\Delta V\right)\right]b^*\right|\\ \leq& \left|\tilde a^{\top}\left[\sum_{i=1}^{n}\frac{1}{n}\left(\langle X_i, \Delta\rangle U^\top X_iV -  U^\top\Delta V\right)\right]\tilde b\right| + \left|\left(a^* - \tilde{a}\right)^{\top}\left[\sum_{i=1}^{n}\frac{1}{n}\left(\langle X_i, \Delta\rangle U^\top X_iV -  U^\top\Delta V\right)\right]{\tilde b}^*\right|\\&+ \left|a^{*\top}\left[\sum_{i=1}^{n}\frac{1}{n}\left(\langle X_i, \Delta\rangle U^\top X_iV -  U^\top\Delta V\right)\right]\left(b^* - \tilde b\right)\right|\\ \leq& \sup_{a \in \calN_1, b \in \calN_2}\left|a^\top\left[\sum_{i=1}^{n}\frac{1}{n}\left(\langle X_i, \Delta\rangle U^\top X_iV -  U^\top\Delta V\right)\right]b\right| + \frac{2}{3}\left\|\sum_{i=1}^{n}\frac{1}{n}\left(\langle X_i, \Delta\rangle U^\top X_iV -  U^\top\Delta V\right)\right\|,
	\end{split}
\end{equation*}
which means that 
\begin{equation}\label{ineq49}
	\left\|\sum_{i=1}^{n}\frac{1}{n}\left(\langle X_i, \Delta\rangle U^\top X_iV -  U^\top\Delta V\right)\right\| \leq 3\sup_{a \in \calN_1, b \in \calN_2}\left|a^\top\left[\sum_{i=1}^{n}\frac{1}{n}\left(\langle X_i, \Delta\rangle U^\top X_iV -  U^\top\Delta V\right)\right]b\right|.
\end{equation}
Combining the previous inequality and \eqref{ineq46}, we have proved the first part.
\item[(2)] For fixed $U_i \in \RR^{p_ir_i}$ and $\calA \in \RR^{p_1 \times p_2 \times p_3}$ satisfying $\|U_i\| \leq 1$ and $\|\calA\|_{\rm F} \leq 1$, by \eqref{ineq47}, we have
\begin{equation*}
	\PP\left(\left\|\frac{1}{n}\sum_{i=1}^{n}\langle \calX_i, \calA\rangle  \calM_1(\calX_i)(U_2 \otimes U_3) - \calM_1(\calA) (U_2 \otimes U_3)\right\| \geq Ct\right) \leq 2\cdot 7^{p_1 + r_2r_3}e^{-C_1\min\{nt^2, nt\}}.
\end{equation*}
By \cite[Lemma 7]{zhang2018tensor}, there exist $1/4$-nets $\bar{\calX}_{p_i, r_i}$ for $\calX_{p_i, r_i} = \{U \in \RR^{p_i \times r_i}: \|U\| \leq 1\}$ with cardinality at most $9^{p_ir_i}$. Therefore, by the union bound, we have
\begin{equation*}
	\begin{split}
	&\PP\left(\sup_{\begin{subarray}{c}
		U_i \in \bar{\calX}_{p_i, r_i}\\\calA \in \bar{\calX}_{\bp, \bar\br} 
		\end{subarray}}\left\|\frac{1}{n}\sum_{i=1}^{n}\langle \calX_i, \calA\rangle U_1^\top \calM_1(\calX_i)(U_2 \otimes U_3) - U_1^\top\calM_1(\calA) (U_2 \otimes U_3)\right\| \geq Ct\right)\\ \leq& 2\cdot 7^{p_1 + r_2r_3}9^{p_2r_2 + p_3r_3}33^{\bar r_1\bar r_2\bar r_3 + \sum_{i=1}^{3}p_i\bar r_i}e^{-C_1\min\{nt^2, nt\}},
	\end{split}
\end{equation*}
where $\bar{\calX}_{\bp, \bar\br}$ is defined in \eqref{ineq48} with $\epsilon = 1/5$.\\
Similarly to \eqref{ineq49}, we have
\begin{equation*}
	\begin{split}
	&\sup_{\begin{subarray}{c}
		U_i \in \RR^{p_i \times r_i}, \|U_i\| \leq 1\\\calA \in \RR^{p_1 \times p_2 \times p_3}, \|\calA\|_{\rm F} \leq 1\\ \rank(\calA) \leq (\bar r_1, \bar r_2, \bar r_3)
		\end{subarray}}\left\|\frac{1}{n}\sum_{i=1}^{n}\langle \calX_i, \calA\rangle U_1^\top \calM_1(\calX_i)(U_2 \otimes U_3) - U_1^\top\calM_1(\calA) (U_2 \otimes U_3)\right\|\\
	\leq& 4\sup_{\begin{subarray}{c}
		U_i \in \bar{\calX}_{p_i, r_i}\\\calA \in \bar{\calX}_{\bp, \br} 
		\end{subarray}}\left\|\frac{1}{n}\sum_{i=1}^{n}\langle \calX_i, \calA\rangle U_1^\top \calM_1(\calX_i)(U_2 \otimes U_3) - U_1^\top\calM_1(\calA) (U_2 \otimes U_3)\right\|.
	\end{split}
\end{equation*}
By combining the two previous inequalities together, we have finished the proof of the second part.
\end{itemize}
\end{proof}

\begin{lemma}\label{lm:sub_Gaussian}
	\begin{itemize}
		\item[(1)] Suppose $X \in \RR^{p_1 \times p_2}$ is a matrix with independent entries satisfying $\EE X_{ij} = 0, \Var(X_{ij}) = 1, \|X_{ij}\|_{\psi_2} \leq C$, and $X_1, \dots, X_n$ are i.i.d. copies of $X$. Suppose $\xi_1, \dots, \xi_n$ are independent zero-mean $C\sigma$-sub-Gaussian random variables. For any fixed $U \in \OO_{p_1, r_1}$ and $V \in \OO_{p_2, r_2}$, we have 
		\begin{equation*}
			\PP\left(\left\|\sum_{i=1}^{n}\xi_iU^\top X_iV\right\| \geq C_2\sqrt{n}\sqrt{r_1 + r_2 + x}\sigma\right) \leq e^{-C_1x} + e^{-c_1n}.
		\end{equation*}
		\begin{equation*}
		\PP\left(\left\|\sum_{i=1}^{n}\xi_iU^\top X_iV\right\|_{\rm F} \geq C_2\sqrt{n}\sqrt{r_1r_2 + x}\sigma\right) \leq e^{-C_1x} + e^{-c_1n}.
		\end{equation*}
		\item[(2)] Suppose $\calX \in \RR^{p_1 \times p_2 \times p_3}$ is a tensor with independent entries satisfying $\EE \calX_{ijk} = 0, \Var(\calX_{ijk}) = 1, \|\calX_{ijk}\|_{\psi_2} \leq C$, and $\calX_1, \dots, \calX_n$ are i.i.d. copies of $\calX$. Suppose $\xi_1, \dots, \xi_n$ are independent zero-mean $C\sigma$-sub-Gaussian random variables. Let $p = \max_{j = 1, 2, 3}p_j$ and $r = \max\{r_1, r_2, r_3\}$. Suppose $r \leq \sqrt{p}$. Then for fixed $V_i \in \OO_{p_i, r_i}$, 
		\begin{equation}\label{ineq:concentration_noise_2}
		\PP\left(\left\|\sum_{i = 1}^{n}\xi_iV_1^\top\calM_1(\calX_i)(V_2 \otimes V_3)\right\| \geq C_2\sqrt{n}\sqrt{r_1 + r_2r_3 + \log(p)}\sigma\right) \leq p^{-C_1} + e^{-c_1n}.
		\end{equation}
		\begin{equation}\label{ineq:concentration_noise_1}
		\begin{split}
		\PP\left(\left\|\sum_{i = 1}^{n}\xi_iV_1^\top\calM_1(\calX_i)(V_2 \otimes V_3)\right\|_{\rm F} \geq C_2\sqrt{n}\sqrt{r_1r_2r_3 + \log(p)}\sigma\right) \leq p^{-C_1} + e^{-c_1n}.
		\end{split}
		\end{equation}
		Moreover,
		\begin{equation}\label{ineq38}
		\begin{split}
		&\PP\left(\sup_{V_i \in \OO_{p_i, r_i}}\left\|\sum_{i = 1}^{n}\xi_iV_1^\top\calM_1(\calX_i)(V_2 \otimes V_3)\right\| \geq C_2\sqrt{npr}\sigma\right) \leq e^{-C_1pr} + e^{-C_1n}.
		\end{split}
		\end{equation}
		\begin{equation}\label{ineq39}
		\begin{split}
		\PP\left(\sup_{V_i \in \OO_{p_i, r_i}}\left\|\sum_{i = 1}^{n}\xi_iV_1^\top\calM_1(\calX_i)(V_2 \otimes V_3)\right\|_{\rm F} \geq C_2\sqrt{npr}\sigma\right) \leq e^{-C_1pr} + e^{-C_1n}.
		\end{split}
		\end{equation}
		\item[(3)] Suppose the conditions in (2) holds and $p_1 \geq r_2r_3 \geq r_1$. Then for fixed $V_2 \in \OO_{p_2, r_2}$ and $V_3 \in \OO_{p_3, r_3}$ and $V_4 \in \OO_{r_2r_3, r_1}$, 
		\begin{equation*}
			\begin{split}
			&\PP\left(\sup_{V_i \in \OO_{p_i, r_i}}\left\|\sum_{i = 1}^{n}\xi_i\calM_1(\calX_i)(V_2 \otimes V_3)V_4\right\| \geq C_2\sqrt{np}\sigma\right) \leq e^{-C_1p} + e^{-C_1n}.
			\end{split}
		\end{equation*}
		\begin{equation*}
			\begin{split}
			&\PP\left(\sup_{V_i \in \OO_{p_i, r_i}}\left\|\sum_{i = 1}^{n}\xi_i\calM_1(\calX_i)(V_2 \otimes V_3)V_4\right\|_{\rm F} \geq C_2\sqrt{npr}\sigma\right) \leq e^{-C_1pr} + e^{-C_1n}.
			\end{split}
		\end{equation*}
		In addition,
		\begin{equation*}
			\begin{split}
			\PP\left(\sup_{V_i \in \OO_{p_i, r_i}, i = 2, 3\atop V_4 \in \OO_{r_2r_3, r_1}}\left\|\sum_{i = 1}^{n}\calM_1(\calX_i)(V_2 \otimes V_3)V_4\right\| \geq C_2\sqrt{npr}\sigma\right) \leq e^{-C_1pr} + e^{-C_1n}.
			\end{split}
		\end{equation*}
		\begin{equation*}
		\begin{split}
		\PP\left(\sup_{V_i \in \OO_{p_i, r_i}, i = 2, 3\atop V_4 \in \OO_{r_2r_3, r_1}}\left\|\sum_{i = 1}^{n}\calM_1(\calX_i)(V_2 \otimes V_3)V_4\right\|_{\rm F} \geq C_2\sqrt{npr}\sigma\right) \leq e^{-C_1pr} + e^{-C_1n}.
		\end{split}
		\end{equation*}
	\end{itemize}
\end{lemma}
\begin{proof}[Proof of Lemma \ref{lm:sub_Gaussian}]
	Without loss of generality, we assume $\sigma = 1$. For any fixed $a = (a_1, \dots, a_n) \in \RR^{n}$, noting that the entries of $\sum_{i = 1}^{n}a_iX_i$ are independent $C\|a\|_2$-sub-Gaussian random variables with mean 0 and variance $\|a\|_2^2$. By Lemma \ref{lm:AZB}, for any $x \geq 0$,
	\begin{equation*}
	\PP\left(\left\|\sum_{i = 1}^{n}a_iU^\top X_iV\right\| \geq C\|a\|_2\sqrt{r_1 + r_2 + x}\right) \leq e^{-C_1x}.
	\end{equation*}
	\begin{equation*}
	\PP\left(\left\|\sum_{i = 1}^{n}a_iU^\top X_iV\right\|_{\rm F} \geq C\|a\|_2\sqrt{r_1r_2 + x}\right) \leq e^{-C_1x}.
	\end{equation*}
	Therefore,
	\begin{equation}\label{ineq35}
	\PP\left(\left\|\sum_{i = 1}^{n}\xi_iU^\top X_iV\right\| \geq C\|\xi\|_2\sqrt{r_1 + r_2 + x}\bigg|\xi_1, \dots, \xi_n\right) \leq e^{-C_1x}.
	\end{equation}
	\begin{equation}\label{ineq36}
	\PP\left(\left\|\sum_{i = 1}^{n}\xi_iU^\top X_iV\right\|_{\rm F} \geq C\|\xi\|_2\sqrt{r_1r_2 + x}\bigg|\xi_1, \dots, \xi_n\right) \leq e^{-C_1x}.
	\end{equation}
	By Bernstein-type inequality (\cite[Proposition 5.16]{vershynin2010introduction}),
	\begin{equation}\label{ineq37}
	\PP\left(\|\xi\|_2 \geq C\sqrt{n}\right) \leq e^{-C_1n}.
	\end{equation}
	Thus we have
	\begin{equation*}
	\PP\left(\left\|\sum_{i = 1}^{n}\xi_iU^\top X_iV\right\| \geq C\sqrt{n}\sqrt{r_1 + r_2 + x}\right) \leq e^{-C_1x} + e^{-C_1n}.
	\end{equation*}
	\begin{equation*}
	\PP\left(\left\|\sum_{i = 1}^{n}\xi_iU^\top X_iV\right\|_{\rm F} \geq C\sqrt{n}\sqrt{r_1r_2r_3 + x}\right) \leq e^{-C_1x} + e^{-C_1n}.
	\end{equation*}
	By setting $x = \log(p), U = U_1, X_i = \calM(\calX_i)$ and $V = U_2 \otimes U_3$ in the previous two inequalities, we have proved \eqref{ineq:concentration_noise_2} and \eqref{ineq:concentration_noise_1}, respectively.
	
	By setting $t = C\sqrt{pr}$ in \eqref{ineq35} and \eqref{ineq36} and using Lemma \ref{lm:epsilon_net}, \eqref{ineq37} and the $\epsilon$-net argument, we have \eqref{ineq38} and \eqref{ineq39}.
	
	Similarly, we can prove the inequalities in Part 3.
\end{proof}

\begin{proof}[Proof of Lemma \ref{lm:orlicz}]
	By the definition of $\psi_{\alpha}$-norm, 
	\begin{equation*}
		\EE \left[\exp\left(\left|\frac{X_i}{K_i}\right|^{\alpha_i}\right)\right] \leq 2, \quad \forall i \in [n].
	\end{equation*}
	By the weighted AM-GM inequality, we have
	\begin{equation*}
		\left|\frac{\prod_{i=1}^n X_i}{\prod_{i=1}^n K_i}\right|^{\frac{1}{\sum_{i=1}^{n}\frac{1}{\alpha_i}}} \leq \sum_{i=1}^{n}\frac{\frac{1}{\alpha_i}}{\sum_{i=1}^{n}\frac{1}{\alpha_i}}\left|\frac{X_i}{K_i}\right|^{\alpha_i},
	\end{equation*}
	and 
	\begin{align*}
		\exp\left(\left|\frac{\prod_{i=1}^n X_i}{\prod_{i=1}^n K_i}\right|^{\frac{1}{\sum_{i=1}^{n}\frac{1}{\alpha_i}}}\right) \leq& \exp\left(\sum_{i=1}^{n}\frac{\frac{1}{\alpha_i}}{\sum_{i=1}^{n}\frac{1}{\alpha_i}}\left|\frac{X_i}{K_i}\right|^{\alpha_i}\right)\\ =& \prod_{i=1}^n\exp\left(\frac{\frac{1}{\alpha_i}}{\sum_{j=1}^{n}\frac{1}{\alpha_j}}\left|\frac{X_i}{K_i}\right|^{\alpha_i}\right)\\
		\leq& \sum_{i=1}^{n}\frac{\frac{1}{\alpha_i}}{\sum_{j=1}^{n}\frac{1}{\alpha_j}}\exp\left(\left|\frac{X_i}{K_i}\right|^{\alpha_i}\right).
	\end{align*}
	By taking expectations on both sides, we have
	\begin{equation*}
		\EE\left[\exp\left(\left|\frac{\prod_{i=1}^n X_i}{\prod_{i=1}^n K_i}\right|^{\frac{1}{\sum_{i=1}^{n}\frac{1}{\alpha_i}}}\right)\right] \leq \sum_{i=1}^{n}\frac{\frac{1}{\alpha_i}}{\sum_{j=1}^{n}\frac{1}{\alpha_j}}\EE\left[\exp\left(\left|\frac{X_i}{K_i}\right|^{\alpha_i}\right)\right] \leq 2,
	\end{equation*}
	which has finished the proof of Lemma \eqref{lm:orlicz}
\end{proof}

\begin{proof}[Proof of Lemma \ref{lm:third_moment_gaussian}]
	Since $\frac{\langle Z_i, M_i\rangle}{\|M_i\|_{\rm F}} \sim N(0, 1)$, $\|\frac{1}{\|M_i\|_{\rm F}}\langle Z_i, M_i\rangle\|_{\psi_2} \leq C$. 
	Notice that $\|Z_i\|_{\rm F}^2 \sim \chi_{pr}^2$ by \cite[Example 2.8]{wainwright2019high}, we have
	\begin{equation*}
	\EE\left(e^{t(\|Z_i\|_{\rm F}^2 - pr)}\right) \leq e^{2prt^2}, \quad \forall |t| < \frac{1}{4}.
	\end{equation*}
	Set $t = \sqrt{\frac{\log(2)}{2pr}} \leq 1/4$, we have
	\begin{equation*}
		\EE\left(e^{\frac{\|Z_i\|_{\rm F}^2 - pr}{\sqrt{2pr/\log(2)}}}\right) \leq 2.
	\end{equation*}
	Therefore, $\left\|\|Z_i\|_{\rm F}^2 - pr\right\|_{\psi_1} \leq C\sqrt{pr}$. By Lemma \ref{lm:orlicz}, 
	\begin{equation*}
		\left\|(\|Z_i\|_{\rm F}^2 - pr)\frac{\langle Z_i, M_i\rangle}{\|M_i\|_{\rm F}}\right\|_{\psi_{2/3}} \leq \left\|\frac{1}{\|M_i\|_{\rm F}}\langle Z_i, M_i\rangle\right\|_{\psi_2}\left\|\|Z_i\|_{\rm F}^2 - pr\right\|_{\psi_1} \leq C\sqrt{pr}.
	\end{equation*}
	Since $$\EE\left(\sum_{i=1}^{n}(\|Z_i\|_{\rm F}^2 - pr)\langle Z_i, M_i\rangle\right) = 0,$$ by \cite[Lemma 8]{hao2020sparse}, with probability at least $1 - p^{-C_1}$,
	\begin{equation*}
		\left|\sum_{i=1}^{n}(\|Z_i\|_{\rm F}^2 - pr)\langle Z_i, M_i\rangle\right| \leq C\sqrt{pr}\left(\left(\sum_{i=1}^{n}\|M_i\|_{\rm F}^2\right)^{1/2}\sqrt{\log(p)} + \left(\max_{1 \leq i \leq n}\|M_i\|_{\rm F}\right)(\log(p))^{3/2}\right).
	\end{equation*}
	In addition, since $\sum_{i=1}^{n}\langle Z_i, M_i\rangle \sim N(0, \sum_{i=1}^{n}\|M_i\|_{\rm F}^2)$, 
	\begin{equation*}
		\PP\left(\left|pr\sum_{i=1}^{n}\langle Z_i, M_i\rangle\right| \geq Cpr\left(\sum_{i=1}^{n}\|M_i\|_{\rm F}^2\right)^{1/2}\sqrt{\log(p)}\right) \leq p^{-C_1}.
	\end{equation*}
	By combining the previous two inequalities, we have finished the proof of Lemma \ref{lm:third_moment_gaussian}.
\end{proof}
	\begin{lemma}\label{lm:asymptotic}
	Suppose $k$ and $d$ are fixed numbers satisfying $1 \leq k \leq d$. $Y_p \stackrel{\rm d.}{\to} N(0, I_d)$ as $p \to \infty$. Then for any deterministic matrix array $\{A_p\}_{p=1}^\infty$ satisfying $A_p \in \OO_{d,k}$, we have
	\begin{equation*}
	A_p^\top Y_p \stackrel{\rm d.}{\to} N(0, I_k) \quad \text{\rm as} \quad p \to \infty.
	\end{equation*}
\end{lemma}
\begin{proof}[Proof of Lemma \ref{lm:asymptotic}]
	By Skorohod's theorem \cite[Theorem 1.8]{shao2003mathematical}, there exist random vectors $\{Z_p\}, Z$ defined on a common probability space such that $Z \sim N(0, I_3)$, $Y_p \stackrel{\rm d.}{=} Z_p$ and $Z_p \stackrel{\rm a.s.}{\to} Z$ as $p \to \infty$.
	Notice that
	\begin{align*}
	A_p^\top Z \sim N(0,1)
	\end{align*}
	and
	\begin{align*}
	&\|A_p^\top (Z_p - Z)\|_2 \leq \|A_p\|\|Z_p - Z\|_2
	= \|Z_p - Z\|_2 \stackrel{\rm a.s.}{\to} 0  \quad \text{\rm as} \quad p \to \infty,
	\end{align*}
	we have
	\begin{align*}
	A_p^\top Y_p \stackrel{\rm d.}{\to} A_p^\top Z_p = A_p^\top Z + A_p^\top(Z_p - Z) \stackrel{\rm d.}{\to} N(0, I_k) \quad \text{\rm as} \quad p \to \infty.
	\end{align*}
\end{proof}

\end{sloppypar}

\end{document}